\documentclass[12pt,reqno]{amsart}
\usepackage{amsmath,amsfonts,amsthm,color,amssymb}
\usepackage[usenames,dvipsnames]{xcolor}
\usepackage[T1]{fontenc}
\usepackage{dsfont}
\usepackage{wasysym}
\usepackage{ulem}
\usepackage{stmaryrd} 
\usepackage{stackrel}
\usepackage[left=1in, right=1in, top=1.1in,bottom=1.1in]{geometry}
\usepackage{enumitem}
\usepackage{hyperref}
\hypersetup{
     colorlinks   = true,
     citecolor    = blue,
     linkcolor=blue
}
\usepackage{cancel}
\usepackage{comment}

\numberwithin{equation}{section}
\theoremstyle{plain}
\newtheorem{thm}{\protect\theoremname}[section]
\theoremstyle{definition}
\newtheorem{defn}[thm]{\protect\definitionname}
\theoremstyle{plain}
\newtheorem{prop}[thm]{\protect\propositionname}
\theoremstyle{remark}
\newtheorem{rem}[thm]{\protect\remarkname}
\theoremstyle{remark}
\newtheorem{notation}[thm]{\protect\notationname}
\theoremstyle{plain}
\newtheorem{lem}[thm]{\protect\lemmaname}
\theoremstyle{plain}
\newtheorem{cor}[thm]{\protect\corollaryname}
\theoremstyle{plain}
\newtheorem{assumption}[thm]{\protect\assumptionname}

\providecommand{\assumptionname}{Assumption}
\providecommand{\definitionname}{Definition}
\providecommand{\lemmaname}{Lemma}
\providecommand{\notationname}{Notation}
\providecommand{\propositionname}{Proposition}
\providecommand{\remarkname}{Remark}
\providecommand{\theoremname}{Theorem}
\providecommand{\corollaryname}{Corollary}

\makeatother

\newcommand{\hb}[1]{\textcolor{blue}{#1}}
\newcommand{\hre}[1]{\textcolor{red}{#1}}

\renewcommand{\emph}{\textit}
\renewcommand{\le}{\leqslant}

\newcommand{\R}{\mathbb R}
\newcommand{\N}{\mathbb N}

\newcommand{\Z}{\mathbb Z}

\newcommand{\ca}{\mathcal A}

\newcommand{\ck}{\mathcal K}
\newcommand{\cl}{\mathcal L}

\newcommand{\cw}{\mathcal W}

\newcommand{\al}{\alpha}
\newcommand{\ep}{\varepsilon}

\newcommand{\si}{\sigma}
\newcommand{\te}{\theta}

\newcommand{\lp}{\left(}
\newcommand{\rp}{\right)}

\newcommand{\lln}{\left|}
\newcommand{\rrn}{\right|}

\begin{document}
\title[Random walks in random environment]
{A Coupling Between Random Walks in Random Environments and Brox's
Diffusion}
\author[X. Geng \and M. Gradinaru \and S. Tindel]
{Xi Geng, Mihai Gradinaru and Samy Tindel}
\date{}

\newcommand{\Addresses}{{% additional braces for segregating \footnotesize
  \bigskip
  \footnotesize

 \noindent
  X.~Geng, \textsc{School of Mathematics and Statistics, The University of Melbourne, Parkville, VIC 3052, Australia.}\\
  \textit{E-mail address}: \texttt{xi.geng@unimelb.edu.au}

  \medskip

  \noindent
 M.~Gradinaru, \textsc{Univ Rennes, CNRS, IRMAR – UMR 6625, F-35000 Rennes, France.}\\
  \textit{E-mail address}: \texttt{mihai.gradinaru@univ-rennes.fr}

  \medskip

  \noindent
 S.~Tindel, \textsc{S. Tindel: Department of Mathematics,
Purdue University, West Lafayette, IN 47907, USA.}\\
  \textit{E-mail address}: \texttt{stindel@purdue.edu}

}}

\maketitle

\begin{abstract}
It has been established in~\cite{Se} that a properly rescaled version of Sinai's random walk converges in distribution to Brox's diffusion. In this article we quantify this convergence by considering a specific coupling between Sinai's walk and Brox's diffusion. Our method relies on convergence results for martingale problems considered in the rough path setting.
\end{abstract}

\tableofcontents{}
\section{Introduction}
Sinai's random walk is a very popular model of random walk in random environment. It is fairly simple to describe at a mathematical level, and yet it exhibits a highly non-trivial asymptotic behaviour as time goes to infinity. Let us describe a simplified version of this object (a more complete version being laid out in Section \ref{sec:prelim-sinai}).

The state space for Sinai's random walk $\{X_{n}^{d}; n\geqslant 1\}$ is $\mathbb Z$. The walk is based on a random environment $\{\omega_{x}^{+}; x\in\mathbb Z\}$ which gives the probability of a right jump when one reaches a level $x\in\mathbb Z$. The i.i.d. sequence $\{\omega_{x}^{+}; x\in\mathbb Z\}$ is defined on a probability space $(\Omega,\mathcal G, \mathbf P)$. In the original model treated by Sinai~\cite{Si}, once $\omega$ is fixed,the random walk $X^d$ is defined through the so-called 
quenched transition probabilities 
\begin{equation}\label{SinaiRW-int}
\mathbb{P}^{\omega}(X_{n+1}^{d}=x+1|X_{n}^{d}=x)=\omega_{x}^{+},\quad
\mathbb{P}^{\omega}(X_{n+1}^{d}=x-1|X_{n}^{d}=x)=\omega_{x}^{-}=
1-\omega_{x}^{+}.
\end{equation}
Then the so-called the following assumption \emph{recurrence assumption} is spelled out as
\begin{equation}\label{SinaiRW-rec-int}
\mathbf{E}[\log(\omega_{x}^{-}/\omega_{x}^{+})]=0,\quad\text{ for all }x\in\mathbb Z. 
\end{equation}
Under hypothesis \eqref{SinaiRW-rec-int} it is shown in \cite{So} that $X^d$ is recurrent (almost surely with respect to the randomness in $\omega$). Moreover (see \cite{Si}) and \cite{Ke76}) the asymptotic behavior of  the path $n\mapsto X^d_n$ exhibits a highly nontrivial behaviour of the form 
\begin{equation}\label{SinaiRW-iterlog-int}
\frac{X^d_n}{\log^2{n}}\stackrel{{\rm (d)}}{\longrightarrow}L, 
\end{equation}
where the limit in distribution is considered with respect to the annealed 
probability $\mathbf P(d\omega)\times\mathbb P^{\omega}$ and where the law of $L$ 
is described in \cite{Ke76}. Fascinating behaviours like \eqref{SinaiRW-iterlog-int} have converted random walks in random environments into very popular objects in discrete probability and mathematical physics. 

\begin{rem}
A typical example of distribution satisfying \eqref{SinaiRW-rec-int} is 
$\omega_{x}^{+}\sim\text{\rm Beta}(a,a)$ distribution for a parameter $a\in(0,\infty)$. In this case it is well known that 
\[\mathbf E\left[\log(\omega_{x}^{+})\right]=\psi(a)-\psi(2a),\]
where $\psi$ is the digamma function. Therefore it is readily checked that $\omega_x$ satisfies the recurrence condition \eqref{SinaiRW-rec-int}. 
\end{rem}

\begin{rem}
Due to some intricate parity issues, we shall in fact work with a lazy version of Sinai's random walk. In order to keep technical details to a minimum in the introduction, we postpone a full description of this lazy random walk to Section~\ref{sec:prelim-sinai}.
\end{rem}

The second object of interest in this article is called Brox diffusion, which has to be considered as a Brownian motion evolving in a Brownian environment. This continuous time process $\{X_t^c;t\geqslant 0\}$ is based on an environment $\{W(x);x\in\mathbb R\}$ which is given as a double-sided Brownian motion. We will define $W$ on the same same space of probability $(\Omega, \mathcal G,\mathbf P)$ as the environment $\omega$ of the Sinai random walk. Then $X^c$ is formally the solution of the following stochastic differential equation: 
\begin{equation}\label{Broxdiff-int}
X_t^c=-\frac{1}{2}\int_{0}^{t}\dot{W}(X_{s}^{c})\,ds+B_{t},
\end{equation} 
where $B$ is a standard Brownian motion independent of $W$. Notice that the drift $\dot{W}$ in \eqref{Broxdiff-int} is a distribution which lies in a Sobolev space of regularity $\alpha=-\frac{1}{2}-\varepsilon$. The roughness rules out the possibility to solve \eqref{Broxdiff-int} in a pathwise sense (see e.g. \cite{BC} or \cite{CG} for optimal results concerning pathwise definitions of stochastic differential equations with distributional drifts). Therefore the process $X^c$ is usually seen (like in the original contribution~\cite{Br}) as the Brownian motion $B$ composed with a properly defined scale function. A behaviour similar to \eqref{SinaiRW-iterlog-int} is proved for $X^c$ in \cite{Br}. 
The similarities pointed out above prompted the community to think that Brox's diffusion $X^c$ is the continuous time equivalent of the random walk $X^d$. This claim has been made precise in the remarkable paper \cite{Se}. In this contribution a version $X^\delta$ of $X^d$, for $\delta>0$, is considered. The process $X^\delta$ is defined through a proper rescaling of time, space and environment $\omega$ (see Section \ref{sec:heuristics} for a detailed description). Then for a fixed time horizon $T$ it is proved in \cite{Se} that 
\begin{equation}\label{SinBro-cv-int}
X^\delta_{|[0,T]}\stackrel{{\rm (d)}}{\longrightarrow}X^c_{|[0,T]}, 
\end{equation}
where the convergence is considered with respect to the Skorokhod metric and also with respect to the annealed probability (as in \eqref{SinaiRW-iterlog-int}). 

The current contribution has to be seen as a progress in the direction of \eqref{SinBro-cv-int}. Indeed, Donsker type theorems like \eqref{SinBro-cv-int}
do not provide any type of knowledge about rates of convergence in distribution. However for diffusion processes or rough differential equations, a good wealth of information is available for weak type convergences of discretisations (see e.g. 
\cite{BT}, \cite{LLT}). In this article we wish to establish a rate of convergence for the limit in distribution \eqref{SinBro-cv-int}. Our main findings can be summarized as follows (a more quantitative version will be spelled out in Section \ref{sec:ConvEst}). 

\begin{thm}\label{SinBro-thm-int}
Let $X^\delta$ be the lazy version of Sinai's random walk as described in Section~\ref{sec:prelim-sinai}, properly rescaled as in Section \ref{rescaled-Sinai-walk}. Consider the weak solution $X^c$ to the equation \eqref{Broxdiff-int}, as well as a function $h\in{\rm C}_{b}^{3}$. Then there exists a coupling $(X^\delta,X^c)_{\delta>0}$ such that for all $t\in[0,T]$ and $\delta\in(0,1)$ we have 
\begin{equation}\label{SinBro-rate-int}
\Big|\mathbb E^\omega\left[h(X_t^c)\right]-\mathbb E^\omega\left[h(X_t^\delta)\right]\Big|\leqslant C_{h,T}(\omega)\delta^{\frac{1}{17}},
\end{equation}
where $C_{h,T}(\omega)$ is a random constant which only depends on $h,T$ and the environment $\omega$. 
\end{thm} 

\noindent
A few comments about Theorem \ref{SinBro-thm-int} are in order:

\noindent
{\bf (1)}\quad 
To the best of our knowledge, equation \eqref{SinBro-rate-int} provides the first quantitative convergence result for the convergence of $X^\delta$ to $X^c$. However, our $\frac{1}{17}$ rate of convergence is not expected to be optimal. Indeed, as the reader will see, our considerations are mostly based on rough path analysis applied in a Brownian context (with H\"older regularity $\alpha=\frac{1}{2}-$). Moreover without anticipating too much on technical details, let us mention that we will use controlled process expansions of order 1 (which lead to a proper rough paths expansion whenever $\alpha>\frac{1}{3}$). We would expect our rate of convergence to be of order
\[
\max_{\alpha\in(0,1/2)}\min\left\{\left(\frac{1}{2}-\alpha\right),\alpha\right\}=\frac{1}{4}.
\]This is seen from (\ref{eq:OpExp}) by taking $\tau\uparrow 1/2-\alpha $, $\beta' \uparrow \alpha$ and $\beta \downarrow0$. However, one has to consider higher order rough path expansions and higher order controlled processes to make this work. Even so, at the moment it is still not clear whether $1/4$ is expected to be the optimal convergence rate.

\noindent
{\bf (2)}\quad Let us also highlight the fact that \eqref{SinBro-rate-int} is a quenched type result. It would be more consistent with \eqref{SinBro-cv-int} to obtain an annealed rate. This would depend on integrability properties of the random constant $C_{h,T}(\omega)$. We anticipate this integrability to be similar to the one obtained in \cite{CLL} for the Jacobian of rough differential 
equation. 

\noindent
{\bf (3)}\quad A rate of convergence for the law of $X^\delta$ in total variation, Wasserstein or other classical distances for probability measure seems to be out of reach at this moment. Indeed, our computations for \eqref{SinBro-rate-int} will involve integrations by parts, for which some derivatives of $g$ are needed. 

\noindent
{\bf (4)}\quad An important aspect of our method is the fact that we are based on our explicit coupling between the rescaled Sinai random walks $X^\delta$ 
and the continuous process $X^c$. While some other type of coupling can already be found in \cite{HS}, we believe our coupling might be interesting in its own right. It certainly allows to transfer some information from $X^\delta$ to $X^c$. 

\noindent
We plan to delve deeper into those aspects in future works. 

In order to prove the main Theorem \ref{SinBro-thm-int}, our analysis hinges on two main ingredients. First, as detailed later in Section \ref{sec:mild-pde-brox}, a weak solution to the Brox diffusion equation \eqref{Broxdiff-int} can be apprehended through its martingale problem. This amounts to consider a family of PDE's of the form 
\begin{equation}\label{BroPDE-int}
\partial_{t}f_{t}(x)-\mathcal L^c f(x)=g_{t}(x),\quad t\in[0,T],\;x\in\mathbb R, 
\end{equation}
where $g$ is a sufficiently smooth function and where the operator $\mathcal L^c$ is given by 
\begin{equation}\label{BroIG-int}
\mathcal L^c f(x)=\frac{1}{2}f''(x)-\frac{1}{2}\dot{W}(x)f'(x).
\end{equation}
Note that equation \eqref{BroPDE-int} is formal at this point, since $\dot{W}$ in \eqref{BroIG-int} is a distribution. An important contribution in \cite{DD} has been to give a pathwise interpretation for a mild form of~\eqref{BroPDE-int}, thanks to rough paths techniques. We are inspired by this approach here. A substantial part of our efforts in the paper consist in considering a discrete version of \eqref{BroPDE-int} and take limits in the discretisation parameter. The intricate technical details are provided in Section \ref{sec:ConvEst}. One should mention at this point the interesting alternative approach in~\cite{HLM} to pathwise interpretations of equation \eqref{Broxdiff-int}. We have decided to stick to~\cite{DD} in our contribution, since it is much more likely to be extended to multidimensional settings. 

The second crucial ingredient in our strategy is related to pathwise approximations in Donsker's theorem. namely the discrete version of the operator  \eqref{BroIG-int} involves a rescaled random walk increment called 
$\dot{U}^\delta$ (see \eqref{eq:def-U-delta} for the definition). Using some classical results by Komlos, Major and Tusnady \cite{KMT76}, we are able to couple $\dot{U}^\delta$ and $\dot{W}$. This coupling will be the backbone of our coupling between $X^\delta$ and $X^c$. Observe that our rough path type approach will force us to state and prove an extension of \cite{KMT76} to H\"older type norms. Also observe that the methodology outlined above is not restricted to the Sinai random walk setting. We plan to explore other applications, like branching processes and super-Brownian motion, in the next future.
 
Our article is structured as follows: in Section~\ref{Sinai-rw} we recall basic facts about Sinai's random walk, we properly define its scaling and we introduce related martingale problems. We also include estimates (some of them new) for the discrete heat kernel  which are essential for the sequel. Section~\ref{BMBE} is dedicated to a definition of Brox's diffusion through martingale problems considered in the rough paths sense. In particular, we will introduce the rough path setting employed throughout the article. In Section~\ref{RPMP} we derive the rough paths estimates ensuring a proper solution to the martingale problem of Section~\ref{BMBE}. While this section is not totally new when compared with respect to~\cite{DD}, our simplified setting yields clearer calculations. Moreover, Section~\ref{RPMP} lays the ground for our convergence analysis. Section~\ref{sec:KMT} delves into the strong Donsker type approximations which are at the heart of our method. Building on the original contribution~\cite{KMT76}, we obtain a strong approximation result in weighted H\"older norms. Eventually, Section~\ref{sec:ConvEst} contains the bulk of our convergence estimates, combining elements contained in the previous sections.

\section{Sinai's random walk}\label{Sinai-rw}

In this section we collect some basic facts about Sinai\textquoteright s
random walk and define its renormalized version on a grid whose mesh
goes to $0$. As anticipated in the introduction, in order to avoid periodicity problems for random walks we will handle a lazy version of Sinai's walk. 

\subsection{Preliminaries on Sinai's walk}\label{sec:prelim-sinai}

In this section we properly define a lazy version of Sinai\textquoteright s random walk
and write some related martingale problems which turn out to be crucial
for our limiting procedure.

\subsubsection{Definition of Sinai's walk}

In order to define our random walk, we first characterize the random
environment under consideration. In our case of interest, it is given
by a sequence of independent random variables and a parameter $\ep\in(0,1)$ quantifying the randomness of the walk.
\begin{defn}
\label{def:environment}Let $\varepsilon\in(0,1)$ be a given small number
which is fixed throughout the rest of the paper. The random environment
is given by a sequence of i.i.d. random variables $\omega^{+}=\{\omega_{x}^{+}:x\in\mathbb{Z}\}$
defined on a probability space $(\Omega,\mathcal{G},\mathbf{P})$
and satisfying the following conditions:

\begin{comment}
\vspace{1mm}\noindent (i) ${\bf P}$-almost surely each $\omega_{x}^{+}$
takes values in $[0,1-\varepsilon]$.\\
(ii) Symmetry assumption: $\mathbb{E}[\omega_{x}^{+}]=\frac{1-\varepsilon}{2}$.
\end{comment}

\vspace{1mm}\noindent (i) [Ellipticity] ${\bf P}$-almost surely each $\omega_{x}^{+}$
takes values in $[\kappa,1-\varepsilon-\kappa]$ where $\kappa$ is some given fixed strictly positive number.\\
(ii) [Recurrence] $\mathbf{E}[\log(\omega_{x}^{-}/\omega_{x}^{+})]=0$ where $\omega_{x}^{-}\triangleq1-\varepsilon-\omega_{x}^{+}$.\\
(iii) [Regularity] The distribution of $\omega_x^+$ has a $\mathcal{C}^1$ density with at most finitely many algebraic singularities.

\end{defn}

\begin{rem}
The third condition is only set for technical convenience in the application of the classical KMT approximation theorem (cf. Remark \ref{rem:DensityCond} for a more detailed discussion). 
\end{rem}

\begin{notation}\label{variance}
In the sequel $1-\varepsilon$ will play the role of a variance parameter. From now on we will thus set $\sigma^{2}=1-\varepsilon$.  
\end{notation}

Having defined our environment $\omega$, we now introduce the random walk itself. 

\begin{comment}
\hre{Why do we impose $\mathbb{E}[\omega_{x}^{+}]=\frac{1-\varepsilon}{2}$ here? When $\ep=0$ we get $\mathbb{E}[\omega_{x}^{+}]=\frac{1}{2}$, which is not the usual $\mathbb{E}[\ln(\frac{1-\omega_{x}^{+}}{\omega_{x}^{+}})]=0$ for recurrence.}
\end{comment}

\begin{defn}
\label{def:SinRWZ}Given the environment of Definition \ref{def:environment},
one can construct a random walk $X^{d}$ (called \textit{Sinai's random
walk}) on another probability space $(\hat{\Omega},\mathcal{F},\mathbb{P}^{\omega})$
in the following way:
\[
\mathbb{P}^{\omega}(X_{n+1}^{d}=y|X_{n}^{d}=x)\triangleq\begin{cases}
\varepsilon, & \text{if }y=x;\\
\omega_{x}^{\pm}, & \text{if }y=x\pm1;\\
0, & \text{otherwise. }
\end{cases}
\]
\end{defn}

Notice that in the above definition and in the sequel, the superscript
$d$ stands for \textit{discrete time parameter}. The probability
$\mathbb{P}^{\omega}$ is usually called the \textit{quenched probability},
for which the randomness of $\omega_{x}^{+}$ is fixed. Otherwise
stated, under $\mathbb{P}^{\omega}$ the process $\{X_{k}:k\geqslant0\}$
is a Markov chain on $\mathbb{Z}$ whose one-step transition matrix
$T^{d}$ can be written as
\begin{equation}
T^{d}f(x)=\omega_{x}^{+}f(x+1)+\omega_{x}^{-}f(x-1)+\varepsilon f(x).\label{eq:DisTranMatZ}
\end{equation}
We call $\mathcal{L}^{d}$ the \textit{discrete generator} of $T^{d},$
defined by
\begin{equation}
\mathcal{L}^{d}=T^{d}-{\rm Id}.\label{eq:DisGenZ}
\end{equation}
Notice that the discrete generator is often defined as ${\rm Id}-T^{d}$
in the literature. However, our notation \eqref{eq:DisGenZ} allows
a better transition to the continuous setting. Furthermore, starting
from (\ref{eq:DisTranMatZ}) one can easily deduce the form of the
discrete generator $T^{d}-{\rm Id}$. Below we state an elementary
proposition expressing this generator as a perturbed discrete Laplace
operator, which will make it easier to relate with its continuous
counterpart.
\begin{prop}
\label{prop:gen-rw-discrete}Let $X^{d}$ be the random walk given
by Definition \ref{def:SinRWZ}. Recall that the discrete generator
$\mathcal{L}^{d}$ of $X^{d}$ is given by (\ref{eq:DisGenZ}). For
$f\in L^{\infty}(\mathbb{Z})$, we set
\begin{equation}
\Delta^{d}f(x)\triangleq f(x+1)+f(x-1)-2f(x),\quad 
\hat{\nabla}f(x)\triangleq\frac{1}{2}[f(x+1)-f(x-1)].\label{eq:discr-laplacian-nabla}
\end{equation}
Also define the discrete potential 
\begin{equation}
\dot{U}(x)\triangleq\omega_{x}^{+}-\omega_{x}^{-}=2\omega_{x}^{+}-(1-\varepsilon)=2\omega_{x}^{+}-\sigma^{2}.\label{eq:def-dot-U}
\end{equation}
Then for $f\in L^{\infty}(\mathbb{Z})$ we have 
\begin{equation}\label{eq:geneLd}
\mathcal{L}^{d}f(x)  =\frac{\sigma^{2}}{2}\cdot\Delta^{d}f(x)+\dot{U}(x)\cdot\hat{\nabla}f(x).
\end{equation}
\end{prop}

\begin{proof}
Recalling that $\mathcal{L}^{d}=T^{d}-{\rm Id}$ we have 
\begin{align*}
\mathcal{L}^{d}f(x) & =\omega_{x}^{+}f(x+1)+\omega_{x}^{-}f(x-1)+\varepsilon f(x)-f(x)\\
 & =\frac{\sigma^{2}}{2}\cdot\left(f(x+1)+f(x-1)-2f(x)\right)+\left(\omega_{x}^{+}-\frac{\sigma^{2}}{2}\right)\left(f(x+1)-f(x-1)\right)\\
 & =\frac{\sigma^{2}}{2}\cdot\Delta^{d}f(x)+(2\omega_{x}^{+}-\sigma^{2})\cdot\hat{\nabla}f(x).
\end{align*}
Taking the definition \eqref{eq:def-dot-U} of $\dot{U}$ into account, this proves our claim \eqref{eq:geneLd}.
\end{proof}
\begin{rem}
The introduction of the ``lazy'' probability $\varepsilon$ is merely
for technical convenience. It avoids parity issues in discrete heat
kernel estimates. Indeed, if $\varepsilon=0$, the second order term
of the generator $\mathcal{L}^{d}$ corresponds to the simple random
walk, whose $n$-step transition probability function is only supported
on either even or odd integer points (depending on the parity of $n$
and the starting point). This creates non-trivial issues when estimating
discrete derivatives of the transition function. Resolution of such
issues requires substantially more technical effort in the analysis.
To reduce technicalities and focus more on the essential parts, we
choose to restrict ourselves to the aperiodic situation. Notice that the generator of the simple lazy random walk with laziness parameter $\ep$ is
\begin{equation}\label{b0}
\bar{\cl}_{x}^{d}f(x)\triangleq \frac{\si^{2}}{2} \lp f(x+1)+f(x-1)-2f(x) \rp,
\end{equation}
where we recall that $\si^{2}=1-\ep$.
\end{rem}

\begin{rem}
According to \cite[Theorem 2.1.2]{Ze}, the reccurrence assumption
in Definition \ref{def:environment} implies that
$X^{d}$ is recurrent. 
\end{rem}

\subsubsection{Some martingale problems}

Our analysis will hinge crucially on asymptotic properties for some
martingale problems related to $X^{d}$. We will thus start by labeling
the basic martingale problem for $X^{d}$.
\begin{prop}
\label{prop:discrete-mg-problem-basic} Let $X^{d}$ be the random
walk given by Definition \ref{def:SinRWZ}, and consider a function
$f\in L^{\infty}(\mathbb{Z})$. For $j\geqslant1$ set 
\begin{equation}
Z_{j}\triangleq f(X_{j}^{d})-T^{d}f(X_{j-1}^{d})=f(X_{j}^{d})-\mathbb{E}_{\omega}\left[f(X_{j}^{d})|\mathcal{F}_{j-1}\right]\quad\text{and}\quad M_{k}\triangleq\sum_{j=1}^{k}Z_{j}.\label{eq:def-Z-and-M}
\end{equation}
Then for $k\geqslant1$ we have 
\begin{equation}
f(X_{k}^{d})-f(X_{0}^{d})-\sum_{j=0}^{k-1}\mathcal{L}^{d}f(X_{j}^{d})=M_{k},\label{eq:discrete-martingale-problem}
\end{equation}
where $\mathcal{L}^{d}$ is the operator introduced in \eqref{eq:geneLd} and the process $M=M^{f}$ is a $\mathbb{P}^{\omega}$-martingale.
\end{prop}

\begin{proof}
We obtain relation (\ref{eq:discrete-martingale-problem}) thanks
to some elementary algebraic manipulations (we refer e.g to \cite[Equation (1.3)]{KLO}
for further details). Indeed, a simple telescopic sum argument reveals
that 
\[
f(X_{k}^{d})-f(X_{0}^{d})=\sum_{j=1}^{k}\left(f(X_{j}^{d})-f(X_{j-1}^{d})\right).
\]
We now insert terms of the form $T^{d}f(X_{j-1}^{d})$ in order to
get 
\[
f(X_{k}^{d})-f(X_{0}^{d})=\sum_{j=1}^{k}\left(f(X_{j}^{d})-T^{d}f(X_{j-1}^{d})\right)+\sum_{j=1}^{k}\left(T^{d}-{\rm Id}\right)f(X_{j-1}^{d}).
\]
Since $\mathcal{L}^{d}=T^{d}-{\rm Id}$, this immediately yields 
\[
f(X_{k}^{d})-f(X_{0}^{d})-\sum_{j=0}^{k-1}\mathcal{L}^{d}f(X_{j}^{d})=\sum_{j=1}^{k}\left(f(X_{j}^{d})-T^{d}f(X_{j-1}^{d})\right).
\]
With our notation (\ref{eq:def-Z-and-M}) in mind, relation~(\ref{eq:discrete-martingale-problem})
is now easily proved. The fact that $M$ is a martingale is also readily
checked.
\end{proof}
In the sequel we will need to introduce some space-time stochastic
equations, for which additional notation has to be introduced.
\begin{notation}
\label{not:gradient-time-or-space} In the remainder of the article
the subscripts in $\nabla_{n}$, $T_{x}^{d}$ etc. denote the variable
concerned by the operator at stake (generally either a time or space
variable). For instance, the discrete time gradient of a function
$f\in L^{\infty}(\mathbb{N}\times\mathbb{Z})$ is written as 
\begin{equation}\label{eq:nablan}
\nabla_{n}f(k,x)=f(k+1,x)-f(k,x).
\end{equation}
\end{notation}

\noindent With Notation \ref{not:gradient-time-or-space} in hand,
we now define a space-time martingale problem related to Sinai's random
walk.
\begin{prop}
\label{prop:discrete-space-time-martingale-problem} Let $X^{d}$
be the random walk given by Definition \ref{def:SinRWZ}, and consider
a function $f\in L^{\infty}(\mathbb{N}\times\mathbb{Z})$. For $j\geqslant1$
we set 
\begin{equation}
Z_{j}\triangleq f(j-1,X_{j}^{d})-T_{x}^{d}f(j-1,X_{j-1}^{d})\quad\text{and}\quad M_{k}\triangleq\sum_{j=1}^{k}Z_{j}.\label{eq:def-ZjN}
\end{equation}
Recall that the discrete generator $\mathcal{L}^{d}$ of $X^{d}$
is given by (\ref{eq:DisGenZ}). Then for $k\geqslant1$ we have 
\begin{equation}
f(k,X_{k}^{d})-f(0,X_{0}^{d})-\sum_{l=0}^{k-1}\left[\mathcal{L}_{x}^{d}f(l,X_{l}^{d})-\nabla_{n}f(l,X_{l+1}^{d})\right]=M_{k}\,,\label{eq:discrete-space-time-martingale-problem}
\end{equation}
and the process $M=M^{f}$ is a $\mathbb{P}^{\omega}$-martingale.
\end{prop}

\begin{proof}
Since $T_{x}^{d}$ defined by (\ref{eq:DisTranMatZ}) is the transition
matrix of our random walk, for any positive integer $k$ we have $T_{x}^{d}f(k,X_{j-1}^{d})=\mathbb{E}_{\omega}[f(k,X_{j}^{d})|\mathcal{F}_{j-1}]$.
Hence the time increment $Z_{j}$ defined by (\ref{eq:def-ZjN}) can
be written as 
\[
Z_{j}=f(j-1,X_{j}^{d})-\mathbb{E}_{\omega}\left[f(j-1,X_{j}^{d})|\mathcal{F}_{j-1}\right],
\]
from which it is easily deduced that $M$ is a martingale. Moreover,
a simple telescoping sum argument allows us to write 
\begin{equation}
f(k,X_{k}^{d})-f(0,X_{0}^{d})=\sum_{j=1}^{k}\left[f(j,X_{j}^{d})-f(j-1,X_{j-1}^{d})\right].\label{c1}
\end{equation}
For any $j\in\{1,\ldots,k\}$, we decompose the terms on the right
hand side of (\ref{c1}) as 
\[
f(j,X_{j}^{d})-f(j-1,X_{j-1}^{d})=\left[f(j,X_{j}^{d})-f(j-1,X_{j}^{d})\right]+\left[f(j-1,X_{j}^{d})-f(j-1,X_{j-1}^{d})\right].
\]
Therefore, recalling Notation \ref{not:gradient-time-or-space} together
with the definitions (\ref{eq:def-ZjN}) and (\ref{eq:DisGenZ}),
we obtain 
\begin{align*}
 & f(j,X_{j}^{d})-f(j-1,X_{j-1}^{d})\\
 & =\nabla_{n}f(j-1,X_{j}^{d})+Z_{j}+\left[T^{d}f(j-1,X_{j-1}^{d})-f(j-1,X_{j-1}^{d})\right]\\
 & =Z_{j}+\mathcal{L}_{x}^{d}f(j-1,X_{j-1}^{d})+\nabla_{n}f(j-1,X_{j}^{d}).
\end{align*}
Plugging this information into (\ref{c1}) and setting $l=j-1$ in
the sum, this yields our claim~(\ref{eq:discrete-space-time-martingale-problem}).
\end{proof}

\subsection{Rescaled version of Sinai's walk}\label{rescaled-Sinai-walk}

In this section, we shall describe the dynamics of Sinai's walk when
this process is accelerated in time and rescaled in space. This will
generate a process which should converge to the Brownian motion in
a Brownian environment. In the sequel, we will simply use $\delta>0$
to denote the generic scaling parameter, and one would like to rescale
$X^{d}$ according to $\delta$. However, a naive approach to this
rescaling procedure yields a convergence to a standard Brownian motion
as $\delta\rightarrow0$ (see e.g. \cite{Na}). In order to get convergence to the Brownian
motion in a Brownian environment, we must also renormalize the environment
$\omega$ in a proper way. The precise renormalization procedure is
specified as follows, starting with some of the heuristic steps in \cite{Se}. Namely going back at least to \cite{Si}, the analysis of $X^{d}$ relies on a potential $V^{d}:\mathbb Z\to\mathbb R$ which can be expressed as 
\[
V(x)=\sum_{j\in\llbracket 0,x\rrbracket}\log(\xi_{j}),
\]
where each random variable $\xi_{x}$ is defined by 
\begin{equation}\label{eq:xi-rw}
\xi_{x}=\log\left(\frac{\sigma^{2}-\omega_{x}^{+}}{\omega_{x}^{+}}\right)
=\log\left(\frac{\omega_{x}^{-}}{\omega_{x}^{+}}\right).
\end{equation}
The potential $V$ is used e.g. to express hitting probabilities for $X^{d}$, see \cite{Ze}. Now we can easily invert \eqref{eq:xi-rw} and write 
\begin{equation}\label{eq:omega-rw}
\omega_{x}^{+}=\frac{\sigma^{2}}{1+e^{\xi_{x}}}=
\frac{\sigma^{2}}{1+e^{\log(\omega_{x}^{-}/\omega_{x}^{+})}}.
\end{equation}
The scaling which is given below is then based on a scaling of $\xi_{x}$ which enables to have each $\omega_{x}^{+}$ close to $\sigma^{2}/2$.
\begin{defn}
\label{def:EnvScale}Let $\omega^{+}=\{\omega_{x}^{+}:x\in\mathbb{Z}\}$
be a given random environment that satisfies Definition \ref{def:environment}.
Recall that $\omega_{x}^{-}=1-\varepsilon-\omega_{x}^{+}=\sigma^{2}-\omega_{x}^{+}$, where $\varepsilon=1-\sigma^{2}$
is the given fixed ``lazy'' probability. For each fixed $\delta>0$,
we define a rescaled version of $\omega^{+}$ on the grid $\delta\mathbb{Z}$
by 
\begin{comment}
\begin{equation}
\omega_{x}^{+,\delta}\triangleq\frac{\sigma^{2}+\sqrt{\delta}(\omega_{x/\delta}^{+}-\omega_{x/\delta}^{-})}{2}\quad\text{and}\quad \omega_{x}^{-,\delta}\triangleq \sigma^{2}-\omega_{x}^{+,\delta},\ \ \ x\in\delta\mathbb{Z}.\label{eq:rescaled-omega-pm}
\end{equation}
%We also set $\omega_{x}^{-,\delta}\triangleq \sigma^{2}-\omega_{x}^{+,\delta}$.
\end{comment}
\begin{equation}
\omega_{x}^{+,\delta}\triangleq\frac{\sigma^{2}}{1+e^{\sqrt{\delta}\log(\omega_{x/\delta}^{-}/\omega_{x/\delta}^{+})}}\quad\text{ and }\quad\omega_{x}^{-,\delta}\triangleq\sigma^{2}-\omega_{x}^{+,\delta},\ \ \text{ for all }\ x\in\delta\mathbb{Z}.\label{eq:rescaled-omega-pm}
\end{equation}
\end{defn}
\begin{comment}
\hre{Why do we have $\omega_{x/\delta}^{+}-\omega_{x/\delta}^{-}$ in the definition of $\omega_{x}^{+,\delta}$? For $\ep=0$ we didn't have that.}
\end{comment}
We now describe a rescaled random walk in the rescaled environment given by 
\eqref{eq:rescaled-omega-pm} in a way which mimics Definition \ref{def:SinRWZ}

\begin{defn}\label{def:rescSinRWZ}
Given the rescaled environment in Definition \ref{def:EnvScale} we define a Sinai type random walk $\hat{X}^{\delta}$ on $\hat{\Omega},\mathcal F,\mathbb P^{\omega})$ with state space $\delta\mathbb Z$ by specifying the 
transitions 
\[
\mathbb{P}^{\omega}(\hat{X}_{n+1}^{\delta}=y|\hat{X}_{n}^{\delta}=x)\triangleq\begin{cases}
\varepsilon, & \text{if }y=x;\\
\omega_{x}^{\pm\delta}, & \text{if }y=x\pm\delta;\\
0, & \text{otherwise. }
\end{cases}
\]
\end{defn}
Our next aim is to obtain a convenient expression for the generator of the walk $\hat{X}^{\delta}$. This is achieved in the lemma below.

\begin{lem}
Let $\hat{X}^{\delta}$ be the walk introduced in Definition \ref{def:rescSinRWZ}. We introduce a rescaled discrete Laplace 
operator $\bar{\mathcal L}_{x}^{\delta}$ and a potential $\dot{U}^{\delta}$ 
on $\delta\mathbb Z$ as 
\begin{equation}
\bar{\mathcal{L}}_{x}^{\delta}f(x)\triangleq\frac{\sigma^{2}}{2\delta^{2}}\left[f(x+\delta)+f(x-\delta)-2f(x)\right].\label{eq:renormalized-discrete-laplace}
\end{equation}
\begin{equation}
\dot{U}^{\delta}(z)\triangleq\omega_{z}^{+,\delta}-\omega_{z}^{-,\delta}=2\omega_{z}^{+,\delta}-\sigma^{2}.\label{eq:def-U-delta}
\end{equation}
Also consider the twisted gradient $\hat{\nabla}_{x}^{\delta}$ defined by 
\begin{equation}
\hat{\nabla}_{x}^{\delta}f(x)\triangleq\frac{1}{2\delta}\left(f(x+\delta)-f(x-\delta)\right)=\frac{1}{2}\left(\nabla_{x}^{\delta}f(x)+\nabla_{x}^{\delta}f(x-\delta)\right),\label{eq:renormalized-discrete-nabla}
\end{equation}
where $\nabla_{x}^{\delta}f$ is defined on $\delta\mathbb{Z}$,  by 
\begin{equation}\label{c2}
\nabla_{x}^{\delta}f(x)=\frac{1}{\delta}\left(f(x+\delta)-f(x)\right).
\end{equation}
Then the generator $\mathcal L^{\delta}$ of $\hat{X}^{\delta}$ admits the 
expression 
\begin{equation}
\mathcal{L}^{\delta}f(x)
=\bar{\mathcal{L}}_{x}^{\delta}f(x)+\frac{1}{\delta}\dot{U}^{\delta}(x)\cdot\hat{\nabla}_{x}^{\delta}f(x).\label{d61}
\end{equation}
\end{lem}

\begin{proof}
The transition operator $T^{\delta}$ for $\hat{X}^{\delta}$ is given, for $f:\delta\mathbb Z\to\mathbb R$ as 
\begin{equation}
T^{\delta}f(x)=\omega_{x}^{+,\delta}f(x+\delta)+\omega_{x}^{-,\delta}f(x-\delta)+\varepsilon f(x).\label{eq:transition-matrix-rescaled}
\end{equation}
Thus, thanks to the relation $\si^{2}=1-\ep$ we have 
\begin{equation}\label{eq:Ldelta}
\mathcal{L}^{\delta}f(x)=\frac{1}{\delta^{2}}\left(T^{\delta}f(x)-f(x)\right)=\frac{1}{\delta^{2}}\left(\omega_{x}^{+,\delta}f(x+\delta)+\omega_{x}^{-,\delta}f(x-\delta)-\sigma^{2}f(x)\right).
\end{equation}
Therefore resorting to relation \eqref{eq:renormalized-discrete-laplace}  we can recast \eqref{eq:Ldelta} as 
\begin{equation*}
\mathcal{L}^{\delta}f(x)
%=\frac{\sigma^{2}}{2}\Delta^{\delta}f(x)+\frac{1}{\delta}\dot{U}^{\delta}(x)\cdot\hat{\nabla}_{x}^{\delta}f(x)
=\bar{\mathcal{L}}_{x}^{\delta}f(x)+\frac{1}{\delta^{2}}\left(\omega_{x}^{+,\delta}-\frac{\sigma^{2}}{2}\right)f(x+\delta)
+\frac{1}{\delta^{2}}\left(\omega_{x}^{-,\delta}-\frac{\sigma^{2}}{2}\right)f(x-\delta).
\end{equation*}
Now recalling from \eqref{eq:rescaled-omega-pm} that $\omega_{x}^{-,\delta}=\sigma^{2}-\omega_{x}^{+,\delta}$ we can write 
\begin{equation*}
\mathcal{L}^{\delta}f(x)
%=\frac{\sigma^{2}}{2}\Delta^{\delta}f(x)+\frac{1}{\delta}\dot{U}^{\delta}(x)\cdot\hat{\nabla}_{x}^{\delta}f(x)
=\bar{\mathcal{L}}_{x}^{\delta}f(x)+\frac{1}{\delta}\left(\omega_{x}^{+,\delta}-\frac{\sigma^{2}}{2}\right)\cdot\frac{1}{\delta}\left(f(x+\delta)
-f(x-\delta)\right).
\end{equation*}
With the definition~\eqref{eq:renormalized-discrete-nabla} of $\hat{\nabla}_{x}^{\delta}f$ in mind, our claim~\eqref{d61} is now easily proved.
\end{proof}

Let us also introduce the time partition and
the related discrete time derivative we will deal with in the sequel.
\begin{defn}
\label{def:partition-gradient} Let $T>0$ be a fixed time horizon.
In the remainder of the paper we consider $\Pi_{\delta}=\{t_{0},\ldots,t_{N}\}$
a division of the interval $[0,T]$ with step size $\delta^{2}$,
where $N=\lfloor T/\delta^{2}\rfloor$. A generic element of this
partition will be denoted by $t_{j}=j\delta^{2}$. We also introduce
the rescaled discrete gradient of a function $f\in L^{\infty}(\mathbb{N}\times\mathbb{Z})$,
given as follows for $t_{j}\in\Pi_{\delta}$ and $x\in\delta\mathbb{Z}$,
\begin{equation}
\nabla_{t}^{\delta}f_{t_{j}}(x)\triangleq\frac{1}{\delta^{2}}(f_{t_{j+1}}(x)-f_{t_{j}}(x)).\label{eq:def-nabla-discrete-renormalized}
\end{equation}
Finally we introduce a notation for discrete intervals, namely for
$s,t\in[0,T]$ we write $t_{j}\in\llbracket s,t\rrbracket$ for $t_{j}\in[s,t]\cap\Pi_{\delta}$.
We will also write $\llparenthesis s,t\rrbracket$ for $t_{j}\in(s,t]\cap\Pi_{\delta}$.
\end{defn}

Let us now define the rescaled random walks which will feature in
our coupling procedure.
\begin{defn}
\label{def:rescaled-rw} Throughout the paper, we designate by $\hat{X}^{\delta}=\{\hat{X}_{j}^{\delta};j\ge1\}$
the random walk on $\delta\mathbb{Z}$ with transition $T^{\delta}$
given by (\ref{eq:transition-matrix-rescaled}). Then the \textit{time-accelerated
random walk} $X^{\delta}$ considered below is given by 
\begin{equation}
X_{t}^{\delta}:=\hat{X}_{\lfloor t/\delta^{2}\rfloor}^{\delta}=\sum_{j\geqslant0}\hat{X}_{j}^{\delta}{\bf 1}_{[t_{j},t_{j+1})}(t),\label{eq:renormalized-process}
\end{equation}
where we have used the notation of Definition~\ref{def:partition-gradient}
for the partition $\Pi_{\delta}=\{t_{0},\ldots,t_{N}\}$. The filtration
related to the process $X^{\delta}$ is then given by $\{\mathcal{F}_{t_{j}}^{\omega,\delta}:t_{j}\in\Pi_{\delta}\}$,
with 
\begin{equation}
\mathcal{F}_{t_{j}}^{\omega,\delta}\triangleq\sigma\big\{\hat{X}_{k}^{\delta}:k\leqslant j\big\}=\sigma\big\{ X_{t_{k}}^{\delta}:k\leqslant j\big\}\label{eq:def-rescaled-filtration}
\end{equation}
and where the superscript $\omega$ in $\mathcal{F}_{t_{j}}^{\omega,\delta}$
means that the random environment $\omega$ is frozen.
\end{defn}

\begin{rem}\label{rk:variance}
The rescaling of the random environment given in Definition \ref{def:EnvScale}
is consistent with the one in \cite{Se} for the corresponding weak convergence result. 
%In that paper, the author introduced the scaling $\xi_{x}^{\delta}\triangleq\sqrt{\delta}\log\frac{1-\omega_{x}^{+}}{\omega_{x}^{+}}$
%and defined the rescaled environment $\omega_{x}^{+,\delta}\triangleq\frac{1}{1+\exp\xi_{x}^{\delta}}$
%accordingly ($\varepsilon=0$ in that context). Our Definition \ref{def:EnvScale}
%is essentially equivalent to this rescaling in the sense that both
Indeed, with relations~\eqref{eq:xi-rw} and 
\eqref{eq:def-U-delta} in mind, it is readily checked that 
\begin{equation}\label{eq:var-udot-delta}
{\rm Var}(\dot{U}_{x}^{\delta})
={\rm Var}\left(\sigma^{2}\left(\frac{2}{1+e^{\sqrt{\delta}\,\xi_{x/\delta}}}-1\right)\right)
=\sigma^{4}\,{\rm Var}(Z_{x}^{\delta}), 
\end{equation}
where we have set 
\[
Z_{x}^{\delta}=\frac{1-e^{\sqrt{\delta}\,\xi_{x/\delta}}}
{1+e^{\sqrt{\delta}\,\xi_{x/\delta}}}\,.
\] 
Now a first order approximation of $Z_{x}^{\delta}$ when $\delta\to 0$ is given by 
\begin{equation}\label{eq:1OrdApXi}
Z_{x}^{\delta}\sim \frac{\sqrt{\delta}}{2}\xi_{x/\delta}.
\end{equation}
Plugging this information into \eqref{eq:var-udot-delta} we obtain 
\[
{\rm Var}(\dot{U}_{x}^{\delta})\propto\delta,
\]
which is a key constraint for proving convergence under the current
perspective (cf. Section~\ref{sec:heuristics} below for a more
detailed explanation on this point).
\end{rem}

In the context of Definition \ref{def:rescaled-rw}, the martingale
problem (\ref{eq:discrete-martingale-problem}) can be rescaled on
the grid $\delta\mathbb{Z}$. We label this property in the proposition
below for further use.
\begin{prop}
For the random environment outlined in Definition \ref{def:EnvScale}
and $\delta\in(0,1)$, let $X^{\delta}$ be the rescaled random walk
introduced in Definition \ref{def:rescaled-rw}. Consider a function
$f\in L^{\infty}(\delta\mathbb{Z})$. Then for all $s,t\in\llbracket0,T\rrbracket$,
we have
\begin{equation}
f(X_{t}^{\delta})-f(X_{s}^{\delta})-\delta^{2}\sum_{t_{j}\in\llbracket s,t\rrparenthesis}\mathcal{L}^{\delta}f(X_{t_{j}}^{\delta})=M_{t}^{\delta}-M_{s}^{\delta},\label{eq:renormalize-discrete-martingale-problem}
\end{equation}
where the process $M^{\delta}=M^{\delta,f}$ is a $\mathbb{P}^{\omega}$-martingale
with respect to the filtration $\{\mathcal{F}_{t_{j}}^{\omega,\delta}\}$
introduced in (\ref{eq:def-rescaled-filtration}) given by 
\[
M_{t}^{\delta}\triangleq\sum_{t_{j}\in\llparenthesis0,t\rrbracket}Z_{t_{j}}^{\delta}\ \ \ \text{with}\ \ \ Z_{t_{j}}^{\delta}\triangleq f(X_{t_{j}}^{\delta})-T^{\delta}f(X_{t_{j-1}}^{\delta}).
\]
\end{prop}

Similarly to Proposition \ref{prop:discrete-space-time-martingale-problem},
we can also define a space-time renormalized discrete martingale problem
related to $X^{\delta}$. Its proof is very similar to the proof of
Proposition \ref{prop:discrete-space-time-martingale-problem} and
is thus omitted for sake of conciseness.
\begin{prop}
Recalling the notation of Definition \ref{def:partition-gradient}
and Definition \ref{def:rescaled-rw}, let $X^{\delta}$ be given
by (\ref{eq:renormalized-process}) and consider $f\in L^{\infty}(\llbracket0,T\rrbracket\times\mathbb{Z})$.
Then for any $s,t\in\llbracket0,T\rrbracket$ and any $\delta\in(0,1)$,
we have 
\begin{equation}
f_{t}(X_{t}^{\delta})-f_{s}(X_{s}^{\delta})-\sum_{t_{j}\in\llbracket s,t\rrparenthesis}\delta^{2}\left[\mathcal{L}^{\delta}f_{t_{j}}(X_{t_{j}}^{\delta})-\nabla_{t}^{\delta}f_{t_{j}}(X_{t_{j+1}}^{\delta})\right]=M_{t}^{\delta}-M_{s}^{\delta},\label{eq:renormalize-time-space-discrete-martingale-problem}
\end{equation}
where $M^{\delta}$ is a $\mathbb{P}^{\omega}$-martingale with respect
to the filtration $\{\mathcal{F}_{t_{j}}^{\omega,\delta}\}$ given
by 
\begin{equation}
M_{t}^{\delta}\triangleq\sum_{t_{j}\in\llparenthesis0,t\rrbracket}Z_{t_{j}}^{\delta}\quad\text{with}\quad Z_{t_{j}}^{\delta}\triangleq f_{t_{j}}(X_{t_{j}}^{\delta})-T^{\delta}f_{t_{j}}(X_{t_{j-1}}^{\delta}).\label{eq:martingale-discrete}
\end{equation}
\end{prop}

We close this section by giving a representation of the martingale
$M^{\delta}$ which will be useful in order to take limits to the
Brownian motion in Brownian environment.
\begin{prop}
\label{fId} Let $\{Z_{t_{j}}^{\delta};j\geqslant1\}$ be the sequence
of random variables defined by (\ref{eq:martingale-discrete}). We
introduce $\{\bar{\zeta}_{j}^{\delta}:j\geqslant1\}$ and $\{\zeta_{j}^{x,\delta}:j\geqslant1\}$
sequences of i.i.d. random variables defined by\begin{equation}
\bar{\zeta}_j^{\delta}\triangleq{\bf 1}_{\{U_j>\varepsilon,V_j\leqslant1/2\}}-{\bf 1}_{\{U_j>\varepsilon,V_j>1/2\}}\label{eq:def-zetas-bar}
\end{equation}
and 
\begin{equation}
\zeta_j^{x,\delta}\triangleq{\bf 1}_{\{U_j>\varepsilon,V_j\leqslant\omega_{x}^{+,\delta}/(1-\varepsilon)\}}-{\bf 1}_{\{U_j>\varepsilon,V_j>\omega_{x}^{+,\delta}/(1-\varepsilon)\}},\label{eq:def-zetas}
\end{equation}
where $\{(U_j,V_j):j\geqslant 1\}$ are independent copies of uniform random variables on $[0,1]$.
Then the following relation holds true in distribution:
\begin{align}\label{eq:Zn-as-stoch-integral}
Z_{t_{j+1}}^{\delta} & =\frac{1}{2}\left(f_{t_{j+1}}(X_{t_{j}}^{\delta}+\delta)-f_{t_{j+1}}(X_{t_{j}}^{\delta}-\delta)\right)\cdot\bar{\zeta}_{j+1}^{\delta} \\
 & \ \ \ +\frac{1}{2}\left(f_{t_{j+1}}(X_{t_{j}}^{\delta}+\delta)-f_{t_{j+1}}(X_{t_{j}}^{\delta}-\delta)\right)\cdot\left(\left(\zeta_{j+1}^{x,\delta}-\bar{\zeta}_{j+1}^{\delta}\right)-\mathbb{E}_{\omega}\left[\zeta_{j+1}^{x,\delta}-\bar{\zeta}_{j+1}^{\delta}|\mathcal{F}_{t_{j}}\right]\right)\nonumber \\
 & \ \ \ +\frac{1}{2}\left(f_{t_{j+1}}(X_{t_{j}}^{\delta}+\delta)-2f_{t_{j+1}}(X_{t_{j}}^{\delta})+f_{t_{j+1}}(X_{t_{j}}^{\delta}-\delta)\right)\cdot\left(\left(\bar{\zeta}_{j+1}^{\delta}\right)^{2}-\mathbb{E}_{\omega}\left[\left(\bar{\zeta}_{j+1}^{\delta}\right)^{2}|\mathcal{F}_{t_{j}}\right]\right), \nonumber
\end{align}
where $x \triangleq X_{t_j}^\delta$ and $\mathcal{F}_{t_{j}}$ is a slight variation of (\ref{eq:def-rescaled-filtration})
defined by 
\[
\mathcal{F}_{t_{j}}\triangleq\sigma\big\{\bar{\zeta}_{k}^{\delta},\zeta_{k}^{x,\delta}:k\leqslant j\big\}.
\]
\end{prop}

\begin{proof}
We invoke a generalization of a discrete It\^o formula stated in \cite{FK}.
Namely if $\{S_{t_{j}};t_{j}=j\delta^{2}\}$ is a random walk on $\delta\mathbb{Z}$
such that $S_{t_{j+1}}=S_{t_{j}}+\xi_{t_{j}}$ with $\xi_{t_{j}}\in\{-\delta,0,\delta\}$,
then for $f:\delta^{2}\mathbb{N}\times\delta\mathbb{Z}\to\mathbb{R}$
we have 
\begin{align}\label{eq:discrete-Ito-formula}
&f_{t_{j+1}}(S_{t_{j+1}})-f_{t_{j}}(S_{t_{j}})=\frac{1}{2\delta}\left(f_{t_{j+1}}(S_{t_{j}}+\delta)-f_{t_{j+1}}(S_{t_{j}}-\delta)\right)(S_{t_{j+1}}-S_{t_{j}})\\
&+\frac{1}{2\delta^{2}}\left(f_{t_{j+1}}(S_{t_{j}}+1)+f_{t_{j+1}}(S_{t_{j}}-1)-2f_{t_{j+1}}(S_{t_{j}})\right)(S_{t_{j+1}}-S_{t_{j}})^{2}+f_{t_{j+1}}(S_{t_{j}})-f_{t_{j}}(S_{t_{j}}). \notag
\end{align}
Moreover, having in mind the sequences $\bar{\zeta}^{\delta}$ and
$\zeta^{x,\delta}$ introduced in (\ref{eq:def-zetas-bar}) and (\ref{eq:def-zetas})
respectively, $X^{\delta}$ admits the following representation: 
\[
X_{t_{j+1}}^{\delta}\stackrel{(d)}{=}X_{t_{j}}^{\delta}+\delta\cdot\zeta_{j+1}^{x,\delta},
\]
Since the difference $X_{t_{j+1}}^{\delta}-X_{t_{j}}^{\delta}$ takes
values in $\{-\delta,0,\delta\}$, one can apply (\ref{eq:discrete-Ito-formula})
in order to get 
\begin{align}
 & f_{t_{j+1}}(X_{t_{j+1}}^{\delta})-f_{t_{j}}(X_{t_{j}}^{\delta})=f_{t_{j+1}}(X_{t_{j}}^{\delta})-f_{t_{j}}(X_{t_{j}}^{\delta})\nonumber \\
 & \hspace{0.5in}+\frac{1}{2\delta}\left(f_{t_{j+1}}(X_{t_{j}}^{\delta}+\delta)-f_{t_{j+1}}(X_{t_{j}}^{\delta}-\delta)\right)
 \left( \delta\cdot\zeta_{j+1}^{x,\delta} \right)\nonumber \\
 & \hspace{0.5in}+\frac{1}{2\delta^{2}}\left(f_{t_{j+1}}(X_{t_{j}}^{\delta}+\delta)-2f_{t_{j+1}}(X_{t_{j}}^{\delta})+f_{t_{j+1}}(X_{t_{j}}^{\delta}-\delta)\right) \left(\delta\cdot\zeta_{j+1}^{x,\delta}\right)^{2}.
 \label{b1}
\end{align}
We now evaluate $T^{\delta}f_{t_{j+1}}(X_{t_{j}}^{\delta})$, for
which we start from the simple decomposition 
\begin{eqnarray*}
T^{\delta}f_{t_{j+1}}(X_{t_{j}}^{\delta}) & = & \mathbb{E}_{\omega}\left[f_{t_{j+1}}(X_{t_{j+1}}^{\delta})|\mathcal{F}_{t_{j}}\right]\\
 & = & \mathbb{E}_{\omega}\left[f_{t_{j+1}}(X_{t_{j+1}}^{\delta})-f_{t_{j}}(X_{t_{j}}^{\delta})|\mathcal{F}_{t_{j}}\right]+f_{t_{j}}(X_{t_{j}}^{\delta}).
\end{eqnarray*}
Therefore, taking expected values in (\ref{b1}) we end up with 
\begin{align}
 & T^{\delta}f_{t_{j+1}}(X_{t_{j}}^{\delta})=\frac{1}{2}\left(f_{t_{j+1}}(X_{t_{j}}^{\delta}+\delta)-f_{t_{j+1}}(X_{t_{j}}^{\delta}-\delta)\right)\,\mathbb{E}_{\omega}\left[\zeta_{j+1}^{x,\delta}|\mathcal{F}_{t_{j}}\right]\nonumber \\
 & \hspace{0.5in}+\frac{1}{2}\left[f_{t_{j+1}}(X_{t_{j}}^{\delta}+\delta)-2f_{t_{j+1}}(X_{t_{j}}^{\delta})+f_{t_{j+1}}(X_{t_{j}}^{\delta}-\delta)\right]\,\mathbb{E}_{\omega}\left[\left(\zeta_{j+1}^{x,\delta}\right)^{2}|\mathcal{F}_{t_{j}}\right]\nonumber \\
 & \hspace{0.5in}+f_{t_{j+1}}(X_{t_{j}}^{\delta})-f_{t_{j}}(X_{t_{j}}^{\delta})+f_{t_{j}}(X_{t_{j}}^{\delta}).\label{b2}
\end{align}
Recall from (\ref{eq:martingale-discrete}) that $Z_{t_{j+1}}^{\delta}=f_{t_{j+1}}(X_{t_{j+1}}^{\delta})-T^{\delta}f_{t_{j+1}}(X_{t_{j}}^{\delta})$.
Hence subtracting (\ref{b2}) from relation (\ref{b1}), we get 
\begin{align}\label{b3}
Z_{t_{j+1}}^{\delta} & =\frac{1}{2}\left(f_{t_{j+1}}(X_{t_{j}}^{\delta}+\delta)-f_{t_{j+1}}(X_{t_{j}}^{\delta}-\delta)\right)\cdot\left(\zeta_{j+1}^{x,\delta}-\mathbb{E}\left[\zeta_{j+1}^{x,\delta}|\mathcal{F}_{t_{j}}\right]\right) \\
 &  +\frac{1}{2}\left(f_{t_{j+1}}(X_{t_{j}}^{\delta}+\delta)-2f_{t_{j+1}}(X_{t_{j}}^{\delta})+f_{t_{j+1}}(X_{t_{j}}^{\delta}-\delta)\right)\cdot\left(\left(\zeta_{j+1}^{x,\delta}\right)^{2}-\mathbb{E}\left[\left(\zeta_{j+1}^{x,\delta}\right)^{2}|\mathcal{F}_{t_{j}}\right]\right).
 \nonumber
\end{align}
In order to conclude, we first observe that 
\begin{equation}\label{eq:ZetSq}
(\zeta_{j+1}^{x,\delta})^{2}=(\bar{\zeta}_{j+1}^{\delta})^{2}.
\end{equation}Indeed, by definition we have 
\begin{align*}
\left(\bar{\zeta}_{j+1}^{\delta}\right)^{2} & ={\bf 1}_{\{U_{j+1}>\varepsilon,V_{j+1}\leqslant1/2\}}+{\bf 1}_{\{U_{j+1}>\varepsilon,V_{j+1}>1/2\}}\\
 & \ \ \ -2{\bf 1}_{\{U_{j+1}>\varepsilon,V_{j+1}\leqslant1/2\}}{\bf 1}_{\{U_{j+1}>\varepsilon,V_{j+1}>1/2\}}={\bf 1}_{\{U_{j+1}>\varepsilon\}}.
\end{align*}The same relation holds for $(\zeta^{x,\delta}_{j+1})^2$. Therefore, (\ref{eq:ZetSq}) holds. By adding and subtracting $\bar{\zeta}_{j+1}^{\delta}-\mathbb{E}_{\omega}[\bar{\zeta}_{j+1}^{\delta} |\mathcal{F}_{t_{j}}]$
to the right hand side of (\ref{b3}), and using the relation \eqref{eq:ZetSq}, we easily obtain our claim (\ref{eq:Zn-as-stoch-integral}).
%\hre{Can we justify this last relation?} {\color{blue}Xi: Updated as above.}
\end{proof}
\begin{rem}
On the right hand side of (\ref{eq:Zn-as-stoch-integral}), we will
see in the next sections that the limit of the sum (over $j$) of
the first term yields a continuous martingale (\ref{eq:continuous-martingale})
(which will be a stochastic integral with respect to a Brownian motion
$B$), while the limit of the sum of the second and third terms is
zero. 
\end{rem}

\subsection{\label{subsec:DiscPDE}Mild formulation of a discrete PDE}

Equation (\ref{eq:renormalize-time-space-discrete-martingale-problem})
involves the operator $\nabla^{\delta}-\mathcal{L}^{\delta}$. Related
to this fact, for a function $g\in L^{\infty}(\llbracket0,T\rrbracket\times\delta\mathbb{Z})$
one would like to give a meaning and solve the following forward discrete
partial differential equation, 
\begin{equation}
\nabla_{t}^{\delta}f_{t_{j}}(x)-\mathcal{L}_{x}^{\delta}f_{t_{j}}(x)=g_{t_{j}}(x).\label{eq: discrete-PDE}
\end{equation}
where we recall that $\nabla_{t}^{\delta}$ is given by (\ref{eq:def-nabla-discrete-renormalized}).
As a preliminary step, we introduce the associated symmetric random
walk and the discrete heat kernel as follows.
\begin{defn}
\label{def:renormalized-semigroups} Recall that $\bar{\mathcal{L}}^{\delta}$
was the operator defined for $x\in\delta\mathbb{Z}$ by \eqref{eq:renormalized-discrete-laplace}.
We consider the renormalized lazy symmetric random walk $Y^{\delta}$,
whose generator is $\bar{\mathcal{L}}^{\delta}$. We will denote the
corresponding renormalized discrete heat kernel by $p^{\delta}$.
Due to the homogeneity of $Y^{\delta},$ the kernel $p^{\delta}$
is defined on $\llbracket0,T\rrbracket\times\delta\mathbb{Z}$ as
the transition function of $Y^{\delta}$, namely 
\begin{equation}
p_{s}^{\delta}(x)\triangleq{\bf P}\left(Y_{t+s}^{\delta}=x|\,Y_{t}^{\delta}=0\right).\label{eq:discrete-heat-kernel}
\end{equation}
The corresponding transition operator on $L^{2}(\delta\mathbb{Z})$
is denoted by $P^{\delta}$. 
\end{defn}

\begin{rem}
\label{rmk:non-rescaled-discrete}In a non-rescaled situation, corresponding
to $\delta=1$ in Definition \ref{def:renormalized-semigroups}, we
shall use a superscript $d$ as in Section \ref{sec:prelim-sinai}.
Therefore we shall consider objects of the form $Y^{d}$, $p^{d}$,
$P^{d},\bar{\mathcal{L}}^{d}$ or $\hat{\nabla}_{x}^{d}$. It should
be noticed that $p^{\delta}$ can be expressed in terms of $p^{d}$
in the following way:
\begin{equation}\label{eq:pd-and-pdelta}
p_{s}^{\delta}(x)=p_{s/\delta^{2}}^{d}(x/\delta),\ \ \ (s,x)\in\delta^{2}\mathbb{Z}\times\delta\mathbb{Z}.
\end{equation}
\end{rem}

We would now like to express the solution of equation (\ref{eq: discrete-PDE})
in terms of a suitable fixed point problem which allows us to compare
with the continuum limit. To this aim, we derive the mild formulation
of our discrete partial differential equation below.
\begin{prop}
\label{prop:discrete-mild-pde} Let $\mathcal{L}^{\delta}$ be the
operator defined by (\ref{eq:transition-matrix-rescaled}) and consider
$g\in L^{\infty}(\llbracket0,T\rrbracket\times\delta\mathbb{Z})$.
Recall that the discrete heat kernel $p^{\delta}$ is introduced in
(\ref{eq:discrete-heat-kernel}). Then the mild form of equation (\ref{eq: discrete-PDE})
can be written for $x\in\delta\mathbb{Z}$ as 
\begin{equation}
f_{t_{j}}(x)=G_{t_{j}}(x)+J_{t_{j}}(x),\label{eq:discrete-mild-pde}
\end{equation}
where the function $G$ is given for $x\in\delta\mathbb{Z}$ by: 
\[
G_{t_{j}}(x)=\sum_{y\in\delta\mathbb{Z}}p_{t_{j}}^{\delta}(x-y)f_{0}(y)+\delta^{2}\sum_{\ell=0}^{j-1}\sum_{y\in\delta\mathbb{Z}}p_{t_{j-1}-t_{\ell}}^{\delta}(x-y)g_{t_{\ell}}(y).
\]
In equation (\ref{eq:discrete-mild-pde}) we also have that $J_{t_{j}}(x)$
can be expressed as below (recall that the twisted gradient $\hat{\nabla}_{x}^{\delta}$
is defined by relation (\ref{eq:renormalized-discrete-nabla})): 
\begin{equation}
J_{t_{j}}(x)=\delta\sum_{\ell=0}^{j-1}\sum_{y\in\delta\mathbb{Z}}p_{t_{j-1}-t_{\ell}}^{\delta}(x-y)\dot{U}^{\delta}(y)\hat{\nabla}_{x}^{\delta}f_{t_{\ell}}(y),%
\label{eq:def-J}
\end{equation}
where we recall that $\dot{U}^{\delta}(y)$ is given by \eqref{eq:def-U-delta}.
\end{prop}

\begin{proof}
We divide the proof into several steps. For the sake of clarity, we
will first derive the mild formulation in a non-rescaled case (namely
$\delta=1$). Specifically, recalling the notation $\dot{U}$ introduced
in (\ref{eq:def-dot-U}) as well as Proposition \ref{prop:gen-rw-discrete}
and Remark \ref{rmk:non-rescaled-discrete}, we can recast equation~(\ref{eq: discrete-PDE}) in the non-rescaled case as 
\begin{equation}
\nabla_{k}f_{k}(x)=\bar{\mathcal{L}}_{x}^{d}f_{k}(x)+\dot{U}(x)\,\hat{\nabla}_{x}f_{k}(x)+g_{k}(x),\label{eq:DiscPDE}
\end{equation}
where $f_{0}$ is given as an initial condition. Note that existence
and uniqueness for equation~(\ref{eq:DiscPDE}) is trivial. In fact,
by the definition of $\nabla_{k},$ one can easily build up the solution
recursively from the initial data $f_{0}.$ We start by deriving the
mild formulation for (\ref{eq:DiscPDE}).

\vspace{2mm}\noindent \textit{Step 1: Mild formulation for the free
equation.} Let us first consider the standard homogeneous free equation:
\begin{equation}
\nabla_{k}^{d}F_{k}(x)=\bar{\mathcal{L}}_{x}^{d}F_{k}(x),\quad\text{with initial condition}\quad F_{0}(x)=f_{0}(x).\label{d1}
\end{equation}
Recalling our Remark \ref{rmk:non-rescaled-discrete} about the symmetric
random walk $Y^{d}$ on $\mathbb{Z}$ and its transition operator
$P^{d}$, the solution to (\ref{d1}) is given by 
\begin{equation}
F_{k}^{H}(x)=P_{k}^{d}f_{0}(x)=\mathbb{E}_{\omega}[f_{0}(Y_{k}^{d})|Y_{0}^{d}=x].\label{d2}
\end{equation}
Although relation (\ref{d2}) is very classical, we now proceed to
its verification due to similar manipulations to be performed later
in the proof. To this aim, let $F^{H}$ be the space-time function
defined by the right hand side of (\ref{d2}). According to the Markov
property, we have 
\begin{align*}
F_{k+1}^{H}(x) & =\sum_{y\in\mathbb{Z}}p_{1}^{d}(x-y)\mathbb{E}_{\omega}\left[f_{0}(Y_{k}^{d})|Y_{0}^{d}=y\right]\\
 & =\varepsilon F_{k}^{H}(x)+\frac{1-\varepsilon}{2}\left(F_{k}^{H}(x+1)+F_{k}^{H}(x-1)\right).
\end{align*}
Therefore, 
\begin{align*}
\nabla_{k}F_{k}^{H}(x) & =F_{k+1}(x)-F_{k}(x)=\frac{\sigma^{2}}{2}\Delta_{x}^{d}F_{k}^{H}(x)=\bar{\mathcal{L}}_{x}^{d}F_{k}^{H}(x).
\end{align*}
\textit{Step 2: Mild formulation of an inhomogeneous }\textsc{pde}\textit{.}
Next we consider a function $g$ defined on $\mathbb{N}\times\mathbb{Z}$
and the standard inhomogeneous equation 
\begin{equation}
\begin{cases}
\nabla_{k}F_{k}(x)=\bar{\mathcal{L}}_{x}^{d}F_{k}(x)+g_{k}(x),\\
F_{0}(x)=f_{0}(x).
\end{cases}\label{eq:InHFree}
\end{equation}
One can derive an expression for $F$ in terms of $F^{H}$ thanks
to Duhamel's principle. Namely we claim that $F$ can be written as
\begin{equation}
F_{k}(x)=F_{k}^{H}(x)+\sum_{j=0}^{k-1}\mathbb{E}_{\omega}[g_{j}(Y_{k-1-j}^{d})|Y_{0}^{d}=x],\label{d3}
\end{equation}
where we have used the convention $\sum_{j=0}^{-1}\triangleq0$. Indeed,
for $F$ defined by (\ref{d3}) we have 
\begin{equation}
\nabla_{k}F_{k}(x)=\nabla_{k}F_{k}^{H}(x)+g_{k}(x)+\sum_{j=0}^{k-1}\mathbb{E}_{\omega}[g_{j}(Y_{k-j}^{d})-g_{j}(Y_{k-1-j}^{d})|Y_{0}^{d}=x].\label{d31}
\end{equation}
Hence recalling the notation for the transition operator $P^{d}$
in Remark \ref{rmk:non-rescaled-discrete} as well as equation~(\ref{d1})
for $F^{H}$, we get 
\begin{equation}
\nabla_{k}F_{k}(x)=\bar{\mathcal{L}}_{x}^{d}F_{k}^{H}(x)+g_{k}(x)+\sum_{j=0}^{k-1}\left(P^{d}-{\rm Id}\right)P_{k-j-1}^{d}g_{j}(x)\label{d32}
\end{equation}
Now thanks to the fact that $P^{d}-{\rm Id}=\bar{\mathcal{L}}^{d}$
and resorting to the definition of $F$ in the right hand side of
(\ref{d3}), we end up with 
\begin{align}
\nabla_{k}F_{k}(x) & =\bar{\mathcal{L}}_{x}^{d}F_{k}^{H}(x)+g_{k}(x)+\sum_{j=0}^{k-1}\bar{\mathcal{L}}_{x}^{d}\mathbb{E}_{\omega}[g_{j}(Y_{k-1-j}^{d})|Y_{0}^{d}=x]\nonumber \\
 & =\bar{\mathcal{L}}_{x}^{d}F_{k}(x)+g_{k}(x).\label{d33}
\end{align}
Our claim (\ref{d3}) is now achieved.

\vspace{2mm}\noindent \textit{Step 3: Mild formulation for the PDE
in a random environment.} Finally, we consider the original equation
(\ref{eq:DiscPDE}), whose solution is called $f$. Also recall that
$F$ designates the solution to (\ref{eq:InHFree}), and we set 
\begin{equation}
\Phi_{k}(x)\triangleq F_{k}(x)-f_{k}(x).\label{d4}
\end{equation}
Subtracting the right hand side of (\ref{eq:InHFree}) from the right
hand side of (\ref{eq:DiscPDE}), it follows that $\Phi_{0}(\cdot)=0$
and 
\[
\nabla_{k}\Phi=\nabla_{k}F-\nabla_{k}f=\bar{\mathcal{L}}_{x}^{d}\Phi-\dot{U}(x)\cdot\hat{\nabla}_{x}f.
\]
The above equation is of the form (\ref{eq:InHFree}) with $g=-\dot{U}\cdot\hat{\nabla}_{x}f$
and $f_{0}=0$. Therefore, applying (\ref{d3}) with $F^{H}=0$ we
get the following expression for $\Phi$: 
\begin{equation}
\Phi_{k}(x)=-\sum_{j=0}^{k-1}\mathbb{E}_{\omega}[\dot{U}(Y_{k-1-j}^{d})\cdot\hat{\nabla}_{x}f_{j}(Y_{k-1-j}^{d})|Y_{0}^{d}=x].\label{d5}
\end{equation}
We now gather relations (\ref{d3}), (\ref{d4}) and (\ref{d5}) in
order to get the following expression for the solution $f$ to (\ref{eq:DiscPDE}),
\begin{align}
f_{k}(x) & =\mathbb{E}_{\omega}[f_{0}(Y_{k}^{d})|Y_{0}^{d}=x]+\sum_{j=0}^{k-1}\mathbb{E}_{\omega}[g_{j}(Y_{k-1-j}^{d})|Y_{0}^{d}=x]\nonumber \\
 & \ \ \ +\sum_{j=0}^{k-1}\mathbb{E}_{\omega}[\dot{U}(Y_{k-1-j}^{d})\cdot\hat{\nabla}_{x}f_{j}(Y_{k-1-j}^{d})|Y_{0}^{d}=x].\label{eq:Solc=00003D1}
\end{align}
Notice that in (\ref{eq:Solc=00003D1}) we have chosen to write the
formula in terms of the symmetric random walk $Y^{d}$. Recalling
Definition \ref{def:renormalized-semigroups} and Remark \ref{rmk:non-rescaled-discrete}
for the symmetric random walk kernel $p^{d}$, one can recast (\ref{eq:Solc=00003D1})
as
\begin{eqnarray}\label{d51}
f_{k}(x)  =\sum_{y\in\mathbb{Z}}p_{k}^{d}(x-y)f_{0}(y)
&+&\sum_{j=0}^{k-1}\sum_{y\in\mathbb{Z}}p_{k-1-j}^{d}(x-y)g_{j}(y) \notag\\
&+&\sum_{j=0}^{k-1}\sum_{y\in\mathbb{Z}}p_{k-1-j}^{d}(x-y)\dot{U}(y)\hat{\nabla}_{x}f_{j}(y).
\end{eqnarray}
Observe that since $p^{d}$ is finitely supported, the sums over $y\in\mathbb{Z}$
in (\ref{d51}) are in fact finite sums. Otherwise stated, we have
proved (\ref{eq:discrete-mild-pde}) for $\delta=1$.

\vspace{2mm}\noindent \textit{Step 4: Rescaling the equation.} Now
let us move to the rescaled picture and write down the solution to
the corresponding discrete \textsc{pde}. We will also illustrate the fact that Definition~\ref{def:EnvScale} is the suitable rescaling for the random environment.
As in Definition \ref{def:EnvScale} and Definition~\ref{def:partition-gradient},
we use $t_{k}$ to denote a generic point in $\delta^{2}\mathbb{N}$
and $x$ to denote a generic point in $\delta\mathbb{Z}$. We also
write $\mathcal{L}^{\delta}$, $\bar{\mathcal{L}}^{\delta}$, ${\hat\nabla}_{x}^{\delta}$  
and $\nabla_{t}^{\delta}$ 
for the rescaled operators respectively given by (\ref{eq:transition-matrix-rescaled}), (\ref{eq:renormalized-discrete-laplace}),  (\ref{eq:renormalized-discrete-nabla}) and \eqref{eq:def-nabla-discrete-renormalized}. Recall that we obtained the relation \eqref{d61} connecting  $\mathcal{L}^{\delta}$, $\bar{\mathcal{L}}^{\delta}$ and ${\hat\nabla}_{x}^{\delta}$ . 

With relation \eqref{d61} in hand and given an initial condition $f_{0}$,
we are now interested in the following rescaled \textsc{pde}: 
\begin{equation}
\nabla_{t_{k}}^{\delta}f_{t_{k}}(x)=\mathcal{L}_{x}^{\delta}f_{t_{k}}(x)+g_{t_{k}}(x).\label{eq:ScalDiscPDE}
\end{equation}
Let us start by mimicking Step 1 and Step 2 in the rescaled picture.
First the solution to the homogeneous free equation is now given by
\[
F_{t_{k}}^{H,\delta}(x)=\mathbb{E}_{\omega}[f_{0}(Y_{t_{k}}^{\delta})|Y_{0}^{\delta}=x],
\]
where $Y^{\delta}$ is the rescaled random walk introduced in Definition~\ref{def:renormalized-semigroups}.
In addition, the solution to the standard inhomogeneous equation corresponding
to~(\ref{eq:ScalDiscPDE}) needs to be adjusted as 
\begin{equation}
F_{t_{k}}^{\delta}(x)=F_{t_{k}}^{H,\delta}(x)+\delta^{2}\sum_{j=0}^{k-1}\mathbb{E}_{\omega}[g_{t_{j}}(Y_{t_{k-1}-t_{j}}^{\delta})|Y_{0}^{\delta}=x].\label{eq:InhomRes}
\end{equation}
With respect to (\ref{d3}), the appearance of the factor $\delta^{2}$
in (\ref{eq:InhomRes}) is easily checked. Indeed, let $F_{t_{k}}^{\delta}(x)$
be defined by equation~(\ref{eq:InhomRes}). Then we have 
\begin{multline}
F_{t_{k+1}}^{\delta}(x)-F_{t_{k}}^{\delta}(x)=F_{t_{k+1}}^{H,\delta}(x)-F_{t_{k}}^{H,\delta}(x)+\delta^{2}g_{t_{k}}(x)\\
+\delta^{2}\sum_{j=0}^{k-1}\left(\mathbb{E}_{\omega}[g_{t_{j}}(Y_{t_{k}-t_{j}}^{\delta})|Y_{0}^{\delta}=x]-\mathbb{E}_{\omega}[g_{t_{j}}(Y_{t_{k-1}-t_{j}}^{\delta})|Y_{0}^{\delta}=x]\right).\label{d7}
\end{multline}
Moreover, similarly to (\ref{d32}) and (\ref{d33}), we obtain 
\begin{equation}
\mathbb{E}_{\omega}[g_{t_{j}}(Y_{t_{k}-t_{j}}^{\delta})|Y_{0}^{\delta}=x]-\mathbb{E}_{\omega}[g_{t_{j}}(Y_{t_{k-1}-t_{j}}^{\delta})|Y_{0}^{\delta}=x]=\delta^{2}\bar{\mathcal{L}}_{x}^{\delta}\mathbb{E}_{\omega}[g_{t_{j}}(Y_{t_{k-1}-t_{j}}^{\delta})|Y_{0}^{\delta}=x].\label{d8}
\end{equation}
Therefore, plugging (\ref{d8}) into (\ref{d7}) we have 
\[
F_{t_{k+1}}^{\delta}(x)-F_{t_{k}}^{\delta}(x)=F_{t_{k+1}}^{H,\delta}(x)-F_{t_{k}}^{H,\delta}(x)+\delta^{2}g_{t_{k}}(x)+\delta^{4}\sum_{j=0}^{k-1}\bar{\mathcal{L}}_{x}^{\delta}\mathbb{E}_{\omega}[g_{t_{j}}(Y_{t_{k-1}-t_{j}}^{\delta})|Y_{0}^{\delta}=x],
\]
Since $\nabla_{t}^{\delta}$ is defined by (\ref{eq:def-nabla-discrete-renormalized}),
we thus get 
\[
\nabla_{t}^{\delta}F_{t_{k}}^{\delta}(x)=\nabla_{t}^{\delta}F_{t_{k}}^{H,\delta}(x)+g_{t_{k}}(x)+\delta^{2}\sum_{j=0}^{k-1}\bar{\mathcal{L}}_{x}^{\delta}\mathbb{E}_{\omega}[g_{t_{j}}(Y_{t_{k-1}-t_{j}}^{\delta})|Y_{0}^{\delta}=x].
\]
On the other hand, if $F^{\delta}$ is given by (\ref{eq:InhomRes})
we also have 
\[
\bar{\mathcal{L}}_{x}^{\delta}F_{t_{k}}^{\delta}(x)=\bar{\mathcal{L}}_{x}^{\delta}F_{t_{k}}^{H,\delta}(x)+\delta^{2}\sum_{j=0}^{k-1}\bar{\mathcal{L}}_{x}^{\delta}\mathbb{E}_{\omega}[g_{t_{j}}(Y_{t_{k-1}-t_{j}}^{\delta})|Y_{0}^{\delta}=x].
\]
Consequently, similarly to (\ref{d33}) we have proved that $F^{\delta}$
defined by~(\ref{eq:InhomRes}) solves the non homogeneous equation
\[
\nabla_{t}^{\delta}F_{t_{k}}^{\delta}(x)=\bar{\mathcal{L}}_{x}^{\delta}F_{t_{k}}^{\delta}(x)+g_{t_{k}}(x).
\]
Starting from the above relation, we let the patient reader check
(along the same lines going from (\ref{d3}) to (\ref{d51})) that
the solution $f$ to (\ref{eq: discrete-PDE}) is expressed in the
mild form by 
\begin{eqnarray*}
f_{t_{k}}(x)=\sum_{y\in\delta\mathbb{Z}}p_{t_{k}}^{\delta}(x-y)f_{0}(y)
&+&\delta^{2}\sum_{j=0}^{k-1}\sum_{y\in\delta\mathbb{Z}}p_{t_{k-1}-t_{j}}^{\delta}(x-y)g_{t_{j}}(y)\\
&+&\delta\sum_{j=0}^{k-1}\sum_{y\in\delta\mathbb{Z}}p_{t_{k-1}-t_{j}}^{\delta}(x-y)\dot{U}^{\delta}(y)\hat{\nabla}_{x}^{\delta}f_{t_{j}}(y).
\end{eqnarray*}
We have thus shown relation (\ref{eq:discrete-mild-pde}), which finishes
the proof.
\end{proof}
In order to connect the mild equation (\ref{eq:discrete-mild-pde})
with its continuous counterpart, it will be beneficial to proceed
to an integration by parts procedure. This is summarized in the next
proposition.
\begin{prop}
\label{prop:discrete-mild-pde-2} Under the setting and conditions
of Proposition \ref{prop:discrete-mild-pde}, let $J$ be the term
defined by (\ref{eq:def-J}). For an arbitrary $a\in\delta\mathbb{Z}$
and $y\in\delta\mathbb{Z}$ also set 
\begin{equation}
\mathcal{I}_{t}^{\delta}(a,y)\triangleq\sum_{z\in\llbracket a+2\delta,y\rrbracket}\dot{U}^{\delta}(z)\,\hat{\nabla}_{x}^{\delta}f_{t}(z),\label{eq:def-Ical}
\end{equation}
where we recall that $\dot{U}^{\delta}$ is defined by (\ref{eq:def-U-delta}).
Then the following holds true: 

\vspace{2mm}\noindent (i) An alternative way to write $J_{t_{j}}(x)$
is given by 
\begin{equation}
J_{t_{j}}(x)=\delta^{2}\sum_{\ell=0}^{j-1}\sum_{y\in\delta\mathbb{Z}}\nabla_{x}^{\delta}p_{t_{j-1}-t_{\ell}}^{\delta}(x-\delta-y)\,\mathcal{I}_{t_{\ell}}^{\delta}(a,y).\label{eq:Jt-with-ibp}
\end{equation}
(ii) Let $f$ be the solution of the mild equation (\ref{eq:discrete-mild-pde}).
The equation governing the derivative $\nabla_{x}^{\delta}f_{t_{j}}(x)$
can be written for $x\in\delta\mathbb{Z}$ as 
\begin{equation}
\nabla_{x}^{\delta}f_{t_{j}}(x)=\nabla_{x}^{\delta}G_{t_{j}}(x)+\delta^{2}\sum_{\ell=0}^{j-1}\sum_{y\in\delta\mathbb{Z}}\nabla_{x}^{2,\delta}p_{t_{j-1}-t_{\ell}}^{\delta}(x-y)\,\mathcal{I}_{t_{\ell}}^{\delta}(a,y),\label{eq:nabla-ft-discrete}
\end{equation}
where for $\varphi\in L^{\infty}(\delta\mathbb{Z})$ we define
\[
\nabla_{x}^{2,\delta}\varphi(x)\triangleq\frac{1}{\delta^{2}}\left(\varphi(x+\delta)+\varphi(x-\delta)-2\varphi(x)\right) \, .
\]
(iii) We also obtain an equation for the twisted derivative $\hat{\nabla}_{x}^{\delta}f$
of $f$, namely: 
\begin{equation}
\hat{\nabla}_{x}^{\delta}f_{t_{j}}(x)=\hat{\nabla}_{x}^{\delta}G_{t_{j}}(x)+\delta^{2}\sum_{\ell=0}^{j-1}\sum_{y\in\delta\mathbb{Z}}\hat{\nabla}_{x}^{\delta}\nabla_{x}^{\delta}p_{t_{j-1}-t_{\ell}}^{\delta}(x-\delta-y)\mathcal{I}_{t_{\ell}}^{\delta}(a,y).\label{eq:TwiDf}
\end{equation}
\end{prop}

\begin{proof}
We start from the expression (\ref{eq:def-J}) for $J$ and perform
a discrete integration by parts. To this aim, we first notice from
(\ref{eq:def-Ical}) that 
\[
\mathcal{I}_{t_{\ell}}^{\delta}(a,y)-\mathcal{I}_{t_{\ell}}^{\delta}(a,y-\delta)=\dot{U}^{\delta}(y)\hat{\nabla}_{x}^{\delta}f_{t_{\ell}}(y).
\]
Hence using the fact that $p_{t}^{\delta}(x,\cdot)$ is finitely supported,
we can rewrite the second sum in~(\ref{eq:def-J}) as 
\begin{eqnarray*}
Q_{j}^{\delta} & = & \sum_{y\in\delta\mathbb{Z}}p_{t_{j-1}-t_{\ell}}^{\delta}(x-y)\,\left(\mathcal{I}_{t_{\ell}}^{\delta}(a,y)-\mathcal{I}_{t_{\ell}}^{\delta}(a,y-\delta)\right)\\
 & = & \sum_{y\in\delta\mathbb{Z}}p_{t_{j-1}-t_{\ell}}^{\delta}(x-y)\mathcal{I}_{t_{\ell}}^{\delta}(a,y)-\sum_{y\in\delta\mathbb{Z}}p_{t_{j-1}-t_{\ell}}^{\delta}(x-y)\mathcal{I}_{t_{\ell}}^{\delta}(a,y-\delta).
\end{eqnarray*}
Therefore resorting to an elementary change of variables in the right
hand side above, we end up with 
\begin{eqnarray*}
Q_{j}^{\delta} & = &\sum_{y\in\delta\mathbb{Z}}\left(p_{t_{j-1}-t_{\ell}}^{\delta}(x-y)-p_{t_{j-1}-t_{\ell}}^{\delta}(x-\delta-y)\right)\,\mathcal{I}_{t_{\ell}}^{\delta}(a,y)\\ 
& = &\delta\sum_{y\in\delta\mathbb{Z}}\nabla_{x}^{\delta}p_{t_{j-1}-t_{\ell}}^{\delta}(x-\delta-y)\,\mathcal{I}_{t_{\ell}}^{\delta}(a,y).
\end{eqnarray*}
Plugging this identity into~(\ref{eq:def-J}) we have proved our
claim (\ref{eq:Jt-with-ibp}). Relations (\ref{eq:nabla-ft-discrete})
and (\ref{eq:TwiDf}) are then easily obtained by taking discrete
derivatives on both sides of (\ref{eq:Jt-with-ibp}), which ends the
proof of Proposition \ref{prop:discrete-mild-pde-2}.
\end{proof}

\subsection{\label{sec:heuristics}Heuristic considerations and rescaling}

With Proposition \ref{prop:discrete-mild-pde-2} in hand, let us briefly
outline the coupling strategy at the heart of our considerations.
Specifically, one can gather Propositions \ref{prop:discrete-mild-pde}
and \ref{prop:discrete-mild-pde-2} and write the mild form of equation
(\ref{eq: discrete-PDE}) as 
\begin{equation}
f_{t_{j}}(x)=G_{t_{j}}^{1}(x)+G_{t_{j}}^{2}(x)+J_{t_{j}}(x),\label{eq:heuristics-a}
\end{equation}
where 
\[
G_{t_{j}}^{1}(x)=\sum_{y\in\delta\mathbb{Z}}p_{t_{j}}^{\delta}(x-y)f_{0}(y),\qquad G_{t_{j}}^{2}(x)=\delta^{2}\sum_{\ell=0}^{j-1}\sum_{y\in\delta\mathbb{Z}}p_{t_{j-1}-t_{\ell}}^{\delta}(x-y)g_{t_{\ell}}(y),
\]
and where the term $J_{t_{j}}(x)$ is defined by (\ref{eq:Jt-with-ibp}).
Let us now figure out heuristically how relation (\ref{eq:heuristics-a})
will converge to a continuous limit. First of all, the local central
limit theorem (see e.g. \cite{LL10}) implies formally that if $(x_{\delta})\in\delta\mathbb{Z}$
converges to $x\in\mathbb{R}$ then we have 
\begin{equation}
\lim_{\delta\to0}\frac{1}{\delta}p_{t}^{\delta}(x_{\delta})=p_{t}(x),\label{eq:heuristics-b}
\end{equation}
where $p^{\delta}$ is introduced in (\ref{eq:discrete-heat-kernel}).
Also note that in (\ref{eq:heuristics-b}), $p$ designates the heat
kernel on $\mathbb{R}$ (with variance $\sigma^{2}\triangleq1-\varepsilon$)
given by 
\begin{equation}
p_{t}(x)=\frac{1}{\sqrt{2\pi\sigma^{2}t}}e^{-\frac{x^{2}}{2\sigma^{2}t}}.\label{eq:heuristics-c}
\end{equation}
Now let us recast the term $G^{1}$ in (\ref{eq:heuristics-a}) as
\begin{equation}
G_{t_{j}}^{1}(x)=\delta\sum_{y\in\delta\mathbb{Z}}\frac{1}{\delta}p_{t_{j}}^{\delta}(x-y)f_{0}(y).\label{eq:heuristics-d}
\end{equation}
Using relation (\ref{eq:heuristics-b}) and convergence of Riemann
sum considerations, it is easily conceived that 
\begin{equation}
\lim_{\delta\to0}G_{t_{j}}^{1}(x)=\int_{\mathbb{R}}p_{t}(x-y)f_{0}(y)dy.\label{eq:heuristics-e}
\end{equation}
In the same way we expect that 
\begin{equation}
\lim_{\delta\to0}G_{t_{j}}^{2}(x)=\int_{0}^{t}\int_{\mathbb{R}}p_{t-s}(x-y)g_{s}(y)dy\,ds.\label{eq:heuristics-f}
\end{equation}
In order to derive the heuristics about the limiting behavior of the
term $J$ in~(\ref{eq:heuristics-a}), let us analyse the sum $\mathcal{I}_{t}^{\delta}(x,y)$
defined by (\ref{eq:def-Ical}). This will also yield an explanation
for our renormalisation of the environment $\omega$. Indeed, the
term $\hat{\nabla}_{x}^{\delta}f_{t}(z)$ in the definition of $\mathcal{I}_{t}^{\delta}(x,y)$
is expected to converge to $\partial_{x}f_{t}(z)$. Thus we also expect
to have (whenever $t_{j}\to t$) 
\begin{equation}
\lim_{\delta\to0}\mathcal{I}_{t_{j}}^{\delta}(a,y)=\lim_{\delta\to0}\sum_{z\in\llbracket a+2\delta,y\rrbracket}\dot{U}^{\delta}(z)\hat{\nabla}_{x}^{\delta}f_{t_{j}}(z)=\int_{a}^{y}\partial_{z}f_{s}(z)dW(z),\label{eq:heuristics-g}
\end{equation}
where $W(z)$ is a Brownian motion on $\mathbb{R}$ with suitable variance $\tau^2$. 
Now if we want (\ref{eq:heuristics-g}) to hold, this imposes that
${\rm Var}(\dot{U}^{\delta}(z))$ is of order $\delta$ for all $z\in\delta\mathbb{Z}$.
According to Remark \ref{rk:variance}, we are led to renormalization of the environment $\omega$
given by Definition \ref{def:EnvScale}. 

Let us summarize our computations so far. If we gather (\ref{eq:heuristics-e}),
(\ref{eq:heuristics-f}) and (\ref{eq:heuristics-g}), we expect the
discrete equation (\ref{eq:heuristics-a}) to converge to the \textsc{PDE}
given in (\ref{eq:forme-mild-edp-ip1}) below. Our aim is to justify
this assertion in the following sections, and provide a rate of convergence
in the limiting procedure thanks to a coupling method.

\subsection{\label{sec:heatkerest}Estimates for the discrete heat kernel}

In Section \ref{subsec:DiscPDE} we have defined the simple random walk $Y^{\delta}$ related to the 
generator $\bar{\mathcal{L}}^{\delta}$ in Definition \ref{def:renormalized-semigroups}. In this section we review 
some notation about the transition kernels for $Y^{\delta}$ and state some useful Gaussian 
type bounds. 

Recall that $Y^{d}=\{Y_{n}^{d}:n\in\mathbb{N}\}$ is the symmetric
random walk with one-step transition density given by 
\[
\mathbb{P}(Y_{n+1}^{d}=l|Y_{n}^{d}=k)=\begin{cases}
\varepsilon, & \text{if }l=k;\\
\frac{1-\varepsilon}{2}, & \text{if }l=k\pm1;\\
0, & \text{otherwise}.
\end{cases}
\]
Noticing that $Y^{d}$ is space time-homogeneous, the $n$-step transition
function of $Y^{d}$ is denoted as 
\begin{equation}
p_{n}^{d}(k)\triangleq\mathbb{P}(Y_{n}^{d}=k|Y_{0}^{d}=0).\label{eq:SRWKernel}
\end{equation}
One can also view $Y_{n}^{d}=Y_{0}^{d}+\xi_{1}+\cdots+\xi_{n}$, where
$\{\xi_{i}:i\geqslant1\}$ is an i.i.d. sequence with distribution
\[
\mathbb{P}(\xi_{1}=0)=\varepsilon,
\quad\text{and}\quad 
\mathbb{P}(\xi_{1}=1)=\mathbb{P}(\xi_{1}=-1)=\frac{1-\varepsilon}{2}.
\]
Note that the variance of $\xi_{1}$ is ${\rm Var}[\xi_{1}]=1-\varepsilon=\sigma^{2}.$
Given $\delta>0,$ we shall consider the following rescaled version
of the kernel $p^{d}$ on the space-time grid $\delta^{2}\mathbb{N}\times\delta\mathbb{Z}$,
defined by 
\begin{equation}\label{eq:hat-p-delta}
\hat{p}_{t}^{\delta}(x)
\triangleq 
\frac{1}{\delta}
p_{t/\delta^{2}}^{d}(x/\delta)
=
\frac{1}{\delta} p_{t}^{\delta}(x),
\quad\text{for}\quad 
(t,x)\in\delta^{2}\mathbb{N}\times\delta\mathbb{Z} \, ,
\end{equation}
where $p^{\delta}$ is introduced in \eqref{eq:pd-and-pdelta}.
We also recall that $p_{t}(x)$ denotes the continuous heat kernel
with variance $\sigma^{2}=1-\varepsilon:$
\begin{equation}
p_{t}(x)\triangleq\frac{1}{\sqrt{2\pi\sigma^{2}t}}e^{-\frac{x^{2}}{2\sigma^{2}t}},\ \ \ t>0,x\in\mathbb{R}.\label{eq:GauKernelSigma}
\end{equation}
Let us summarize some notation about discrete spatial gradients on the grid $\delta\Z$ for $\delta>0$. First we have set, for $x\in\delta\Z$,
\begin{equation}\label{d9}
\nabla_{x}^{\delta} f(x)
=\frac{1}{\delta} \lp  f(x+\delta) -f(x) \rp,
\quad\text{and}\quad
\hat{\nabla}_{x}^{\delta} f(x)
=\frac{1}{2\delta} \lp  f(x+\delta) -f(x-\delta) \rp .
\end{equation}
For higher order rescaled discrete gradients, we have used the following conventions:  
%\hb{(Samy: I have replaced $\widetilde{\nabla}_{x}^{4,\delta}$ by ${\nabla}_{x^{4,\delta}$ below)}:
\begin{eqnarray}
\nabla_{x}^{2,\delta} f(x)
&=&
\frac{1}{\delta^{2}} \left(  f(x+\delta) + f(x-\delta) -2 f(x) \right) \label{d10}\\
\widetilde{\nabla}_{x}^{3,\delta}f(x)&=&
\nabla_{x}^{\delta} \nabla_{x}^{2,\delta} f(x)
=
\frac{1}{\delta} \left(  \nabla_{x}^{2,\delta}f(x+\delta) -\nabla_{x}^{2,\delta}f(x) \right) \label{d11}\\
%\widetilde
{\nabla}_{x}^{4,\delta}f(x)&=&
\frac{1}{\delta^{4}}\left(  f(x+2\delta)-4f(x+\delta)+
6f(x)-4f(x-\delta)+f(x-2\delta)\right) \label{d12}
\end{eqnarray}
Also recall from \eqref{eq:def-nabla-discrete-renormalized} that the rescaled discrete time gradient $\nabla_{t}^{\delta}$ is given on $\delta^{2}\mathbb{N}\times\delta\mathbb{Z}$ by 
\[
\nabla_{t}^{\delta}f_{t_{j}}(x)\triangleq\frac{1}{\delta^{2}}(f_{t_{j+1}}(x)-f_{t_{j}}(x)).
\]
With this series of notation in hand, let us highlight the fact that the discrete heat kernel $\hat{p}^{\delta}$ displayed in \eqref{eq:hat-p-delta} satisfies the discrete heat equation 
\begin{equation}\label{eq:discrete-heat-equation}
\nabla_{t}^{\delta}\hat{p}^{\delta}_{t_{j}}(x)
=\frac{\sigma^{2}}{2}
\nabla^{2,\delta}\hat{p}^{\delta}_{t_{j}}(x)
=
\bar{\mathcal{L}}_{x}^{\delta}\hat{p}^{\delta}_{t_{j}}(x),
\end{equation}
where the operator $\bar{\mathcal{L}}_{x}^{\delta}$ is introduced in \eqref{eq:renormalized-discrete-laplace}.
We now state a quantitative local central limit theorem (CLT) as well as a uniform bound on the discrete kernel $\hat{p}^{\delta}$.

\begin{thm}[cf. \cite{LL10}]\label{thm:LCLT} Recall that $p_{n}^{d}(k)$
is the discrete kernel of $Y^{d}$ defined by (\ref{eq:SRWKernel}),
and $p_{t}(x)$ is the Gaussian kernel with variance $\sigma^{2}=1-\varepsilon$
defined by (\ref{eq:GauKernelSigma}). Then for $m=2,4$ there exists a constant $C_{m,\varepsilon}>0$ depending only
on $m$ and $\varepsilon,$ such that 
\begin{equation}
\big|\nabla_{x}^{m}p_{n}^{d}(x)-\nabla_{x}^{m}p_{n}(x)\big|\leqslant\frac{C_{m,\varepsilon}}{n^{(m+3)/2}} ,
\quad\text{for all}\quad n\geqslant1\;\text{ and }\;k\in\mathbb{Z}.\label{eq:LCLT}
\end{equation}

\end{thm}

The following uniform Gaussian upper bound for the discrete kernel $\hat{p}^\delta$ is classical. The proof for the cases when $m=0,1$ is essentially contained in \cite[Theorem 5.1]{HS93}. However, to our best knowledge the case for higher derivatives does not seem to be easily available in the literature and the extension of the argument in \cite{HS93} is not entirely straightforward. We thus provide an independent proof based on the local CLT and the Markov property, which works for any order of spatial derivatives $\nabla^m$.

\begin{prop}
\label{prop:DiscUnifGau} For $\delta>0$, let $\hat{p}^{\delta}$ the kernel defined by \eqref{eq:hat-p-delta}.  
Consider the gradient $\nabla_{x}^{2,\delta}$ given by~\eqref{d10}, as well as the gradient $\nabla_{x}^{4,\delta}$ introduced in~\eqref{d12}. Then for $m=2,4$ there exist
two universal constant $C_{1},C_{2}>0$ depending only on $m$ and the lazy parameter $\varepsilon$,
such that 
\begin{equation}
\big|\nabla^{m,\delta}\hat{p}_{t}^{\delta}(u)\big|\leqslant\frac{C_{1}}{t^{(m+1)/2}}e^{-C_{2}u^{2}/t},\label{eq:DiscUnifGau>0}
\end{equation}
for all $(t,u)\in\delta^{2}\mathbb{N}\times\delta\mathbb{Z}$ with
$t>0$. For $t=0$, the bound \eqref{eq:DiscUnifGau>0} becomes 
\begin{equation}\label{eq:DiscUnifGau=0}
\big|\nabla^{m,\delta}\hat{p}_{0}^{\delta}(u)\big|
\leqslant\frac{C_{2}}{\delta^{m+1}} \, \mathds{1}_{\{|u|\leqslant\delta\}}.
\end{equation}
\end{prop}

\begin{proof}
For $t=0$ we have $\hat{p}_{0}^{\delta}(x)=\delta^{-1}\mathds{1}_{\{x=0\}}$. Hence the bound \eqref{eq:DiscUnifGau=0} derives immediately from the definitions \eqref{d10}-\eqref{d12}. In what follows we thus focus on proving \eqref{eq:DiscUnifGau>0}. 

\noindent
\emph{Step 1: Strategy.} In order to prove \eqref{eq:DiscUnifGau>0}, fix $m=2$ for sake of clarity (the case $m=4$ is treated very similarly). Writing explicitly the definition \eqref{d10} of $\nabla^{m,\delta}$ and resorting to the expression \eqref{eq:hat-p-delta} for $\hat{p}_{t}^{\delta}$, we let the reader check that \eqref{eq:DiscUnifGau>0} is equivalent to the non rescaled version 
\begin{equation}
\big|\nabla^{m}p_{n}^{d}(k)\big|\leqslant\frac{a}{n^{(m+1)/2}}e^{-bk^{2}/n},\quad\text{for all}\quad n \geqslant1\;\text{ and }\;k\in\mathbb{Z} ,
\label{eq:DiscGaussEst}
\end{equation}
where $\nabla^{m}$ is defined by \eqref{d10} with $\delta=1$ and $p^{d}$ is defined in \eqref{eq:SRWKernel}. In \eqref{eq:DiscGaussEst}, the numbers 
$a,b$ are universal constants depending only on $m$ and $\varepsilon$.
The choice of $a,b$ will be clear in the course of the proof. We
are going to prove (\ref{eq:DiscGaussEst}) by induction on $n$.
The main idea is that the quantitative local CLT (cf. Theorem \ref{thm:LCLT})
easily yields (\ref{eq:DiscGaussEst}) for the regime $\{k:|k|^{2}\lesssim n\}$.
The other regime is then handled by induction and the Markov property.
In what follows, the notation $C_{\text{subscript}}$ denotes a constant
depending only on the parameters specified in the subscript whose
value may change from line to line.

\noindent
\emph{Step 2: Case $|k^{2}|\leqslant\Lambda n$.} To begin with, by applying the triangle inequality to (\ref{eq:LCLT}) we obtain
that 
\begin{equation}
\big|\nabla^{m}p_{n}^{d}(k)\big|\leqslant\big|\nabla^{m}p_{n}(k)\big|+\frac{C_{m,\varepsilon}}{n^{(m+3)/2}}\quad\text{for all}\quad n \geqslant1\;\text{ and }\;k\in\mathbb{Z}.\label{eq:DiscGPf1}
\end{equation}
By the explicit expression \eqref{eq:GauKernelSigma} of $p_{n}(k)$, it is easy to show that
\begin{equation}
\big|\nabla^{m}p_{n}(k)\big|\leqslant\frac{C_{m,\varepsilon}}{n^{(m+1)/2}}e^{-\frac{k^{2}}{2\sigma^{2}n}}\quad\text{for all}\quad n \geqslant1\;\text{ and }\;k\in\mathbb{Z}.\label{eq:DiscGPf2}
\end{equation}
Now let $\Lambda$ be a positive universal constant whose value will
be specified later on. According to (\ref{eq:DiscGPf1}) and (\ref{eq:DiscGPf2}),
we arrive at the following estimate:
\begin{equation}
\big|\nabla^{m}p_{n}^{d}(k)\big|\leqslant\frac{C_{m,\varepsilon,\Lambda}}{n^{(m+1)/2}}e^{-\frac{k^{2}}{2\sigma^{2}n}}\quad\text{ for all }\; n,k\;\text{ such that }\;|k^{2}|\leqslant\Lambda n.\label{eq:DiscGPf3}
\end{equation}
This proves \eqref{eq:DiscGaussEst} for $|k^{2}|\leqslant\Lambda n$.

\noindent
\emph{Step 3: Inductive procedure.}
Let $C_{m,\varepsilon,\Lambda}$ and $\si^{2}$ be the constants in~\eqref{eq:DiscGPf3}.
Next, we are going to prove (\ref{eq:DiscGaussEst}) by induction
on $n$ with suitably chosen constants 
\begin{equation}
a>C_{m,\varepsilon,\Lambda},\quad\text{ and }\quad\ b<\frac{1}{2\sigma^{2}}.\label{eq:DiscGPf4}
\end{equation}
Note that if the estimate (\ref{eq:DiscGaussEst}) is valid, it remains
true for larger $a$ and smaller $b$. We first fix $a_{0},b_{0}$
to be such that (\ref{eq:DiscGaussEst}) holds (with $a=a_{0},b=b_{0}$)
when $n=1$ for all $k\in\mathbb{Z}$. Since $\nabla^{m}p_{1}^{d}$
is finitely supported, the existence of such numbers is obvious. We
will define 
\[
a\triangleq a_{0}\vee C_{m,\varepsilon,\Lambda},
\]
where $C_{m,\varepsilon,\Lambda}$ is the constant appearing in (\ref{eq:DiscGPf3})
and the choices of $b,\Lambda$ will be clear later on. Now suppose
that the estimate (\ref{eq:DiscGaussEst}) holds for fixed $n$ and
all $k\in\mathbb{Z}$. To establish the induction step, we only need
to consider the case when $|k|^{2}>\Lambda(n+1)$ (the other case
is proved in (\ref{eq:DiscGPf3}) with the presumed constraint (\ref{eq:DiscGPf4})
for $a,b$). By applying $\nabla^{m}$ to the Markov property, we
see that

\[
\nabla^{m}p_{n+1}(k)=\frac{1-\varepsilon}{2}\left(\nabla^{m}p_{n}(k-1)+\nabla^{m}p_{n}(k+1)\right)+\varepsilon\nabla^{m}p_{n}(k).
\]
According to the induction hypothesis, we have 
\begin{align*}
\big|\nabla^{m}p_{n+1}(k)\big| & \leqslant\frac{1-\varepsilon}{2}\left(\frac{a}{n^{(m+1)/2}}e^{-\frac{b(k-1)^{2}}{n}}+\frac{a}{n^{(m+1)/2}}e^{-\frac{b(k+1)^{2}}{n}}\right)+\varepsilon\cdot\frac{a}{n^{(m+1)/2}}e^{-\frac{bk^{2}}{n}}\\
 & =\frac{a}{(n+1)^{(m+1)/2}}\cdot\left(1+\frac{1}{n}\right)^{\frac{m+1}{2}}\left(\frac{1-\varepsilon}{2}\left(e^{-\frac{b(k-1)^{2}}{n}}+e^{-\frac{b(k+1)^{2}}{n}}\right)+\varepsilon e^{-\frac{bk^{2}}{n}}\right).
\end{align*}
Multiplying both sides above by $e^{b k^{2}/n+1}$ and performing some elementary manipulations on the exponential functions, it follows that 
\[
\big|\nabla^{m}p_{n+1}(k)\big|\cdot e^{\frac{bk^{2}}{n+1}}\leqslant\frac{a}{(n+1)^{(m+1)/2}}\cdot\mathcal D_{n} \, ,
\]
where we have set
\begin{equation}\label{eq:matcal D}
\mathcal D_{n}=\left(1+\frac{1}{n}\right)^{\frac{m+1}{2}}\cdot e^{-\frac{bk^{2}}{n(n+1)}}\left((1-\varepsilon)e^{-\frac{b}{n}}\cosh\left(\frac{2bk}{n}\right)+\varepsilon\right).
\end{equation}
In order to complete the induction step, the crucial point is thus to prove 
that one can choose $b$ and $\Lambda$ so that the
factor $\mathcal D_{n}$ satisfies
\begin{equation}\label{eq:mathcal D<=1}
\mathcal D_{n}\leqslant 1,\quad\text{ for all }\; |k|^{2}>\Lambda(n+1).
\end{equation}
We now proceed to prove \eqref{eq:mathcal D<=1}.

\noindent
\emph{Step 4: Bounding $\mathcal D_{n}$.} In order to prove \eqref{eq:mathcal D<=1}, it is important to note that we are working
with the regime $|k|^{2}>\Lambda(n+1)$. Another basic observation
is that $\nabla^{m}p_{n+1}(k)$ is only non-zero when $|k|\leqslant C_{m}n$
(since the one-step transition probabilities are finitely supported).
As a result, the effective $k$-region under consideration is 
\begin{equation}
\Lambda(n+1)\leqslant k^{2}\leqslant C_{m}^{2}n^{2}.\label{eq:DiscGPf9}
\end{equation}
We also recall the following elementary inequalities to be used later
on:
\begin{equation}
\left(1+\frac{1}{n}\right)^{\frac{m+1}{2}}\leqslant1+\frac{C'_{m}}{n}\ \ \ \forall n\geqslant1;\quad  e^{-x}\leqslant1-\frac{1}{2}x;\quad\cosh x\leqslant1+x^{2}\ \ \ \forall x\in[0,\eta] , \label{eq:DiscGPf6}
\end{equation}
where $\eta$ is some universal constant.
Now we choose $b$ to satisfy
\begin{equation}
b<b_{0}\wedge\frac{1}{2\sigma^{2}}\wedge\frac{\eta}{C_{m}^{2}\wedge2C_{m}}\wedge\frac{1}{16},\label{eq:DiscGPf5}
\end{equation}
where $C_{m}$ is the constant appearing in (\ref{eq:DiscGPf9}) and $a_0, b_0$ are defined after \eqref{eq:DiscGPf4}.
Then taking relations~\eqref{eq:DiscGPf4} and~\eqref{eq:DiscGPf5} into account we know that 
\[
\frac{bk^{2}}{n(n+1)}\vee\frac{2b|k|}{n}<\eta,\quad\text{ for all }\;|k|\leqslant C_{m}n.
\]
Therefore, the inequalities
in (\ref{eq:DiscGPf6}) imply that 
\begin{eqnarray}
  \mathcal D_{n}
 &\leqslant&\left(1+\frac{1}{n}\right)^{\frac{m+1}{2}}\cdot e^{-\frac{bk^{2}}{n(n+1)}}\cosh\left(\frac{2bk}{n}\right)\nonumber \\
 & \leqslant&\left(1+\frac{C'_{m}}{n}\right)\left(1-\frac{1}{2}\frac{bk^{2}}{n(n+1)}\right)\left(1+\frac{4b^{2}k^{2}}{n^{2}}\right) \, ,\label{eq:DiscGPf7}
\end{eqnarray}
whenever $k$ satisfies $|k|\leqslant C_{m}n.$ To analyze the last
expression, we first note that 
\begin{equation}\label{eq:ineg-elem-1}
\left(1-\frac{1}{2}\frac{bk^{2}}{n(n+1)}\right)\left(1+\frac{4b^{2}k^{2}}{n^{2}}\right)\leqslant1-b\left(\frac{1}{2}\frac{k^{2}}{n(n+1)}-4b\frac{k^{2}}{n^{2}}\right).
\end{equation}
Since $b<1/16$, it is elementary to see that 
\begin{equation}\label{eq:ineg-elem-2}
4b\frac{k^{2}}{n^{2}}<\frac{1}{4}\frac{k^{2}}{n(n+1)}.
\end{equation}
Plugging \eqref{eq:ineg-elem-2} back to \eqref{eq:ineg-elem-1} and then into \eqref{eq:DiscGPf7}, we end up with 
\begin{equation}\label{eq:mathcalD-first-upper}
\mathcal D_{n}\leqslant\left(1+\frac{C_{m}'}{n}\right)\left(1-\frac{1}{4}\frac{bk^{2}}{n(n+1)}\right).
\end{equation}
In addition, we have $k^{2}>\Lambda(n+1)$. Therefore we can further bound the right hand side of~\eqref{eq:mathcalD-first-upper} as 
\[
\mathcal D_{n}\leqslant\left(1+\frac{C_{m}'}{n}\right)\left(1-\frac{b\Lambda}{4n}\right)\leqslant1-\left(\frac{b\Lambda}{4}-C_{m}'\right)\frac{1}{n}.
\]
It is now clear that we have 
\begin{equation}
\mathcal D_{n}\leqslant 1,\quad\text{ as long as }\quad\frac{b\Lambda}{4}>C_{m}'.\label{eq:DiscGPf8}
\end{equation}
Consequently, with such choices of $b,\Lambda$ we conclude that our claim 
\eqref{eq:mathcal D<=1} holds true for $n+1$ if we assume that it holds for $n\geqslant 1$.

\noindent
\emph{Step 4: Conclusion.} The bound \eqref{eq:mathcal D<=1} we have just proved allows to complete the induction step for \eqref{eq:DiscUnifGau>0}. To summarize the above procedure, we highlight the following facts: 

\noindent
(i) We first choose $a_{0},b_{0}$ to
obtain the initial step $n=1$. Next, we choose $b$ to satisfy (\ref{eq:DiscGPf5}).
Then we choose $\Lambda$ to satisfy (\ref{eq:DiscGPf8}). Finally,
we choose $a=a_{0}\vee C_{m,\varepsilon,\Lambda}.$ 

\noindent
(ii) With the above constants fixed, for $k^{2}\leqslant\Lambda(n+1)$ our desired inequality \eqref{eq:DiscUnifGau>0} is ensured by \eqref{eq:DiscGPf2}.

\noindent
(iii) For $k^{2}>\Lambda(n+1)$, our induction procedure yields 
\eqref{eq:DiscUnifGau>0}.

\noindent
The proof of the proposition is thus complete.
\end{proof}

\begin{rem}
Proposition \ref{prop:DiscUnifGau} and Theorem \ref{thm:LCLT}
of course hold for other type of discrete derivatives (e.g. forward
differences) for all $m\in\mathbb{N}$. Here we have chosen a formulation
that is directly applicable to our situation.
\end{rem}

\section{Brownian motion in Brownian environment}\label{BMBE}

In this section, we recall some basic facts about the continuous analogue
of the random walk $X^{d}$ given by Definition (\ref{def:SinRWZ}).
This classical object is called \textit{Brox diffusion} (or Brownian
motion with Brownian potential), denoted here by $X^{c}$ (in the
sequel, the superscript $c$ stands for \emph{continuous time parameter}).
We will recall some important definitions concerning $X^{c}$ in Section
\ref{sec:def-brox}. Then we will introduce a martingale problem related to $X^{c}$ in Section \ref{sec:mild-pde-brox}. Eventually we introduce some rough paths metric which allows to get a pathwise meaning to the martingale problem in 
Section \ref{sec:rough-path-struct}. Throughout this section we will use the following notation. 
\begin{notation}\label{nabla-delta-continuous}
For a smooth enough function $f:
\mathbb R\to\mathbb R$, we set 
\[\nabla^{c}f\equiv\partial_{x}f,\quad\text{ and }\quad\Delta^{c}f\equiv\partial^{2}_{xx}f.\]
\end{notation}

\subsection{Definition of the Brox diffusion}\label{sec:def-brox} 

Let $W$ be a one-dimensional, two-sided Brownian motion defined on
a probability space $(\Omega,\mathcal{G},{\bf P})$ with suitable variance $\tau^2$ (the exact value of $\tau^2$ is specified by (\ref{eq:WVar}) below). The process $X^{c}$
can be seen as the solution to the following formal stochastic differential
equation:
\begin{equation}
dX_{t}^{c}=-\frac{1}{2}\dot{W}(X_{t}^{c})\,dt+dB_{t},\label{eq:bm-in-br-environment}
\end{equation}
where $B$ is a one-dimensional standard Brownian motion independent
of $W$. To be consistent with the discrete Definition \ref{def:SinRWZ}, we will assume that $B$ is defined on a probability space $(\hat{\Omega}, \mathcal F,\mathbb P)$ which is independent of $(\Omega, \mathcal G, \mathbf P)$.  For our later purpose of comparison with the discrete case,
throughout the rest we always assume that $B$ has variance $\sigma^{2}\triangleq1-\varepsilon$
(i.e. with generator $\frac{\sigma^{2}}{2}\Delta^{c}\equiv\frac{\sigma^{2}}{2}\partial^{2}_{xx}$), where $\varepsilon$
is the parameter introduced in Definition \ref{def:environment}.

Since the drift $\dot{W}$ in (\ref{eq:bm-in-br-environment}) is
a distribution, this equation does not admit a strong solution. Therefore
a more standard way to introduce the process $X^{c}$ is to define
it as a Feller diffusion with the following generator: 
\begin{equation}
\mathcal{L}^{c}f(x)=\frac{\sigma^{2}}{2}e^{W(x)/\sigma^{2}}\partial_{x}\left(e^{-W(x)/\sigma^{2}}\partial_{x}f\right)(x)=\frac{\sigma^{2}}{2}\Delta^{c}f(x)-\frac{1}{2}\dot{W}(x)\nabla^{c}f(x),\label{eq:generator-continuous}
\end{equation}
where we observe that the second expression is still formal. Also
observe that $X^{c}$ can be seen as a Markov process whose Dirichlet
form on ${\rm L}^{2}(\mathbb R,\mu)$, with $\mu(dx)=e^{-W(x)/\sigma^{2}}dx$,   is given by 
\begin{equation}\label{eq:feler-gene-cont}
\mathcal{E}^{c}(f)=\si^{2}\int_{\mathbb{R}}e^{-W(x)/\sigma^{2}}|\partial_{x}f(x)|^{2}\,dx,
\end{equation}
and we notice that expressions like \eqref{eq:feler-gene-cont} now make sense even if $W$ is not differentiable. 

The definition of $X^{c}$ as a Markov process is a well established
fact. Below we summarize some classical results in this direction,
which can be found e.g in \cite{GH} (also see references therein).
\begin{prop}
\label{prop:very-weak-brox} Let $W$ be a one-dimensional two-sided
Brownian motion defined on a probability space $(\Omega,\mathcal{G},{\bf P})$.
We consider the operator $\mathcal{L}^{c}$ defined by (\ref{eq:generator-continuous}).
Then we have:

\vspace{2mm}\noindent (i) The domain of $\mathcal{L}^{c}$ is dense
in the space $C_{0}(\mathbb{R})$ of continuous functions vanishing
at infinity.\\
(ii) On a probability space $(\hat{\Omega},\mathcal{F},\mathbb{P})$,
one can construct a Markov process $X^{c}$ whose generator is given
by $\mathcal{L}^{c}$.
\end{prop}

While Proposition \ref{prop:very-weak-brox} is certainly a substantial progress in the understanding of equation~\eqref{eq:bm-in-br-environment}, its main shortcoming is that it yields a very weak notion of solution. A considerable amount of effort has been devoted in the recent past to a more pathwise definition of $X^{c}$. The first important contribution in this direction is \cite{HLM}, which hinges on the It\^o-McKean representation (see \cite{Br, Sh}) of $X^{c}$, as well as a thorough analysis of local times. The second main work on pathwise interpretations of \eqref{eq:bm-in-br-environment} can be found in \cite{DD}. It relies on a pathwise interpretation of the martingale problem related to \eqref{eq:bm-in-br-environment} thanks to rough paths techniques. In the current contribution we will
mostly stick to the setting of \cite{DD}, since it might be easier
to generalize to higher dimensional contexts. As in the introduction, one should also mention the articles \cite{BC,CG}, which show that equations like \eqref{eq:bm-in-br-environment} admit a strong solution even in cases where the drift 
$\dot{W}$ is a distribution. However \cite{BC,CG} fall short of handling the case of a drift $\dot{W}\in {\rm C}^{-1/2-\varepsilon}$ like a white noise on $\mathbb R$. 

In conclusion, our interpretation of~\eqref{eq:bm-in-br-environment} will rely on the martingale problem and related pathwise rough PDEs developed in \cite{DD}. We now proceed to give a heuristic derivation of the martingale problem framework.

\subsection{\label{sec:mild-pde-brox} Heuristics about the \textsc{PDE} problem
related to Brox diffusion}

In order to understand the nature of the family of mild \textsc{PDE}s related
to equation (\ref{eq:bm-in-br-environment}), let us first consider
the following smoothened version of the white noise $\dot{W}$ defined
for $\eta>0$: 
\begin{equation}
\dot{W}^{\eta}=\dot{W}*p_{\eta},\label{e1}
\end{equation}
where $p_{t}$ denotes the continuous heat kernel in $\mathbb{R}$,
defined as in \eqref{eq:heuristics-c} by 
\begin{equation}\label{heat_kernel}
p_{t}(x)=\frac{1}{\sqrt{2\pi\sigma^{2}t}}e^{-\frac{x^{2}}{2\sigma^{2}t}} .
\end{equation}
Also denote by $P_{t}$ the heat semigroup associated to the generator
$\frac{\sigma^{2}}{2}\Delta^{c}\equiv\frac{\sigma^{2}}{2}\partial_{xx}^{2}$. In order to alleviate notation,
we will still write $\mathcal{L}^{c}$ for the operator defined by
(\ref{eq:generator-continuous}) with $\dot{W}$ replaced by $\dot{W}^{\eta}$.
Namely we set 
\begin{equation}\label{eq:reset-infgene}
\mathcal{L}^{c}f(x)=\mathcal{L}^{c,\eta}f(x)=\frac{\sigma^{2}}{2}\Delta^{c}f(x)-\frac{1}{2}\dot{W}^{\eta}(x)\nabla^{c}f(x),
\end{equation}
where we recall that we have set $\nabla^{c}f\equiv\partial_{x}f$.

The martingale problem for equation \eqref{eq:bm-in-br-environment} relies on a family of PDEs. Namely for a generic  function $g\in{\rm C}(\mathbb{R}_{+}\times\mathbb{R})$, we consider the solution $f$ of the following equation: 
\begin{equation}
\partial_{t}f_{t}(x)-\mathcal{L}_{x}^{c}f_{t}(x)=g_{t}(x),\quad t\in[0,\tau],\,x\in\mathbb{R}.\label{eq:continous-PDE}
\end{equation}
Since $\dot{W}^{\eta}$ in \eqref{e1} is smooth, equation \eqref{eq:continous-PDE} can be solved in a strong sense. However, in order to take limits as $\eta\to 0$, we will first write \eqref{eq:continous-PDE} in a mild form which echoes Proposition~\ref{prop:discrete-mild-pde}. This is the content of the following proposition.
\begin{prop}
\label{prop:CtsPDE} Fix $\eta>0$ and recall that $\dot{W}^{\eta}$
is defined by~(\ref{e1}). Then if $g\in{\rm C}(\mathbb{R}_{+}\times\mathbb{R})$,
the mild form of equation (\ref{eq:continous-PDE}) can be written
as 
\begin{equation}
f_{t}(x)=P_{t}f_{0}(x)+\int_{0}^{t}\int_{\mathbb{R}}p_{t-s}(x-y)g_{s}(y)dyds+\int_{0}^{t}\int_{\mathbb{R}}p_{t-s}(x-y)\dot{W}^{\eta}(y)\nabla^{c}f_{s}(y)dyds.\label{eq:forme-mild-edp}
\end{equation}
Furthermore, an alternative way to write equation~(\ref{eq:forme-mild-edp})
is 
\begin{multline}
f_{t}(x)=P_{t}f_{0}(x)
+\int_{0}^{t}ds\int_{\mathbb{R}}p_{t-s}(x-y)g_{s}(y)\,dy \\
-\frac{1}{2}\int_{0}^{t}ds\int_{\mathbb{R}}dy\,\partial_{x}p_{t-s}(x-y)\int_{a}^{y}\partial_{x}f_{s}(z)dW^{\eta}(z),
\label{eq:forme-mild-edp-ip1}
\end{multline}
where $a$ is an arbitrary real number. In addition, an equation is
also available for the derivative of $f_{t}$ with respect to the
space variable: 
\begin{multline}
\partial_{x}f_{t}(x)=\partial_{x}P_{t}f_{0}(x)+\int_{0}^{t}\int_{\mathbb{R}}\partial_{x}p_{t-s}(x-y)g_{s}(y)dyds\\
-\frac{1}{2}\int_{0}^{t}\int_{\mathbb{R}}\partial_{xx}^{2}p_{t-s}(x-y)\left(\int_{a}^{y}\partial_{z}f_{s}(z)dW^{\eta}(z)\right)dyds.\label{eq-forme-mild-edp-ip2}
\end{multline}
\end{prop}

\begin{proof}
We proceed as in the proof of Proposition \ref{prop:discrete-mild-pde}.
Namely start with the simple \textsc{PDE}
\[
\partial_{t} F=\frac{\sigma^{2}}{2}\Delta^{c}F+g,\qquad F_{0}=f_{0},
\]
whose solution can be written as 
\begin{equation}
F_{t}(x)=P_{t}f_{0}(x)+\int_{0}^{t}\int_{\mathbb{R}}p_{t-s}(x-y)g_{s}(y)dyds.\label{e2}
\end{equation}
Also set $J\triangleq F-f$, where $f$ solves (\ref{eq:continous-PDE}).
Then $J$ satisfies $J_{0}=0$ and 
\begin{align*}
\partial_{t}J & =\partial_{t} F-\partial_{t} f=\frac{\sigma^{2}}{2}\Delta^{c}F+g-\frac{\sigma^{2}}{2}\Delta^{c}f-\frac{1}{2}\dot{W}^{\eta}\,\nabla^{c}f_{t}-g=\frac{\sigma^{2}}{2}\Delta^{c}J-\frac{1}{2}\dot{W}^{\eta}\,\nabla^{c}f_{t}.
\end{align*}
Therefore the function $J$ can be written in mild form as 
\begin{equation}
J_{t}(x)=-\frac{1}{2}\int_{0}^{t}\int_{\mathbb{R}^{n}}p_{t-s}(x-y)\dot{W}^{\eta}(y)\,\nabla^{c}f_{s}(y)\,dyds.\label{e3}
\end{equation}
Substracting (\ref{e3}) from (\ref{e2}), it follows that $f=F-J$
satisfies equation (\ref{eq:forme-mild-edp-ip1}).

As for equation (\ref{eq:TwiDf}), we go from (\ref{eq:forme-mild-edp})
to (\ref{eq:forme-mild-edp-ip1}) thanks to an integration by parts
procedure. More specifically, due to the vanishing properties of the
heat kernel at infinity one can write 
\begin{multline}
\int_{\mathbb{R}}p_{t-s}(x-y)\dot{W}^{\eta}(y)\,\nabla^{c}f_{s}(y)dy=-\int_{\mathbb{R}}p_{t-s}(x-y)d\left(\int_{a}^{y}\partial_{z}f_{s}(z)dW^{\eta}(z)\right)\\
=-\int_{\mathbb{R}}\partial_{x}p_{t-s}(x-y)\left(\int_{a}^{y}\partial_{z}f_{s}(z)dW^{\eta}(z)\right)dy.\label{eq:integration-parts}
\end{multline}
Plugging this relation into \eqref{eq:forme-mild-edp}, our claim \eqref{eq:forme-mild-edp-ip1}
is easily proved.

Our last argument is as follows: exactly as in the proof of Proposition
\ref{prop:discrete-mild-pde}, relation (\ref{eq-forme-mild-edp-ip2})
is obtained from \eqref{eq:forme-mild-edp-ip1} by differentiating
with respect to $x$ on both sides of the relation. This finishes
our proof.
\end{proof}
\begin{rem}\label{rem:limiting-proc-eta}
As mentioned earlier, Proposition \ref{prop:CtsPDE} is stated for
a regularized version of the noise $\dot{W}$. The pathwise interpretation
of the martingale problem related to equation (\ref{eq:bm-in-br-environment})
can then be reduced to a limiting procedure in equation (\ref{eq-forme-mild-edp-ip2}). We give some hints about this procedure in the current section, the technical details being deferred to Section 
\ref{RPMP}.
\end{rem}

%%%%%%%%%%%%%%%%%%%%%%%%%
%%%%%%%%%%%%%%%%%%%%%%%%%
%%%%%%%%%%%%%%%%%%%%%%%%%

At the heart of the approach in \cite{DD} is the fact that one can
solve equation~(\ref{eq:continous-PDE}), or better said equation~(\ref{eq-forme-mild-edp-ip2}),
in a pathwise way. As mentioned in Remark \ref{rem:limiting-proc-eta}, this is achieved by taking limits in equation \eqref{eq:integration-parts}. We summarize
this result in the following theorem, which is stated here quite informally.
\begin{thm}
\label{thm:rough-pde-informal} Let $f_{0},g$ be two given $C_{b}^{2}$-functions,
and consider an arbitrary time horizon $T>0$. Then there exists
a unique function $f$ in a proper space of controlled process with
respect to $W$, satisfying the following equation in the rough paths
sense on $[0,T]\times\mathbb{R}$: 
\begin{multline}
\partial_{x}f_{t}(x)=\int_{\mathbb{R}}\partial_{x}p_{t}(x-y)f_{0}(y)dy+\int_{0}^{t}\int_{\mathbb{R}}\partial_{x}p_{t-s}(x-y)g_{s}(y)dyds\\
-\frac{1}{2}\int_{0}^{t}\int_{\mathbb{R}}\partial_{xx}^{2}p_{t-s}(x-y)\left(\int_{x}^{y}\partial_{z}f_{s}(z)dW(z)\right)dyds,\label{eq:FPP}
\end{multline}
\end{thm}

Our aim in this section is to recall the main setting allowing to
properly state and prove Theorem~\ref{thm:rough-pde-informal}. Before
getting into the computational details, let us recall that this theorem
has to be seen as the main building block in order to set up the martingale
problem for equation~(\ref{eq:bm-in-br-environment}). Namely our
ultimate goal is to get the result below.
\begin{thm}
\label{thm:rough-pde-informal-martingale} Let $f_{0},g$ be two given
$\mathcal{C}_{b}^{2}$-functions, and consider an arbitrary time horizon
$T>0$. Let $f$ be the unique solution to equation~(\ref{eq:FPP}).
Then there exists a probability $\mathbb{P}^{W}$ and a canonical
process $X^{c}$ on a filtered space $(\hat{\Omega},(\mathcal{F}_{t})_{t\in[0,T]})$
such that under $\mathbb{P}^{W}$ the process $M=M^{f}$ is a martingale,
where 
\begin{equation}
M_{t}=f_{t}(X_{t}^{c})-f_{0}(X_{0}^{c})-\int_{0}^{t}\left(\mathcal{L}_{x}^{c}f_{s}(X_{s}^{c})-\partial_{t}f_{s}(X_{s}^{c})\right)\,ds.\label{eq:continuous-time-space-martingale-problem}
\end{equation}
In addition, there exists a Brownian motion $B$ defined on 
$(\hat{\Omega},(\mathcal{F}_{t})_{t\in[0,T]})$ such that 
the martingale $M$ admits the representation 
\begin{equation}
M_{t}=\int_{0}^{t}\partial_{x}f_{s}(X_{s}^{c})\,dB_{s}.\label{eq:continuous-martingale}
\end{equation}
\end{thm}

Section  \ref{sec:rough-path-struct} below is devoted to specify the setting under which we will achieve Theorem~\ref{thm:rough-pde-informal} and Theorem \ref{thm:rough-pde-informal-martingale}.

\subsection{\label{sec:rough-path-struct}The rough path structure for the fixed point problem}

In order to get a better grasp on the rough path setting employed below, let us first vaguely summarize the main philosophy
invoked in \cite{DD} in order to get Theorem~\ref{thm:rough-pde-informal}:

\vspace{2mm}\noindent (i) We think of $\partial_{x}f_{t}(x)$ in~(\ref{eq:FPP})
as the unknown object, denoted as $v_{t}(x).$ For each fixed $t$,
the function $x\mapsto v_{t}(x)$ is regarded as a rough path controlled
by the Brownian motion $W.$ Therefore, the integral $\int_{x}^{y}v_{s}(z)dW(z)$
is well-defined in the rough path sense. Notice that the rough paths point of view is needed here. This is due to the fact that $W$ is a double sided Brownian motion, and therefore $z\mapsto v_{s}(z)$ cannot be thought of as an adapted process.

\noindent
(ii) We introduce an essential transformation 
\begin{equation}
\mathcal{M}:v\mapsto(\mathcal{M}v)_{t}(x)\triangleq\int_{0}^{t}\int_{\mathbb{R}}\partial_{xx}^{2}p_{t-s}(x-y)\left(\int_{x}^{y}v_{s}(z)dW(z)\right)dyds.\label{e4}
\end{equation}
This is a transformation on a suitable space of controlled rough paths. Equation \eqref{eq:FPP} can thus be written as $v=\mathcal{M}v$. Namely $v$ is a fixed point for the mapping $\mathcal{M}$. 

\noindent
(iii) In order to solve the fixed point problem, we must expect that
$\mathcal{M}$ is a contraction. The major technical challenge here
is to define the rough path metric in a delicate way, so that the
transformation $\mathcal{M}$ is indeed a contraction. This is a highly
non-trivial point, as there is no a priori evidence about why $\mathcal{M}$
needs to be a contraction at all (it is not surprising that $\mathcal{M}$
is a bounded linear transformation though).

\vspace{2mm} Notice that in the setting of \cite{DD},
an auxiliary component $Z_{t}(x)$ is introduced to form a two-dimensional
rough path $(W_{t}(x),Z_{t}(x))$. This is needed due to the time-inhomogeneity
of the diffusion equation therein. However, the generator in our situation
is time-homogeneous. This largely simplifies the rough path viewpoint.
In particular, it allows us to treat $W$ as a one-dimensional rough
path with the obvious lifting, and the analysis is also simplified
accordingly. This is why we have included a self contained version
of the computations in Section \ref{RPMP}. 
We now introduce the rough path setting which will be used in order to solve equation \eqref{eq:FPP}.

\subsubsection{Metric on the Brownian rough path $W$}

Our regularities will be quantified in terms of two parameters $1/3<\beta<\alpha<1/2$.
The parameter $\alpha$ will be the H\"older-exponent for $W,$ while
$\beta$ will be used for the H\"older-exponent of the solution path
$v$.

Recall that $x\mapsto W(x)$ is a two-sided Brownian motion. On
each compact interval $[-a,a],$ it can be viewed as a one-dimensional
$\alpha$-H\"older rough path in the obvious way, that is ${\bf W}=(W^{1},W^{2})$
with 
\begin{equation}
W^{1}(x,y)=W(y)-W(x),\quad\text{and}\quad W^{2}(x,y)=\frac{1}{2}\left(W(y)-W(x)\right)^{2}.\label{eq:def-w1-w2}
\end{equation}
It is important for our purposes to consider $W$ on the entire real
line. Thus we need to introduce a suitable metric on ${\bf W}$ to
take into account its growth at infinity.
\begin{lem}
\label{lem:GrowthWHol} Let $W^{1}$ be defined by (\ref{eq:def-w1-w2})
and recall that $\alpha\in(1/3,1/2)$. For each $a\geqslant1$, we
define a random variable $H_{a}$ by 
\begin{equation}
H_{a}\triangleq\|W^{1}\|_{\alpha}^{[-a,a]}\triangleq\sup_{-a\leqslant x\neq y\leqslant a}\frac{|W^{1}(x,y)|}{|y-x|^{\alpha}}.\label{eq:def-Ha}
\end{equation}
Then for any $\chi>1/2-\alpha$, we have $\sup_{a\geqslant1}\frac{H_{a}}{a^{\chi}}<\infty$
almost surely.
\end{lem}

\begin{proof}
The Brownian scaling property shows that $\frac{H_{a}}{a^{\chi}}\stackrel{\text{law}}{=}\frac{1}{a^{\chi-(1/2-\alpha)}}\cdot H_{1}.$
As a result, we have 
\[
\sum_{n=1}^{\infty}\mathbb{P}\left(\frac{H_{n+1}}{n^{\chi}}>\varepsilon\right)=\sum_{n=1}^{\infty}\mathbb{P}\left(H_{1}>\varepsilon(n+1)^{\chi-(1/2-\alpha)}\cdot\left(\frac{n}{n+1}\right)^{\chi}\right).
\]
The series is clearly convergent since $H_{1}$ has Gaussian tail
(according to Fernique's lemma). Therefore a standard application
of Borel-Cantelli shows that $\frac{H_{n+1}}{n^{\chi}}\rightarrow0$
almost surely, further implying that 
\[
\sup_{a\geqslant1}\frac{H_{a}}{a^{\chi}}\leqslant\sup_{n\geqslant1}\frac{H_{n+1}}{n^{\chi}}<\infty\ \text{a.s.}
\]
This finishes our proof.
\end{proof}
In view of Lemma \ref{lem:GrowthWHol}, we will first label our assumptions
on $\alpha,\beta$ and $\chi$: 
\begin{equation}
1/3<\beta<\alpha<1/2,\quad\text{and}\quad\frac{1}{2}-\alpha<\chi<\frac{\beta}{2}\label{f1}
\end{equation}
We let the reader check that the hypothesis on $\chi$ in (\ref{f1})
can be met whenever $1/3<\alpha<1/2$. Next, given $\alpha,\beta$
and $\chi$ satisfying (\ref{f1}), we define a random variable $\kappa_{\alpha,\chi}({\bf W})$
by 
\begin{equation}
\kappa_{\alpha,\chi}({\bf W})\triangleq\sup_{a\geqslant1}\left(\frac{\|W^{1}\|_{\alpha}^{[-a,a]}}{a^{\chi}}+\frac{\|W^{2}\|_{2\alpha}^{[-a,a]}}{a^{2\chi}}\right),\label{f11}
\end{equation}
where the norm $\|W^{2}\|_{2\alpha}^{[-a,a]}$ is defined similarly
to (\ref{eq:def-Ha}). The quantity $\kappa_{\alpha,\chi}({\bf W})$
is finite almost surely, as seen from Lemma \ref{lem:GrowthWHol}.

\subsubsection{The solution space and the corresponding rough path metric}

The definition of the solution space for equation \eqref{e4} and the corresponding metric
is much more involved. This is largely due to the need of obtaining
a contraction property for the transformation $\mathcal{M}$ defined
by (\ref{e4}).

We shall identify a Banach space $\mathcal{B}$ where the solution
path $v$ for the fixed point problem given by (\ref{eq:FPP}) or
(\ref{e4}) lives. Let us first describe what the object $v$ looks
like. Recall that, in the fixed point problem, we want to think of
the integral $\int_{x}^{y}v_{t}(z)dW(z)$ as a rough path integral.
This requires viewing $z\mapsto v_{t}(z)$ (for each fixed $t$) as
a rough path controlled by ${\bf W}$ (similar to \cite[Definition 1]{Gu}).
Therefore, the object $v_{t}(\cdot)$ must also come with a Gubinelli
derivative path $\partial_{W}v_{t}(\cdot)$. In other words, the a
priori shape of the object $v$ should be given by a pair 
\begin{equation}\label{f12}
\mathcal{V}=(v,\partial_{W}v),
\end{equation}
where 
\[
v:[0,T]\times\mathbb{R}\rightarrow\mathbb{R},\quad\text{and}\quad\partial_{W}v:[0,T]\times\mathbb{R}\rightarrow\mathbb{R}.
\]
To make sense of the rough integral, we require that, for each fixed
$t\in[0,T]$, the path 
\[
z\mapsto\mathcal{V}_{t}(z)=(v_{t}(z),\partial_{W}v_{t}(z))
\]
is a $\beta$-H\"older path controlled by ${\bf W}.$ Namely, both $v_{t}(\cdot)$
and $\partial_{W}v_{t}(\cdot)$ are $\beta$-H\"older continuous on
compact intervals, and the remainder 
\begin{equation}
\mathcal{R}^{\mathcal{V}_{t}}(x,y)\triangleq v_{t}(y)-v_{t}(x)-\partial_{W}v_{t}(x)\cdot(W(y)-W(x))\label{f2}
\end{equation}
has $2\beta$-H\"older regularity in the sense that for each $a\geqslant1$
we have 
\begin{equation}\label{eq:two-beta-holder}
\|\mathcal{R}^{\mathcal{V}_{t}}\|_{2\beta}^{[-a,a]}\triangleq\sup_{-a\leqslant x\neq y\leqslant a}\frac{|\mathcal{R}^{\mathcal{V}_{t}}(x,y)|}{|y-x|^{2\beta}}<\infty.
\end{equation}

Our next task is to define a rough path metric on $\mathcal{V}$ quantitatively.
As mentioned above, the main effort here is to tune this metric in
a careful way so that we end up with a contraction for the transformation
$\mathcal{M}$ described by (\ref{e4}). We first introduce the following
general notation. For a given function $f:[0,T]\times\mathbb{R}\rightarrow\mathbb{R}$,
we define: 
\begin{equation}
\llbracket f\rrbracket_{\beta/2,\beta}^{[0,t]\times[-a,a]}\triangleq\|f\|_{\infty}^{[0,t]\times[-a,a]}+a^{-\beta/2}\|f\|_{\beta/2,\beta}^{[0,t]\times[-a,a]},\label{eq:GenHolNorm}
\end{equation}
where 
\[
\|f\|_{\beta/2,\beta}^{[0,t]\times[-a,a]}\triangleq\sup_{(s,x)\neq(s',x')\in[0,t]\times[-a,a]}\frac{|f_{s'}(x')-f_{s}(x)|}{|s'-s|^{\beta/2}+|x'-x|^{\beta}}.
\]
We now define some weight functions in order to take into account
the fact that we are considering controlled processes on the real
line $\mathbb{R}$.
\begin{defn}
\label{def:weights} Recall that $\alpha,\beta,\chi$ satisfy relation~(\ref{f1}).
Consider another set of parameters $\lambda,\theta>1$. Then for $a\geqslant1$
and $t\geqslant0$ we set 
\begin{equation}
E^{\theta,\lambda}(a,t)\triangleq e^{\lambda t+\theta a+\theta at},\quad\text{and}\quad Q(a,t):=a^{\chi}\cdot\left(a^{\beta/2}+t^{-\beta/2}\right),\label{eq:EQDef}
\end{equation}
with the convention $Q(a,t)^{-1}\triangleq0$ if $t=0$.
\end{defn}

With Definition \ref{def:weights} in hand, we can now 
introduce the proper notion of controlled processes we wish to consider
in this article. Namely the space $\mathcal{B}^{\theta,\lambda}$
below will be the underlying space on which the fixed point problem \eqref{e4}
is solved.
\begin{defn}
\label{def:controlled-process} Let the notation of Definition \ref{def:weights}
prevail, and consider a controlled process $\mathcal{V}=(v,\partial_{W}v)$
such that \eqref{eq:two-beta-holder} holds true. We introduce a new parameter $\gamma\triangleq\frac{\alpha-\beta}{4}>0$.
Then we define the norm of $\mathcal{V}$ in the following way: 
\begin{multline}
\Theta^{\theta,\lambda}(\mathcal{V})\triangleq\sup_{t\in[0,T],a\geqslant1}E^{\theta,\lambda}(a,t)^{-1}\\
\times\left(\llbracket v\rrbracket_{\beta/2,\beta}^{[0,t]\times[-a,a]}+\lambda^{-\gamma}\llbracket\partial_{W}v\rrbracket_{\beta/2,\beta}^{[0,t]\times[-a,a]}+\lambda^{-\gamma}Q(a,t)^{-1}\|\mathcal{R}^{\mathcal{V}_{t}}\|_{2\beta}^{[-a,a]}\right).\label{eq:def_normTheta}
\end{multline}
We denote by $\mathcal{B}^{\theta,\lambda}$ the space of those $\mathcal{V}$'s
such that $\Theta^{\theta,\lambda}(\mathcal{V})<\infty$. It follows
that $\mathcal{B}^{\theta,\lambda}$ is a Banach space under the norm
$\Theta^{\theta,\lambda}$.
\end{defn}

\begin{rem}
\label{rem:MonBThe} Notice that the choice $\gamma=\frac{\alpha-\beta}{4}$
in Definition~\ref{def:controlled-process} might seem obscure at
this point. The reason for this will be made clear in the following
sections. Also observe that the norm $\Theta^{\theta,\lambda}$ is
decreasing in $\lambda.$ As a result, we have $\mathcal{B}^{\theta,\lambda}\subseteq\mathcal{B}^{\theta,\lambda'}$
if $\lambda<\lambda'$.
\end{rem}

Before getting into any technical estimates, it would be helpful to
point out why the norm $\Theta^{\theta,\lambda}$ is designed in such
a way. Unfortunately the explanation can only be vague at this stage. 

\vspace{2mm}\noindent 
(i) The weight $e^{\theta a}$ in $E^{\theta,\lambda}(a,t)$
allows the solution to grow at most exponentially in space at infinity.\\
(ii) The weights $a^{-\beta/2}$ in (\ref{eq:GenHolNorm}) and the
term $a^{\chi+\beta/2}$ in $Q(a,t)$ are used to absorb the polynomial
factors in $a$ which come out when estimating the transformation
$\mathcal{M}\mathcal{V}$. The appearance of polynomial factors is
not surprising due to the polynomial growth of $\|W\|_{\alpha}^{[-a,a]}$
as $a\rightarrow\infty$.\\
(iii) The weights $e^{\lambda t+\theta at}$ and $\lambda^{-\gamma}$
are used to ensure that $\mathcal{M}$ is a contraction. As we will
see, when estimating the norm of the transformation $\mathcal{M}$,
a negative power of $\lambda$ comes out due to the introduction of
these weights. As a result, if we make the a priori choice of $\lambda$
to be large enough, $\mathcal{M}$ becomes a contraction. This is
the most magical part of the analysis.\\
(iv) The term $t^{-\beta/2}$ in $Q(a,t)$ accounts for the singularity
of the heat kernel at the origin. Such singularity will appear naturally
when estimating a remainder term in $\mathcal{M}\mathcal{V}$. This
is also related to our one-dimensional rough path viewpoint.

\subsubsection{\label{Mtransformation}The transformation $\mathcal{M}$}

Having introduced the Banach space $(\mathcal{B}^{\theta,\lambda},\Theta^{\theta,\lambda}(\cdot))$,
we will now define the key transformation $\mathcal{M}$ as a rough
path in $\mathcal{B}^{\theta,\lambda}$. Namely given $\mathcal{V}=(v,\partial_{W}v)$,
we define 
\begin{equation}
\mathcal{M}\mathcal{V}\triangleq[(t,x)\mapsto((\mathcal{M}\mathcal{V})_{t}(x),\partial_{W}(\mathcal{M}\mathcal{V})_{t}(x))],\label{eq:def-trsf-M}
\end{equation}
where $\mathcal{M}\mathcal{V}$ is given by (\ref{e4}). As far as
$\partial_{W}(\mathcal{M}\mathcal{V})$ is concerned, we will see
in Section~\ref{subsec:RemEst} that for $(t,x)\in[0,T]\times\mathbb{R}$
we have 
\begin{equation}
\partial_{W}(\mathcal{M}\mathcal{V})_{t}(x)\triangleq-\frac{2}{\sigma^{2}} \, v_{t}(x),\label{eq:deriv-Gubi-M}
\end{equation}
where we recall that $\sigma^{2}$ is introduced in Notation \ref{variance}.
Furthermore, in view of the formal equation (\ref{eq:FPP}), it is
natural to consider the following transformation on $(\mathcal{B}^{\theta,\lambda},\Theta^{\theta,\lambda}(\cdot))$
\begin{equation}
\hat{\mathcal{M}}:\mathcal{V}\mapsto\hat{\mathcal{M}}\mathcal{V}\triangleq\psi^{1}+\psi^{2}+\mathcal{M}\mathcal{V},\label{f3}
\end{equation}
where $\psi^{1},\psi^{2}$ are functions respectively defined by 
\begin{equation}
\psi_{t}^{1}(x)\triangleq\int_{\mathbb{R}}\partial_{x}p_{t}(x-y)f_{0}(y)dy,\quad\text{and}\quad\psi_{t}^{2}(x)\triangleq\int_{0}^{t}\int_{\mathbb{R}}\partial_{x}p_{t-s}(x-y)g_{s}(y)dyds.\label{f4}
\end{equation}
Our goal in Section \ref{RPMP} will be to show that $\hat{\mathcal{M}}$ admits
a unique fixed point in the Banach space $\mathcal{B}^{\theta,\lambda}$.
More precisely, this is the content of the following result.
\begin{thm}
\label{thm:FPT} Let $T$ be a finite time horizon and consider $\alpha,\beta,\chi,\theta$
as in Definitions~\ref{def:weights} and~\ref{def:controlled-process}.
We also consider the quantity $\kappa_{\alpha,\chi}({\bf W}$) introduced
in~(\ref{f11}). Recall that the transformation $\hat{\mathcal{M}}$
is defined by (\ref{f3}). Then there exists $\Lambda=\Lambda_{\alpha,\beta,\chi,\theta,T,\kappa}>0$,
such that for any $\lambda>\Lambda$ and $f_{0},g\in C_{b}^{2},$
we have:

\noindent
\textnormal{(i)}
The functions $\psi^{1},\psi^{2}\in\mathcal{B}^{\theta,\lambda}$
defined by (\ref{f4}) are elements of $\mathcal{B}^{\theta,\lambda}$.

\noindent
\textnormal{(ii)}
There exists a unique fixed point for $\hat{\mathcal{M}}$, that
is a unique element $\mathcal{V}$ of the space $\mathcal{B}^{\theta,\lambda}$
satisfying 
\[
\mathcal{V}=\psi^{1}+\psi^{2}+\mathcal{M}\mathcal{V}.
\]
Otherwise stated one can solve equation \eqref{eq:continous-PDE} for $\eta=0$, under its mild form \eqref{eq:FPP}. The unique solution sits in the space $\mathcal{B}^{\theta,\lambda}$.

\noindent
\textnormal{(iii)}
The controlled process $\mathcal{V}$ satisfies 
\[
\Theta^{\theta,\lambda}(\mathcal{V})\leqslant c,
\]
for a constant $c$ depending on $f_{0}$ and $g$ only.
\end{thm}

To put Theorem \ref{thm:FPT} into perspective, let us mention again that this result allows to solve equation \eqref{eq:continous-PDE} in a pathwise sense for all continuous functions $g$. Hence there is a unique solution to the martingale problem related to equation \eqref{eq:bm-in-br-environment}, that is Theorem \ref{thm:rough-pde-informal} holds true. Going from Theorem \ref{thm:FPT} to Theorem \ref{thm:rough-pde-informal} is a 
matter of standard considerations, for which we refer e.g. to 
\cite{KS,SV}. Let us mention that the martingale $M=M^{f}$ in \eqref{eq:continuous-martingale} can now be written as 
\[
M_{t}=f_{t}(X_{t}^{c})-f_{0}(X_{0}^{c})-\int_{0}^{t}g_{s}(X_{s}^{c})ds, 
\]
where $f$ is the solution to the following rough PDE on 
$[0,\tau]\times\mathbb R$: 
\begin{equation}\label{eq:rough-pde}
f_{t}(x)=P_{\tau-t}g_{\tau}(x)+\int_{t}^{\tau}\int_{\mathbb R}p_{s-t}(x-y)\partial_{x}f_{s}(y)W(dy)\,ds. 
\end{equation}
Notice that \eqref{eq:rough-pde} is a backward version of the mild equation 
\eqref{eq:forme-mild-edp}, which can be solved exactly in the same way. In our paper we have chosen to deal with the forward equation~\eqref{eq:forme-mild-edp} for notational convenience. Let us also mention that with Theorem \ref{thm:rough-pde-informal} in hand, the existence-uniqueness of a weak solution to equation \eqref{eq:bm-in-br-environment} is again a matter of standard arguments. We summarize this in the following theorem. 
\begin{thm}
Fix a realization of the process $W$ in $(\Omega,\mathcal G,\mathbf P)$. Then there exists a probability $\mathbb P^{W}$ on the canonical space 
$\left({\rm C}(\mathbb R_{+}),\mathcal B({\rm C}(\mathbb R_{+}))\right)$ and a Brownian motion $B$ defined on that canonical space, such that $(X,B)$ 
satisfies the following integral form of equation \eqref{eq:bm-in-br-environment}:
\[
X_{t}^{c}=-\frac{1}{2}\int_{0}^{t}\dot{W}(X_{s}^{c})ds+B_{t},
\quad\text{ for all }\quad t\geqslant 0. 
\]
\end{thm}

We close this section by labeling some notation which will prevail
for our computations below.
\begin{notation}
We use $C_{\text{subscript}}$ to denote a constant depending only
on the parameters specified in the subscript but nothing else. We
also use $C$ to denote a constant depending only on $\alpha,\beta,\chi,\theta,T,\varepsilon$
but nothing else. Careful inspection on the analysis will reveal that
the dependence of $C$ on $\theta,T$ is at most $e^{\theta^{2}(1+T)^{3}}.$
The value of these constants can change from line to line. We also
use $P_{u}(w)$ to denote universal polynomials of $|w|$ which only
depend on $\alpha,\beta,\chi$ and $T$ (the dependence on $T$ only
appears in the coefficients and is polynomial).
\end{notation}

\section{\label{RPMP} The rough path formulation of the martingale problem}

As mentioned in Section \ref{BMBE}, Theorem \ref{thm:FPT} is the key to establish existence and uniqueness for the martingale problem related to equation \eqref{eq:bm-in-br-environment}. Arguably, our Theorem \ref{thm:FPT} is contained in \cite{DD} Theorem 5. However we have decided to include a detailed proof of the main estimates here for two reasons: first our context where $W$ only depends on the space variable lead to simpler considerations than in \cite{DD}. In addition, the bounds presented in this section will 
play a prominent role in our analysis of the convergence in Section \ref{sec:ConvEst}. We will thus handle the rough path type estimates for $\mathcal{M}v$ in Sections 
\ref{sec:RIP}-\ref{sec:KEHKRI}-\ref{sec:EstMV}, and complete the proof of Theorem \ref{thm:FPT} in Sections \ref{sec:estimate-deter-funct}-\ref{sec:FPT}. 

\subsection{Rough integral estimates}\label{sec:RIP}

As mentionned before, the integral $\int_{x}^{y}v_{s}(z)dW(z)$ in~\eqref{e4} 
is defined in the rough paths sense. With respect to the usual rough paths setting of \cite{Gu}, those integrals in $\mathbb R$ have to involve weighted norms in space. In the current section we develop some basic tools for integrals of processes in spaces of the form  $\mathcal{B}^{\theta,\lambda}$ 
(see Definition~\ref{def:controlled-process}). Let us start with a bound on those integrals, depending on an interval $[-a,a]\subset\mathbb R$. 

\begin{prop}
\label{prop:RIEstNoTime} Let ${\bf W}=(W^{1},W^{2})$ be the one-dimensional rough path given by (\ref{eq:def-w1-w2}). We consider
$\alpha,\beta$ satisfying \eqref{f1} and a controlled path $x\mapsto\mathcal{V}(x)$ with $\mathcal{V}(x)=(v(x),\partial_{W}v(x))$
as introduced in~\eqref{f12}-\eqref{f2}. Then for
any $a\geqslant1$ and $x,y\in[-a,a]$ the integral $\int_{x}^{y}v(z)dW(z)$
is well-defined, and we have 
\begin{multline}
\Big|\int_{x}^{y}v(z)dW(z)-v(x)W^{1}(x,y)-\partial_{W}v(x)W^{2}(x,y)\Big|\\
\leqslant C_{\beta}\cdot\left(\|\partial_{W}v\|_{\beta}^{[-a,a]}\cdot\|W^{2}\|_{2\alpha}^{[-a,a]}\cdot|y-x|^{2\alpha+\beta}+\|\mathcal{R}^{\mathcal{V}}\|_{2\beta}^{[-a,a]}\cdot\|W\|_{\alpha}^{[-a,a]}\cdot|y-x|^{\alpha+2\beta}\right),\label{eq:RI}
\end{multline}
where the H\"older norms above are understood as in~(\ref{eq:def-Ha}).
\end{prop}

\begin{proof}
This is a mere elaboration of \cite[Theorem 1]{Gu}, whose proof is omitted for sake of conciseness.
\end{proof}

Notice that we have stated Proposition~\ref{prop:RIEstNoTime} for
a usual controlled process, since we are only using the function $x\mapsto\mathcal{V}(x)$
in our estimate. Therefore here $v$ is just considered as one single
rough path rather than a family of rough paths parametrized by time.
However, in the sequel we shall need a corollary stated for processes
in $\mathcal{B}^{\theta,\lambda}$. This is summarized in the lemma
below.
\begin{lem}
\label{eq:RIEst} Consider $\alpha,\beta,\chi$ such that \eqref{f1} is fulfilled. Let $\mathcal{V}=[(t,x)\mapsto\mathcal{V}_{t}(x)=(v_{t}(x),\partial_{W}v_{t}(x))]$
be an element in $\mathcal{B}^{\theta,\lambda}$, as given in Definition~\ref{def:controlled-process}.
To simplify notation, let us set 
\begin{equation}
\kappa\triangleq\kappa_{\alpha,\chi}({\bf W}),\quad\Theta\triangleq\Theta^{\theta,\lambda}(\mathcal{V}),\quad E\triangleq E^{\theta,\lambda}(a,t),\label{eq: kappa-Theta-E}
\end{equation}
and 
\begin{equation}
D(a,t,z)\triangleq a^{2\chi}|z|^{2\alpha}+a^{2\chi+\beta/2}|z|^{2\alpha+\beta}+a^{2\chi}(a^{\beta/2}+t^{-\beta/2})|z|^{\alpha+2\beta}.\label{eq:DFunction}
\end{equation}
Then for any $t\in(0,T]$, $a\geqslant1$ and $x,y\in[-a,a],$ we
have 
\begin{equation}
\Big|\int_{x}^{y}v_{t}(z)dW(z)-v_{t}(x)W^{1}(x,y)\Big|\leqslant C_{\beta}\,\kappa\,\Theta E\,\lambda^{\gamma}\,D(a,t,y-x)\label{eq:RI0}
\end{equation}
and 
\begin{equation}
\Big|\int_{x}^{y}v_{t}(z)dW(z)\Big|\leqslant C_{\beta}\,\kappa\,\Theta E\,\lambda^{\gamma}\,\left(a^{\chi}|y-x|^{\alpha}+D(a,t,y-x)\right).\label{eq:RI1}
\end{equation}
\end{lem}

\begin{proof}
In order to get the first estimate we resort to relation (\ref{eq:RI}),
which yields 
\begin{multline*}
\Big|\int_{x}^{y}v_{t}(z)dW(z)-v_{t}(x)W^{1}(x,y)\Big|\leqslant\Big|\partial_{W}v_{t}(x)\Big|\cdot\big|W^{2}(x,y)\big|\\
+C_{\beta}\cdot\left(\|\partial_{W}v_{t}\|_{\beta}^{[-a,a]}\cdot\|W^{2}\|_{2\alpha}^{[-a,a]}\cdot|y-x|^{2\alpha+\beta}+\|\mathcal{R}^{\mathcal{V}_{t}}\|_{2\beta}^{[-a,a]}\cdot\|W\|_{\alpha}^{[-a,a]}\cdot|y-x|^{\alpha+2\beta}\right).
\end{multline*}
Next we bound $\big|\partial_{W}v_{t}(x)\big|$ and $\big|W^{2}(x,y)\big|$
above thanks to their H\"older and supremum norms. We get 
\begin{multline*}
\Big|\int_{x}^{y}v_{t}(z)dW(z)-v_{t}(x)W^{1}(x,y)\Big|\leqslant\|\partial_{W}v\|_{\infty}^{[0,t]\times[-a,a]}\cdot\|W^{2}\|_{2\alpha}^{[-a,a]}\cdot|y-x|^{2\alpha}\\
+C_{\beta}\cdot\left(\|\partial_{W}v\|_{\beta/2,\beta}^{[0,t]\times[-a,a]}\cdot\|W^{2}\|_{2\alpha}^{[-a,a]}\cdot|y-x|^{2\alpha+\beta}+\|\mathcal{R}^{\mathcal{V}_{t}}\|_{2\beta}^{[-a,a]}\cdot\|W\|_{\alpha}^{[-a,a]}\cdot|y-x|^{\alpha+2\beta}\right).
\end{multline*}
We now use the definition (\ref{f11}) of $\kappa_{\alpha,\chi}({\bf W})$
as well as Definition \ref{def:controlled-process} in order to obtain
\begin{align}
 & \Big|\int_{x}^{y}v_{t}(z)dW(z)-v_{t}(x)W^{1}(x,y)\Big|\nonumber\\
 & \leqslant C_{\beta}\cdot\kappa\cdot\Theta E\cdot\lambda^{\gamma}\cdot\left(a^{2\chi}|y-x|^{2\alpha}+a^{2\chi+\beta/2}|y-x|^{2\alpha+\beta}+a^{2\chi}(a^{\beta/2}+t^{-\beta/2})|y-x|^{\alpha+2\beta}\right),\label{eq:estim-fourterms}
\end{align}
from which our estimate (\ref{eq:RI0}) is easily achieved. Our second
claim (\ref{eq:RI1}) is also deduced from (\ref{eq:RI0}) plus
the trivial bound 
\[
|v_{t}(x)W^{1}(x,y)|\leqslant\|v\|_{\infty}^{[0,t]\times[-a,a]}\cdot a^{\chi}\cdot\kappa\cdot|y-x|^{\alpha}\leqslant\kappa\Theta E\cdot a^{\chi}|y-x|^{\alpha}.
\]
The proof is now complete.
\end{proof}

\subsection{A key estimate for heat kernel convolutions with rough integrals}\label{sec:KEHKRI}

In view of the expression (\ref{e4}) of $(\mathcal{M}\mathcal{V})_{t}(x)$,
an essential ingredient in the analysis of $\mathcal{M}\mathcal{V}$
is to estimate heat kernel convolutions with rough path integrals
carefully. This is summarized in the following key lemma.
\begin{lem}
\label{lem:KeyLem} Consider $\alpha,\beta,\chi$ such that (\ref{f1})
is fulfilled. Let $\mathcal{V}=(v,\partial_{W}v)$ be a controlled
path in the space $\mathcal{B}^{\theta,\lambda}$ introduced in Definition
\ref{def:controlled-process}. Recall that the heat kernel on $\mathbb{R}$
is given by (\ref{heat_kernel}) and the constants $\kappa,\Theta,E$
are spelled out in (\ref{eq: kappa-Theta-E}). Then the following
two estimates hold true:

\vspace{2mm}\noindent (i) Let $0\leqslant\gamma_{1}\leqslant\gamma_{2}\leqslant\beta/2$
and $k\geqslant1.$ For any $a\geqslant1,$ $x\in[-a,a]$ and $\tau_{1}\leqslant\tau_{2}\in(0,T],$
we have 
\begin{equation}
\int_{\mathbb{R}}\int_{0}^{\tau_{1}}\frac{|\partial_{x\cdots x}^{k}p_{1}(w)|}{s^{1+\gamma_{1}}}\left|\int_{x}^{x+\sqrt{s}w}v_{\tau_{2}-s}(z)dW(z)\right|dwds\leqslant C\cdot\kappa\Theta E\cdot\lambda^{\gamma-\frac{\alpha-\beta}{2}}a^{\gamma_{2}}\cdot\tau_{1}^{\gamma_{2}-\gamma_{1}},\label{eq:Key1}
\end{equation}
(ii) Let $0\leqslant\sigma_{1}\leqslant\beta$ and $k\geqslant1.$
For any $a\geqslant1,$ $x\in[-a,a]$ and $\tau_{1}\leqslant\tau_{2}\in(0,T],$
we have 
\begin{align}
 & \int_{\mathbb{R}}\int_{0}^{\tau_{1}}\frac{|\partial_{x\cdots x}^{k}p_{1}(w)|}{s^{1+\sigma_{1}}}\left|\int_{x}^{x+\sqrt{s}w}(v_{\tau_{2}-s}(z)-v_{\tau_{2}-s}(x))dW(z)\right|dwds\nonumber \\
 & \ \ \ \leqslant C\cdot\kappa\Theta E\cdot\lambda^{\gamma-\frac{\alpha-\beta}{4}}a^{\chi}\left(a^{\beta/2}+\tau_{2}^{-\beta/2}\right)\cdot\tau_{1}^{\beta-\sigma_{1}}.\label{eq:Key2}
\end{align}

\end{lem}

\begin{rem}
\label{rem:ChoGam}The factors $\lambda^{\gamma-\frac{\alpha-\beta}{2}}$
and $\lambda^{\gamma-\frac{\alpha-\beta}{4}}$ appearing in the above
two estimates are crucial. In fact, we will choose $\gamma=\frac{\alpha-\beta}{4}$
in the a priori definition of $\gamma.$ This implies that the factor
$\lambda^{\gamma-\frac{\alpha-\beta}{2}}=\lambda^{-\frac{\alpha-\beta}{4}}$
is decaying. This
point is critical for obtaining the contraction property of $\mathcal{M}\mathcal{V}.$
\end{rem}

\begin{proof}[Proof of Lemma \ref{lem:KeyLem}] We will divide
this proof into several steps.

\noindent \textit{Step 1: Estimate for the rough integral.} Set 
\begin{equation}
{\mathcal{A}}_{x,s,w,\tau_{2}}=\int_{x}^{x+\sqrt{s}w}v_{\tau_{2}-s}(z)dW(z).\label{f40}
\end{equation}
Then according to (\ref{eq:RI1}), for any $x\in[-a,a]$ we can upper
bound $|{\mathcal{A}}_{x,s,w,\tau_{2}}|$ by 
\begin{align*}
 & C_{\beta}\cdot\kappa\Theta\cdot E(a+\sqrt{s}|w|,\tau_{2}-s)\cdot\lambda^{\gamma}\cdot\left((a+\sqrt{s}|w|)^{\chi}s^{\alpha/2}|w|^{\alpha}+(a+\sqrt{s}|w|)^{2\chi}s^{\alpha}|w|^{2\alpha}\right.\\
 & \ \ \ \left.+(a+\sqrt{s}|w|)^{2\chi+\beta/2}s^{\alpha+\beta/2}|w|^{2\alpha+\beta}+(a+\sqrt{s}|w|)^{2\chi+\beta/2}s^{\alpha/2+\beta}|w|^{\alpha+2\beta}\right.\\
 & \ \ \ \left.+(a+\sqrt{s}|w|)^{2\chi}(\tau_{2}-s)^{-\beta/2}s^{\alpha/2+\beta}|w|^{\alpha+2\beta}\right).
\end{align*}
Hence recalling that we write $P_{u}(w)$ for any polynomial in $w$
and that our time horizon is $T$, we get that $|{\mathcal{A}}_{x,s,w,\tau_{2}}|$
can be upper bounded as follows: 
\begin{align}
 & C_{\beta}\kappa\Theta\lambda^{\gamma}\cdot P_{u}(w)\cdot E(a+\sqrt{s}|w|,\tau_{2}-s)\cdot\left((a+\sqrt{T}|w|)^{\chi}s^{\alpha/2}+(a+\sqrt{T}|w|)^{2\chi}s^{\alpha}\right.\label{eq:esti}\\
 & \left.+(a+\sqrt{T}|w|)^{2\chi+\beta/2}s^{\alpha+\beta/2}+(a+\sqrt{T}|w|)^{2\chi+\beta/2}s^{\alpha/2+\beta}+(a+\sqrt{T}|w|)^{2\chi}(\tau_{2}-s)^{-\beta/2}s^{\alpha/2+\beta}\right).\nonumber 
\end{align}
In addition, still invoking the fact that our time horizon is $T$,
we can write 
\begin{align}
E(a+\sqrt{s}|w|,\tau_{2}-s) & =e^{\lambda(\tau_{2}-s)+\theta(a+\sqrt{s}|w|)+\theta(a+\sqrt{s}|w|)(\tau_{2}-s)}\\
 & \leqslant E(a,\tau_{2})\cdot e^{\theta(1+T)\sqrt{T}|w|}\cdot e^{-\left(\lambda+\theta(a+\sqrt{T}|w|)\right)s}.\label{estimateEE}
\end{align}
Consider now $\gamma_{1}\leqslant\gamma_{2}\leqslant\beta/2$ as in the statement of our current lemma.
Setting $\rho\triangleq a+\sqrt{T}|w|$ and plugging \eqref{estimateEE} in (\ref{eq:esti}), we obtain that 
\begin{equation}
s^{-(1+\gamma_{1})}\big|{\mathcal{A}}_{x,s,w,\tau_{2}}\big|\leqslant C_{\beta}\kappa\Theta\lambda^{\gamma}P_{u}(w)E(a,\tau_{2})e^{\theta(1+T)\sqrt{T}|w|}\cdot e^{-(\lambda+\theta\rho)s}\phi_{1}(\rho,s),\label{eq:Key1sInt}
\end{equation}
where we define two functions $\phi_{i}$, $i=1,2$, as follows:
\begin{multline}
\phi_{i}(\rho,s)=\left(\rho^{\chi}s^{\alpha/2-1-\gamma_{i}}+\rho^{2\chi}s^{\alpha-1-\gamma_{i}}\right.\\
\ \ \ \left.+\rho^{2\chi+\beta/2}s^{\alpha+\beta/2-1-\gamma_{i}}+\rho^{2\chi+\beta/2}s^{\alpha/2+\beta-1-\gamma_{i}}+\rho^{2\chi}(\tau_{2}-s)^{-\beta/2}s^{\alpha/2+\beta-1-\gamma_{i}}\right),\label{eq:phifunction}
\end{multline}
and we notice that $\phi_{2}$ will be used in \eqref{eq:JancientI}.

\noindent \textit{Step 2: Decomposition of a time integral.} With
a proper estimate for ${\mathcal{A}}_{x,s,w,\tau_{2}}$ in hand, we
turn to an estimate of the double integral~(\ref{eq:Key1}). To this
aim, we first consider the $s$-integral and we set 
\begin{equation}
{\mathcal{D}}_{\tau_{1}}(w)=\int_{0}^{\tau_{1}}s^{-(1+\gamma_{1})}{\mathcal{A}}_{x,s,w,\tau_{2}}ds.\label{f41}
\end{equation}
Then taking (\ref{eq:Key1sInt}) and (\ref{eq:phifunction}) into
account we get 
\[
\big|{\mathcal{D}}_{\tau_{1}}(w)\big|\leqslant C_{\beta}\kappa\Theta\lambda^{\gamma}P_{u}(w)E(a,\tau_{2})e^{\theta(1+T)\sqrt{T}|w|}\int_{0}^{\tau_{1}}e^{-(\lambda+\theta\rho)s}\phi_{1}(\rho,s)ds.
\]
We will now gain a small factor $\tau_{1}^{\gamma_{2}-\gamma_{1}}$
by playing with the fact that $\gamma_{2}\geqslant\gamma_{1}$. Namely we
easily obtain 
\begin{equation}
\big|{\mathcal{D}}_{\tau_{1}}(w)\big|
\leqslant
C_{\beta}\tau_{1}^{\gamma_{2}-\gamma_{1}}\kappa\Theta\lambda^{\gamma}P_{u}(w)E(a,\tau_{2})e^{\theta(1+T)\sqrt{T}|w|}{\mathcal{J}}_{\tau_{1}},\label{f42}
\end{equation}
where we have set 
\begin{equation}
{\mathcal{J}}_{\tau_{1}}=\int_{0}^{\tau_{1}}e^{-(\lambda+\theta\rho)s}\phi_{2}(\rho,s)ds,
%\left(\rho^{\chi}s^{\alpha/2-1-\gamma_{2}}+\rho^{2\chi}s^{\alpha-1-\gamma_{2}}+\rho^{2\chi+\beta/2}s^{\alpha+\beta/2-1-\gamma_{2}}
\label{eq:JancientI}
\end{equation}
and where we recall that $\phi_{2}$ is introduced in \eqref{eq:phifunction}.
In addition, we shall decompose the integral $\mathcal{J}$ into ${\mathcal{J}}={\mathcal{J}}^{1}+{\mathcal{J}}^{2}$,
where we define 
\begin{multline*}
{\mathcal{J}}_{\tau_{1}}^{1}\triangleq\int_{0}^{\tau_{1}}e^{-(\lambda+\theta\rho)s}\cdot\left(\rho^{\chi}s^{\alpha/2-1-\gamma_{2}}+\rho^{2\chi}s^{\alpha-1-\gamma_{2}}+\rho^{2\chi+\beta/2}s^{\alpha+\beta/2-1-\gamma_{2}}\right.\\
\ \ \ \left.+\rho^{2\chi+\beta/2}s^{\alpha/2+\beta-1-\gamma_{2}}\right)ds
\end{multline*}
and 
\begin{equation}
{\mathcal{J}}_{\tau_{1}}^{2}\triangleq\int_{0}^{\tau_{1}}e^{-(\lambda+\theta\rho)s}\rho^{2\chi}(\tau_{2}-s)^{-\beta/2}s^{\alpha/2+\beta-1-\gamma_{2}}ds.\label{f5}
\end{equation}
The methods in order to estimate ${\mathcal{J}}^{1}$ and ${\mathcal{J}}^{2}$
are slightly different, and thus those two tasks will be carried out
separately. The main non-trivial effort here is to see that 
\begin{equation}
{\mathcal{J}}_{\tau_{1}}^{i}\leqslant C_{\alpha,\beta}\rho^{\gamma_{2}}\cdot\lambda^{\frac{\beta-\alpha}{2}},\quad\text{for}\quad i=1,2.\label{eq:Key1I_12}
\end{equation}
This is the critical point where a negative power of $\lambda$ appears
(recall that $\beta<\alpha$), which in turn contributes to the contraction
property of the transformation $\mathcal{M}.$

\noindent \textit{Step 3: Estimation of }${\mathcal{J}}^{1}$. By
applying a change of variables $r=(\lambda+\theta\rho)s,$ we obtain
\begin{multline*}
\rho^{-\gamma_{2}}{\mathcal{J}}_{\tau_{1}}^{1}=\int_{0}^{(\lambda+\theta\rho)\tau_{1}}e^{-r}\left(\frac{\rho^{\chi-\gamma_{2}}}{(\lambda+\theta\rho)^{\alpha/2-\gamma_{2}}}r^{\alpha/2-\gamma_{2}-1}+\frac{\rho^{2\chi-\gamma_{2}}}{(\lambda+\theta\rho)^{\alpha-\gamma_{2}}}r^{\alpha-\gamma_{2}-1}\right.\\
\left.+\frac{\rho^{2\chi+\beta/2-\gamma_{2}}}{(\lambda+\theta\rho)^{\alpha+\beta/2-\gamma_{2}}}r^{\alpha+\beta/2-\gamma_{2}-1}+\frac{\rho^{2\chi+\beta/2-\gamma_{2}}}{(\lambda+\theta\rho)^{\alpha/2+\beta-\gamma_{2}}}r^{\alpha/2+\beta-\gamma_{2}-1}\right)dr.
\end{multline*}
Therefore collecting the terms of the form $\rho^{\alpha_{1}}/(\lambda+\theta\,\rho)^{\alpha_{2}}$
and computing the integrals $\int_{0}^{\infty}e^{-r}r^{\alpha_{3}}\,dr$,
the reader can easily check that 
\begin{equation}
\rho^{-\gamma_{2}}{\mathcal{J}}_{\tau_{1}}^{1}\leqslant C_{\alpha,\beta}\cdot\left(\frac{\rho^{\chi-\gamma_{2}}}{(\lambda+\theta\rho)^{\alpha/2-\gamma_{2}}}+\frac{\rho^{2\chi-\gamma_{2}}}{(\lambda+\theta\rho)^{\alpha-\gamma_{2}}}+\frac{\rho^{2\chi+\beta/2-\gamma_{2}}}{(\lambda+\theta\rho)^{\alpha+\beta/2-\gamma_{2}}}+\frac{\rho^{2\chi+\beta/2-\gamma_{2}}}{(\lambda+\theta\rho)^{\alpha/2+\beta-\gamma_{2}}}\right).\label{f51}
\end{equation}
Hence to reach the estimate (\ref{eq:Key1I_12}), we need to show
that 
\begin{equation}
\frac{\rho^{\chi-\gamma_{2}}}{(\lambda+\theta\rho)^{\alpha/2-\gamma_{2}}}\vee\frac{\rho^{2\chi-\gamma_{2}}}{(\lambda+\theta\rho)^{\alpha-\gamma_{2}}}\vee\frac{\rho^{2\chi+\beta/2-\gamma_{2}}}{(\lambda+\theta\rho)^{\alpha+\beta/2-\gamma_{2}}}\vee\frac{\rho^{2\chi+\beta/2-\gamma_{2}}}{(\lambda+\theta\rho)^{\alpha/2+\beta-\gamma_{2}}}\leqslant\lambda^{\frac{\beta-\alpha}{2}}.\label{eq:QuoEstP1}
\end{equation}

In order to achieve (\ref{eq:QuoEstP1}), we first observe that the
above four quotients all have the form 
\begin{equation}
\frac{\rho^{b-\gamma_{2}}}{(\lambda+\theta\rho)^{a-\gamma_{2}}},\qquad\ \ \ \text{with }0<b<a\ \text{and }a>\gamma_{2}.\label{eq:GenFormQuo}
\end{equation}
More precisely, we recall from (\ref{f1}) that $\beta<\alpha$ and
$\chi<1/2$, and we have also assumed that $\gamma_{1}\le\gamma_{2}\le\beta/2$.
Thus it is readily checked that the exponents in~(\ref{eq:QuoEstP1})
satisfy~(\ref{eq:GenFormQuo}), by respectively considering the following
values for $a,b$: 
\begin{equation}
\text{(i) }\begin{cases}
b=\chi\\
a=\alpha/2
\end{cases},\ \text{(ii) }\begin{cases}
b=2\chi\\
a=\alpha
\end{cases},\ \text{(iii) }\begin{cases}
b=2\chi+\beta/2\\
a=\alpha+\beta/2
\end{cases},\ \text{(iv) }\begin{cases}
b=2\chi+\beta/2\\
a=\alpha/2+\beta
\end{cases}\label{f6}
\end{equation}
With (\ref{eq:QuoEstP1}) in mind, we now claim that for all $\lambda,\rho\geqslant1$
we have: 
\begin{equation}
\frac{\rho^{b-\gamma_{2}}}{(\lambda+\theta\rho)^{a-\gamma_{2}}}\leqslant\lambda^{b\vee\gamma_{2}-a}\ \ \ \forall\lambda,\rho\geqslant1.\label{eq:LamEstGen}
\end{equation}
Indeed, if $\rho\geqslant\lambda,$ we have (recall $\lambda,\theta>1$
and $b<a$) 
\[
\frac{\rho^{b-\gamma_{2}}}{(\lambda+\theta\rho)^{a-\gamma_{2}}}\leqslant\rho^{b-a}\leqslant\lambda^{b-a}\leqslant\lambda^{b\vee\gamma_{2}-a}.
\]
On the other hand, if $1\leqslant\rho<\lambda$ and $b>\gamma_{2}$,
we have 
\[
\frac{\rho^{b-\gamma_{2}}}{(\lambda+\theta\rho)^{a-\gamma_{2}}}\leqslant\frac{\rho^{b-\gamma_{2}}}{\lambda^{a-\gamma_{2}}}\leqslant\lambda^{b-a}\leqslant\lambda^{b\vee\gamma_{2}-a}.
\]
Eventually, if $1\leqslant\rho<\lambda$ and $b\leqslant\gamma_{2}$,
then 
\[
\frac{\rho^{b-\gamma_{2}}}{(\lambda+\theta\rho)^{a-\gamma_{2}}}=\frac{1}{\rho^{\gamma_{2}-b}(\lambda+\theta\rho)^{a-\gamma_{2}}}\leqslant\lambda^{\gamma_{2}-a}\leqslant\lambda^{b\vee\gamma_{2}-a}.
\]
Therefore, the claim (\ref{eq:LamEstGen}) follows. The estimate (\ref{eq:Key1I_12})
for ${\mathcal{J}}^{1}$ is then established by reporting~(\ref{eq:LamEstGen})
into the cases (i)-(iv) in~(\ref{f6}). Notice that we have 
\[
a-b\vee\gamma_{2}=\begin{cases}
\text{(i) }\frac{\alpha}{2}-\chi\vee\gamma_{2}\geqslant (\alpha-\beta)/2\\
\text{(ii) }\alpha-2\chi\vee\gamma_{2}\geqslant\alpha-\beta\\
\text{(iii) }\left(\alpha+\frac{\beta}{2}\right)-\left(2\chi+\frac{\beta}{2}\right)\vee\gamma_{2}=\alpha-2\chi\geqslant\alpha-\beta\\
\text{(iv) }\left(\frac{\alpha}{2}+\beta\right)-\left(2\chi+\frac{\beta}{2}\right)\vee\gamma_{2}=\frac{\alpha+\beta}{2}-2\chi
\end{cases}\geqslant\frac{\alpha-\beta}{2},
\]
in each of the four cases respectively. This achieves~(\ref{eq:QuoEstP1}),
and recall that plugging~(\ref{eq:QuoEstP1}) into (\ref{f51}) we
obtain 
\begin{equation}
\rho^{-\gamma_{2}}\mathcal{J}_{\tau_{1}}^{1}\le C_{\alpha,\beta}\,\lambda^{\frac{\beta-\alpha}{2}},\label{f61}
\end{equation}
that is relation (\ref{eq:Key1I_12}) for ${\mathcal{J}}^{1}$.

\noindent \textit{Step 4: Estimation of} ${\mathcal{J}}^{2}$. We
first recall an elementary estimate, which holds true for all $r,c>0$
\begin{equation}
e^{-r}\leqslant c^{c}r^{-c}.\label{eq:EleEst}
\end{equation}
Observe that (\ref{eq:EleEst}) can be seen from $x\leqslant e^{x}$
for $x>0$ and taking $x=r/c$. Now we recall that ${\mathcal{J}}^{2}$
is defined by~(\ref{f5}). We use~(\ref{eq:EleEst}) in this definition
and choose $r\triangleq(\lambda+\theta\rho)s$, with a parameter $c>0$
to be specified later on. We end up with 
\[
\rho^{-\gamma_{2}}{\mathcal{J}}_{\tau_{1}}^{2}\leqslant\frac{c^{c}\rho^{2\chi-\gamma_{2}}}{(\lambda+\theta\rho)^{c}}\int_{0}^{\tau_{1}}s^{\alpha/2+\beta-\gamma_{2}-c-1}(\tau_{2}-s)^{-\beta/2}ds.
\]
Therefore the elementary change of variable $s\triangleq\tau_{1}u$
yields 
\[
\rho^{-\gamma_{2}}{\mathcal{J}}_{\tau_{1}}^{2}\leqslant\frac{c^{c}\rho^{2\chi-\gamma_{2}}\tau_{1}^{\alpha/2+\beta-\gamma_{2}-c}}{(\lambda+\theta\rho)^{c}}\int_{0}^{1}u^{\alpha/2+\beta-\gamma_{2}-c-1}(\tau_{2}-\tau_{1}u)^{-\beta/2}du,
\]
and since we have assumed $\tau_{1}\le\tau_{2}$ we can write $(\tau_{2}-\tau_{1}u)^{-\beta/2}\le\tau_{1}^{-\beta/2}(1-u)^{-\beta/2}$.
We get 
\begin{equation}
\rho^{-\gamma_{2}}{\mathcal{J}}_{\tau_{1}}^{2}\leqslant\frac{c^{c}\rho^{2\chi-\gamma_{2}}\tau_{1}^{\alpha/2+\beta/2-\gamma_{2}-c}}{(\lambda+\theta\rho)^{c}}\int_{0}^{1}u^{\alpha/2+\beta-\gamma_{2}-c-1}(1-u)^{-\beta/2}du.\label{f7}
\end{equation}
We now pick a convenient value for the parameter $c$ above. Indeed,
one can choose $c=\frac{\alpha}{2}+\frac{\beta}{2}-\gamma_{2}$, so
that $\tau_{1}^{\alpha/2+\beta/2-\gamma_{2}-c}=1$. With this value
of $c$, we let the reader check that the exponent $\alpha/2+\beta-\gamma_{2}-c-1$
in the right hand side of (\ref{f7}) is larger than $-1$. It follows
that 
\begin{equation}
\rho^{-\gamma_{2}}{\mathcal{J}}_{\tau_{1}}^{2}\leqslant C_{\alpha,\beta}\cdot\frac{\rho^{2\chi-\gamma_{2}}}{(\lambda+\theta\rho)^{\alpha/2+\beta/2-\gamma_{2}}}.\label{f8}
\end{equation}
Notice that the ratio in (\ref{f8}) has the form of (\ref{eq:GenFormQuo}),
with $b=2\chi$ and $a=(\alpha+\beta)/2$. In addition, our relation
(\ref{f1}) on $\alpha,\beta,\chi$ ensures that $b<a$. Thus according
to (\ref{eq:LamEstGen}), we obtain that 
\begin{equation}
\rho^{-\gamma_{2}}\mathcal{J}_{\tau_{1}}^{2}\leqslant C_{\alpha,\beta}\cdot\lambda^{(2\chi)\vee\gamma_{2}-\frac{\alpha+\beta}{2}}\leqslant C_{\alpha,\beta}\cdot\lambda^{\frac{\beta-\alpha}{2}},\label{f9}
\end{equation}
where we have resorted to the fact that $(2\chi)\vee\gamma_{2}\le\beta$
for the second inequality above. Therefore, we have established the
estimate (\ref{eq:Key1I_12}) for ${\mathcal{J}}^{2}$.

Summarizing our considerations so far, we have proved~(\ref{eq:Key1I_12})
by gathering (\ref{f9}) and~(\ref{f61}). Therefore recalling that
we have set ${\mathcal{J}}={\mathcal{J}}^{1}+{\mathcal{J}}^{2}$ and
relation~(\ref{f42}), we get that the term $\mathcal{D}_{\tau_{1}}$
defined by~(\ref{f41}) satisfies 
\begin{equation}
\big|{\mathcal{D}}_{\tau_{1}}(w)\big|\leqslant C_{\alpha,\beta}\cdot\tau_{1}^{\gamma_{1}-\gamma_{2}}\kappa\Theta P_{u}(w)E(a,\tau_{2})e^{\theta(1+T)\sqrt{T}|w|}%{\mathcalJ}_{\tau_{1}}
\cdot\lambda^{\gamma+\frac{\beta-\alpha}{2}}\cdot\rho^{\gamma_{2}}.\label{f10}
\end{equation}

\noindent \textit{Step 5: Conclusion}. Recall that our aim is to achieve
the upper bound~(\ref{eq:Key1}). With~(\ref{f40}) and~(\ref{f41})
in hand, the left hand side of~(\ref{eq:Key1}) can be written as
\[
\int_{\mathbb{R}}|\partial_{x\cdots x}^{k}p_{1}(w)|\,{\mathcal{D}}_{\tau_{1}}(w)\,dw.
\]Therefore plugging~(\ref{f10}) into the above expression, we get
that (recall $\rho=a+\sqrt{T}|w|$): 
\begin{multline*}
\int_{\mathbb{R}}|\partial_{x\cdots x}^{k}p_{1}(w)|\,{\mathcal{D}}_{\tau_{1}}(w)\,dw\leqslant C_{\alpha,\beta}\cdot\kappa\Theta E(a,\tau_{2})\lambda^{\gamma+\frac{\beta-\alpha}{2}}\tau_{1}^{\gamma_{2}-\gamma_{1}}\\
\times\int_{\mathbb{R}}|\partial_{x\cdots x}^{k}p_{1}(w)|\cdot P_{u}(w)\cdot e^{\theta(1+T)\sqrt{T}|w|}(a+\sqrt{T}|w|)^{\gamma_{2}}dw \, .
\end{multline*}
Hence some elementary heat kernel estimates entail
\begin{eqnarray}
\int_{\mathbb{R}}|\partial_{x\cdots x}^{k}p_{1}(w)|\,{\mathcal{D}}_{\tau_{1}}(w)\,dw & \leqslant & C_{\alpha,\beta,k,\sigma}e^{\theta^{2}(1+T)^{3}}\cdot\kappa\Theta E(a,\tau_{2})\lambda^{\gamma+\frac{\beta-\alpha}{2}}\tau_{1}^{\gamma_{2}-\gamma_{1}}\cdot a^{\gamma_{2}}\nonumber \\
 & = & C_{\alpha,\beta,k,\theta,T,\sigma}\cdot\kappa\Theta E(a,\tau_{2})\cdot\lambda^{\gamma+\frac{\beta-\alpha}{2}}a^{\gamma_{2}}\tau_{1}^{\gamma_{2}-\gamma_{1}}.\label{eq:elem-heatk-estim}
\end{eqnarray}
This finishes the proof of~(\ref{eq:Key1}). For sake of conciseness
we leave the proof of~(\ref{eq:Key2}) to the patient reader. It
is based on the same kind of arguments as for~(\ref{eq:Key1}), also
taking into account the extra regularity brought in by the increments
$v_{\tau_{2}-s}(z)-v_{\tau_{2}-s}(x)$. \end{proof}

\noindent For our future computations it will be useful to extend
the estimate (\ref{eq:Key2}) to a context with time increments of
$v$. This is the content of the following lemma.
\begin{lem}
\noindent \label{lem:ExtKey2} We assume that the hypothesis of Lemma
\ref{lem:KeyLem} hold true. In particular, we consider a process
$v$ in $\mathcal{B}^{\theta,\lambda}$ and we recall that the parameter
$\gamma$ in Definition \ref{def:controlled-process} satisfies $\gamma=\frac{\alpha-\beta}{4}$.
Then the following bound holds true: 
\begin{multline}
\int_{\mathbb{R}}\int_{0}^{\tau_{1}}\frac{|\partial_{x\cdots x}^{k}p_{1}(w)|}{s^{1+\sigma_{1}}} \,
\Big|\int_{x}^{x+\sqrt{s}w}(v_{\tau_{2}-s}(z)-v_{\tau_{2}}(x))dW(z)\Big|dwds\\
\leqslant C\cdot\kappa\Theta E(a,\tau_{2})a^{\chi}\left(a^{\beta/2}+\tau_{2}^{-\beta/2}\right)\cdot\tau_{1}^{\beta-\sigma_{1}}.\label{eq:ExtKey2}
\end{multline}
\end{lem}

\begin{proof}
\noindent We decompose the increment $v_{\tau_{2}-s}(z)-v_{\tau_{2}}(x)$
into 
\[
v_{\tau_{2}-s}(z)-v_{\tau_{2}}(x)=\left[v_{\tau_{2}-s}(z)-v_{\tau_{2}-s}(x)\right]+\left[v_{\tau_{2}-s}(x)-v_{\tau_{2}}(x)\right] \, .
\]
Then notice that the integral (\ref{eq:ExtKey2}) corresponding to
$v_{\tau_{2}-s}(z)-v_{\tau_{2}-s}(x)$ has been handled in (\ref{eq:Key2}).
We will thus focus on the following term: 
\begin{equation}
{\mathcal{E}}_{x,\tau_{1},\tau_{2}}=\int_{\mathbb{R}}\int_{0}^{\tau_{1}}\frac{|\partial_{x\cdots x}^{k}p_{1}(w)|}{s^{1+\sigma_{1}}}\big|\int_{x}^{x+\sqrt{s}w}(v_{\tau_{2}-s}(x)-v_{\tau_{2}}(x))dW(z)\big|dwds.\label{eq:defEcal}
\end{equation}
In order to bound ${\mathcal{E}}_{x,\tau_{1},\tau_{2}}$ we observe
that the increment $v_{\tau_{2}-s}(x)-v_{\tau_{2}}(x)$ does not depend
on the variable $z$. Therefore we have 
\begin{equation}\label{eq:Ecal}
{\mathcal{E}}_{x,\tau_{1},\tau_{2}}=\int_{\mathbb{R}}\int_{0}^{\tau_{1}}\frac{|\partial_{x\cdots x}^{k}p_{1}(w)|}{s^{1+\sigma_{1}}}|v_{\tau_{2}-s}(x)-v_{\tau_{2}}(x)|\cdot|W^{1}(x,x+\sqrt{s}w)|dwds \, . 
\end{equation}
Next we bound the increments of $v$ resorting to the bound on $\llbracket v\rrbracket_{\beta/2,\beta}$
ensured by Definition~\ref{def:controlled-process}, together with
the definition (\ref{eq:GenHolNorm}) of $\llbracket\,\cdot\,\rrbracket_{\beta/2,\beta}$.
This yields 
\begin{equation}
\big|v_{\tau_{2}-s}(x)-v_{\tau_{2}}(x)\big|\leqslant a^{\beta/2}E(a,\tau_{2})\Theta\cdot s^{\beta/2},\label{eq:bound-increment-v}
\end{equation}
The increments of $W$ can be estimated thanks to the fact that $\kappa_{\alpha,\chi}({\bf W})$
(see equation (\ref{f11})) is a finite quantity. We get 
\begin{equation}
\big|W^{1}(x,x+\sqrt{s}w)\big|\leqslant\kappa\cdot(a+\sqrt{T}|w|)^{\chi}\cdot s^{\alpha/2}\cdot|w|^{\alpha}\leqslant C_{\chi,T}\cdot\kappa\cdot P_{u}(w)\cdot a^{\chi}s^{\alpha/2},\label{eq:bound-increment-W}
\end{equation}
where we recall that $P_{u}(w)$ designates any polynomial in the
$w$ variable. Now plugging (\ref{eq:bound-increment-v}) and (\ref{eq:bound-increment-W})
into (\ref{eq:Ecal}) we obtain 
\[
{\mathcal{E}}_{x,\tau_{1},\tau_{2}}\leqslant C\cdot\kappa a^{\chi+\beta/2}E(a,\tau_{2})\Theta\cdot\int_{0}^{\tau_{1}}s^{\frac{\alpha+\beta}{2}-1-\sigma_{1}}ds=C\cdot\kappa a^{\chi+\beta/2}E(a,\tau_{2})\Theta\cdot\tau_{1}^{\frac{\alpha+\beta}{2}-\sigma_{1}}.
\]
Then we trivially bound $\tau_{1}^{\frac{\alpha+\beta}{2}}$ by $C\tau_{1}^{\beta}$
(recall that $\tau_{1}$ will be chosen as a small constant) and $a^{\beta/2}$
by $a^{\beta/2}+t^{-\beta/2}$. We end up with 
\[
{\mathcal{E}}_{x,\tau_{1},\tau_{2}}\leqslant C\cdot\kappa E(a,\tau_{2})\Theta\cdot\tau_{1}^{\beta-\sigma_{1}}a^{\chi}\left(a^{\beta/2}+t^{-\beta/2}\right),
\]
which means that we have achieved the bound (\ref{eq:ExtKey2}) for
${\mathcal{E}}_{x,\tau_{1},\tau_{2}}$. Our claim (\ref{eq:ExtKey2})
is thus obtained thanks to the considerations at the beginning of
our proof.
\end{proof}

\subsection{\label{sec:EstMV} Estimating the norm of $\mathcal{M}\mathcal{V}$}

Recall that our main object of interest is the transformation $\mathcal{M}\mathcal{V}$
introduced in Section \ref{Mtransformation}. For convenience recall
that this transformation is defined in (\ref{eq:def-trsf-M}) as 
\begin{equation}
\mathcal{M}\mathcal{V}\triangleq[(t,x)\mapsto((\mathcal{M}\mathcal{V})_{t}(x),\partial_{W}(\mathcal{M}\mathcal{V})_{t}(x))],\label{eq:recall-def-calM}
\end{equation}
where the first term $(\mathcal{M}v)_{t}(x)$ is defined by (\ref{e4})
and where the Gubinelli derivative $\partial_{W}(\mathcal{M}\mathcal{V})_{t}$
is defined in (\ref{eq:deriv-Gubi-M}). One of our main steps toward
the proof of Theorem \ref{thm:FPT} is the following contraction property
for the norm $\Theta^{\theta,\lambda}(\mathcal{M}\mathcal{V})$.
\begin{prop}
\label{prop:EstM} Consider a set of parameters $\alpha,\beta,\gamma,\chi,\theta,\lambda$
satisfying (\ref{f1}). Let $\mathcal{M}$ be the map whose definition
is recalled in (\ref{eq:recall-def-calM}). Then $\mathcal{M}$ is
a well-defined linear transformation on the Banach space $\mathcal{B}^{\theta,\lambda}$
introduced in Definition \ref{def:controlled-process}. It satisfies
the following estimate: 
\begin{equation}
\Theta^{\theta,\lambda}(\mathcal{M}\mathcal{V})\leqslant C\cdot(1+\kappa_{\alpha,\chi}({\bf W}))\cdot\lambda^{-\frac{\alpha-\beta}{4}}\cdot\Theta^{\theta,\lambda}(\mathcal{V}),\label{eq:estim-ThetaMV}
\end{equation}
where $C$ is a constant depending only on all the parameters $\alpha,\beta,\gamma,\chi,\theta,\lambda$,
but not on ${\bf W}$ and $\mathcal{V}$. In particular, by choosing
$\lambda$ to be large enough, we can ensure that 
\begin{equation}
\Theta^{\theta,\lambda}(\mathcal{M}\mathcal{V})\leqslant\frac{1}{2}\Theta^{\theta,\lambda}(\mathcal{V}).\label{eq:normeMVnormev}
\end{equation}
Therefore $\mathcal{M}$ is a contraction on $\mathcal{B}^{\theta,\lambda}$.
\end{prop}

Going back to the Definition \ref{def:controlled-process} of $\mathcal{B}^{\theta,\lambda}$
and its norm $\Theta^{\theta,\lambda}$, the estimate for $\mathcal{M}$
can be split in four main terms (i) a uniform estimate for $\mathcal{M}\mathcal{V}$,
(ii) an estimate for the time fluctuations of $\mathcal{M}\mathcal{V}$,
(iii) a bound on the spatial fluctuations of $\mathcal{M}\mathcal{V}$
and (iv) a control on the remainder $\mathcal{R}^{\mathcal{M}\mathcal{V}}$.
We develop these four estimates in the following sections. At this
point it is helpful to recall the following change of variables, which
will be extensively used in the sequel and is valid for $s>0$, $x,y\in\mathbb{R}$,
with $y=x+\sqrt{s}w$, and $k\geq1$: 
\begin{equation}
\left.\partial_{x\cdots x}^{k}p_{s}(x-y)\right|_{y=x+\sqrt{s}w}=C_{\sigma}\cdot s^{-\frac{k+1}{2}}\partial_{x\cdots x}^{k}p_{1}(w).\label{eq:deriv-heat-kernel}
\end{equation}
Unless otherwise stated, we always shorten our notation as in (\ref{eq: kappa-Theta-E})
in order to avoid lengthy expressions.

\subsubsection{The uniform estimate for $\mathcal{M}\mathcal{V}$}

We begin the proof of Proposition \ref{prop:EstM} with the following
simple uniform estimate for $(\mathcal{M}\mathcal{V})_{t}(x).$
\begin{lem}
\label{lem:UniformEst} Under the conditions of Proposition \ref{prop:EstM},
consider $t\in[0,T]$ and $x\in[-a,a]$. Then we have 
\begin{equation}
\big|(\mathcal{M}\mathcal{V})_{t}(x)\big|\leqslant C\cdot\kappa\Theta E(a,t)\cdot\lambda^{-\frac{\alpha-\beta}{4}}.\label{eq:unif-estiM}
\end{equation}
\end{lem}

\begin{proof}
Starting from the expression (\ref{e4}) for $\mathcal{M}v$, writing
$y=x+\sqrt{s}w$ and invoking relation~(\ref{eq:deriv-heat-kernel})
we get 
\begin{equation}
(\mathcal{M}\mathcal{V})_{t}(x)=\int_{0}^{t}\int_{\mathbb{R}}\frac{\partial_{xx}^{2}p_{1}(w)}{s}\left(\int_{x}^{x+\sqrt{s}w}v_{t-s}(z)dW(z)\right)dw\,ds.\label{eq:MV-with-A}
\end{equation}
Our claim (\ref{eq:unif-estiM}) is thus a direct consequence of Lemma
\ref{lem:KeyLem} (i) with $k=2$, $\gamma_1=\gamma_2=0$ and $\tau_{1}=\tau_{2}=t$.
\end{proof}

\subsubsection{\label{seq:TimVarMV}The time variation estimate for $\mathcal{M}\mathcal{V}$}

In this section we investigate the time fluctuations of $\mathcal{M}\mathcal{V}$,
which is another step toward the proof of Proposition \ref{prop:EstM}.
\begin{lem}
\label{lem:TimEst}We work under the conditions of Proposition \ref{prop:EstM}.
In particular, recall that $\alpha,\beta$ satisfy (\ref{f1}).
Then for any elements $x\in[-a,a]$ and $t_{1},t_{2}\in[0,T]$ such
that $t_{1}<t_{2}$, we have 
\begin{equation}
\big|(\mathcal{M}\mathcal{V})_{t_{2}}(x)-(\mathcal{M}\mathcal{V})_{t_{1}}(x)\big|\leqslant C\kappa\Theta E(a,t_{2})a^{\beta/2}\lambda^{-\frac{\alpha-\beta}{4}}|t_{2}-t_{1}|^{\beta/2}.\label{eq:estim-time-flucM}
\end{equation}
\end{lem}

\begin{proof} Invoking equation (\ref{e4}), it is readily checked
that the time increments of $\mathcal{M}\mathcal{V}$ can be written
as 
\begin{multline*}
(\mathcal{M}\mathcal{V})_{t_{2}}(x)-(\mathcal{M}\mathcal{V})_{t_{1}}(x)=\int_{t_{1}}^{t_{2}}\int_{\mathbb{R}}\partial_{xx}^{2}p_{t_{2}-s}(x-y)\int_{x}^{y}v_{s}(z)dW(z)dyds\\
+\int_{0}^{t_{1}}\int_{\mathbb{R}}\left(\partial_{xx}^{2}p_{t_{2}-s}(x-y)-\partial_{xx}^{2}p_{t_{1}-s}(x-y)\right)\int_{x}^{y}v_{s}(z)dW(z)dyds.
\end{multline*}
Hence, writing the increment $\partial_{xx}^{2}p_{t_{2}-s}(x-y)-\partial_{xx}^{2}p_{t_{1}-s}(x-y)$
in terms of a time derivative and owing to the fact that $p$ satisfies
$\partial_{t}p=\frac{\sigma^{2}}{2}\partial_{xx}^{2}p$, we get a
decomposition of the form 
\begin{equation}
(\mathcal{M}\mathcal{V})_{t_{2}}(x)-(\mathcal{M}\mathcal{V})_{t_{1}}(x)=\frac{\sigma^{2}}{2}\mathcal{T}_{1}+\mathcal{T}_{2},\label{eq:decompT1T2}
\end{equation}
where the terms $\mathcal{T}_{1}$ and $\mathcal{T}_{2}$ are respectively
defined by 
\begin{eqnarray}
\mathcal{T}_{1}&=&
C_{
\sigma}\int_{0}^{t_{1}}ds\int_{\mathbb{R}}dy\int_{t_{1}-s}^{t_{2}-s}\partial_{x}^{4}p_{u}(x-y)du\int_{x}^{y}v_{s}(z)dW(z) \notag \\
\mathcal{T}_{2}&=&
\int_{t_{1}}^{t_{2}}\int_{\mathbb{R}}\partial_{xx}^{2}p_{t_{2}-s}(x-y)\int_{x}^{y}v_{s}(z)dW(z)dyds, \label{eq:cal-T2}
\end{eqnarray}
where we write $\partial_{x}^{4}p$ instead of $\partial_{xxxx}^{4}p$
for notational sake. We now estimate $\mathcal{T}_{1}$ and $\mathcal{T}_{2}$
separately.

\vspace{2mm}
\noindent
\textit{Estimation of} $\mathcal{T}_{1}$: By a change of variables
$w=(y-x)/\sqrt{u},$ and owing to relation (\ref{eq:deriv-heat-kernel})
we can write $\mathcal{T}_{1}$ as 
\begin{equation}
\mathcal{T}_{1}=C_{\sigma}\int_{\mathbb{R}}\partial_{x}^{4}p_{1}(w)dw\int_{0}^{t_{1}}ds\int_{t_{1}-s}^{t_{2}-s}u^{-2}du\int_{x}^{x+\sqrt{u}w}v_{s}(z)dW(z).\label{eq:cal-T1}
\end{equation}
We first look at the triple integral inside the $w$-integral, that
is 
\begin{equation}
\mathcal{T}_{1,x,w,t_{1},t_{2}}=\int_{0}^{t_{1}}ds\int_{t_{1}-s}^{t_{2}-s}u^{-2}du\int_{x}^{x+\sqrt{u}w}v_{s}(z)dW(z).\label{eq:cal-T1xwt1t2}
\end{equation}
In order to upper bound $\mathcal{T}_{1,x,w,t_{1},t_{2}}$ we perform
the change of variable $u=t_{1}+r-s$ so that we get 
\[
\mathcal{T}_{1,x,w,t_{1},t_{2}}=\int_{0}^{t_{1}}ds\int_{0}^{t_{2}-t_{1}}\frac{dr}{(t_{1}+r-s)^{2}}\int_{x}^{x+\sqrt{t_{1}+r-s}w}v_{s}(z)dW(z).
\]
Then switching the order of integration in $r$ and $s$ and setting
$\rho=t_{1}+r-s$ we obtain 
\[
\mathcal{T}_{1,x,w,t_{1},t_{2}}=\int_{0}^{t_{2}-t_{1}}dr\int_{r}^{t_{1}+r}\frac{d\rho}{\rho^{1+\beta/2+1-\beta/2}}\int_{x}^{x+\sqrt{\rho}w}v_{t_{1}+r-\rho}(z)dW(z).
\]
This can be easily bounded as 
\begin{equation}
\Big|\mathcal{T}_{1,x,w,t_{1},t_{2}}\big|\leqslant\int_{0}^{t_{2}-t_{1}}r^{\beta/2-1}dr\int_{0}^{t_{1}+r}\frac{d\rho}{\rho^{1+\beta/2}}\left|\int_{x}^{x+\sqrt{\rho}w}v_{t_{1}+r-\rho}(z)dW(z)\right|.\label{eq:estim-cal-T1xwt1t2}
\end{equation}
Therefore, plugging (\ref{eq:estim-cal-T1xwt1t2}) into (\ref{eq:cal-T1xwt1t2})
and then (\ref{eq:cal-T1}), we end up with 
\begin{align*}
\mathcal{T}_{1} & \leqslant\int_{\mathbb{R}}\big|\partial_{x}^{4}p_{1}(w)\big|dw\int_{0}^{t_{2}-t_{1}}r^{\beta/2-1}dr\int_{0}^{t_{1}+r}\frac{d\rho}{\rho^{1+\beta/2}}\Big|\int_{x}^{x+\sqrt{\rho}w}v_{t_{1}+r-\rho}(z)dW(z)\Big|\\
 & =\int_{0}^{t_{2}-t_{1}}r^{\beta/2-1}dr\cdot\int_{\mathbb{R}}\big|\partial_{x}^{4}p_{1}(w)\big|dw\int_{0}^{t_{1}+r}\frac{d\rho}{\rho^{1+\beta/2}}\Big|\int_{x}^{x+\sqrt{\rho}w}v_{t_{1}+r-\rho}(z)dW(z)\Big|.
\end{align*}
We can now resort to Lemma \ref{lem:KeyLem} (i) with $\gamma_{1}=\gamma_{2}=\beta/2$.
This yields 
\begin{equation}
\mathcal{T}_{1}\leqslant C\kappa\Theta E(a,t_{2})a^{\beta/2}\lambda^{\gamma+\frac{\beta-\alpha}{2}}|t_{2}-t_{1}|^{\beta/2}=C\kappa\Theta E(a,t_{2})a^{\beta/2}\lambda^{-\frac{\alpha-\beta}{4}}|t_{2}-t_{1}|^{\beta/2},\label{eq:uppbdT1}
\end{equation}
where in the second identity we have used the fact that $\gamma=\frac{\alpha-\beta}{4}$.

 \vspace{2mm}
\noindent\textit{Estimation of} $\mathcal{T}_{2}$: The upper bound for $\mathcal{T}_{2}$
is obtained similarly to $\mathcal{T}_{1}$. Namely starting from
expression (\ref{eq:cal-T2}) some elementary change of variables
in the space and time variables yield 
\[
\mathcal{T}_{2}=\int_{0}^{t_{2}-t_{1}}\int_{\mathbb{R}}\partial_{xx}^{2}p_{r}(x-y)\int_{x}^{y}v_{t_{2}-r}(z)dW(z)dydr.
\]
Therefore setting $y=x+\sqrt{r}w$ and resorting to (\ref{eq:deriv-heat-kernel})
in order to replace $p_{r}$ by $p_{1}$, we get 
\[
\mathcal{T}_{2}=C_{\sigma}\int_{0}^{t_{2}-t_{1}}\frac{dr}{r}\int_{\mathbb{R}}\partial_{xx}^{2}p_{1}(w)\int_{x}^{x+\sqrt{r}w}v_{t_{2}-r}(z)dW(z)dw.
\]
The above identity enables the application of (\ref{eq:Key1}), where
we choose $\gamma_{1}=\gamma_{2}=0$. We get 
\begin{equation}
\mathcal{T}_{2}\leqslant C\kappa\Theta E(a,t_{2})a^{\beta/2}\lambda^{-\frac{\alpha-\beta}{4}}|t_{2}-t_{1}|^{\beta/2}.\label{eq:uppbdT2}
\end{equation}
Hence reporting \eqref{eq:uppbdT1} and \eqref{eq:uppbdT2} into the decomposition \eqref{eq:decompT1T2}, the proof of 
our claim \eqref{eq:unif-estiM} is easily achieved. 
\end{proof}

\subsubsection{\label{sec:SpVarEstRP} The space variation estimate for $\mathcal{M}\mathcal{V}$}

This section is devoted to the third ingredient in our global strategy,
namely the upper bound on the $\beta$-H\"older norm for the spatial
increments of $\mathcal{M}\mathcal{V}$. 
In order to ease notation, we will often use the convention
\begin{equation}\label{eq:notation-spacediff}
f(x,x')\equiv f(x)-f(x'),
\end{equation}
valid for any function $f$ of a spatial variable $x$ in $\mathbb{R}$ or 
$\delta\mathbb{Z}$. This notation will prevail for the remainder of the article.
Our main aim is to prove
the lemma below.
\begin{lem}
\label{lem:SpaEst} Let the assumptions and notations of Proposition
\ref{prop:EstM} prevail. Then for any $t\in(0,T]$ and $x,x'\in[-a,a],$
we have 
\begin{equation}
\big|(\mathcal{M}\mathcal{V})_{t}(x,x')\big|
\leqslant 
C\cdot\kappa\Theta E(a,t)\cdot\lambda^{-\frac{\alpha-\beta}{4}}a^{\beta/2}|x'-x|^{\beta}.\label{g1}
\end{equation}
\end{lem}

\begin{proof} Before starting our discussion, recall our identity (\ref{eq:integration-parts}),
for which we had chosen an arbitrary $a\in\mathbb{R}$. Differentiating
this relation with respect to $x$ and taking limits as $\eta\to0$
in the rough path sense, we get that the integral 
\begin{equation}
\int_{0}^{t}\int_{\mathbb{R}}\partial_{xx}^{2}p_{s}(x-y)\int_{a}^{y}v_{t-s}(z)dW(z)dyds\label{eq:integral_indep_a}
\end{equation}
does not depend on $a$. This fact will be used without further mention
in the remainder of the proof.

As an application of \eqref{eq:integral_indep_a},  let $t\in(0,T]$ and $x,x'\in[-a,a].$ Starting from
the expression~(\ref{e4}) for $\mathcal{M}\mathcal{V}$ and owing
to the fact that the lower bound $x$ can be chosen arbitrarily in
$\int_{x}^{y}v_{s}(z)dW(z)$, it is readily checked that 
\begin{equation}
(\mathcal{M}\mathcal{V})_{t}(x,x')
=
\int_{0}^{t}\int_{\mathbb{R}}\left(\partial_{xx}^{2}p_{s}(x'-y)-\partial_{xx}^{2}p_{s}(x-y)\right)\int_{x}^{y}v_{t-s}(z)dW(z)dyds.
\label{eq:space-increments}
\end{equation}
We now divide our analysis in two cases
according to the respective value of $x'-x$ and $t$.

\noindent \textit{Case }(i): $|x'-x|^{2}\leqslant t.$ With identity
(\ref{eq:space-increments}) in mind, let us split the interval $[0,t]$
into $[0,|x'-x|]\cup[|x'-x|,t]$ in order to get the following decomposition
\begin{equation}
(\mathcal{M}\mathcal{V})_{t}(x')-(\mathcal{M}\mathcal{V})_{t}(x)=\mathcal{I}_{1}(x')-\mathcal{I}_{1}(x)+\mathcal{I}_{2},\label{g2}
\end{equation}
where for $\xi\in\mathbb{R}$ we have set 
\begin{equation}
\mathcal{I}_{1}(\xi)=\int_{0}^{|x'-x|^{2}}\int_{\mathbb{R}}\left(\partial_{xx}^{2}p_{s}(\xi-y)\int_{\xi}^{y}v_{t-s}(z)dW(z)\right)dyds,\label{eq:cal-I1}
\end{equation}
and where the term $\mathcal{I}_{2}$ is given by 
\begin{equation}
\mathcal{I}_{2}=\int_{|x'-x|^{2}}^{t}ds\int_{\mathbb{R}}dy\int_{x}^{x'}\partial_{x}^{3}p_{s}(u-y)du\int_{x}^{y}v_{t-s}(z)dW(z),\label{eq:cal-I2}
\end{equation}
where we use again $\partial_{x}^{3}p$ instead of $\partial_{xxx}^{3}p$
for notational sake. We proceed to estimate the term $\mathcal{I}_{1}(\xi)$
in (\ref{eq:cal-I1}). To this aim we resort to the same type of change
of variables as in Lemma \ref{lem:TimEst} and invoke the identity
(\ref{eq:deriv-heat-kernel}) once more. We get 
\[
|\mathcal{I}_{1}(\xi)|\leqslant\int_{0}^{|x'-x|^{2}}\frac{ds}{s}\int_{\mathbb{R}}\big|\partial_{xx}^{2}p_{1}(w)\big|dw\left|\int_{\xi}^{\xi+\sqrt{s}w}v_{t-s}(z)dW(z)\right|.
\]
We are thus in position to apply Lemma \ref{lem:KeyLem} (i) with
$\tau_{1}=|x'-x|^{2}$, $\gamma_{2}=\beta/2$ and $\gamma_{1}=0$.
Also we recall that we have chosen $\gamma=\frac{\alpha-\beta}{4}$.
This implies 
\begin{equation}
|\mathcal{I}_{1}(\xi)|\leqslant C\cdot\kappa\Theta E(a,t)\cdot\lambda^{-\frac{\alpha-\beta}{4}}a^{\beta/2}|x'-x|^{\beta}.\label{g3}
\end{equation}
Let us turn to an upper bound on the term $\mathcal{I}_{2}$ in (\ref{eq:cal-I2}).
Switching the order of integration and thanks to (\ref{eq:deriv-heat-kernel}),
one can write 
\begin{align*}
|\mathcal{I}_{2}| & \leqslant\int_{x}^{x'}du\int_{|x'-x|^{2}}^{t}\frac{ds}{s^{1/2-\beta/2+1+\beta/2}}\int_{\mathbb{R}}\big|\partial_{x}^{3}p_{1}(w)\big|dw\left|\int_{u}^{u+\sqrt{s}w}v_{t-s}(z)dW(z)\right|\\
 & \leqslant|x'-x|^{\beta-1}\int_{x}^{x'}du\int_{0}^{t}\frac{ds}{s^{1+\beta/2}}\int_{\mathbb{R}}\big|\partial_{x}^{3}p_{1}(w)\big|dw\left|\int_{u}^{u+\sqrt{s}w}v_{t-s}(z)dW(z)\right|,
\end{align*}
where we have enlarged the term $\frac{1}{s^{1/2-\beta/2}}$ to $|x'-x|^{\beta-1}$
and then enlarged the $s$-integral over the region $[0,t].$ By applying
Lemma \ref{lem:KeyLem} (i) to the inner $(s,w,z)$-integral with
$\gamma_{1}=\gamma_{2}=\beta/2$, we get 
\begin{equation}
|\mathcal{I}_{2}|\leqslant C\cdot\kappa\Theta E(a,t)\cdot\lambda^{-\frac{\alpha-\beta}{4}}a^{\beta/2}|x'-x|^{\beta}.\label{g4}
\end{equation}
Reporting (\ref{g3}) and (\ref{g4}) into (\ref{g2}) we have 
\[
\big|\mathcal{M}\mathcal{V}_{t}(x,x')\big|\leqslant C\cdot\kappa\Theta E(a,t)\cdot\lambda^{-\frac{\alpha-\beta}{4}}a^{\beta/2}|x'-x|^{\beta}.
\]
Therefore our claim (\ref{g1}) is proved whenever $|x'-x|^{2}\leqslant t$.

\noindent \textit{Case }(ii): $|x'-x|^{2}>t.$ This case is essentially
contained in the above $\mathcal{I}_{1}(\xi)$-estimate. Indeed, starting
directly from (\ref{eq:space-increments}) we can write 
\[
\mathcal{MV}_{t}(x,x')
=
\mathcal{I}'_{1}(x')-\mathcal{I}'_{1}(x),
\]
where 
\[
\mathcal{I}_{1}'(\xi)\triangleq\int_{0}^{t}\int_{\mathbb{R}}\partial_{xx}^{2}p_{s}(\xi-y)\int_{\xi}^{y}v_{t-s}(z)dW(z)dyds,\ \ \ \xi\in\mathbb{R}.
\]
The only difference between $\mathcal{I}_{1}'(\xi)$ and $\mathcal{I}_{1}(\xi)$
introduced in Case (i) is that the region for the $s$-integral becomes
$[0,t]$ rather than $[0,|x'-x|^{2}]$. By exactly the same argument
as before, the estimate (\ref{g3}) becomes 
\[
|\mathcal{I}_{1}'(\xi)|\leqslant C\cdot\kappa\Theta E(a,t)\cdot\lambda^{-\frac{\alpha-\beta}{4}}a^{\beta/2}t^{\beta/2}.
\]
Since $|x'-x|^{2}>t$, it follows that 
\begin{align*}
\big|\mathcal{MV}_{t}(x,x')\big| & \leqslant\big|\mathcal{I}_{1}'(x')\big|+\big|\mathcal{I}_{1}'(x)\big|\\
 & \leqslant C\cdot\kappa\Theta E(a,t)\cdot\lambda^{-\frac{\alpha-\beta}{4}}a^{\beta/2}|x'-x|^{\beta}.
\end{align*}
This achieves (\ref{g1}) for $|x'-x|^{2}>t$ and finishes our proof.
\end{proof}

\subsubsection{\label{subsec:RemEst}The remainder estimate for $\mathcal{M}\mathcal{V}$}

Recall that the remainder of a controlled process is introduced in
(\ref{f2}). Moreover, we defined the Gubinelli derivative of $\mathcal{M}\mathcal{V}$
as $-2v/\sigma^{2}$ in~(\ref{eq:deriv-Gubi-M}). Gathering those
two relations and recalling notation~\eqref{eq:notation-spacediff}, we get that the remainder of $(\mathcal{M}\mathcal{V})_{t}$
is given as 
\begin{equation}
\mathcal{R}_{t}^{\mathcal{M}\mathcal{V}}(x,x')=
\mathcal{M}\mathcal{V}_{t}(x,x')+\frac{2}{\sigma^{2}}v_{t}(x)W^{1}(x,x').\label{eq:ReMV}
\end{equation}
We will now prove that $\mathcal{R}_{t}^{\mathcal{M}\mathcal{V}}$
is indeed a $2\beta$-H\"older increment.
\begin{lem}
\label{lem:RemEst} Let $\alpha,\beta,\chi,\gamma,\theta,\lambda$
be parameters such that (\ref{f1}) is fulfilled. Recall that the
functions $E$ and $Q$ are introduced in Definition \ref{def:weights},
and that the norm $\Theta$ is given in Definition~\ref{def:controlled-process}.
Then the remainder $\mathcal{R}^{\mathcal{M}\mathcal{V}}$ defined
by (\ref{eq:ReMV}) satisfies 
\begin{equation}
\big|\mathcal{R}_{t}^{\mathcal{M}\mathcal{V}}(x,x')\big|\leqslant C\cdot\kappa E(a,t)\Theta\cdot Q(a,t)\cdot|x'-x|^{2\beta},\label{eq:estimRMV}
\end{equation}
for all $t\in(0,T]$ and $x,x'\in[-a,a]$.
\end{lem}

A crucial point for proving Lemma \ref{lem:RemEst} is to make use
of the decomposition in the following lemma. This decomposition, together with
the estimates developed in Lemma \ref{lem:R0Est} and Lemma \ref{lem:SGEst}
below, will explain why the Gubinelli derivative of $(\mathcal{M}\mathcal{V})_{t}$
should be defined to be $-2v_{t}/\sigma^{2}$.
\begin{lem}
\label{lem:RemDec}As in Lemma \ref{lem:RemEst}, we consider the
remainder $\mathcal{R}^{\mathcal{M}\mathcal{V}}$ given by (\ref{eq:ReMV}).
Then for all $t\in[0,T]$ and $x,x'\in[-a,a]$ we have 
\begin{equation}
\mathcal{R}_{t}^{\mathcal{M}\mathcal{V}}(x,x')=\mathcal{R}_{t}^{0}(x,x')+\frac{2}{\sigma^{2}}v_{t}(x)\left(P_{t}W(x')-P_{t}W(x)\right),\label{eq:ReMVbis}
\end{equation}
where the increment $\mathcal{R}^{0}$ is defined by 
\begin{equation}
\mathcal{R}_{t}^{0}(x,x')=\int_{0}^{t}\int_{\mathbb{R}}\left(\partial_{xx}^{2}p_{s}(x'-y)-\partial_{xx}^{2}p_{s}(x-y)\right)
\left(\int_{x}^{y}\left(v_{t-s}(z)-v_{t}(x)\right)dW(z)\right)dyds\label{eq:ReMV0}
\end{equation}
and where $P_{t}$ has been defined after (\ref{heat_kernel}) as
the heat semigroup in $\mathbb{R}$ with variance $\sigma^{2}$.
\end{lem}

\begin{proof} Recall that starting from the expression (\ref{e4})
for $\mathcal{M}\mathcal{V}$ we already obtained the space increment
(\ref{eq:space-increments}). We now add and substract $v_{t}(x)$
in the integral on the right hand side of~(\ref{eq:space-increments}).
With the definition (\ref{eq:ReMV0}) of $\mathcal{R}^{0}$ in mind,
this yields 
\begin{equation}
\mathcal{MV}_{t}(x,x')=\mathcal{R}_{t}^{0}(x,x')
+v_{t}(x)\int_{0}^{t}\int_{\mathbb{R}}\left(\partial_{xx}^{2}p_{s}(x'-y)-\partial_{xx}^{2}p_{s}(x-y)\right)\,W^{1}(x,y)\,dyds.\label{eq:ReMV01}
\end{equation}
Furthermore, the decaying properties of the heat kernel $p$ imply
that $\int_{\mathbb{R}}\partial_{xx}^{2}p_{s}(x'-y)dy=0$ for all
$x'\in\mathbb{R}$. Therefore one can insert $\pm W(x')$ in (\ref{eq:ReMV01})
and one gets 
\begin{equation}
\mathcal{MV}_{t}(x,x')=\mathcal{R}_{t}^{0}(x,x')+v_{t}(x)\left(\mathcal{B}_{t}(x')-\mathcal{B}_{t}(x)\right),\label{eq:calB}
\end{equation}
where the term $\mathcal{B}_{t}(\xi)$ is defined for $\xi\in\mathbb{R}$
by 
\[
\mathcal{B}_{t}(\xi)=\int_{0}^{t}\int_{\mathbb{R}}\partial_{xx}^{2}p_{s}(\xi-y)\,W^{1}(\xi,y)dyds.
\]
Next we simplify the expression for $\mathcal{B}_{t}(\xi)$ by involving
the relation $\partial_{xx}^{2}p_{s}(y)=\frac{2}{\sigma^{2}}\partial_{s}p_{s}(y)$.
Interchanging integral and differentiation, we obtain 
\begin{equation}
\mathcal{B}_{t}(\xi)=\frac{2}{\sigma^{2}}\int_{0}^{t}\partial_{s}\left(\int_{\mathbb{R}}p_{s}(\xi-y)\,W^{1}(\xi,y)dy\right)ds=\frac{2}{\sigma^{2}}\int_{0}^{t}\partial_{s}\left(P_{s}W(\xi)-W(\xi)\right)ds,\label{eq:calB-bis}
\end{equation}
where we recall that $P_{s}$ designates the heat semigroup and where
we have used the fact that $\int_{\mathbb{R}}p_{s}(z)dz=1$. We now
simply evaluate the time integral in (\ref{eq:calB-bis}) in order
to get 
\[
\mathcal{B}_{t}(\xi)=\frac{2}{\sigma^{2}}\left(P_{t}W(\xi)-W(\xi)\right).
\]
Plugging this identity in (\ref{eq:calB}) we end up with 
\begin{eqnarray*}
  \mathcal{MV}_{t}(x,x')
  &=&\mathcal{R}_{t}^{0}(x,x')+\frac{2}{\sigma^{2}}v_{t}(x)\left(\left(P_{t}W(x')-W(x')\right)-\left(P_{t}W(x)-W(x)\right)\right),\\
 & =&\mathcal{R}_{t}^{0}(x,x')+\frac{2}{\sigma^{2}}v_{t}(x)\left(P_{t}W(x')-P_{t}W(x)\right)-\frac{2}{\sigma^{2}}v_{t}(x)\left(W(x')-W(x)\right),
\end{eqnarray*}
from which our claim (\ref{eq:ReMVbis}) is easily deduced. \end{proof}

In view of Lemma \ref{lem:RemDec}, the estimate of $\mathcal{R}^{(\mathcal{M}\mathcal{V})_{t}}(x,x')$
contains two ingredients: the $\mathcal{R}_{t}^{0}(x,x')$-estimate
and the heat semigroup variation estimate. We develop these two ingredients
separately, starting with an estimate for $\mathcal{R}_{t}^{0}(x,x')$.
\begin{lem}
\label{lem:R0Est}We consider the setting of Lemma \ref{lem:RemEst},
and let $\mathcal{R}^{0}$ the increment defined by~(\ref{eq:ReMV0}).
Then for any $t\in(0,T]$ and $x,x'\in[-a,a],$ we have 
\begin{equation}
\big|\mathcal{R}_{t}^{0}(x,x')\big|\leqslant C\cdot\kappa E(a,t)\Theta\cdot Q(a,t)\cdot|x'-x|^{2\beta}.\label{eq:estimate_R0t}
\end{equation}
\end{lem}

\begin{proof} As in the proof of Lemma \ref{lem:SpaEst} we divide
our discussion into two cases, according to the value of $|x'-x|$.

\vspace{2mm}
\noindent
\textit{Case (i):} $|x'-x|^{2}\leqslant t.$ In this case, we start
with the following decomposition: 
\begin{equation}
\mathcal{R}_{t}^{0}(x,x')=\mathcal{J}+\mathcal{K},\label{eq:R0eJpK}
\end{equation}
where the terms $\mathcal{J}$ and $\mathcal{K}$ are respectively
defined by 
\begin{align}
 & \mathcal{J}:=\int_{0}^{|x'-x|^{2}}ds\int_{\mathbb{R}}\left(\partial_{xx}^{2}p_{s}(x'-y)-\partial_{xx}^{2}p_{s}(x-y)\right)\left(\int_{x}^{y}(v_{t-s}(z)-v_{t}(x))dW(z)\right)dy\label{eq:calJ}\\
 & \mathcal{K}:=\int_{|x'-x|^{2}}^{t}ds\int_{\mathbb{R}}\left(\partial_{xx}^{2}p_{s}(x'-y)-\partial_{xx}^{2}p_{s}(x-y)\right)\left(\int_{x}^{y}(v_{t-s}(z)-v_{t}(x))dW(z)\right)dy.\label{eq:calK}
\end{align}
Next we further divide the integral $\mathcal{J}$ in (\ref{eq:calJ})
as 
\begin{equation}
\mathcal{J}=\mathcal{J}_{1}(x')-\mathcal{J}_{1}(x)+\mathcal{J}_{2},\label{eq:JeJ1J1pJ2}
\end{equation}
where the terms $\mathcal{J}_{1}(\xi)$ and $\mathcal{J}_{2}$ are
given by 
\begin{eqnarray}
\mathcal{J}_{1}(\xi) & = & \int_{0}^{|x'-x|^{2}}ds\int_{\mathbb{R}}\partial_{xx}^{2}p_{s}(\xi-y)\left(\int_{\xi}^{y}(v_{t-s}(z)-v_{t}(\xi))dW(z)\right)dy,\label{eq:calJ1}\\
\mathcal{J}_{2} & = & \int_{0}^{|x'-x|^{2}}ds\int_{\mathbb{R}}\partial_{xx}^{2}p_{s}(x'-y)\left((v_{t}(x')-v_{t}(x))\cdot W^{1}(x',y)\right)dy.\label{eq:calJ2}
\end{eqnarray}
We now proceed to estimate the terms in (\ref{eq:calJ})-(\ref{eq:calJ2}).
Let us start with a bound on the term $\mathcal{J}_{1}(\xi)$ defined
by (\ref{eq:calJ1}). To this aim, we first resort to the same kind
of change of variables as for (\ref{eq:cal-T1}). We get 
\[
|\mathcal{J}_{1}(\xi)|\leqslant\int_{0}^{|x'-x|^{2}}ds\int_{\mathbb{R}}\frac{|\partial_{xx}^{2}p_{1}(w)|}{s}\cdot\left|\int_{\xi}^{\xi+\sqrt{s}w}(v_{t-s}(z)-v_{t}(\xi))dW(z)\right|dw.
\]
We are now in a position to apply Lemma \ref{lem:ExtKey2} directly (with $\sigma_{1}=0$),
which yields 
\begin{equation}
|\mathcal{J}_{1}(\xi)|\leqslant C\cdot\kappa\Theta E(a,t)a^{\chi}\cdot Q(a,t)\cdot|x'-x|^{2\beta}.\label{eq:estimJ1}
\end{equation}
As far as the $\mathcal{J}_{2}$ integral is concerned, we perform
our usual change of variable in space. We obtain 
\begin{equation}
\mathcal{J}_{2}=\int_{0}^{|x'-x|^{2}}\frac{ds}{s}\int_{\mathbb{R}}\partial_{xx}^{2}p_{1}(w)\cdot(v_{t}(x')-v_{t}(x))\cdot W^{1}(x',x'+\sqrt{s}w)\,dw.\label{eq:J2Z}
\end{equation}
Then we estimate the right hand side of (\ref{eq:J2Z}) thanks to
the fact that $v$ verifies Definition~\ref{def:controlled-process}
and $\kappa_{\alpha,\chi}(w)$, as introduced in (\ref{f11}), is
finite. We obtain 
\begin{eqnarray}
|v_{t}(x')-v_{t}(x)| & \leqslant & E(a,t)\Theta a^{\beta/2}|x'-x|^{\beta},\label{eq:vVarZ}\\
\big|W^{1}(x',x'+\sqrt{s}w)\big| & \leqslant & C\kappa P_{u}(w)a^{\chi}s^{\alpha/2}.\label{eq:wVarZ}
\end{eqnarray}
By substituting (\ref{eq:vVarZ})-(\ref{eq:wVarZ}) into (\ref{eq:J2Z})
and performing the $w$-integral, we arrive at 
\begin{multline}
|\mathcal{J}_{2}|\leqslant C\cdot\kappa E(a,t)\Theta a^{\chi+\beta/2}|x'-x|^{\beta}\cdot\int_{0}^{|x'-x|^{2}}s^{\alpha/2-1}ds\\
\leqslant C\cdot\kappa E(a,t)\Theta a^{\chi+\beta/2}|x'-x|^{\beta}\cdot|x'-x|^{\alpha}\leqslant C\cdot\kappa\Theta E(a,t)\cdot Q(a,t)\cdot|x'-x|^{2\beta}.\label{eq:estimJ2}
\end{multline}

Next we estimate the $\mathcal{K}$-integral given by (\ref{eq:calK}).
We start by writing $\partial_{xx}^{2}p_{s}(x'-y)-\partial_{xx}^{2}p_{s}(x-y)$
in terms of a spatial derivative. This yields 
\[
\mathcal{K}=\int_{|x'-x|^{2}}^{t}ds\int_{\mathbb{R}}dy\int_{x}^{x'}\partial_{x}^{3}p_{s}(u-y)du\int_{u}^{y}(v_{t-s}(z)-v_{t}(x))dW(z)dy,
\]
where we have used the fact that integrals like \eqref{eq:integral_indep_a} do not depend on the parameter $a$. 
Then we insert $\pm v_{t}(u)$ in the integral above in order to get
the decomposition 
\begin{equation}
\mathcal{K}=\mathcal{K}_{1}+\mathcal{K}_{2},\label{eq:KeK1pK2}
\end{equation}
where $\mathcal{K}_{1}$ and $\mathcal{K}_{2}$ are defined by 
\begin{eqnarray}
\mathcal{K}_{1} & = & \int_{|x'-x|^{2}}^{t}ds\int_{\mathbb{R}}dy\int_{x}^{x'}\partial_{x}^{3}p_{s}(u-y)du\int_{u}^{y}(v_{t-s}(z)-v_{t}(u))dW(z)dy,\label{eq:K1}\\
\mathcal{K}_{2} & = & \int_{|x'-x|^{2}}^{t}ds\int_{\mathbb{R}}dy\int_{x}^{x'}\partial_{x}^{3}p_{s}(u-y)du(v_{t}(u)-v_{t}(x))W^{1}(u,y)dy.\label{eq:K2}
\end{eqnarray}
We now devote our efforts to upper bound $\mathcal{K}_{1}$ and $\mathcal{K}_{2}$.

In order to estimate the $\mathcal{K}_{1}$-integral, we
set again $y=u+\sqrt{s}w$. We obtain 
\[
\mathcal{K}_{1}=\int_{x}^{x'}du\int_{|x'-x|^{2}}^{t}\frac{ds}{s^{3/2}}\int_{\mathbb{R}}\partial_{x}^{3}p_{1}(w)\int_{u}^{u+\sqrt{s}w}(v_{t-s}(z)-v_{t}(u))dW(z)dw.
\]
Along the same lines as for the estimation of $\mathcal{I}_{2}$ in
the proof of Lemma \ref{lem:SpaEst}, we bound the term $s^{-(1/2-\beta)}$
by $|x'-x|^{\beta-1}$ and we enlarge the $s$-integral to the interval
$[0,t]$. We get 
\[
\mathcal{K}_{1}\leqslant|x'-x|^{1-2\beta}\cdot\int_{x}^{x'}du\int_{0}^{t}\frac{ds}{s^{1+\beta}}\int_{\mathbb{R}}\big|\partial_{x}^{3}p_{1}(w)\big|\cdot\left|\int_{u}^{u+\sqrt{s}w}(v_{t-s}(z)-v_{t}(u))dW(z)\right|dw.
\]
We are now able to apply Lemma \ref{lem:ExtKey2} with $\sigma_{1}=\beta$
and we obtain 
\begin{equation}
\mathcal{K}_{1}\leqslant C\cdot\kappa\Theta E(a,t)\cdot Q(a,t)\cdot|x'-x|^{2\beta}.\label{eq:estimK1}
\end{equation}
Let us bound the term $\mathcal{K}_{2}$ defined by (\ref{eq:K2}).
The change of variable $y=u+\sqrt{s}w$ together with (\ref{eq:deriv-heat-kernel})
yields 
\begin{equation}
\mathcal{K}_{2}=\int_{x}^{x'}du\int_{|x'-x|^{2}}^{t}\frac{ds}{s^{3/2}}\int_{\mathbb{R}}\partial_{x}^{3}p_{1}(w)(v_{t}(u)-v_{t}(x))W^{1}(u,u+\sqrt{s}w)dw.\label{eq:K2bis}
\end{equation}
In addition, owing to (\ref{f11}) and the fact that $v$ fulfills
the assumptions of Definition \ref{def:controlled-process}, we have
\[
|v_{t}(u)-v_{t}(x)|\leqslant E(a,t)\Theta\cdot a^{\beta/2}|x'-x|^{\beta}\quad\mbox{ and }\quad\big|W^{1}(u,u+\sqrt{s}w)\big|\leqslant C\kappa P_{u}(w)a^{\chi}s^{\alpha/2}.
\]
Plugging this information into (\ref{eq:K2bis}) we get 
\begin{equation}
\mathcal{K}_{2}\leqslant C\cdot\kappa\Theta E(a,t)\cdot a^{\chi+\beta/2}\cdot|x'-x|^{1+\beta}\cdot\int_{|x'-x|^{2}}^{t}s^{\alpha/2-3/2}ds.\label{eq:estimK2}
\end{equation}
Performing the $s$-integral and resorting to the relation $|x'-x|^{2}\leqslant t$,
we end up with 
\[
\mathcal{K}_{2}\leqslant C\cdot\kappa\Theta E(a,t)\cdot a^{\chi+\beta/2}\cdot|x'-x|^{\alpha+\beta}\leqslant C\cdot\kappa\Theta E(a,t)\cdot Q(a,t)\cdot|x'-x|^{2\beta}.
\]
Let us summarize our considerations so far: we plug (\ref{eq:estimJ1})
and (\ref{eq:estimJ2}) into the decomposition (\ref{eq:JeJ1J1pJ2})
for $\mathcal{J}$. We also gather (\ref{eq:estimK1}) and (\ref{eq:estimK2})
into the decomposition (\ref{eq:KeK1pK2}) for $\mathcal{K}$. Then
we report those estimates into the main decomposition (\ref{eq:R0eJpK}).
This achieves our claim (\ref{eq:estimate_R0t}) for $|x'-x|^{2}\leqslant t$.

\vspace{2mm}
\noindent
\textit{Case (ii):} $|x'-x|^{2}>t.$ This case is essentially contained
in the estimate of the term $\mathcal{J}$ defined by (\ref{eq:calJ}).
That is, going back to the expression (\ref{eq:ReMV0}) for $\mathcal{R}_{t}^{0}$
and using the fact that the expression \eqref{eq:integral_indep_a} does
not depend on $a$, one can write we have 
\begin{equation}
\mathcal{R}_{t}^{0}(x,x')=\mathcal{L}_{1}(x)+\mathcal{L}_{1}(x')+\mathcal{L}_{2},\label{eq:R0eL1L1pL2}
\end{equation}
where the quantities $\mathcal{L}_{1}(\xi)$ and $\mathcal{L}_{2}$
are given by 
\begin{eqnarray}
\mathcal{L}_{1}(\xi) & = & \int_{0}^{t}ds\int_{\mathbb{R}}\partial_{xx}^{2}p_{s}(\xi-y)\int_{\xi}^{y}(v_{t-s}(z)-v_{t}(\xi))dW(z)dy,\\
\mathcal{L}_{2} & = & \int_{0}^{t}ds\int_{\mathbb{R}}\partial_{xx}^{2}p_{s}(x'-y)\cdot\left(v_{t}(x')-v_{t}(x)\right)W^{1}(x',y)dy.\label{eq:L2}
\end{eqnarray}
Then invoking the same kind of arguments as for the term $\mathcal{I}_{1}(\xi)$
in (\ref{eq:calJ1}), plus the fact that $|x'-x|^{2}>t$, we obtain
\[
|\mathcal{L}_{1}(\xi)|\leqslant C\cdot\kappa\Theta E(a,t)\cdot Q(a,t)\cdot t^{\beta}\leqslant C\cdot\kappa\Theta E(a,t)\cdot Q(a,t)\cdot|x'-x|^{2\beta}.
\]
Furthermore, thanks to the space regularity of $v$ and owing to the
decay of the heat kernel $p$ (see (\ref{eq:vVarZ})-(\ref{eq:wVarZ})
for similar arguments) we get 
\begin{multline}
|\mathcal{L}_{2}|\leqslant C\cdot\kappa E(a,t)\Theta a^{\chi+\beta/2}|x'-x|^{\beta}\cdot\int_{0}^{t}s^{\alpha/2-1}ds\\
\leqslant C\cdot\kappa E(a,t)\Theta a^{\chi+\beta/2}|x'-x|^{\beta}\cdot t^{\alpha/2}\leqslant C\cdot\kappa E(a,t)\Theta\cdot Q(a,t)\cdot|x'-x|^{2\beta}.\label{eq:estimL2}
\end{multline}
Gathering (\ref{eq:L2}) and (\ref{eq:estimL2}) into (\ref{eq:R0eL1L1pL2}),
we have shown (\ref{eq:estimate_R0t}) for $|x'-x|^{2}>t$. Now the
proof of Lemma \ref{lem:R0Est} is complete. \end{proof}

\subsubsection{Estimation of the heat semigroup variation}

Remember that we have obtained the decomposition (\ref{eq:ReMVbis})
for $\mathcal{R}^{\mathcal{MV}}$. This decomposition involves two
main terms, namely $\mathcal{R}_{t}^{0}(x,x')$ and the increment
$P_{t}W(x')-P_{t}W(x)$. In this section we handle the latter increment.
Our main result is summarized in the following lemma.
\begin{lem}
\label{lem:SGEst} We stick to the notation of Lemma \ref{lem:RemEst}
and consider the double-sided Brownian motion $W$ of equation (\ref{eq:bm-in-br-environment}).
Then for any $t\in(0,T]$ and $x,x'\in[-a,a],$ we have 
\begin{equation}
\big|P_{t}W(x')-P_{t}W(x)\big|\leqslant C\cdot\kappa\cdot a^{\chi}\cdot t^{-\beta/2}\cdot|x'-x|^{2\beta}.\label{eq:estim-incr-PtW}
\end{equation}
\end{lem}

\begin{proof} As in Lemmas \ref{lem:SpaEst} and \ref{lem:R0Est},
we will separate the cases $|x'-x|^{2}\leqslant t$ and $|x'-x|^{2}>t$.

\vspace{2mm}
\noindent
\textit{Case (i)}: $|x'-x|^{2}\leqslant t$. We start from the from
the very definition of $P_{t}W(x)$ and write the increment $p_{t}(x'-y)-p_{t}(x-y)$
as a spatial derivative in order to get 
\[
P_{t}W(x')-P_{t}W(x)=\int_{\mathbb{R}}(p_{t}(x'-y)-p_{t}(x-y))W(y)dy=\int_{x}^{x'}du\int_{\mathbb{R}}\partial_{x}p_{t}(u-y)W(y)dy.
\]
Then perform the change of variable $y=u+\sqrt{t}w$, invoke relation
(\ref{eq:deriv-heat-kernel}) and use the decay properties of $\partial_{x}p_{1}$.
This yields 
\begin{eqnarray}
P_{t}W(x')-P_{t}W(x) & = & \frac{1}{\sqrt{t}}\int_{x}^{x'}du\int_{\mathbb{R}}\partial_{x}p_{1}(w)W(u+\sqrt{t}w) \, dw
\nonumber \\
 & = & \frac{1}{\sqrt{t}}\int_{x}^{x'}du\int_{\mathbb{R}}\partial_{x}p_{1}(w)W^{1}(u,u+\sqrt{t}w)dw.\label{eq:increm_PtW}
\end{eqnarray}
Note that in (\ref{eq:increm_PtW}) the variable $u$ is an element
of $[x,x']\subseteq[-a,a].$ Therefore, recalling the definition (\ref{f11})
of $\kappa$ one can bound the increment $W^{1}(u,u+\sqrt{t}w)$ as
\[
|W^{1}(u,u+\sqrt{t}w)|\leqslant\kappa\cdot(a+\sqrt{t}|w|)^{\chi}\cdot t^{\alpha/2}|w|^{\alpha}\leqslant C\kappa P_{u}(w)a^{\chi}t^{\alpha/2},
\]
where we recall that $P_{u}(w)$ designates a generic polynomial in
$w$. Plugging this inequality in~(\ref{eq:increm_PtW}) and performing
the integral with respect to the variable $w$, we end up with 
\begin{eqnarray}
\big|P_{t}W(x')-P_{t}W(x)\big| & \leqslant & C\cdot\kappa a^{\chi}t^{\alpha/2-1/2}|x'-x|\nonumber \\
 & = & C\cdot\kappa a^{\chi}|x'-x|^{2\beta}\cdot|x'-x|^{1-2\beta}t^{-\beta/2}t^{\frac{\alpha-\beta}{2}}\cdot\left(\frac{|x'-x|^{2}}{t}\right)^{1/2-\beta}\nonumber \\
 & \leqslant & C\cdot\kappa a^{\chi}t^{-\beta/2}\cdot|x'-x|^{2\beta},\label{eq:estim-incr-PtW-bis}
\end{eqnarray}
where we have resorted to the fact that $|x'-x|^{2}\leqslant t$ for
the last step. The estimate (\ref{eq:estim-incr-PtW-bis}) proves
our claim (\ref{eq:estim-incr-PtW}) when $|x'-x|^{2}\leqslant t$.

\vspace{2mm}
\noindent
\textit{Case (ii): $|x'-x|^{2}>t$.} In this case, we simply write
\[
P_{t}W(x')-P_{t}W(x)=\int_{\mathbb{R}}p_{t}(y)W^{1}(x-y,x'-y)dy.
\]
Then we invoke (\ref{f11}) and the change of variable $y=u+\sqrt{t}w$
again in order to get 
\[
\big|P_{t}W(x')-P_{t}W(x)\big|\leqslant\kappa|x'-x|^{\alpha}\int_{\mathbb{R}}(a+|y|)^{\chi}p_{t}(y)dy=\kappa|x'-x|^{\alpha}\int_{\mathbb{R}}(a+\sqrt{t}|w|)^{\chi}p_{1}(w)dw
\]
Observe that whenever $t<|x'-x|^{2}$, we also have $\sqrt{t}\leqslant a$.
Hence we obtain 
\begin{align}
\big|P_{t}W(x')-P_{t}W(x)\big| & \leqslant C\cdot\kappa a^{\chi}\cdot|x'-x|^{\alpha}=C\cdot\kappa a^{\chi}\cdot t^{-\beta/2}|x'-x|^{2\beta}\cdot t^{\beta/2}|x'-x|^{\alpha-2\beta}\nonumber \\
 & \leqslant C\cdot\kappa a^{\chi}\cdot t^{-\beta/2}|x'-x|^{2\beta}\cdot t^{\frac{\alpha-\beta}{2}}\leqslant C'\cdot\kappa a^{\chi}\cdot t^{-\beta/2}|x'-x|^{2\beta}.\label{eq:estim-incr-PtW-ter}
\end{align}
We have now proved our claim (\ref{eq:estim-incr-PtW}) when $|x'-x|^{2}>t$.
Hence gathering (\ref{eq:estim-incr-PtW-bis}) and (\ref{eq:estim-incr-PtW-ter})
the desired estimate (\ref{eq:estim-incr-PtW}) is achieved. \end{proof}

\noindent We now gather the previous lemmas and prove our main estimate
for the remainder $\mathcal{R}^{\mathcal{MV}}$.

\begin{proof}[Proof of Lemma \ref{lem:RemEst}] The remainder $\mathcal{R}^{\mathcal{MV}}$
has been decomposed in (\ref{eq:ReMVbis}). Then for an estimate of
the term $\mathcal{R}_{t}^{0}(x,x')$ in (\ref{eq:ReMVbis}), we simply
refer to Lemma \ref{lem:R0Est}. Next we deal with the other term
in (\ref{eq:ReMVbis}), namely 
\[
\frac{2}{\sigma^{2}}v_{t}(x)\left(P_{t}W(x')-P_{t}W(x)\right).
\]
To this aim, we combine the fact that according to Definition \ref{def:weights}
we have 
\[
\llbracket v\rrbracket_{\beta/2,\beta}^{[0,t]\times[-a,a]}\leqslant E^{\theta,\lambda}(a,t),
\]
together with Lemma \ref{lem:SGEst}. We obtain 
\[
\big|v_{t}(x)\big|\cdot\big|P_{t}W(x')-P_{t}W(x)\big|\leqslant C\cdot\kappa\Theta E(a,t)\cdot a^{\chi}t^{-\beta/2}\cdot|x'-x|^{2\beta}.
\]
Gathering this upper bound with our estimate of $\mathcal{R}_{t}^{0}(x,x')$,
the proof of Lemma \ref{lem:RemEst} is now complete. \end{proof}

\subsubsection{\label{Put-toghether} Putting all the estimates together}

Having all the previous variational estimates at hand, we can now
complete the proof of Proposition \ref{prop:EstM}.

First of all, Lemma \ref{lem:UniformEst} gives 
\begin{equation}
\|\mathcal{M}\mathcal{V}\|_{\infty}^{[0,t]\times[-a,a]}\leqslant C\cdot\kappa\Theta E(a,t)\cdot\lambda^{-\frac{\alpha-\beta}{4}}.\label{eq:estim-normMVinfty}
\end{equation}
Next, Lemma \ref{lem:TimEst} and Lemma \ref{lem:SpaEst} yield 
\begin{equation}
\|\mathcal{M}\mathcal{V}\|_{\beta/2,\beta}^{[0,t]\times[-a,a]}\leqslant C\cdot\kappa\Theta E(a,t)\cdot\lambda^{-\frac{\alpha-\beta}{4}}a^{\beta/2}.\label{eq:estim-normMVbeta}
\end{equation}
In addition, since we defined $\partial_{W}(\mathcal{MV})$ as $-2v/\sigma^{2}$
in (\ref{eq:deriv-Gubi-M}) and $\mathcal{V}$ is a controlled process
satisfying Definition \ref{def:weights}, we get 
\begin{equation}
\llbracket\partial_{W}(\mathcal{M}\mathcal{V})\rrbracket_{\beta/2,\beta}^{[0,t]\times[-a,a]}=\frac{2}{\sigma^{2}}\cdot\llbracket v\rrbracket_{\beta/2,\beta}^{[0,t]\times[-a,a]}\leqslant\frac{2}{\sigma^{2}}\cdot\Theta E(a,t).\label{eq:estim-norm-derivMV}
\end{equation}
Finally, Lemma \ref{lem:RemEst} yields 
\begin{equation}
\|\mathcal{R}^{(\mathcal{M}\mathcal{V})_{t}}\|_{2\beta}\leqslant C\cdot\kappa\Theta E(a,t)\cdot Q(a,t),\label{eq:estim-normRMV}
\end{equation}
where the functions $E$ and $Q$ are introduced in (\ref{eq:EQDef}).
As a result, plugging (\ref{eq:estim-normMVinfty})-(\ref{eq:estim-normRMV})
into the definition (\ref{eq:def_normTheta}) of the norm $\Theta=\Theta^{\theta,\lambda}(\mathcal{V})$
and recalling that we have chosen $\gamma=\frac{\alpha-\beta}{4}$,
we end up with 
\begin{multline*}
\Theta^{\theta,\lambda}(\mathcal{M}\mathcal{V})=\sup_{a\geqslant1,t\in[0,T]}E(a,t)^{-1}\cdot\left(\llbracket\mathcal{M}\mathcal{V}\rrbracket_{\beta/2,\beta}^{[0,t]\times[-a,a]}+\lambda^{-\gamma}\llbracket\partial_{W}(\mathcal{M}\mathcal{V})\rrbracket_{\beta/2,\beta}^{[0,t]\times[-a,a]}\right.\\
\left.+\lambda^{-\gamma}Q(a,t)^{-1}\|\mathcal{R}^{(\mathcal{M}\mathcal{V})_{t}}\|_{2\beta}\right)\leqslant C(\kappa+1)\lambda^{-\frac{\alpha-\beta}{4}}\Theta.
\end{multline*}
Therefore, the proof of Proposition \ref{prop:EstM} is complete.

\subsection{\label{sec:estimate-deter-funct} Estimating the deterministic input
functions}

Recall that the main fixed point problem is associated with the transformation
$\hat{\mathcal{M}}$ defined by 
\[
\hat{\mathcal{M}}:\mathcal{V}\mapsto\hat{\mathcal{M}}\mathcal{V}\triangleq\psi^{1}+\psi^{2}+\mathcal{M}\mathcal{V},
\]
where 
\begin{equation}
\psi_{t}^{1}(x)=\int_{\mathbb{R}}\partial_{x}p_{t}(x-y)f_{0}(y)dy,\quad\text{and}\quad\psi_{t}^{2}(x)=\int_{0}^{t}\int_{\mathbb{R}}\partial_{x}p_{t-s}(x-y)g_{s}(y)dyds.\label{eq:defpsi1psi2}
\end{equation}
In order to complete the proof of Theorem \ref{thm:FPT}, we now need
to show that the functions $\psi^{1},\psi^{2}$ can be considered
as controlled processes in the space $\mathcal{B}^{\theta,\lambda}$
introduced in Definition \ref{def:weights}. In order, to achieve
this goal, we assume in this section that $f_{0}\in C_{b}^{2}(\mathbb{R})$
and $g\in C_{b}^{2}([0,T]\times\mathbb{R}).$

\subsubsection{\label{psi1} Estimation of $\psi^{1}$}

In this section we focus on the term $\psi^{1}$ and prove that this
function can be considered as a controlled process. We summarize our
conclusions in the following lemma.
\begin{lem}
\label{lem:psi1} Let $f_{0}$ be a function in $C_{b}^{2}(\mathbb{R})$,
and consider $\psi^{1}$ defined by (\ref{eq:defpsi1psi2}). Then
$\psi^{1}$ is an element of $\mathcal{B}^{\theta,\lambda}$, whose
norm $\Theta$ (see relation (\ref{eq:def_normTheta})) can be bounded
as 
\[
\Theta^{\theta,\lambda}(\psi^{1})\leqslant C\|f_{0}''\|_{\infty}.
\]
\end{lem}

\begin{proof} We first look at the time variations of $\psi^{1}$.
Namely, let $t_{1}<t_{2}.$ By definition, we have 
\[
\psi_{t_{2}}^{1}(x)-\psi_{t_{1}}^{1}(x)=\int_{\mathbb{R}}(\partial_{x}p_{t_{2}}(x-y)-\partial_{x}p_{t_{1}}(x-y))f_{0}(y)dy.
\]
Therefore writing $\partial_{x}p_{t_{2}}-\partial_{x}p_{t_{1}}$ in
terms of a time derivative, resorting to the relation $\partial_{s}p_{s}=\frac{\sigma^{2}}{2}\partial_{xx}^{2}p_{s}$
and setting $y=\sqrt{s}w$ we get 
\begin{equation}
\psi_{t_{2}}^{1}(x)-\psi_{t_{1}}^{1}(x)=\frac{\sigma^{2}}{2}\int_{t_{1}}^{t_{2}}\frac{ds}{s^{3/2}}\int_{\mathbb{R}}\partial_{x}^{3}p_{1}(w)f_{0}(x+\sqrt{s}w)dw.\label{eq:timeincre-psi1}
\end{equation}
Furthermore, thanks to the decaying properties of the heat kernel
$p_{t}$, a straightforward integration by parts procedure shows that
\[
\int_{\mathbb{R}}\partial_{x}^{3}p_{t}(w)wdw=-\int_{\mathbb{R}}\partial_{xx}^{2}p_{t}(w)dw=-2\partial_{t}\int_{\mathbb{R}}p_{t}(w)dw=0,
\]
where we have invoked the relation $\int_{\mathbb{R}}p_{t}(w)dw=1$
for the last step. In the same way, we also have $\int_{\mathbb{R}}\partial_{x}^{3}p_{t}(w)dw=0$
for all $t>0$. Hence we can recast (\ref{eq:timeincre-psi1})
as 
\begin{equation}
\psi_{t_{2}}^{1}(x)-\psi_{t_{1}}^{1}(x)=\frac{\sigma^{2}}{2}\int_{t_{1}}^{t_{2}}\frac{ds}{s^{3/2}}\int_{\mathbb{R}}\partial_{x}^{3}p_{1}(w)\left(f_{0}(x+\sqrt{s}w)-f_{0}(x)-f_{0}'(x)\sqrt{s}w\right)dw.\label{eq:timeincre-psi1-bis}
\end{equation}
The Taylor expansion in (\ref{eq:timeincre-psi1-bis}), together with
the fact that $f_{0}\in C_{b}^{2}(\mathbb{R})$, allow to get rid
of the singularity $s^{-3/2}$ at $s=0$. We obtain 
\begin{equation}
|\psi_{t_{2}}^{1}(x)-\psi_{t_{1}}^{1}(x)|\leqslant C\|f''_{0}\|_{\infty}(\sqrt{t_{2}}-\sqrt{t_{1}})\leqslant C\|f''_{0}\|_{\infty}|t_{2}-t_{1}|^{1/2}.\label{eq:Psi1TVar}
\end{equation}

Next, we consider the space variation of $\psi^{1}$. Namely, let
$t\in(0,T]$ and $x,x'$ be elements of $[-a,a].$ Then we have 
\[
\psi_{t}^{1}(x')-\psi_{t}^{1}(x)=\int_{\mathbb{R}}(\partial_{x}p_{t}(x'-y)-\partial_{x}p_{t}(x-y))f_{0}(y)dy.
\]
Performing the same kind of manipulations as in (\ref{eq:timeincre-psi1})
and (\ref{eq:timeincre-psi1-bis}), we easily get 
\begin{equation}
\psi_{t}^{1}(x')-\psi_{t}^{1}(x)=\frac{1}{t}\int_{x}^{x'}du\int_{\mathbb{R}}\partial_{xx}^{2}p_{1}(w)\left(f_{0}(u+\sqrt{t}w)-f_{0}(u)-f'_{0}(u)\sqrt{t}w\right)dw.\label{eq:timeincre-psi1-ter}
\end{equation}
Invoking the $C^{2}$-regularity of $f_{0}$, we thus obtain 
\begin{equation}
|\psi_{t}^{1}(x')-\psi_{t}^{1}(x)|\leqslant C\cdot\|f_{0}''\|_{\infty}|x'-x|.\label{eq:Psi1SVar}
\end{equation}
Taking limits, notice that the estimate (\ref{eq:Psi1SVar}) is also
valid for $t=0$. Owing to (\ref{eq:Psi1TVar})-(\ref{eq:Psi1SVar}),
one can thus take $\partial_{W}\psi^{1}\equiv0$ in relation (\ref{f2}),
and write 
\[
\psi_{t}^{1}(y)-\psi_{t}^{1}(x)=\mathcal{R}^{\psi_{t}^{1}}(x,y),
\]
with a remainder $\mathcal{R}^{\psi_{t}^{1}}$ enjoying a H\"older regularity
of order $2\beta$ (recall that $2/3<2\beta<1$ according to (\ref{f1})).
Plugging this information in the definition (\ref{eq:def_normTheta})
of $\Theta^{\theta,\lambda}(\psi^{1})$, we thus get that $\psi^{1}\in\mathcal{B}^{\theta,\lambda}$
with $\partial_{W}\psi_{t}^{1}\equiv0$. 
\end{proof} 

\subsubsection{\label{psi2} Estimation of $\psi^{2}$}

The term $\psi^{2}$ can be estimated similarly to $\psi^{1}$ in
Section \ref{psi1}. We label our results in the following lemma.
\begin{lem}
\label{lem:psi2} Let $g$ be a function in $C_{b}^{2}([0,T]\times\mathbb{R})$,
and consider $\psi^{2}$ defined by (\ref{eq:defpsi1psi2}). Then
$\psi^{2}$ sits in the space $\mathcal{B}^{\theta,\lambda}$ given
in Definition \ref{def:weights}, and we have 
\[
\Theta^{\theta,\lambda}(\psi^{2})\leqslant C\sup_{0\leqslant t\leqslant T}\left(\|\partial_{x}g_{t}\|_{\infty}+\|\partial_{xx}^{2}g_{t}\|_{\infty}\right).
\]
\end{lem}

\begin{proof} The proof is very similar to the proof of Lemma \ref{lem:psi1},
and we omit some details for the sake of conciseness. Let us start
with the time variations of $\psi^{2}$. Some easy algebraic manipulations
show that they can be decomposed as 
\begin{equation}
\psi_{t_{2}}^{2}(x)-\psi_{t_{1}}^{2}(x)=\mathcal{I}_{1}+\mathcal{I}_{2},\label{eq:psi2eI1pI2}
\end{equation}
where $\mathcal{I}_{1}$ and $\mathcal{I}_{2}$ are respectively defined
by 
\begin{eqnarray}
\mathcal{I}_{1} & = & \int_{0}^{t_{1}}ds\int_{\mathbb{R}}\left(\partial_{x}p_{t_{2}-s}(x-y)-\partial_{x}p_{t_{1}-s}(x-y)\right)g_{s}(y)dy\label{eq:psi2I1}\\
\mathcal{I}_{2} & = & \int_{t_{1}}^{t_{2}}\int_{\mathbb{R}}\partial_{x}p_{t_{2}-s}(x-y)g_{s}(y)dyds.\label{eq:psi2I2}
\end{eqnarray}
we will now treat those two terms separately.

In order to estimate $\mathcal{I}_{1}$ above we proceed as in (\ref{eq:timeincre-psi1}).
Namely we express $\partial_{x}p_{t_{2}-s}-\partial_{x}p_{t_{1}-s}$
in terms of a time derivative, invoke the relation $\partial_{s}p_{s}=\frac{\sigma^{2}}{2}\partial_{x}^{2}p_{s}$,
and set $y=\sqrt{s}w$. We end up with 
\[
\mathcal{I}_{1}=\frac{\sigma^{2}}{2}\int_{0}^{t_{1}}ds\int_{t_{1}-s}^{t_{2}-s}\frac{du}{u^{3/2}}\int_{\mathbb{R}}\partial_{x}^{3}p_{1}(w)g_{s}(x+\sqrt{u}w)dw.
\]
Next we successively set $u=t_{1}+r-s$ and $\rho=t_{1}+r-s$ in the
above integral. We get 
\begin{equation}
\mathcal{I}_{1}=\frac{\sigma^{2}}{2}\int_{0}^{t_{2}-t_{1}}dr\int_{r}^{t_{1}+r}\frac{d\rho}{\rho^{3/2}}\int_{\mathbb{R}}\partial_{x}^{3}p_{1}(w)g_{t_{1}+r-\rho}(x+\sqrt{\rho}w)dw.\label{eq:psi2I1-bis}
\end{equation}
As in (\ref{eq:timeincre-psi1-bis}), we now cancel the singularity
$\rho^{-3/2}$ by means of a Taylor expansion for $g$ and take advantage
of the decaying properties of the heat kernel $p$. We obtain 
\begin{equation}
|\mathcal{I}_{1}|\leqslant C\sup_{0\leqslant s\leqslant t_{1}}\|\partial_{xx}^{2}g_{s}\|_{\infty}\cdot\int_{0}^{t_{2}-t_{1}}dr\int_{r}^{t_{1}+r}\frac{d\rho}{\sqrt{\rho}}\leqslant C\sqrt{T}\sup_{0\leqslant s\leqslant t_{1}}\|\partial_{xx}^{2}g_{s}\|_{\infty}\cdot|t_{2}-t_{1}|.\label{eq:estimI1psi2}
\end{equation}
Let us now turn to the time variations in $\mathcal{I}_{2}$. That
is setting $u=t_{2}-s$ in the integral defining $\mathcal{I}_{2}$
and proceeding as in (\ref{eq:timeincre-psi1-ter}), we get 
\[
\mathcal{I}_{2}=\int_{0}^{t_{2}-t_{1}}\frac{du}{u^{1/2}}\int_{\mathbb{R}}\partial_{x}p_{1}(w)\left(g_{t_{2}-u}(x+\sqrt{u}w)-g_{t_{2}-u}(x)\right)dw.
\]
It follows that 
\begin{equation}
|\mathcal{I}_{2}|\leqslant C\sup_{0\leqslant s\leqslant t_{2}}\|\partial_{x}g_{s}\|_{\infty}\cdot|t_{2}-t_{1}|.\label{eq:estimI2psi2}
\end{equation}
To summarize, plugging (\ref{eq:estimI1psi2}) and (\ref{eq:estimI2psi2})
into (\ref{eq:psi2eI1pI2}), we can bound the time variations of $\psi^{2}$
as follows 
\begin{equation}
\big|\psi_{t_{2}}^{2}(x)-\psi_{t_{1}}^{2}(x)\big|\leqslant C\cdot\sup_{0\leqslant s\leqslant t_{2}}\|\partial_{x}g_{s}\|_{\infty}\cdot|t_{2}-t_{1}|.\label{eq:Psi2SVar}
\end{equation}
We now handle the spatial variations of $\psi^{2}$. Namely consider
$t\in[0,T]$ and $x,x'\in[-a,a].$ Then we have 
\[
\psi_{t}^{2}(x')-\psi_{t}^{2}(x)=\int_{0}^{t}ds\int_{\mathbb{R}}\left(\partial_{x}p_{s}(x'-y)-\partial_{x}p_{s}(x-y)\right)g_{t-s}(y)dy.
\]
Along the same lines as for (\ref{eq:timeincre-psi1-bis}), (\ref{eq:timeincre-psi1-ter})
and (\ref{eq:estimI2psi2}), we obtain 
\[
\psi_{t}^{2}(x')-\psi_{t}^{2}(x)=\int_{x}^{x'}du\int_{0}^{t}\frac{ds}{s}\int_{\mathbb{R}}\partial_{xx}^{2}p_{1}(w)\left(g_{t-s}(u+\sqrt{s}w)-g_{t-s}(u)-\partial_{x}g_{t-s}(u)\sqrt{s}w\right)dw.
\]
Owing to the $C^{2}-$regularity of $g$, this yields 
\begin{equation}
|\psi_{t}^{2}(x')-\psi_{t}^{2}(x)|\leqslant C\cdot\sup_{0\leqslant s\leqslant t}\|\partial^{2}g_{s}\|_{\infty}\cdot t\cdot|x'-x|.\label{eq:spacePsi2SVar}
\end{equation}
Gathering (\ref{eq:Psi2SVar}) and (\ref{eq:spacePsi2SVar}) we can
now conclude as in Lemma \ref{lem:psi1} that $\psi^{2}\in\mathcal{B}^{\theta,\lambda}$
with $\partial_{W}\psi_{t}^{2}\equiv0.$ \end{proof}

\subsection{Completing the proof of Theorem \ref{thm:FPT}}\label{sec:FPT}

We will use a classical fixed point argument. Specifically, we will
prove that if $\lambda$ is large enough we have 
\begin{equation}
\Theta^{\theta,\lambda}\left(\hat{\mathcal{M}}\mathcal{V}_{2}-\hat{\mathcal{M}}\mathcal{V}_{1}\right)\leqslant\frac{1}{2}\Theta^{\theta,\lambda}(\mathcal{V}_{2}-\mathcal{V}_{1}),\label{eq:contraction}
\end{equation}
which is enough to conclude for the existence of a unique $\mathcal{V}^{*}\in\mathcal{B}^{\theta,\lambda}$
such that $\hat{\mathcal{M}}\mathcal{V}^{*}=\mathcal{V}^{*}.$ In
order to achieve (\ref{eq:contraction}), observe that due to our
definition (\ref{f3}), we have 
\begin{equation}
\hat{\mathcal{M}}\mathcal{V}_{2}-\hat{\mathcal{M}}\mathcal{V}_{1}=\mathcal{MV}_{2}-\mathcal{MV}_{1}=\mathcal{M}(\mathcal{V}_{2}-\mathcal{V}_{1}),\label{eq:diffMV}
\end{equation}
where the second equality stems from the fact that $\mathcal{M}$
is linear. Next, owing to inequality~(\ref{eq:normeMVnormev}) in
Proposition \ref{prop:EstM} (whose proof is spelled out in Section
\ref{Put-toghether}), we have 
\[
\Theta^{\theta,\lambda}\left(\mathcal{M}(\mathcal{V}_{2}-\mathcal{V}_{1})\right)\leqslant\frac{1}{2}\Theta^{\theta,\lambda}\left(\mathcal{V}_{2}-\mathcal{V}_{1}\right),
\]
for $\lambda$ large enough. Together with (\ref{eq:diffMV}), this
entails (\ref{eq:contraction}) and finishes our proof of Theorem
\ref{thm:FPT} item (ii). Item (iii) is then achieved taking into
account the fact that we initiate the Picard iterations from ${\mathcal{V}}^{0}=\psi^{1}+\psi^{2}$,
where we recall that $\psi^{1}$ and $\psi^{2}$ are defined by relation
(\ref{f4}).
\begin{rem}
The above choice of $\lambda$ relies on the rough path norm $\kappa_{\alpha,\chi}({\bf W})$.
As a result, the space $(\mathcal{B}^{\theta,\lambda},\Theta^{\theta,\lambda})$
also implicitly relies on $\kappa_{\alpha,\chi}({\bf W}).$ Nonetheless,
we have mentioned in Remark \ref{rem:MonBThe} that the space $\mathcal{B}^{\theta,\lambda}$
is increasing in $\lambda$. As a result, we have a canonical notion
of existence and uniqueness in the space $\cup_{\lambda>1}\mathcal{B}^{\theta,\lambda}$
which is independent of the quantity $\kappa_{\alpha,\chi}({\bf W}).$
\end{rem}

\section{\label{sec:KMT} Strong approximation of Brownian motion under H\"older
metric}

In this section, we construct a coupling between the discrete and continuous environments and estimate their (rough path) distance with respect to the type of norm introduced in (\ref{f11}). This is a necessary ingredient for establishing the continuity estimate between Sinai's random walk and the Brox diffusion in Section \ref{sec:ConvEst}.

Recall that Sinai's random walk is defined through an environment 
$\{\omega_{x}^{+}:x\in\mathbb{Z}\}$ which is given in Definition \ref{def:environment}. In order to compare the discrete and continuous environments properly, we shall first define a continuously interpolated process from the discrete environment.  

\begin{notation}\label{not:h20}
For $\delta>0$, the rescaled environment $\{\omega_{x}^{+,\delta}:x\in\delta\mathbb{Z}\}$ is introduced in \eqref{eq:rescaled-omega-pm}. We have defined an approximate white noise $\{\dot{U}^{\delta}(x):x\in\delta\mathbb{Z}\}$ in \eqref{eq:def-U-delta}. We now introduce a modified approximate white noise $\{\bar{U}^{\delta}(x):x\in\delta\mathbb{Z}\}$ by setting
\begin{equation}\label{eq:approx-white-noise}
\bar{U}^{\delta}=-2\dot{U}^{\delta},\quad\text{where}\quad\dot{U}^{\delta}(x)=2\omega^{+,\delta}_{x}-\sigma^{2}=\omega^{+,\delta}_{x}-\omega^{-,\delta}_{x}.
\end{equation}
This noise will be mapped into a field with continuous spatial parameter $\{\hat{U}^{\delta}(x):x\in\mathbb{R}\}$ by setting 
\begin{equation}\label{eq:hatUdelta-linear-interpolation}
\hat{U}^{\delta}\equiv\text{linear interpolation obtained from }\bar{U}^{\delta}.
\end{equation}
More specifically, we set $\hat{U}^{\delta}(0)=0$ and for $x\in\delta\mathbb{Z}$ we define $\hat{U}^{\delta}(x)$ recursively by 
\[
\hat{U}^{\delta}(x-\delta,x)\equiv\hat{U}^{\delta}(x)-\hat{U}^{\delta}(x-\delta)=\bar{U}^{\delta}(x).
\]
Then we require that $\hat{U}^\delta$ is linear on each subinterval $[x-\delta,x]$.
\end{notation}

With this notation in mind we now emphasize a technical point in our considerations. Namely, according to Sinai's assumptions for recurrence, we have that $\mathbf{E}[\log(\omega_{x}^{-}/\omega_{x}^{+})]=0$, and therefore the random variables $\dot{U}^{\delta}(x)$ defined by \eqref{eq:approx-white-noise} are not necessarily centered. Nevertheless, we get an asymptotic centering. Let us label this property for further use. Let $\sigma_1^2$ denote the variance of the variable $\xi_x\triangleq \log(\omega_x^- / \omega_x^+)$. 

\begin{lem}\label{lem:U=U1+U2}
For $\delta>0$ and $x\in\delta\mathbb{Z}$, let $\bar{U}^{\delta}(x)$ be defined by \eqref{eq:approx-white-noise}. We decompose  $\bar{U}^{\delta}(x)$
as 
\begin{equation}{\label{eq:BarUDecomp}}
\bar{U}^{\delta}(x)=\bar{U}_{1}^{\delta}(x)+\bar{U}_{2}^{\delta}(x),
\end{equation}
where $\bar{U}_{1}^{\delta},\bar{U}_{2}^{\delta}$ are respectively defined by 
\begin{equation}\label{eq:BarU1Decomp}
\bar{U}_{1}^{\delta}(x)\triangleq\bar{U}^{\delta}(x)-\mathbf{E}\left[\bar{U}^{\delta}(x)\right],
\quad\text{and}\quad 
\bar{U}_{2}^{\delta}(x)\triangleq\mathbf{E}\left[\bar{U}^{\delta}(x)\right].
\end{equation}
Then the following holds true:

\vspace{2mm}\noindent 
\textnormal{(i)}\: the random variable  $\bar{U}_{1}^{\delta}(x)$ is centered  with variance ${\rm Var}\left(\bar{U}_{1}^{\delta}(x)\right)=\sigma^4\sigma_1^2 \delta + o(\delta)$.\\
\textnormal{(ii)}\: $\bar{U}_{2}^{\delta}(x)$ is of order $O(\delta^{3/2})$. 
\end{lem}
\begin{proof}
We have already argued that (i) holds true in Remark \ref{rk:variance} (cf. (\ref{eq:var-udot-delta}) and (\ref{eq:1OrdApXi})). We will thus focus on item (ii). Next let us recast \eqref{eq:rescaled-omega-pm} by 
setting 
\begin{equation}\label{eq:xixdelta}
\xi_{x}^{\delta}=\log(\omega_{x/\delta}^{-}/\omega_{x/\delta}^{+}).
\end{equation}
Thanks to relation \eqref{eq:xi-rw} we get 
\begin{equation}
\omega_{x}^{+,\delta} =\frac{\sigma^{2}}{2}\cdot\frac{2}{1+e^{\sqrt{\delta}\xi_{x}^{\delta}}}
=\frac{\sigma^{2}}{2}\Bigg(1-\tanh\lp\frac{\sqrt{\delta}\xi_{x}^{\delta}}{2}\rp\Bigg)\label{eq:OmeDel}
\end{equation}
Hence owing to expression \eqref{eq:approx-white-noise} we get 
\[
\bar{U}_{2}^{\delta}(x)=2\sigma^{2}\mathbf{E}\left[\tanh\left(\sqrt{\delta}\,\frac{\xi_{x}^{\delta}}{2}\right)\right].
\]
Resorting to a Taylor expansion of the function $\tanh$ and an easy dominated convergence argument, we thus obtain 
\begin{equation}\label{k1}
\bar{U}_{2}^{\delta}(x)=\sigma^{2}\Bigg(\frac{\delta^{1/2}}{2}\mathbf{E}[\xi_{x}^{\delta}]
+\frac{\delta^{3/2}}{3}\mathbf{E}\left[\left(\frac{\xi_{x}^{\delta}}{2}\right)^{3}\right]+o(\delta^{3/2})\Bigg)
=c_{\sigma}\delta^{3/2}+o(\delta^{3/2}),
\end{equation}
where we have invoked the recurrence hypothesis $\mathbf{E}[\xi_{x}^{\delta}]=0$ for the last identity. This finishes the proof. 
\end{proof}

From the above discussion, it is natural to view the centred process $\hat{U}^\delta_{1}$ as a discrete approximation of $W$ (which will be coupled to $W$ according to Theorem \ref{thm:StrongApprox} below) and regard $\hat{U}^\delta_2$ as a remainder. According to Lemma \ref{lem:U=U1+U2} (ii), the exact variance of the Brownian motion $W$ should be given by
\begin{equation}\label{eq:WVar}
\tau^2 \triangleq \sigma^4 \sigma_1^2.
\end{equation}Here we do \textit{not} want to view $\hat{U}^\delta$ as an approximation of $W$ (under rough path metric) since the straight line $x\mapsto\hat{U}^\delta_2(x)$ does not have the correct spatial growth $a^\chi$ encoded in the definition (\ref{f11}). With those preliminaries in hand, our main result in this section is stated as follows.
\begin{thm}\label{thm:StrongApprox} Let $W=\{W(x):x\in\mathbb{R}\}$ be a given two-sided Brownian motion with variance $\tau^2$ defined on some probability space $(\Omega,\mathcal{F},\mathbf{P})$. Let $\alpha,\chi$ be given fixed parameters
satisfying the constraint (\ref{f1}). We denote by $F$ the common distribution of the random variables $\xi_{x}=\log(\omega_{x}^{-}/\omega_{x}^{+})$ defined by \eqref{eq:xi-rw}. Then for each $\delta\in(0,1)$, one can construct a stochastic process $\bar{U}_{1}^{\delta}=\{\bar{U}_{1}^{\delta}(x):x\in\delta\mathbb{Z}\}$ on $(\Omega,\mathcal{F},\mathbf{P})$ such that the following properties hold true.

\vspace{2mm}\noindent \textnormal{(i)} The random variables $\{\bar{U}_{1}^{\delta}(x):x\in\delta\mathbb{Z}\}$
are independent and identically distributed, with a distribution determined by the following relation:
\[
\bar{U}_{1}^{\delta}(x)\stackrel{{\rm law}}{=}
-4\left(\frac{\sigma^{2}}{1+e^{\sqrt{\delta}Z}}-\mathbf{E}\left[\frac{\sigma^{2}}{1+e^{\sqrt{\delta}Z}}\right]\right),\quad\text{ with }\;\; Z\stackrel{{\rm law}}{=}F.
\]
\textnormal{(ii)} Let $\hat{U}_{1}^{\delta}$ be the interpolated process constructed from $\bar{U}_1^\delta$ over the grid $\delta\mathbb{Z}$. Recall that the rough path lifting ${\bf W}$ of $W$ is
defined by (\ref{eq:def-w1-w2}). Correspondingly, we also define $\hat{{\bf U}}_{1}^{\delta}=(\hat{U}_{1}^{\delta},\hat{U}_1^{\delta;2})$
where $\hat{U}_1^{\delta;2}(x,y)\triangleq\hat{U}_{1}^{\delta}(x,y){}^{2}/2$.
Consider the rough path distance between $\hat{{\bf U}}_{1}^{\delta}$
and ${\bf W}$ defined by
\begin{equation}\label{eq:RhoUW}
\rho_{\alpha,\chi}\left(\hat{{\bf U}}_{1}^{\delta},{\bf W}\right)\triangleq
\sup_{a\geqslant1}\Bigg(\frac{\|\hat{U}_{1}^{\delta}-W\|_{\alpha}^{[-a,a]}}{a^{\chi}}+\frac{\|\hat{U}_{1}^{\delta;2}-W^{2}\|_{\alpha}^{[-a,a]}}{a^{2\chi}}\Bigg).
\end{equation}
Then for any given $\eta\in(0,1/2-\alpha)$ and $q>(1/2-\alpha-\eta)^{-1}$,
there exists a constant $C=C_{\alpha,\chi,\eta,q}>0$ such that
\begin{equation}\label{eq:RhoLqEst}
\big\|\rho_{\alpha,\chi}\left(\hat{{\bf U}}_{1}^{\delta},{\bf W}\right)\big\|_{L^{q}}\leqslant C\delta^{\eta}
\end{equation}
for all $\delta\in(0,1)$. In particular, by the Borel-Cantelli lemma,
for any $\tau\in(0,\eta)$ we have the following a.s. estimate:
\begin{equation}\label{eq:KMTAS}
\rho_{\alpha,\chi}\left(\hat{{\bf U}}_{1}^{\delta},{\bf W}\right)\leqslant\Xi\delta^{\tau}
\end{equation}
where $\Xi$ is some a.s. finite random variable that is independent
of $\delta.$

\end{thm}

The rest of this section is devoted to the proof of Theorem \ref{thm:StrongApprox}. The first ingredient is a classical approximation theorem proved by
Koml\'os-Major-Tusn\'ady \cite{KMT76},  which we now recall.

\begin{thm}[Classical KMT Approximation]\label{thm:ClassicKMT} Let
$G$ be a distribution function on $\mathbb{R}$ with mean zero and
unit variance. We denote by $R$ the moment generating function of $G$, considered as a function of $z\in\mathbb{C}$:
\[
R(z)=\int_{\mathbb{R}}e^{zx}G(dx), 
\]
whenever the right hand side above is properly defined. Then we suppose that 
$G$ satisfies the following assumptions. 

\vspace{2mm}\noindent 
(i) The function $R$
is well-defined in some neighbourhood $(-t_{0},t_{0})$ of the origin.\\
(ii) Either $G$ is lattice-valued or there exists $p>1$ such that $R$ 
verifies
\begin{equation}
\int_{\mathbb{R}}|R(t+iu)|^{p}du<\infty,\;\text{ for all $t$ such that }\;|t|<t_{0}.\label{eq:KMTAssump}
\end{equation}
Then given a sequence $\{Y_{n}:n\geqslant1\}$
of i.i.d. standard normal random variables on some probability space
$(\Omega,\mathcal{F},\mathbb{P})$, there exist a sequence $\{f_{n}:n\geqslant 1\}$ of Borel functions from 
$\mathbb{R}^{2n}$to  $\mathbb{R}$ such
that 
\[
X_{n}\triangleq f_{n}(Y_{1},\ldots,Y_{2n})
\]
defines an i.i.d. sequence with distribution $G$ and the following
estimate holds true:

\begin{equation}\label{eq:ClassicKMTEst}
\mathbf{P}\left(\max_{1\leqslant k\leqslant n}|S_{k}-T_{k}|>C\log n+x\right)\leqslant Ke^{-\lambda x}\ \ \ \forall n\geqslant1,x>0.
\end{equation}
Here the random variables $S_{k}$ and $T_{k}$ are respectively defined by $S_{k}\triangleq X_{1}+\cdots+X_{k},$ $T_{k}\triangleq Y_{1}+\cdots+Y_{k}$,
and $C,K,\lambda$ are constants depending only on $G$.

\end{thm}

\begin{rem}\label{rem:DensityCond}
The condition (\ref{eq:KMTAssump}) is satisfied if $G$ has a $\mathcal{C}^{1}$
density function with at most finitely many algebraic singularities. 
%Condition (ii) in Theorem \ref{thm:ClassicKMT} is just 
It should be seen as a mere technical assumption. The coupling between the two sequences $\{X_n\},\{Y_n\}$ still exists without this assumption. However, the price to pay is that both sequences need to be constructed together in this case (one can no longer first fix $Y$ and then construct $X$ from $Y$ to satisfy the estimate (\ref{eq:ClassicKMTEst})). Correspondingly, without Condition (iii) in Definition \ref{def:environment} of the discrete environment, one cannot fix $W$ in advance; the coupling between $W$ and $\hat{U}_1^\delta$ (in particular, $W$ itself) will also depend on $\delta$. This will only lead to the sacrifice of an arbitrarily small power of $\delta$ in the final result. To reduce technicalities, we therefore decide to impose Condition (iii) in Definition~\ref{def:environment}.

\end{rem}

Our next ingredient is a rough path convergence result based on a Kolmogorov type estimate. Although it will be applied to a (trivial) $\mathbb{R}$-valued 
rough path, we state the result in a general finite dimensional vector space $V$. This does not affect the difficulty of our proof and might be interesting in its own right. Let $X$ be a $V$-valued path, whose increments are written as $X_{s,t}^{1}$. We assume that $X$ can be enhanced as a second order rough path 
\begin{equation}\label{eq:enhanced-path}
{\bf X}_{s,t}=(X_{s,t}^{1},X_{s,t}^{2})\in V\oplus V^{\otimes2},\ \ \ -a\leqslant s\leqslant t\leqslant a
\end{equation}
where the rough path notation is borrowed from \cite{FH14}. For $\alpha\in(1/3,1/2)$
we define 
\begin{equation}\label{eq:Holder-norms-enhanced}
\|X^{1}\|_{\alpha}\triangleq\sup_{-a\leqslant s<t\leqslant a}\frac{\|X_{s,t}^{1}\|}{|t-s|^{\alpha}},\ \|X^{2}\|_{2\alpha}\triangleq\sup_{-a\leqslant s<t\leqslant a}\frac{\|X_{s,t}^{2}\|}{|t-s|^{2\alpha}}.
\end{equation}
Recall that the rough path property also means that $\mathbf{X}$ satisfies Chen's relation. Namely, for $-a\leqslant s<u<t\leqslant a$ we have
\begin{equation}
X^2_{s,t} - X^2_{s,u} - X^2_{u,t} = X^1_{s,u}\otimes X^1_{u,t}.\label{eq:Chen_relation}
\end{equation}
We now state our main rough path convergence theorem based on Kolmogorov type estimates. 

\begin{thm}
\label{thm:KolCty}Let ${\bf X}_{s,t}=(X_{s,t}^{1},X_{s,t}^{2})$
and ${\bf Y}_{s,t}=(Y_{s,t}^{1},Y_{s,t}^{2})$ be two random second-order
rough paths over the time horizon $[-a,a]$ in some finite dimensional vector space $V$ as defined in \eqref{eq:enhanced-path}-\eqref{eq:Holder-norms-enhanced}. Let $C,\rho,\varepsilon,\alpha,\nu,q$ be
given positive parameters such that $q>2,$ $\nu>1/q$ and $\alpha\in(0,\nu-1/q)$.
Suppose that the following moment bounds hold true:
\begin{equation}
\|X_{s,t}^{1}\|_{L^{q}}\leqslant C\rho|t-s|^{\nu},\quad \|Y_{s,t}^{1}\|_{L^{q}}\leqslant C\rho|t-s|^{\nu}\label{eq:MEstOrg}
\end{equation}
and 
\begin{equation}
\|X_{s,t}^{1}-Y_{s,t}^{1}\|_{L^{q}}\leqslant C\rho\varepsilon|t-s|^{\nu},\quad \|X_{s,t}^{2}-Y_{s,t}^{2}\|_{L^{q/2}}\leqslant C\rho^{2}\varepsilon|t-s|^{2\nu}.\label{eq:MEstErr}
\end{equation}
Then there exist positive random variables $K^{1}\in L^{q},$ $K^{2}\in L^{q/2}$,
such that
\begin{equation}\label{eq:bounds_XY12}
\|X^{1}-Y^{1}\|_{\alpha}\leqslant K^{1},\quad \|X^{2}-Y^{2}\|_{2\alpha}\leqslant K^{2}
\end{equation}
and 
\[
\|K^{1}\|_{L^{q}}\leqslant C'\rho a^{\nu-\alpha}\varepsilon,\quad \|K^{2}\|_{L^{q/2}}\leqslant C'\rho^{2}a^{2(\nu-\alpha)}\varepsilon,
\]
where $C'$ is a constant depending only on $C, \alpha,\nu,q$.
\end{thm}

\begin{proof} We divide this proof in several steps. Here we use the notation
$\lesssim$ to denote an estimate up to a multiplicative constant
depending only on $C,\alpha,\nu,q$. 

\noindent \textit{Step 1: Inequality on successive dyadic points.}
Let $D_{n}\triangleq\{ka/2^{n}:-2^n\leqslant k\leqslant 2^n\}$ be the $n-$th order
dyadic partition of $[-a,a]$ and $D=\cup_{n=0}^{\infty}D_{n}.$ For
$i=1,2$, we define 
\begin{equation}
K_{n}^{i}=\max_{-2^n+1\leqslant k\leqslant2^{n}}\left|X_{\frac{k-1}{2^{n}}a,\frac{k}{2^{n}}a}^{i}-Y_{\frac{k-1}{2^{n}}a,\frac{k}{2^{n}}a}^{i}\right|.\label{eq:Kn12}
\end{equation}
Then for $i=1,2$ we trivially have 
\[
|K_{n}^{i}|^{q}\leqslant\sum_{k=-2^n+1}^{2^{n}}\left|X_{\frac{k-1}{2^{n}}a,\frac{k}{2^{n}}a}^{i}-Y_{\frac{k-1}{2^{n}}a,\frac{k}{2^{n}}a}^{i}\right|^{q}.
\]
%Let us focus on estimating the first-level components.
Owing to assumption (\ref{eq:MEstErr}) we get 
\begin{equation*}
\mathbf{E}[|K_{n}^{1}|^{q}]\leqslant2^{n+1}\cdot C^{q}\rho^{q}\varepsilon^{q}\left(\frac{a}{2^{n}}\right)^{\nu q}=2C^{q}\rho^{q}\varepsilon^{q}a^{\nu q}2^{-n(\nu q-1)}.
\end{equation*}
This yields the following inequality:
\begin{equation}
\|K_{n}^{1}\|_{L^{q}}\lesssim\varepsilon\rho a^{\nu}2^{-n(\nu-1/q)}.\label{eq:estim_Kn1}
\end{equation}Similarly, for the second level component we also have 
\begin{equation}
\|K_{n}^{2}\|_{L^{q/2}}\lesssim\varepsilon\rho^{2} a^{2\nu}2^{-2n(\nu-1/q)}.\label{eq:estim_Kn2}
\end{equation}

\noindent \textit{Step 2: Decomposition of dyadic intervals.} Let
$s<t$ be two dyadic points. There is a unique $m\geqslant -1$ such
that 
\begin{equation}
2^{-(m+1)}a<t-s\leqslant2^{-m}a.\label{eq:mdyadic}
\end{equation}
We claim that there exists a finite partition $s=\tau_{0}<\tau_{1}<\cdots<\tau_{L}=t,$
such that:

\vspace{2mm}\noindent (i) For each $i,$ $[\tau_{i},\tau_{i+1}]$
is a dyadic sub-interval of order $n$ for some $n\geqslant m$. 

\smallskip

\noindent
(ii) For each $n\geqslant m$, there are at most two intervals among
the collection $\{[\tau_{i},\tau_{i+1}];\,0\le i\le L-1\}$ that belong
to the dyadic partition $D_{n}$.

 To prove the claim, we start from
the point $s$. Denote $\rho_{0}\triangleq s$. By uniquely writing
$s=\frac{k}{2^{n}}a$ 
with $k$ being an odd number, we set $\rho_{1}\triangleq\frac{k+1}{2^{n}}a$ .
We then reduce $\rho_{1}$ to the unique form of $\frac{k'}{2^{n'}}a$
where $k'$ is odd, and iterate this construction. In this way we
produce a family $(\rho_{i})$ such that if we write $\rho_{i}-\rho_{i-1}=2^{-n_{i}}a$,
then we have $n_{i+1}<n_{i}$. Denote by $M$ the last index such
that $\rho_{M}\leqslant t$. Since $n_{0}=n$, $n_{i+1}<n_{i}$,
$\rho_{i}-\rho_{i-1}=2^{-n_{i}}a$ and we are working on the interval
$[0,1]$, the quantity $M$ is clearly finite. Therefore we have obtained
a family $\{\rho_{1},\ldots,\rho_{M+1}\}$ such that 
\[
s<\rho_{1}<\rho_{2}<\cdots<\rho_{M}\leqslant t<\rho_{M+1}.
\]
As mentioned above, the indices $n_{i}$ satisfy $n_{i+1}<n_{i}$.
Thus the intervals $[\rho_{i},\rho_{i+1}]$ correspond to different
orders of dyadic partitions $D_{n_{i}}$. In addition, $[\rho_{i-1},\rho_{i}]$
is always a dyadic interval of some order $n\geqslant m$ since 
\[
\rho_{i}-\rho_{i-1}\leqslant t-s\leqslant2^{-m}a.
\]
If $\rho_{M}$ happens to be equal to $t$, the $\rho_{i}$'s
give the desired partition and we are done. If $\rho_{M}<t$, we
start from the $t$-end and propagate towards the left direction in
exactly the same way as above to obtain 
\[
\mu_{N+1}<\rho_{M}\leqslant\mu_{N}<\mu_{N-1}<\cdots<\mu_{1}<\mu_{0}\triangleq t,
\]
where $\mu_{N}$ is the last point that is not smaller than $\rho_{M}$.
We claim that $\rho_{M}=\mu_{N}.$ Indeed, suppose on the contrary
that $\rho_{M}\neq\mu_{N}.$ Then 
\[
\mu_{N+1}<\rho_{M}<\mu_{N}\leqslant t<\rho_{M+1}.
\]
Let us write 
\[
[\rho_{M},\rho_{M+1}]=\left[\frac{k}{2^{p}}a,\frac{k+1}{2^{p}}a\right],\quad\text{and}\quad[\mu_{N+1},\mu_{N}]=\left[\frac{l-1}{2^{q}}a,\frac{l}{2^{q}}a\right].
\]
Since $\mu_{N}\in(\rho_{M},\rho_{M+1})$ and $\rho_{M}$, $\rho_{M+1}$
are adjacent dyadic points, $\mu_{N}$ has to be a finer dyadic point,
i.e. $q>p$. Similarly, since $\rho_{M}\in(\mu_{N+1},\mu_{N})$
we must also have $p>q$. This clearly gives a contradiction. As a
result, we have $\rho_{M}=\mu_{N}$. Summarizing our considerations
for this step, the partition 
\[
\{\rho_{0},\rho_{1},\ldots,\rho_{M}=\mu_{N},\mu_{N-1},\ldots,\mu_{1},\mu_{0}\}
\]
satisfies our two claims (i) and (ii), with $L=M+N$.

\noindent \textit{Step 3: Proof for the first level of the rough path.}
We now take the results of Step 2 for granted and turn to the proof
of (\ref{eq:bounds_XY12}) for $X^{1}-Y^{1}$. We start by introducing
another piece of notation. Namely, for $i=1,2$ and $-a\leqslant s\leqslant t\leqslant a$
we write 
\begin{equation}
Z_{s,t}^{i}:=X_{s,t}^{i}-Y_{s,t}^{i}.\label{eq:ZiXiYi}
\end{equation}
As in Step 2, consider $s,t$ satisfying (\ref{eq:mdyadic}) for a
given $m\geq-1$. From the properties of the partition $\{\tau_{0},\ldots,\tau_{L}\}$
established in Step 2, we get 
\[
|Z_{s,t}^{1}|\leqslant\sum_{j=1}^{M+N}|Z_{\tau_{j-1},\tau_{j}}^{1}|\leqslant2\sum_{n=m}^{\infty}K_{n}^{1},
\]
where we recall that $K_{n}^{1}$ is defined by (\ref{eq:Kn12}).
Thus, thanks to the fact that $t-s\geqslant 2^{-(m+1)}a$, we obtain
\begin{equation}
\frac{|Z_{s,t}^{1}|}{|t-s|^{\alpha}}\leqslant2\sum_{n=m}^{\infty}\frac{K_{n}^{1}}{|t-s|^{\alpha}}\lesssim\sum_{n=1}^{\infty}\frac{K_{n}^{1}}{a^{\alpha}2^{-n\alpha}}:=K^{1}.\label{h1}
\end{equation}
Note that we have established (\ref{h1}) for dyadic points only.
However, due to the fact that we have assumed $Z^{1}$ to be continuous,
one can easily extend (\ref{h1}) to any couple $-a\leqslant s\leqslant t\leqslant a$
by taking limits over dyadic points. Hence we proved that 
\begin{equation}
\|Z^{1}\|_{\alpha}\leqslant K^{1}.\label{eq:estimZ1}
\end{equation}
Moreover, due to the constraint $\alpha\in(0,\nu-1/q),$ the random variable
$K^{1}$ has finite $q$-moments. Indeed, owing to (\ref{eq:estim_Kn1})
one has
\begin{equation}
\|K^{1}\|_{L^{q}}\leqslant a^{-\alpha}\sum_{n=1}^{\infty}2^{n\alpha}\|K_{n}^{1}\|_{L^{q}}\lesssim\rho\varepsilon a^{\nu-\alpha}\sum_{n=1}^{\infty}2^{-n(\nu-1/q-\alpha)}\lesssim\rho\varepsilon a^{\nu-\alpha}.\label{eq:K1Est}
\end{equation}
Gathering (\ref{eq:estimZ1}) and (\ref{eq:K1Est}), we have
thus proved (\ref{eq:bounds_XY12}) for the first component $X^{1}-Y^{1}$.

\noindent \textit{Step 4. Proof for the second level of the rough
path.} Recall that $Z^{2}$ is defined by (\ref{eq:ZiXiYi}). Then
invoking Chen's relation (\ref{eq:Chen_relation}), it is readily
checked that for $s\leqslant u\leqslant t$ we have 
\begin{equation}
Z_{s,t}^{2}=Z_{s,u}^{2}+Z_{u,t}^{2}+X_{s,u}^{1}\otimes(X_{u,t}^{1}-Y_{u,t}^{1})+(X_{s,u}^{1}-Y_{s,u}^{1})\otimes Y_{u,t}^{1}.\label{eq:Z2}
\end{equation}
Considering two dyadic points $s<t$ such that (\ref{eq:mdyadic})
holds true and iterating relation (\ref{eq:Z2}) over the partition
$\{\tau_{0},\ldots,\tau_{L}\}$ constructed in Step 2, we end up with
\begin{equation}
Z_{s,t}^{2}=\sum_{l=1}^{N}Z_{\tau_{l-1},\tau_{l}}^{2}+\sum_{l=2}^{N}Z_{\tau_{0},\tau_{l-1}}^{1}\otimes X_{\tau_{l-1},\tau_{l}}^{1}+\sum_{l=2}^{N}Y_{\tau_{0},\tau_{l-1}}^{1}\otimes Z_{\tau_{l-1},\tau_{l}}^{1}.\label{eq:2ndErr}
\end{equation}
Therefore we will bound $Z_{s,t}^{2}$ in the following way: 
\begin{equation}
|Z_{s,t}^{2}|\leqslant\sum_{l=1}^{N}|Z_{\tau_{l-1},\tau_{l}}^{2}|+\sum_{l=2}^{N}|Z_{\tau_{0},\tau_{l-1}}^{1}|\cdot|X_{\tau_{l-1},\tau_{l}}^{1}|+\sum_{l=2}^{N}|Y_{\tau_{0},\tau_{l-1}}^{1}|\cdot|Z_{\tau_{l-1},\tau_{l}}^{1}|.\label{eq:estimZ2}
\end{equation}
Next we recall that $K_{n}^{1},K_{n}^{2}$ are defined by (\ref{eq:Kn12}),
and that we have obtained their moment estimates in (\ref{eq:estim_Kn1}) and (\ref{eq:estim_Kn2}) respectively.
\begin{comment}
Thanks to (\ref{eq:MEstErr}), it is straightforward to check that
a similar inequality holds for $K_{n}^{2}$. Therefore we obtain that
for $n\geq1$ and $\beta>q^{-1}$ we have 
\begin{equation}
\|K_{n}^{1}\|_{L^{q}}\vee\|K_{n}^{2}\|_{L^{q/2}}\leqslant C\varepsilon2^{-n(\beta-q^{-1})}.\label{eq:estimKn1Kn2}
\end{equation}
\end{comment}
In addition to the random variables $K_{n}^{1},K_{n}^{2}$, we define
two more random quantities $M_{n}$ and $N_{n}$ as follows: 
\[
M_{n}=\max_{-2^n+1\leqslant k\leqslant 2^n}\left|X_{\frac{k-1}{2^{n}}a,\frac{k}{2^{n}}a}^{1}\right|,\quad\text{ and }\quad N_{n}=\max\left|Y_{\frac{k-1}{2^{n}}a,\frac{k}{2^{n}}a}^{1}\right|.
\]
Similarly to (\ref{eq:estim_Kn1}), it is easily seen that 
\begin{equation}
\|M_{n}\|_{L^{q}}\vee\|N_{n}\|_{L^{q}}\lesssim\rho a^{\nu}2^{-n(\nu-1/q)}.\label{eq:estimMnNn}
\end{equation}
Plugging this information into (\ref{eq:estimZ2}) we get that 
\[
|Z_{s,t}^{2}|\leqslant2\sum_{n=m+1}^{\infty}K_{n}^{2}+4\left(\sum_{n=m+1}^{\infty}K_{n}^{1}\right)\cdot\left(\sum_{n=m+1}^{\infty}(M_{n}+N_{n})\right),
\]
from which we obtain 
\begin{equation}
\frac{|Z_{s,t}^{2}|}{|t-s|^{2\alpha}}\lesssim L+K^{1}\cdot R,\quad\text{ where }\quad L\triangleq\sum_{n=1}^{\infty}\frac{K_{n}^{2}}{2^{-2n\alpha}a^{2\alpha}},\quad R\triangleq\sum_{n=1}^{\infty}\frac{M_{n}+N_{n}}{2^{-n\alpha}a^\alpha}.\label{eq:LR}
\end{equation}
We can now estimate the terms $L$ and $R$ in (\ref{eq:LR}). Indeed,
the inequality~(\ref{eq:estim_Kn2})
yields
\begin{equation}
\|L\|_{L^{q/2}}\leqslant a^{-2\alpha}\sum_{n=1}^{\infty}2^{2n\alpha}\|K_{n}^{2}\|_{L^{q/2}}\lesssim a^{2(\nu-\alpha)}\rho^{2}\varepsilon\sum_{n=1}^{\infty}2^{2n\alpha}2^{-2n(\nu-1/q)}\lesssim a^{2(\nu-\alpha)}\rho^{2}\varepsilon.\label{eq:LEst}
\end{equation}Similarly, the inequality~\eqref{eq:estimMnNn} yields 
\begin{equation}
\|R\|_{L^{q}}\lesssim\rho a^{\nu-\alpha}.\label{eq:REst}
\end{equation}
Reporting these relations into (\ref{eq:LR}) and using an approximation
procedure along dyadics as we did in Step 3, we end up with \[
\|Z^{2}\|_{2\alpha}\lesssim L+K^{1}R =:K^2.
\] 
Gathering the estimates (\ref{eq:K1Est}), (\ref{eq:LEst}) and (\ref{eq:REst}), we arrive at \[
\|K^{2}\|_{L^{q/2}}\leqslant\|L\|_{L^{q/2}}+\|K^{1}R\|_{L^{q/2}}\leqslant\|L\|_{L^{q/2}}+\|K^{1}\|_{L^{q}}\|R\|_{L^{q}}\lesssim\varepsilon\rho^{2}a^{2(\nu-\alpha)}.
\]This gives the desired estimate for the second level difference $X^2-Y^2$ and finishes the proof. 
\end{proof}

\noindent Now we are in a position to prove Theorem \ref{thm:StrongApprox}.

\begin{proof}[Proof of Theorem \ref{thm:StrongApprox}] Recall that $\bar{U}_{1}^{\delta}$ is defined by \eqref{eq:BarU1Decomp}. Before introducing the coupling, we  make one more technical observation on the structure of this random field. In view of the Taylor expansion \eqref{k1}, let us further write 
\begin{equation}\label{eq:U1=U11+U12}
\bar{U}_{1}^{\delta}(x)=\bar{U}_{1,1}^{\delta}(x)+\bar{U}_{1,2}^{\delta}(x),
\end{equation}
where we recall that $\xi_{x}^{\delta}$ is defined by \eqref{eq:xixdelta} and 
\begin{equation}\label{eq:u11-u12-delta}
\bar{U}_{1,1}^{\delta}(x)\triangleq\sqrt{\delta}\sigma^{2}\xi_{x}^{\delta},\qquad \bar{U}_{1,2}^{\delta}(x)\triangleq\bar{U}_{1}^{\delta}(x)-\sqrt{\delta}\sigma^{2}\xi_{x}^{\delta}.
\end{equation}
In the decomposition \eqref{eq:u11-u12-delta}, the term $\bar{U}_{1,2}^\delta(x)$ is easily handled. Indeed, notice that $\{\bar{U}_{1,2}^\delta(x):x\in\delta\mathbb{Z}\}$ is a sequence of i.i.d. centered random variables. Next the ellipticity assumption (i) in Definition \ref{def:environment} together with our definition \eqref{eq:xixdelta} imply that $\xi_{x}^{\delta}$ 
is uniformly bounded in $(x,\delta)$. Owing to the Taylor expansion 
(\ref{k1}), we thus get 
\begin{equation}\label{eq:uppboU12}
\big|\bar{U}_{1,2}^\delta (x)\big|\leqslant c\delta^{3/2},  
\end{equation}
for a universal constant $c$. The term   $\bar{U}_{1,2}^{\delta}(x)$ in 
\eqref{eq:u11-u12-delta} can thus be treated as asymptotically null, and we will mostly focus our attention on $\bar{U}_{1,1}^{\delta}(x)$. 

Let us now describe the coupling between $W$ and $\bar{U}_{1}^{\delta}$. Namely let $W$ be a two-sided Brownian motion with variance $\tau^{2}\triangleq\sigma^{4}{\rm Var}(F)$ defined on some probability space $(\Omega,\mathcal{F},\mathbf{P})$, where we recall that $F$ is the common distribution of the random variables $\xi_{x}$.
For each fixed $\delta\in(0,1),$ we set $W_{t}^{\delta}\triangleq\delta^{-1/2}W_{\delta t}.$
Note that $W^{\delta}$ is again a two-sided Brownian motion with
variance $\tau^2$. With obvious adaptation of constants,
one can apply Theorem \ref{thm:ClassicKMT} to construct an i.i.d.
family $\{X_{m}^{\delta}:m\in\mathbb{Z}\}$ on the same probability
space $(\Omega,\mathcal{F},\mathbf{P})$, such that 
\begin{equation}
X_{m}^{\delta}\stackrel{{\rm law}}{=}
\sigma^{2}\log\lp\frac{\omega_{0}^{-}}{\omega_{0}^{+}}\rp=\sigma^{2}\xi_{0},\label{eq:DistKMT}
\end{equation}
where the second identity stems from \eqref{eq:xi-rw}, and the following estimate holds true:
\begin{equation}
\mathbf{P}\left(\max_{-k\leqslant l<m\leqslant k}\big|S_{l,m}^{\delta}-W_{l,m}^{\delta}\big|>C\log k+x\right)\leqslant Ke^{-\lambda x}\ \ \ \forall k\geqslant1,x>0,\label{eq:KMTEnv}
\end{equation}
where $S_{l,m}^{\delta}\triangleq X_{l+1}^{\delta}+\cdots+X_{m}^{\delta}$
and $C,K,\lambda$ are constants depending only on the distribution
of (\ref{eq:DistKMT}).

We now proceed to define an approximate white noise $\hat{U}_{1}^{\delta}$. Namely, in view of (\ref{eq:OmeDel}), we define 
\[
\omega_{x}^{+,\delta}\triangleq\frac{\sigma^{2}}{2}\cdot\frac{2}{1+e^{\sigma^{-2}\sqrt{\delta}X_{x/\delta}^{\delta}}},\ \ \ x\in\delta\mathbb{Z}.
\]
With \eqref{eq:BarU1Decomp} and \eqref{eq:u11-u12-delta} in mind, we then set 
\begin{equation}\label{eq:barU1-barU11-barU12}
\bar{U}_{1}^{\delta}(x)\triangleq-4\left(\omega_{x}^{+,\delta}-\mathbf{E}[\omega_{x}^{+,\delta}]\right),\quad \bar{U}_{1,1}^{\delta}(x)\triangleq\sqrt{\delta}X_{x/\delta}^{\delta},\quad \bar{U}_{1,2}^{\delta}(x)\triangleq\bar{U}_{1}^{\delta}(x)-\bar{U}_{1,1}^{\delta}(x)
\end{equation}
and define the linearly interpolated processes $\hat{U}_{1}^{\delta},\hat{U}_{1,1}^{\delta},\hat{U}_{1,2}^{\delta}$
accordingly. We also introduce the trivial rough path lifting 
\begin{equation}\label{eq:hatU1-hatU2}
\hat{{\bf U}}_{1}^{\delta}\triangleq(\hat{U}_{1}^{\delta},\hat{U}^{\delta;2}),\quad\text{ where }\quad \hat{U}_{1}^{\delta;2}(x,y)\triangleq\frac{\hat{U}_{1}^{\delta}(x,y)^{2}}{2}
\end{equation}
as required in the theorem (and recall the lifting ${\bf W}$ of $W$
defined by \eqref{eq:def-w1-w2}). Note that all these processes have
the correct distributions associated with the a priori distribution
of the discrete environment and they are all defined on the same probability
space $(\Omega,\mathcal{F},\mathbf{P})$ where the Brownian motion
$W$ is given.

For now we fix $a\geqslant1$ and we want to apply Theorem \ref{thm:KolCty}
to the random rough paths ${\bf X}=\hat{{\bf U}}_{1}^{\delta}$ and
${\bf Y}={\bf W}$ over $[-a,a]$. To this end, we need to first establish
moment estimates in the form of (\ref{eq:MEstOrg}) and (\ref{eq:MEstErr}).
Recall that $\alpha,\chi$ are given fixed parameters satisfying~(\ref{f1})
and we also fix $\eta,q$ to be such that $\eta\in(0,1/2-\alpha)$
and $q>(1/2-\alpha-\eta)^{-1}\vee2.$ Set $\nu\triangleq1/2-\eta$.
It is easily checked that the parameters $\alpha,q,\nu$ satisfy
the constraints in Theorem \ref{thm:KolCty}. Moreover, we will prove in Lemma~\ref{lem:CtyMoment} below that for all $x,y\in[-a,a]$ we have
\begin{eqnarray}
\|\hat{U}_{1}^{\delta}(x,y)\|_{L^{q}}
&\leqslant&
 Ca^{\eta}|y-x|^{\nu} \label{eq:KolCtyEst1}\\
 \|W(x,y)\|_{L^{q}}
 &\leqslant& 
 Ca^{\eta}|y-x|^{\nu}\label{eq:KolCtyEst1-2} \\
 \|\hat{U}_{1}^{\delta}(x,y)-W(x,y)\|_{L^{q}}
 &\leqslant& 
 Ca^{\eta}\cdot\delta^{\eta}\left(1+\log(a\delta^{-1})\right)\cdot|y-x|^{\nu},\label{eq:KolCtyEst2} \\
 \|\hat{U}_{1}^{\delta;2}(x,y)-W^{2}(x,y)\|_{L^{q/2}}
 &\leqslant& 
 Ca^{2\eta}\cdot\delta^{\eta}\left(1+\log(a\delta^{-1})\right)\cdot|y-x|^{2\nu},\label{eq:KolCtyEst3}
\end{eqnarray}
Applying the above inequalities and Theorem \ref{thm:KolCty}
with 
\[
\rho=a^{\eta},\ \varepsilon=\delta^{\eta}\left(1+\log a\delta^{-1}\right),
\]
we conclude that 
\begin{equation}
\|\hat{U}_{1}^{\delta}-W\|_{\alpha}^{[-a,a]}\leqslant K^{1}(a),\ \|\hat{U}_{1}^{\delta;2}-W^{2}\|_{2\alpha}\leqslant K^{2}(a),\label{eq:e}
\end{equation}
where $K^{1}(a)\in L^{q},K^{2}(a)\in L^{q/2}$ are random variables
such that 
\begin{equation}
\|K^{1}(a)\|_{L^{q}}\leqslant Ca^{1/2-\alpha}\delta^{\eta}\left(1+\log a\delta^{-1}\right),\ \|K^{2}(a)\|_{L^{q/2}}\leqslant Ca^{1-2\alpha}\delta^{\eta}\left(1+\log a\delta^{-1}\right),\label{eq:K12Bound}
\end{equation}
with some constant $C$ depending only on $\alpha,\eta,q.$ Our goal
is to estimate the distance $\rho_{\alpha,\chi}(\hat{{\bf U}}_{1}^{\delta},{\bf W})$
defined by (\ref{eq:RhoUW}). Now plugging (\ref{eq:e}) into the definition \eqref{eq:RhoUW} we easily get
\[
\rho_{\alpha,\chi}(\hat{{\bf U}}_{1}^{\delta},{\bf W})\leqslant K^{1}+K^{2},
\]
where the random constants $K^{1}$, $K^{2}$ are respectively given by 
\[
K^{1}\triangleq\sup_{a\geqslant1}\frac{K^{1}(a)}{a^{\chi}},
\quad\text{ and }\quad 
K^{2}\triangleq\sup_{a\geqslant1}\frac{K^{2}(a)}{a^{2\chi}}.
\]
Since $\alpha,\chi$ satisfy the constraint (\ref{f1}), by slightly
adjusting $\eta$ if necessary, we can rewrite the bound (\ref{eq:K12Bound})
by 
\[
\big\|\frac{K^{1}(a)}{a^{\chi}}\big\|_{L^{p}}\leqslant Ca^{-\zeta}\delta^{\eta},
\quad\text{ and }\quad  
\big\|\frac{K^{2}(a)}{a^{2\chi}}\big\|_{L^{p/2}}\leqslant Ca^{-2\zeta}\delta^{\eta},
\]
where $\zeta$ is any fixed number in $(0,\chi-1/2+\alpha)$ and $C$
depends on $\zeta,\alpha,\eta,p$. Note that we have used a different
index $p$ (instead of $q$ in \eqref{eq:K12Bound}) here.

To proceed further, we first estimate 
\[
\mathbb{P}\left(K^{1}>y\right)=\mathbb{P}\left(\sup_{a\geqslant1}\frac{K^{1}(a)}{a^{\chi}}>y\right)\leqslant\sum_{a=1}^{\infty}\mathbb{P}\left(\frac{K^{1}(a)}{a^{\chi}}>y\right)\leqslant C_{p}y^{-p}\delta^{p\eta}\sum_{a=1}^{\infty}a^{-p\zeta}.
\]
We fix a large $p$ such that the last series is finite. It follows
that 
\begin{align*}
\mathbb{E}[(K^{1})^{q}] & =q\int_{0}^{\delta^{\eta}}y^{q-1}\mathbb{P}\left(K^{1}>y\right)dy+q\int_{\delta^{\eta}}^{\infty}y^{q-1}\mathbb{P}\left(K^{1}>y\right)dy\\
 & \leqslant\delta^{q\eta}+C_{p,\zeta}\delta^{p\eta}q\int_{\delta^{\eta}}^{\infty}y^{q-p-1}dy\leqslant C_{p,q,\zeta}\delta^{q\eta}.
\end{align*}
Therefore, we arrive at 
\[
\|K^{1}\|_{L^{q}}\leqslant C_{p,q,\zeta}\delta^{\eta},
\]
and a similar estimate is proved along the same lines for $K^{2}$. This yields the desired
inequality (\ref{eq:RhoLqEst}) and thus completes the proof of Theorem
\ref{thm:StrongApprox}.
\end{proof}

We close this section by giving the technical lemma which is used in the proof of Theorem~\ref{thm:StrongApprox}.

\begin{lem}\label{lem:CtyMoment}
Recall that the coefficients $\alpha,\beta,\eta,\chi,q$ satisfy $1/3<\beta<\alpha<1/2$, $0<\eta<1/2-\alpha<\chi<\beta/2$ and 
$q>(1/2-\alpha-\eta)^{-1}$. We have also set $1/2-\alpha=\nu$. The processes 
$(\hat{U}_{1}^{\delta},\hat{U}_{1}^{\delta;2})$ have been introduced in 
\eqref{eq:hatU1-hatU2}. Then the following moment estimates hold true, for all $x,y\in[-a,a]$:
\begin{eqnarray}
\|\hat{U}_{1}^{\delta}(x,y)\|_{L^{q}}
&\leqslant&
 Ca^{\eta}|y-x|^{\nu} \label{eq:KolCtyEst1-a}\\
 \|W(x,y)\|_{L^{q}}
 &\leqslant& 
 Ca^{\eta}|y-x|^{\nu}\label{eq:KolCtyEst1-2-a} \\
 \|\hat{U}_{1}^{\delta}(x,y)-W(x,y)\|_{L^{q}}
 &\leqslant& 
 Ca^{\eta}\cdot\delta^{\eta}\left(1+\log(a\delta^{-1})\right)\cdot|y-x|^{\nu},\label{eq:KolCtyEst2-a} \\
 \|\hat{U}_{1}^{\delta;2}(x,y)-W^{2}(x,y)\|_{L^{q/2}}
 &\leqslant& 
 Ca^{2\eta}\cdot\delta^{\eta}\left(1+\log(a\delta^{-1})\right)\cdot|y-x|^{2\nu},\label{eq:KolCtyEst3-a}
\end{eqnarray}
where $C$ is a constant depending only on $q$ and $\eta$.
\end{lem}

\begin{proof}
Consider the partition $\mathcal{P}_{a}^{\delta}\triangleq[-a,a]\cap\delta\mathbb{Z}$.
Let $x<y\in[-a,a].$ We divide this proof in several steps. 

\noindent
\textit{Step 1: Proof of \eqref{eq:KolCtyEst1-a} when $x,y$ belong to the same interval}.  We assume for this step that both $x$ and $y$ sit in 
$I_{m}^{\delta}=[(m-1)\delta,m\delta]$ with some $m\in\mathbb{Z}$.
In this case, by the definition of $\hat{U}_{1,1}^{\delta}$ we have
\[
\hat{U}_{1,1}^{\delta}(x,y)=\frac{y-x}{\delta}\cdot\bar{U}_{1,1}^{\delta}(m\delta)=\frac{y-x}{\sqrt{\delta}}\cdot X_{m}^{\delta},
\]
where $X_{m}^{\delta}$ has been defined in \eqref{eq:DistKMT} and where we recall from \eqref{eq:barU1-barU11-barU12} that $\hat{U}_{1}^{\delta}$ is obtained through the relation $\bar{U}_{1,1}^{\delta}(x)=\sqrt{\delta}X_{x/\delta}^{\delta}$. Since $X_{m}^{\delta}$ has a given fixed distribution and since $\eta$ satisfies $1/2+\eta+\nu=1+\eta-\alpha<1$ we obtain that 
\begin{equation}\label{eq:LqboundU11}
\|\hat{U}_{1,1}^{\delta}(x,y)\|_{L^{q}}=\frac{|y-x|}{\sqrt{\delta}}\|X_{m}^{\delta}\|_{L^{q}}\leqslant C_{q}\frac{|y-x|^{1/2}}{\sqrt{\delta}}|y-x|^{\eta+\nu}\leqslant C_{q}|y-x|^{\nu},
\end{equation}
where the last inequality follows from the fact that $|y-x|\leqslant\delta\leqslant1$.
As far as the process $\hat{U}_{1,2}^{\delta}$ defined by \eqref{eq:u11-u12-delta} is concerned, we invoke the almost sure uniform bound \eqref{eq:uppboU12} to write
\begin{equation}\label{eq:LqboundU12}
\|\hat{U}_{1,2}^{\delta}(x,y)\|_{L^{q}}=\frac{|x-y|}{\delta}\|\bar{U}_{1,2}^{\delta}(m\delta)\|_{L^{q}}\leqslant C\sqrt{\delta}|y-x|\leqslant C|y-x|^{\nu}.
\end{equation}
Gathering \eqref{eq:LqboundU11} and \eqref{eq:LqboundU12}, we have thus proved \eqref{eq:KolCtyEst1-a} for $x,y$ in a generic interval $I_{m}^{\delta}$. In order to check \eqref{eq:KolCtyEst1-2-a} in this interval, we just use the Brownian distribution of $W$. we trivially get 
\begin{align}
\|W(x,y)\|_{L^{q}} & =\mathbb{E}\left[\left|\frac{W(y)-W(x)}{\sqrt{y-x}}\right|^{q}\right]^{1/q}\times\sqrt{y-x}=C_{q}\sqrt{y-x}\nonumber \\
 & =C_{q}|y-x|^{\eta}|y-x|^{\nu}\leqslant C_{q}2^{\eta}a^{\eta}|y-x|^{\nu}\, ,
 \label{eq:Kol2Pf}
\end{align}
which proves \eqref{eq:KolCtyEst1-2-a}. Eventually let us check 
\eqref{eq:KolCtyEst2-a} on the interval $I_{m}^{\delta}$ (inequality 
\eqref{eq:KolCtyEst3-a} is left to the patient reader). To this aim we put together \eqref{eq:KolCtyEst1-a} and \eqref{eq:KolCtyEst1-2-a}. This yields 
\[
\|\hat{U}_{1}^{\delta}(x,y)-W(x,y)\|_{L^{q}}\leqslant\|\hat{U}_{1}^{\delta}(x,y)\|_{L^{q}}+\|W(x,y)\|_{L^{q}}\leqslant C_{q,\eta}\sqrt{y-x}\leqslant C_{q,\eta}\delta^{\eta}|y-x|^{\nu},
\]
from which we conclude \eqref{eq:KolCtyEst2-a}.

\noindent
{\it Step 2: Proof of \eqref{eq:KolCtyEst1-a} when $x,y$ belong to different subintervals.} For simplicity, in what follows we just assume that $x=l\delta$ and
$y=m\delta$ for some $l<m$ (otherwise we just add an extra error
term that falls into the first case). Therefore the quantities of the form $\hat{U}^{\delta}$ will be equal to the terms $\bar{U}^{\delta}$.  As in Step 1, we will divide the estimates into an estimate for $U_{1,1}^{\delta}$ and $U_{1,2}^{\delta}$ separately, according to \eqref{eq:U1=U11+U12}. Now in order to prove \eqref{eq:KolCtyEst1-a} for $\hat{U}_{1,1}^{\delta}$ we go back to~\eqref{eq:u11-u12-delta} and~\eqref{eq:barU1-barU11-barU12}. This allows to write 
\begin{equation}\label{eq:U11hat-rewrite}
\hat{U}_{1,1}^{\delta}(x,y)=\sqrt{\delta}\sum_{j=l+1}^{m}X_{j}^{\delta}=\sqrt{\delta(m-l)}\times\frac{\sum_{j=l+1}^{m}X_{j}^{\delta}}{\sqrt{m-l}}=\sqrt{y-x}\times\frac{\sum_{j=l+1}^{m}X_{j}^{\delta}}{\sqrt{m-l}}.
\end{equation}
Moreover, since $\{X_{j}^{\delta}:j\in\mathbb{Z}\}$ are i.i.d. with a fixed
distribution, according to the central limit theorem we know that
\[
\left\|\frac{\sum_{j=l+1}^{m}X_{j}^{\delta}}{\sqrt{m-l}}\right\|_{L^{q}}\leqslant C_{q},\quad\text{ for all }l<m.
\]
As a result, we have 
\begin{equation}
\|\hat{U}_{1,1}^{\delta}(x,y)\|_{L^{q}}\leqslant C_{q}|y-x|^{\eta}|y-x|^{\nu}\leqslant C_{q,\eta}a^{\eta}|y-x|^{\nu}.\label{eq:Lqnorm-againU11}
\end{equation}
Let us turn now to the process  $\hat{U}_{1,2}^{\delta}$ in \eqref{eq:u11-u12-delta} and \eqref{eq:barU1-barU11-barU12}. Here the main observation is the following: by writing 
\begin{equation}\label{eq:U12-as-sum}
\hat{U}_{1,2}^{\delta}(x,y)=\sum_{j=l+1}^{m}\bar{U}_{1,2}^{\delta}(j\delta)
\end{equation}
it is easily deduced from \eqref{eq:barU1-barU11-barU12} that the right hand side above above is a sum of  i.i.d. centered random variables. Owing to \eqref{eq:uppboU12} we end up with 
\[
\mathbf{E}\left[\left(\hat{U}_{1,2}^{\delta}(x,y)\right)^{2}\right]
%=\sum_{j=l+1}^{m}\mathbb{E}\left[\left(\bar{U}_{1,2}^{\delta}(j\delta
%\right)^{2}\right]
\leqslant C(m-l)\delta^{3}=C\delta^{2}(y-x).
\]
Considering the $L^{2}$-norm and recalling that $x,y\in[-a,a]$ this yields 
\begin{equation}
\|\hat{U}_{1,2}^{\delta}(x,y)\|_{L^{2}}\leqslant C\delta\sqrt{y-x}\leqslant C_{\eta}\delta a^{\eta}|y-x|^{\nu}.\label{eq:L2norm-againU12}
\end{equation}
In addition, as easy  BDG type argument in the martingale increment 
\eqref{eq:U12-as-sum} enables to improve \eqref{eq:L2norm-againU12} into a
$L^{q}$ bound of the form
\begin{equation}
\|\hat{U}_{1,2}^{\delta}(x,y)\|_{L^{q}}\leqslant C_{\eta,q}\delta a^{\eta}|y-x|^{\nu}.\label{eq:Lqnorm-againU12}
\end{equation}
Now combining 
(\ref{eq:Lqnorm-againU11}) and \eqref{eq:Lqnorm-againU12} into \eqref{eq:u11-u12-delta}, we have proved inequality \eqref{eq:KolCtyEst1-a}.

\noindent
{\it Step 3: Proof of \eqref{eq:KolCtyEst2-a}.} For the sake of conciseness, we will only treat the case $x=l\delta$ and $y=m\delta$ for $l<m$ as in 
Step 2. Let us separate again the study of $\hat{U}_{1,1}^{\delta}$ and 
$\hat{U}_{1,2}^{\delta}$. For the term $\hat{U}_{1,1}^{\delta}$ we recast 
relation \eqref{eq:U11hat-rewrite} as 
\[
\hat{U}_{1,1}^{\delta}(x,y)=\sqrt{\delta}S_{l,m},\quad\text{ with }
S_{l,m}=\sum_{j=l+1}^{m}X_{j}^{\delta}. 
\]
Then adopting the notation of \eqref{eq:KMTEnv} we have 
\[
\hat{U}_{1,1}^{\delta}(x,y)-W(x,y)=\sqrt{\delta}\left(S_{l,m}^{\delta}-W_{l,m}^{\delta}\right).
\]
Therefore with $k\triangleq|l|\vee|m|$ we see that 
\[
\big|\hat{U}_{1,1}^{\delta}(x,y)-W(x,y)\big|\leqslant\sqrt{\delta}\max_{-k\leqslant i<j\leqslant k}\big|S_{i,j}^{\delta}-W_{i,j}^{\delta}\big|=:\sqrt{\delta}L^{\delta}(k).
\]
We now estimate $L^{\delta}(k)$ in the following way. First write 
\begin{multline*}
\mathbf{E}\left[(L^{\delta}(k))^{q}\right]  =q\int_{0}^{\infty}y^{q-1}\mathbb{P}\left(L^{\delta}(k)>y\right)dy\\
  =q\int_{0}^{C\log k}y^{q-1}\mathbb{P}\left(L^{\delta}(k)>y\right)dy+q\int_{0}^{\infty}(C\log k+x)^{q-1}\mathbb{P}\left(L^{\delta}(k)>C\log k+x\right)dx \, .
 \end{multline*}
In the first integral above we trivially bound the probability by 1. For the second integral we invoke the tail estimate \eqref{eq:KMTEnv}. We obtain
\[
\mathbf{E}\left[(L^{\delta}(k))^{q}\right] 
\leqslant(C\log k)^{q}+q\int_{0}^{\infty}(C\log k+x)^{q-1}Ke^{-\lambda x}dx
\leqslant C_{q}(1+\log k)^{q}.
\]
It follows that
\begin{equation}
\big\|\hat{U}_{1,1}^{\delta}(x,y)-W(x,y)\big\|_{L^{q}}\leqslant C_{q}\sqrt{\delta}\log k\leqslant C_{q}\delta^{\eta}\left(1+\log a\delta^{-1}\right)|y-x|^{\nu},\label{eq:Lqnorm-U11W}
\end{equation}
where we have used our assumption that $|x-y|\geqslant\delta$ for the second inequality. It remains to estimate the term $\hat{U}_{1,2}^{\delta}$, 
which is easily done as in \eqref{eq:U12-as-sum}-\eqref{eq:Lqnorm-againU12}. We get 
\begin{equation}
\|\hat{U}_{1,2}^{\delta}(x,y)\|_{L^{q}}\leqslant C\delta|y-x|^{\eta}|y-x|^{\nu}\leqslant C_{\eta}\delta a^{\eta}|y-x|^{\nu}.\label{eq:Lqnorm-again-U12}
\end{equation}
Combining (\ref{eq:Lqnorm-U11W}) and (\ref{eq:Lqnorm-again-U12}), we obtain the inequality
(\ref{eq:KolCtyEst2-a}).

\noindent
{\it Step 4. Proof of \eqref{eq:KolCtyEst3-a}.} Here we will limit our analysis to the term $\hat{U}_{1}^{\delta;2}$, the other term being easily handled. Now, for \eqref{eq:KolCtyEst3-a} we simply observe that 
\begin{align*}
\big\|\hat{U}_{1}^{\delta;2}(x,y)-W^{2}(x,y)\big\|_{L^{q/2}} & =\big\|\frac{1}{2}\left(\hat{U}_{1}^{\delta}(x,y)-W(x,y)\right)\left(\hat{U}_{1}^{\delta}(x,y)+W(x,y)\right)\big\|_{L^{q/2}}\\
 & \leqslant\frac{1}{2}\big\|\hat{U}_{1}^{\delta}(x,y)-W(x,y)\big\|_{L^{q}}\left(\big\|\hat{U}_{1}^{\delta}(x,y)\big\|_{L^{q}}+\big\| W(x,y)\big\|_{L^{q}}\right)\\
 & \leqslant C_{q}\cdot a^{\eta}\delta^{\eta}\left(1+\log a\delta^{-1}\right)|y-x|^{\nu}\cdot a^{\eta}|y-x|^{\nu}\\
 & =C_{q}a^{2\eta}\delta^{\eta}\left(1+\log a\delta^{-1}\right)|y-x|^{2\nu}.
\end{align*}
This gives the estimate \eqref{eq:KolCtyEst3-a} and completes the proof
of the lemma.
\end{proof}

\section{Convergence estimates}\label{sec:ConvEst}

This section is devoted to the main steps towards a convergence of
the discrete martingale problems related to Sinai's random walk to
those related to the Brox diffusion. We will first explain our global
strategy in Section \ref{sec:strategy-convergence} and then delve
into the details of computations.

\subsection{\label{sec:strategy-convergence} The main convergence estimates and global strategy}

Recall that in Proposition \ref{prop:discrete-mild-pde-2} we have
obtained a discrete PDE for $v^{\delta}:=\nabla^{\delta}f$. With
this notation in hand, it is readily checked that the discrete PDE in Proposition \ref{prop:discrete-mild-pde-2}
can be spelled out as 
\begin{equation}
v_{t_{k}}^{\delta}(x)=\eta_{t_{k}}^{\delta}(x)-\frac{1}{2}\delta^{3}\cdot\sum_{j=0}^{k-1}\sum_{y\in\delta\mathbb{Z}}\nabla_{x}^{2,\delta}\hat{p}_{t_{j}}^{\delta}(x-y)\mathcal{I}_{t_{k-1}-t_{j}}^{\delta}(x,y),\label{eq:disc-PDE-vdelta}
\end{equation}
where the quantity $\eta_{t_{k}}^{\delta}(x)$ is defined by 
\begin{equation}
\eta_{t_{k}}^{\delta}(x):=\delta\sum_{y\in\delta\mathbb{Z}}\nabla_{x}^{1,\delta}\hat{p}_{t_{k}}^{\delta}(x-y)f_{0}(y)+\delta^{3}\sum_{j=0}^{k-1}\sum_{y\in\delta\mathbb{Z}}\nabla_{x}^{1,\delta}\hat{p}_{t_{j}}^{\delta}(x-y)g_{t_{k-1}-t_{j}}(y).\label{eq:disc-eta-delta}
\end{equation}
In (\ref{eq:disc-PDE-vdelta})-(\ref{eq:disc-eta-delta}), recall that  
 the kernel
$\hat{p}^{\delta}$ is given by \eqref{eq:hat-p-delta}.
Namely, we have 
\begin{equation}
\hat{p}_{t}^{\delta}(x-y)=\frac{1}{\delta}p_{t/\delta^{2}}^{d}\left(\frac{x-y}{\delta}\right),\label{eq:HatPDel}
\end{equation}
where $p^{d}$ is introduced in Remark \ref{rmk:non-rescaled-discrete}.
We have also used the notation 
\begin{equation}
\mathcal{I}_{t}^{\delta}(x,y)=\sum_{z=x}^{y}\frac{1}{2}\,\left(v_{t}^{\delta}(z)+v_{t}^{\delta}(z-\delta)\right)\bar{U}^{\delta}(z),\quad\mbox{ with }\quad\bar{U}^{\delta}(z)\triangleq-2\dot{U}^{\delta}(z),\label{eq:disc-Ideltaxy}
\end{equation}
where $\bar{U}^{\delta}$ and  $\dot{U}^{\delta}$ have been expressed in \eqref{eq:approx-white-noise} and
where the expression $\frac{1}{2}\,(v_{t}^{\delta}(z)+v_{t}^{\delta}(z-\delta))$
stems from our formula \eqref{eq:discr-laplacian-nabla} for $\hat{\nabla}f$.
Indeed, notice that \eqref{eq:discr-laplacian-nabla} can be read
as $\hat{\nabla}f(x)=\frac{1}{2}(\nabla f(x)+\nabla f(x-1))$. Also
observe again that the change from $\dot{U}^{\delta}$ to $\bar{U}^{\delta}$
is for matching the $(-1/2)$-factor in the continuous equation \eqref{eq:bm-in-br-environment}.

Next recall that in Section \ref{RPMP} we have analyzed the corresponding
continuous space-time equation \eqref{eq:FPP}, which can be recast
similarly to \eqref{eq:disc-PDE-vdelta} as 
\begin{eqnarray}
v_{t}(x)&=&
\eta_{t}(x)-\frac{1}{2}\int_{0}^{t}\int_{\mathbb{R}}\partial_{x}^{2}p_{s}(x-y)\left(\int_{x}^{y}v_{t-s}(z)dW(z)\right)dyds  \notag\\
&=&\eta_{t}(x)-\frac{1}{2}(\mathcal{M}\mathcal{V})_{t}(x),\label{eq:CtsFPPStra}
\end{eqnarray}
where the function $\mathcal{M}\mathcal{V}$ is defined by \eqref{e4} in the rough path sense and $\eta$ is given by 
\begin{equation}\label{eq:EtaTerm}
\eta_{t}(x):=\int_{\mathbb{R}}\partial_{x}p_{t}(x-y)f_{0}(y)dy+\int_{0}^{t}\int_{\mathbb{R}}\partial_{x}p_{s}(x-y)g_{t-s}(y)dyds.
\end{equation}
In this section in order to ease notations we will also write 
\begin{equation}\label{eq:W=MV}
\mathcal{W}_{t}(x)=(\mathcal{M}\mathcal{V})_{t}(x).
\end{equation}
Our global aim is to quantify the convergence of $v^{\delta}$ to
$v$ as $\delta\to0$, provided that $\bar{U}^{\delta}$ in (\ref{eq:disc-PDE-vdelta})-(\ref{eq:disc-Ideltaxy})
converges to $\dot{W}$ in (\ref{eq:CtsFPPStra}). More specifically,
the result we wish to achieve is the following.

\begin{prop}
\label{prop:exponent-1/10} [Preliminary version] Let $T$ be a fixed positive time horizon.
Recall that we consider some exponents $\alpha,\beta,\chi,\theta$
fulfilling condition (\ref{f1}) and we have chosen $\gamma=\frac{\alpha-\beta}{4}$.
We consider the norms $\kappa_{\alpha,\chi}$ as in (\ref{f11}) and
$\Theta^{\theta,\lambda}$ as in (\ref{eq:def_normTheta}). Pick $\lambda$
such that 
\begin{equation}\label{eq:choice-lambda-kappa}
\lambda^{-\gamma}\left(\sup_{\delta\in(0,1]}\kappa_{\alpha,\chi}(\hat{{\bf U}}_1^{\delta})+\kappa_{\alpha,\chi}({\bf W})\right)=\frac{1}{4}.
\end{equation}
Then for the processes $v^{\delta},v$ defined respectively by (\ref{eq:disc-PDE-vdelta})
and (\ref{eq:CtsFPPStra}) we have 
\[
d(\tilde{v}^{\delta},v)\lesssim\delta^{1/17}
\]for all sufficiently small $\delta$.
Here $\tilde{v}^{\delta}$ denotes a proper linear interpolation
of $v^{\delta}$ and $d(\tilde{v}^{\delta},v)$ is a suitable rough
path distance between $\tilde{v}^{\delta}$ and $v$.
\end{prop}

Notice that Proposition \ref{prop:exponent-1/10} is stated here rather
informally. The definitions of $\tilde{v}^{\delta}$, $d(\tilde{v}^{\delta},v)$
as well as a more precise formulation of this proposition will be
given in Proposition \ref{prop:MainEst} below. We now elaborate on the global strategy employed to prove
this basic result.

\smallskip

\noindent 
\textit{Step 0: Setting.} Let us start with the discrete process $v^{\delta}$
defined by the equation~(\ref{eq:disc-PDE-vdelta}). In order to
compare it with the continuous rough path $v$, we shall consider
a suitable linear interpolation of $v^{\delta}$, which is denoted
as $\tilde{v}^{\delta}$ and is defined precisely in the following
way. 
\begin{defn}\label{def:linear-interpolation-rough-path-v}
For $t\in\delta^{2}\mathbb{N}_{+}$ and $x\in\mathbb{R}$, we
set 
\begin{equation}\label{eq:vtilde-linear-interpolation}
\tilde{v}_{t}^{\delta}(x)=\frac{x_{2}-x}{\delta}v_{t}^{\delta}(x_{1})+\frac{x-x_{1}}{\delta}v_{t}^{\delta}(x_{2}),
\end{equation}
where $x_{1},x_{2}$ are the adjacent grid points in $\delta\mathbb{Z}$
such that $x\in(x_{1},x_{2}]$. For general $t\geqslant0$, we further
set 
\begin{equation}
\tilde{v}_{t}^{\delta}(x)=\begin{cases}
\frac{t_{2}-t}{\delta^{2}}\tilde{v}_{t_{1}}^{\delta}(x)+\frac{t-t_{1}}{\delta^{2}}\tilde{v}_{t_{2}}^{\delta}(x), & t\geqslant\delta^{2};\\
\tilde{v}_{\delta^{2}}^{\delta}(x), & t<\delta^{2},
\end{cases}\label{eq:PLInterp}
\end{equation}
where $t_{1},t_{2}$ are adjacent grid points in $\delta^{2}\mathbb{N}$
such that $t\in(t_{1},t_{2}]$. 
\end{defn}

We also recall our Notation \ref{not:h20} for the linearly interpolated discrete environment, in particular, the decomposition $\hat{U}^{\delta}=\hat{U}_{1}^{\delta}+\hat{U}_{2}^{\delta}$ where $\hat{U}^\delta_1$ is coupled with the Brownian environment $W$ according to Theorem \ref{thm:StrongApprox} and $\hat{U}_2^\delta$ is merely a deterministic straight line $x\mapsto c_\delta x$ where $c_\delta = O(\delta^{3/2})$ as $\delta\rightarrow0$ (see Lemma \ref{lem:U=U1+U2}).
With this notation and decomposition in hand,
the path $\tilde{v}_{t}^{\delta}(\cdot)$
should now be viewed as a continuous rough path controlled by $\hat{{\bf U}}^{\delta}_1$.
Since $\partial_{W}v_{t}=-2v_{t}$ in the
continuous case, in order to expect a suitable convergence estimate
one should naturally define the Gubinelli derivative of $\tilde{v}_{t}^{\delta}$
as 
\[
\partial_{\hat{{\bf U}}_1^{\delta}}\tilde{v}_{t}^{\delta}\triangleq-2\tilde{v}_{t}^{\delta}.
\]
In addition, let us define $\kappa_{\alpha,\chi}(\hat{{\bf U}}_1^{\delta})$
and $\Theta^{\theta,\lambda}(\tilde{v}^{\delta})$ in exactly the
same way as in the Brownian case by replacing ${\bf W}$ with $\hat{{\bf U}}_1^{\delta}$
and $v$ with $\tilde{v}^{\delta}$ in relevant places. We also recall that the rough
path distance between $\hat{{\bf U}}_1^{\delta}$ and $\mathbf{W}$ is defined by (\ref{eq:RhoUW}).
\begin{comment}
as 
\begin{equation}
\rho(\hat{{\bf U}}_1^{\delta},{\bf W})\triangleq\sup_{a\geqslant1}\left(\frac{\|\hat{U}_1^{\delta}-W\|_{\alpha}^{[-a,a]}}{a^{\chi}}+\frac{\|\hat{U}_1^{\delta;2}-W^{2}\|_{2\alpha}^{[-a,a]}}{a^{2\chi}}\right),\label{def:rho-hatUdelta-W}
\end{equation}
where $\hat{U}_1^{\delta;2}(x,y)\triangleq\frac{1}{2}\hat{U}_1^{\delta}(x,y)^{2}.$
\end{comment}
The \textit{controlled distance} between $\tilde{v}^{\delta}$ and
$v$ is then defined by the following function: 
\begin{multline}
d_{\hat{{\bf U}}_1^{\delta},{\bf W}}(\tilde{v}^{\delta},v)\triangleq\sup_{t\in[0,T]}E^{\theta,\lambda}(a,t)^{-1}\times\left(\llbracket\tilde{v}^{\delta}-v\rrbracket_{\beta/2,\beta}^{[0,t]\times[-a,a]}+\lambda^{-\gamma}\llbracket\partial_{\hat{U}_1^{\delta}}\tilde{v}^{\delta}-\partial_{W}v\rrbracket_{\beta/2,\beta}^{[0,t]\times[-a,a]}\right.\\
\left.+\lambda^{-\gamma}Q(a,t)^{-1}\big\|\mathcal{R}_{\hat{U}_1^{\delta}}^{\tilde{v}_{t}^{\delta}}-\mathcal{R}_{W}^{v_{t}}\big\|_{2\beta}^{[-a,a]}\right).\label{eq:ConDistVTilV}
\end{multline}
The estimation of $d_{\hat{{\bf U}}_1^{\delta},{\bf W}}(\tilde{v}^{\delta},v)$
is the main goal of Proposition~\ref{prop:exponent-1/10}.

To describe the next few steps, let us introduce a slightly nonstandard notation for integer parts. This notation will be useful to analyze  the discretization procedure.

\begin{notation}\label{h2}
Let $t\in\R_{+}$ and $x\in\R$. We denote by $\lfloor t\rfloor$ the largest grid point in $\delta^{2}\N$ such that $\lfloor t\rfloor\le t$. In the same way,  $\lfloor x\rfloor$ is the largest point in $\delta\Z$ such that  $\lfloor x\rfloor \le x$. Writing  $\lfloor \cdot\rfloor_{\N}$ and $\lfloor \cdot\rfloor_{\Z}$ for the usual integer parts, we have
\begin{equation*}
\lfloor t\rfloor =\left\lfloor \frac{t}{\delta^{2}}\right\rfloor_{\N}  \delta^{2},
\quad\text{and}\quad
\lfloor x\rfloor = \left\lfloor \frac{x}{\delta}\right\rfloor_{\Z}  \delta \, .
\end{equation*}
\end{notation}

We now somewhat artificially transform the discrete equation \eqref{eq: discrete-PDE} for $v^{\delta}$ into a continuous equation for 
$\tilde{v}^{\delta}$. That is we write
\begin{equation}\label{eq:JDel}
\tilde{v}_{t}^{\delta}(x)=\mathcal{J}_{t}^{\delta}(x)
-\frac{1}{2}\int_{0}^{t}\int_{\mathbb{R}}\nabla_{x}^{2,\delta}\hat{p}_{\lfloor s\rfloor}^{\delta}(\lfloor x-y\rfloor)
\left(\int_{x}^{y}\tilde{v}_{t-s}^{\delta}(z)d\hat{U}^{\delta}(z)\right)dyds,
\end{equation}
where we recall that $\hat{p}$ is given by \eqref{eq:HatPDel}, and where $\mathcal{J}_{t}^{\delta}(x)$ is \textit{defined} by the difference
of the two sides of the above equation. Also notice that the quantities $\lfloor\cdot\rfloor$ in~\eqref{eq:JDel} are introduced in Notation~\ref{h2}. We obtain the following decomposition 
\begin{equation}
\tilde{v}_{t}^{\delta}(x)-v_{t}(x)=(\mathcal{J}_{t}^{\delta}(x)-\eta_{t}(x))-\frac{1}{2}\left(\mathcal{W}_{t}^{\delta}(x)-\mathcal{W}_{t}(x)\right),\label{eq:disc-vdelta-wdelta-w}
\end{equation}
where $\mathcal{W}^{\delta},\mathcal{W}$ are processes given by (recall our notation \eqref{eq:W=MV})
\begin{eqnarray}
\mathcal{W}_{t}^{\delta}(x)
&\triangleq&
\int_{0}^{t}\int_{\mathbb{R}}\nabla_{x}^{2,\delta}\hat{p}_{\lfloor s\rfloor}^{\delta}(\lfloor x-y\rfloor)\left(\int_{x}^{y}\tilde{v}_{t-s}^{\delta}(z)d\hat{U}^{\delta}(z)\right)dyds,\label{eq:disc-wdelta} \\
\mathcal{W}_{t}(x)
&\triangleq&
\int_{0}^{t}\int_{\mathbb{R}}\partial_{x}^{2}p_{s}(x-y)\left(\int_{x}^{y}v_{t-s}(z)dW(z)\right)dyds=(\mathcal{M}\mathcal{V})_{t}(x).
\label{eq:cont-w}
\end{eqnarray}
Notice that $\mathcal{W}^{\delta},\mathcal{W}$ will be controlled rough paths thanks
to our a priori estimates from Theorem \ref{thm:FPT} and Step 1 below. In the sequel we will also use the following notation for the stochastic integrals in \eqref{eq:disc-vdelta-wdelta-w}-\eqref{eq:disc-wdelta}: 
\begin{equation}\label{eq:stochintegvv}
\mathcal{I}_{s}^{\delta}(x,y)\triangleq\int_{x}^{y}\tilde{v}_{s}^{\delta}(z)d\hat{U}^{\delta}(z),\quad\text{and}\quad
\mathcal{I}_{s}(x,y)\triangleq\int_{x}^{y}v_{s}(z)dW(z).
\end{equation}

With those notational preliminaries in hand, our strategy in order to get the convergence of \eqref{eq:JDel} to \eqref{eq:CtsFPPStra} can be summarized as follows. 

\noindent 
\textit{Step 1: A priori estimates} \textit{for the discrete equation}.
Based on the discrete equation (\ref{eq:disc-PDE-vdelta}), we shall
obtain a uniform estimate on the discrete path $v_{t_{k}}^{\delta}(x)$
under a discrete rough path metric. This will be summarized in Proposition~\ref{prop:trapezoidal}
below. This result shall allow us to write down a corresponding uniform
estimate on $\Theta^{\theta,\lambda}(\tilde{v}^{\delta})$ with respect
to $\delta$. More specifically, one expects that 
\begin{equation}
\sup_{\delta>0}\Theta^{\theta,\lambda}\left(\tilde{v}^{\delta}\right)\leqslant M\label{eq:S1Const}
\end{equation}
where $M$ is a constant depending on $\kappa_{\alpha,\chi}({\bf W})$ as well as the initial data $f_0,g$.

\noindent 
\textit{Step }2: \textit{Estimate} $d_{\hat{{\bf U}}_1^{\delta},\mathbf{W}}(\mathcal{J}^{\delta},\eta)$.
Our next task will be reduced to an analysis of the difference $\mathcal{J}^{\delta}-\eta^{\delta}$,
where $\eta^{\delta}$ is introduced in (\ref{eq:disc-eta-delta}).
Since $\eta_{t}^{\delta}(x)$ and $\eta_{t}(x)$ are close to each
other, it is reasonable to expect that 
\begin{equation}
d_{\hat{{\bf U}}_1^{\delta},{\bf W}}(\mathcal{J}^{\delta},\eta)\leqslant C_{M}\delta^{\gamma},\label{eq:Step2Det}
\end{equation}
with some exponent $\gamma>0$, where $C_{M}$ is a constant depending
on $M$ arising from (\ref{eq:S1Const}).

\noindent 
\textit{Step }3: \textit{Estimate} $d_{\hat{{\bf U}}_1^{\delta},\mathbf{W}}(\cw^{\delta},\cw)$.
For the processes $\cw^{\delta}$ and $\cw$ respectively defined by (\ref{eq:disc-wdelta})
and (\ref{eq:cont-w}), we will express this distance in terms of
$d_{\hat{{\bf U}}_1^{\delta},\mathbf{W}}(\tilde{v}^{\delta},v)$, $\rho(\hat{{\bf U}}_1^{\delta},{\bf W})$
as well as errors coming from the local CLT and $\hat{U}^\delta_2$. In this step, it is important to obtain
a contraction factor (which can be made $<1$) in front of the quantity $d_{\hat{{\bf U}}_1^{\delta},\mathbf{W}}(\tilde{v}^{\delta},v)$. 

We are now in a position to formulate our main convergence estimate for the distance $d_{\hat{{\bf U}}_1^{\delta},\mathbf{W}}(\tilde{v}^{\delta},v)$. Since the fixed point problems also involve an initial data  $f_0$ and an inhomogeneous term $g$, we first introduce the suitable spaces where these functions are assumed to live. 

\begin{defn}\label{def:spaceCrL}
Let $L>0$ and $r\in\mathbb{N}$. We say that a function $f_{0}$ of the spatial variable $x$ is an element of 
$\mathcal{C}_{L}^{r}$ if 
\[
\|f_{0}\|_{\mathcal{C}_{L}^{r}}\triangleq\sup_{a\geqslant1}a^{-L}\sup_{k\leqslant r,\,x\in[-a,a]}\big|D_{x}^{k}f_{0}(x)\big|<\infty.
\]
Similarly, we say that a function $g$ of the space-time parameter $(t,x)$ is an element of $\mathcal{C}_{L}^{r}$ if, denoting by  $D_{t,x}^{ij}g$ the derivative $\frac{\partial^{i+j}g}{\partial t^{i}\partial x^{j}}$, we have
\[
\|g\|_{\mathcal{C}_{L}^{r}}\triangleq\sup_{a\geqslant1}a^{-L}\sup_{\substack{i+j\leqslant r\\
t\in[0,T],\,x\in[-a,a]
}
}\big|D_{t,x}^{ij}g(t,x)\big|<\infty.
\]
\end{defn}

Our main convergence estimate for $d_{\hat{{\bf U}}_1^{\delta},\mathbf{W}}(\tilde{v}^{\delta},v)$ is stated as follows.

\begin{prop}\label{prop:MainEst}
Let $f_0,g\in\mathcal{C}_L^r$ be given functions. Let $\alpha,\beta,\chi,\theta$ be given parameters satisfying (\ref{f1}) and 
let $\beta'\in(\beta,\alpha)$ be also given
fixed. Recall that the distance $d_{\hat{\mathbf{U}}_{1}^{\delta},\mathbf{W}}(\tilde{v}^{\delta},v)$ is defined by (\ref{eq:ConDistVTilV}) with respect to the parameters $\alpha,\beta,\theta$ and any given $\lambda>1$. We introduce two additional quenched variables related to the paths $\hat{\mathbf{U}}_{1}^{\delta}$ and $W$, for a given constant $C_{3}$:
\begin{eqnarray}
\bar{\kappa}(\omega)
&\triangleq&
\sup_{\delta\in(0,1]}\kappa_{\alpha,\chi}(\hat{\mathbf{U}}_{1}^{\delta})+\kappa_{\alpha,\chi}(\mathbf{W})\label{eq:BarKap}\\
\bar{\kappa}'(\omega)
&\triangleq&
\min\big\{\left(4C_{3}\right)^{-2},\left(4C_{3}\right)^{-\frac{3}{\alpha-\beta}}\cdot\bar{\kappa}(\omega)^{-\frac{1}{\alpha-\beta}}\big\}.\label{eq:BarKap'}
\end{eqnarray}
Then there exist positive constants $C_{1},C_{2},C_3$ depending only on these parameters and $T$, such that the following quenched estimate holds true for all $\delta\in(0,\bar{\kappa}'(\omega)]$:
\begin{align}
d_{\hat{\mathbf{U}}_{1}^{\delta},\mathbf{W}}(\tilde{v}^{\delta},v) & \leqslant C_{1}e^{C_{2}\lambda}\left(\|f_{0}\|_{\mathcal{C}_{L}^{3}}+\|g\|_{\mathcal{C}_{L}^{3}}\right)\times \nonumber\\
 & \ \ \ \left[\delta^{\beta}+\left(1+\kappa_{\alpha,\chi}(\hat{\mathbf{U}}_{1}^{\delta})+\kappa_{\alpha,\chi}(\mathbf{W})\right)^2\left(\rho_{\alpha,\chi}(\hat{\mathbf{U}}_{1}^{\delta},\mathbf{W})+\delta^{\frac{\beta'(\beta'-\beta)}{\beta'+\beta}}\right)\right],\label{eq:MainEst}
\end{align}
where $\lambda=\lambda_{\omega}$ is chosen by the relation below:
\begin{equation}\label{eq:Lamb}
C_3 \lambda^{-\frac{\alpha-\beta'}{4}}\bar{\kappa}(\omega)=\frac{1}{4}.
\end{equation}
\end{prop}

Finally, we give the precise formulation of Theorem \ref{SinBro-thm-int} that was stated in the introduction.  

\begin{thm}\label{thm:Main} Let $X^\delta$ be the Sinai type random walk whose law is described by \eqref{SinaiRW-int}, properly rescaled as in Section \ref{rescaled-Sinai-walk}. Respectively, consider the weak solution $X^c$ to the equation \eqref{Broxdiff-int}. Let $\alpha,\beta,\chi,\theta,\beta'$ be given fixed exponents as in Proposition \ref{prop:MainEst}. Then there exists a coupling $(X^\delta,X^c)_{\delta>0}$, such that $X_0^\delta=X^c_0=0$ and for all $t\in[0,T]$, $\delta\in(0,\bar{\kappa}'(\omega)]$, $h\in\mathcal{C}_{L}^{3}$ we have 
\begin{align}
\left|\mathbb{E}^{\omega}\left[h(X_{t}^{c})\right]-\mathbb{E}^{\omega}\left[h(X_{t}^{\delta})\right]\right| & \leqslant C_{1}\exp\left(C_{2}\bar{\kappa}(\omega)^{\frac{4}{\alpha-\beta'}}\right)\|h\|_{\mathcal{C}_{L}^{3}}\nonumber\\
 & \ \ \ \times\left(1+\bar{\kappa}(\omega)\right)^2\times\left(\rho_{\alpha,\chi}(\hat{\mathbf{U}}_{1}^{\delta},\mathbf{W})+\delta^{\frac{\beta'(\beta'-\beta)}{\beta'+\beta}}\right).\label{eq:Main}
\end{align}
where $C_1, C_2$ are constants depending only on the exponents $\alpha,\beta,\chi,\theta,\beta',T$ and  $\bar{\kappa}(\omega)$, $\bar{\kappa}'(\omega)$ are defined by (\ref{eq:BarKap}), (\ref{eq:BarKap'}) respectively.

\end{thm} 

The corollary below completes the link between Theorem \ref{SinBro-thm-int} and Theorem \ref{thm:Main}. 

\begin{cor}
Under the conditions of Theorem \ref{thm:Main}, Theorem \ref{SinBro-thm-int} holds true. Namely for $h\in{\rm C}_{L}^{3}$ (with ${\rm C}_{L}^{3}$ given in Definition \ref{def:spaceCrL}) we have 
\[
\Big|\mathbf{E}_{\omega}\big[h(X_{t}^{c})\big]-\mathbf{E}_{\omega}\big[h(X_{t}^{\delta})\big]\Big|\leqslant C_{h,T}(\omega)\,\delta^{\zeta},
\]
with $\zeta=\frac{9-\sqrt{57}}{24}\simeq\frac{1}{17}\simeq 0.06$. 
\end{cor}
\begin{proof}
To see how Theorem \ref{thm:Main} implies Theorem \ref{SinBro-thm-int}, we first recall from (\ref{eq:KMTAS}) that \[
\rho_{\alpha,\chi}(\hat{\mathbf{U}}_{1}^{\delta},\mathbf{W})\leqslant\Xi(\omega)\delta^{\tau}\ \ \ \forall \delta\in(0,1]
\]for any given fixed $\tau\in (0,1/2-\alpha)$, where $\Xi(\omega)$ is an a.s. finite random variable depending on $\alpha,\chi,\tau$. As a result, the convergence rate in \eqref{eq:Main} is given by $\delta^\zeta$ with
\begin{equation}\label{eq:OpExp}
\zeta=\sup\Big\{\min\big(\tau,Q(\beta,\beta')\big); \, (\tau,\beta,\beta,\alpha)\in\mathcal{D}\Big\} \, ,
\end{equation} 
where the quantity $Q(\beta,\beta')$ is given by 
\[
Q(\beta,\beta')=\frac{\beta'(\beta'-\beta)}{\beta'+\beta},
\]
and where the domain $\mathcal{D}$ is defined by 
\[
\mathcal{D}=\left\{(\tau,\beta,\beta',\alpha)\in\mathbb{R}^{4}:0<\tau<\frac{1}{2}-\alpha,\ \frac{1}{3}<\beta<\beta'<\alpha<\frac{1}{2}\right\}.
\]
The elementary and tedious procedure leading to the computation of \eqref{eq:OpExp} can be summarized as follows. We first optimize the function $Q$, for which we start with a change of variables $\beta'=\alpha-k$, $\beta=\alpha-k-\ell$. This leads to slightly simpler maximisation problem. Namely we wish to find 
\[
\mathcal{A}(\alpha)\equiv\sup\{\hat{Q}(k,\ell):k,\ell\in\hat{\mathcal{D}}_{\alpha}\},
\]
where $\hat{Q}$ and $\hat{\mathcal{D}}_{\alpha}$ are respectively defined by 
\[
\hat{Q}(k,\ell)=\frac{(\alpha-k)\ell}{2(\alpha-k)-\ell},\quad\text{and}\quad
\hat{\mathcal{D}}_{\alpha}=\Big\{(k,\ell)\in\mathbb{R}_{+}^{2}:k+\ell\leq\alpha-\frac{1}{3}\Big\}.
\]
Now an analysis of $\nabla\hat{Q}$ reveals that the supremum of $\hat{Q}$ is attained on the boundary of $\hat{\mathcal{D}}_{\alpha}$. Looking at the values of $\hat{Q}$ on the boundary, we deduce that the supremum is in fact reached at $(k,\ell)=(0,\alpha-1/3)$. hence 
\[
\mathcal{A}(\alpha)=\hat{Q}\big(0,\alpha-\frac{1}{3}\big)=\frac{\alpha\big(\alpha-\frac{1}{3}\big)}{\alpha+\frac{1}{3}}.
\]
Plugging this information back into \eqref{eq:OpExp}, we thus get 
\begin{equation}\label{eq:OpExp-bis}
\zeta=\max\Big\{\min\big(\tau, \mathcal{A}(\alpha)\big) : \, 0<\tau<\frac{1}{2}-\alpha, \, \frac{1}{3}<\alpha<\frac{1}{2}\Big\}. 
\end{equation}
Then one can easily see that the maximum in \eqref{eq:OpExp-bis} is obtained for $\alpha^{*}$ such that $1/2-\alpha^{*}=\mathcal{A}(\alpha^{*})$. This is given by $\alpha^{*}=(3+\sqrt{57})/24\approx 0.42$. The corresponding value of $\zeta$ is 
\[
\zeta=\frac{9-\sqrt{57}}{24}\approx 0.06.
\]
This finishes the proof.
\end{proof}

\begin{proof}[Proof of Theorem \ref{thm:Main}]

This is a rather straightforward corollary of Proposition \ref{prop:MainEst}. Indeed, it is enough to consider $t=T$. Then with $f_0 = h$ and $g = 0$ we have  \[
\mathbb{E}^{\omega}[h(X_{T}^{c})]=v(0),\quad\text{and}\quad \mathbb{E}^{\omega}[h(X_{T}^{\delta})]=\tilde{v}^{\delta}(0).
\]According to the definition of the controlled distance $d_{\hat{\mathbf{U}}_{1}^{\delta},\mathbf{W}}(\tilde{v}^{\delta},v)$, we have 
\begin{align*}
\big|\mathbb{E}^{\omega}[h(X_{T}^{c})]-\mathbb{E}^{\omega}[h(X_{T}^{\delta})]\big| & =\big|\tilde{v}^{\delta}(0)-v(0)\big|\leqslant E^{\theta,\lambda}(a,T)\times d_{\hat{\mathbf{U}}_{1}^{\delta},\mathbf{W}}(\tilde{v}^{\delta},v).
\end{align*}
Note that the above estimate holds for all $a>0$ (since $0\in[-a,a]$) and we can essentially take $a=0$. This makes $E^{\theta,\lambda}(a,T)=e^{\lambda T}$. With the choice of $\lambda$ given by (\ref{eq:Lamb}),  it follows from Proposition \ref{prop:MainEst} that 
\begin{align*}
\big|\mathbb{E}^{\omega}[h(X_{T}^{c})]-\mathbb{E}^{\omega}[h(X_{T}^{\delta})]\big| & \leqslant e^{\lambda T}\times C_{1}e^{C_{2}\lambda}\left(\|h\|_{\mathcal{C}_{L}^{3}}\delta^{\beta}+(1+\bar{\kappa}(\omega))\left(\rho_{\alpha,\chi}(\hat{\mathbf{U}}_{1}^{\delta},\mathbf{W})+\delta^{\frac{\beta'(\beta'-\beta)}{\beta'+\beta}}\right)\right)\\
 & \leqslant C_{3}e^{C_{4}\bar{\kappa}(\omega)^{\frac{4}{\alpha-\beta'}}}(1+\bar{\kappa}(\omega)+\|h\|_{\mathcal{C}_{L}^{3}})\left(\rho_{\alpha,\chi}(\hat{\mathbf{U}}_{1}^{\delta},\mathbf{W})+\delta^{\frac{\beta'(\beta'-\beta)}{\beta'+\beta}}\right).
\end{align*}This completes the proof of Theorem \ref{thm:Main}.
\end{proof}

\subsection{Preliminary notions: discrete controlled
rough paths}\label{sec:discrete-CRP} 

In this section we introduce some notation on discrete rough paths
which will be useful in order to handle the convergence results outlined
in Section \ref{sec:strategy-convergence}. We start by introducing
some space-time partitions. Namely for a given time horizon $T$,
a generic partition $\mathcal{P}_{t}$ will be of the form 
\[
\mathcal{P}_{t}=\{0=t_{0}<t_{1}<\cdots<t_{N-1}<t_{N}=T\}.
\]
Similarly, if $a\in\mathbb{R}_{+}$ a generic partition of $[-a,a]$
can be written as 
\[
\mathcal{P}_{x}=\{-a=x_{0}<x_{1}<\cdots<x_{N}=a\}.
\]
Define $\Delta_{T}^{\mathcal{P}}\triangleq\{(s,t):s\leqslant t,\ s,t\in\mathcal{P}\}$
and $\Delta_{a}^{\mathcal{P}}$ accordingly. Given $(s,t)\in\Delta_{T}^{\mathcal{P}},$
we set $\llbracket s,t\rrbracket\triangleq[s,t]\cap\mathcal{P}.$
Given the above partitions, we will define discrete H\"older norms similarly
to (\ref{eq:def-Ha}) or (\ref{eq:GenHolNorm}). As an example we
set 
\begin{equation}
\|f\|_{\alpha,\beta}^{\mathcal{P}_{t},\mathcal{P}_{a}}\:=\sup\left\{ \frac{|f_{s'}(x')-f_{s}(x)|}{|s'-s|^{\alpha}+|x'-x|^{\beta}};\: s,s'\in\mathcal{P}_{t},\: x,x'\in\mathcal{P}_{x},\: (s,x)\neq(s',x')\right\} .\label{eq:discGenHolNorm}
\end{equation}
Of course the quantity above is always finite, since $\mathcal{P}$
is a finite set. We will omit the superscript $\mathcal{P}$ if the
context is clearly discrete. One of our main tasks will be to bound
quantities like (\ref{eq:discGenHolNorm}) uniformly over a sequence
of partitions whose mesh goes to $0$.
\begin{assumption}
\noindent As in Section \ref{rescaled-Sinai-walk}, we are working here on uniform grids
constructed according to a parabolic scaling. Namely for a discretization
parameter $\delta$ we consider $t_{j}=j\delta^{2}$
and $x_{k}\in\delta\mathbb{Z}$.
\end{assumption}

We now recall some notation concerning discrete rough paths. In the
sequel we consider two H\"older exponents $\alpha,\beta$ such that
$1/3<\beta<\alpha\leqslant1/2$. Next we introduce a discrete augmented
path taking values in $\mathbb{R}$.
\begin{defn}
\label{def:discr-rough-path} Let $X$ be a path defined on $\mathcal{P}_{x}$.
For $s,t\in\Delta_{a}^{\mathcal{P}_{x}}$ we set 
\[
X^{1}(x,y):=X(y)-X(x),\;\mbox{ and }\;X^{2}(x,y):=\frac{1}{2}\left(X(y)-X(x)\right)^{2}.
\]
The norms $\|X^{1}\|_{\alpha}^{[-a,a]}$ and $\|X^{2}\|_{2\alpha}^{[-a,a]}$
are defined similarly to (\ref{eq:def-Ha}). Then the 1D discrete
rough path above $X$ is given by 
\[
{\bf X}(x,y)=(1,X^{1}(x,y),X^{2}(x,y)),\;\mbox{ for }\;(x,y)\in\Delta_{a}^{\mathcal{P}}.
\]
We also mimic expression (\ref{f2}) for controlled paths. Namely
for a path $\mathcal{Y}_{x}=(Y,\partial_{x}Y)$ defined on $\mathcal{P}_{x}$,
we set 
\begin{equation}
\mathcal{R}^{\mathcal{Y}_{x}}:=Y(y)-Y(x)-\partial_{x}Y(x)X^{1}(x,y).\label{eq:RYx}
\end{equation}
Then quantities of interest in order to describe $\mathcal{Y}$ as
a discrete controlled path will be $\|\partial_{x}Y\|_{\beta}^{\mathcal{P}_{x}}$
and $\|\mathcal{R}^{\mathcal{Y}_{x}}\|_{2\beta}^{\mathcal{P}_{x}}$.
\end{defn}

Next, we turn to a definition of discrete spatial rough integrals
suited to our context.
\begin{defn}
\label{def:discr-rough-integr} Let $\mathcal{Y}_{x}=(Y,\partial_{x}Y)$
be a discrete controlled path defined on a grid $\mathcal{P}_{x}$.
Consider $a,b\in\mathbb{R}$ with $a<b$. The discrete rough integral
of $\mathcal{Y}$ with respect to $X$ is a discrete rough path $\mathcal{Z}=(Z,\partial_{x}Z)$,
where for $(x,y)$ in $\Delta_{a,b}^{\mathcal{P}}$ we have 
\begin{equation}
Z(x,y)\triangleq\int_{x}^{y}\mathcal{Y}(u)d{\bf X}(u):=\sum_{\llbracket u,v\rrbracket\in\llbracket x,y\rrbracket}\left(Y(u)X^{1}(u,v)+\partial_{x}Y(u)X^{2}(u,v)\right),\ x,y\in\mathcal{P}_{x}.\label{eq:disc-rough-Z}
\end{equation}
and $\partial_{x}Z(x)=Y(x)$. Observe that in (\ref{eq:disc-rough-Z})
the notation $\llbracket u,v\rrbracket\in\llbracket x,y\rrbracket$
stands for 
\begin{equation}
\{(t_{k},t_{k+1})\in\mathcal{P}_{x}^{2};x\leqslant t_{k}<y\}\label{eq:uv-in-xy}
\end{equation}
\end{defn}

In order to bound discrete integrals like (\ref{eq:disc-rough-Z})
we will resort to a discrete sewing lemma borrowed from \cite[Lemma 2.5]{LT}.
This lemma is recalled here for the sake of completeness.
\begin{lem}
\label{lem:discr-sewing} Let $\mathcal{R}(x,y)$ be an increment
defined on the grid $\mathcal{P}_{x}=\{x_{0}<x_{1}<\cdots<x_{n}\}$.
We assume that $\mathcal{R}(x_{i},x_{i+1})=0$ for all $i=0,\ldots,n-1$.
We also suppose that $\|\delta\mathcal{R}\|_{\mu}^{\mathcal{P}_{x}}<\infty$
for a given $\mu>1$. Then there exists a constant $c_{\mu}$ such
that 
\[
\|\mathcal{R}\|_{\mu}^{\mathcal{P}_{x}}\leqslant c_{\mu}\|\delta\mathcal{R}\|_{\mu}^{\mathcal{P}_{x}}
\]
\end{lem}

We now state a basic upper bound for discrete spatial integrals which
will be used in subsequent computations.
\begin{prop}
\label{prop:DisRIGen} Let $X$ be a discrete rough path as introduced
in Definition \ref{def:discr-rough-path}. Consider a controlled path
$\mathcal{Y}$ whose remainder is given by (\ref{eq:RYx}). The corresponding
integral $\int\mathcal{Y}d{\bf X}$ is expressed in Definition \ref{def:discr-rough-integr}.
We assume that for $1/3<\beta<\alpha\leqslant1/2$ we have 
\[
\|X^{1}\|_{\alpha}+\|X^{2}\|_{2\alpha}+\|\partial_{x}Y\|_{\beta}+\|\mathcal{R}^{\mathcal{Y}_{x}}\|_{2\beta}\leqslant M,
\]
where we drop the superscripts $[-a,a]$ in the norms above for notational
sake. Then for $(x,y)\in\Delta_{a}^{p}$ the following estimate holds
true: 
\begin{multline}
\Big|\int_{x}^{y}\mathcal{Y}(u)d{\bf X}(u)-Y(x)X^{1}(x,y)-\partial_{X}Y(x)X^{2}(x,y)\Big| \\
\leqslant 
C_{\alpha,\beta}\left(\|\partial_{X}Y\|_{\beta}\cdot\|X^{2}\|_{2\alpha}\cdot|y-x|^{2\alpha+\beta}
+\|\mathcal{R}^{\mathcal{Y}_{x}}\|_{2\beta}\cdot\|X^{1}\|_{\alpha}\cdot|y-x|^{\alpha+2\beta}\right),\label{i1}
\end{multline}
where the constant $C_{\alpha,\beta}$ depends only on $\alpha,\beta$
only.
\end{prop}

\begin{proof}
We define two increments $\Xi$ and $\mathcal{I}\Xi$ on ${\mathcal{P}}_{x}$
in the following way: 
\[
\Xi(x,y)\triangleq Y(x)X^{1}(x,y)+\partial_{X}Y(x)X^{2}(x,y),
\]
\[
\mathcal{I}\Xi(x,y)\triangleq\sum_{\llbracket u,v\rrbracket\in\llbracket x,y\rrbracket}\left(Y(u)X^{1}(u,v)+\partial_{X}Y(u)X^{2}(u,v)\right)=\int_{x}^{y}{\mathcal{Y}}(u)d{\mathbf{X}}(u).
\]
Then the left hand side of (\ref{i1}) can be written as a remainder
increment $\mathcal{R}$ of the form 
\[
R(x,y)\triangleq\mathcal{I}\Xi(x,y)-\Xi(x,y).
\]
In this context, recall that $\delta\mathcal{R}$ is defined as a
function of three variables, by 
\[
\delta R(x,u,y)\triangleq R(x,y)-R(x,u)-R(u,y),\ \ \ x<u<y\in\Delta^{\mathcal{P}_{x}}.
\]
Using \eqref{eq:RYx}, it is plain algebra to check that 
\[
|\delta R(x,u,y)|=|\delta\Xi(x,u,y)|\leqslant\|\mathcal{R}^{\mathcal{Y}_{x}}\|_{2\beta}\cdot\|X^{1}\|_{\alpha}\cdot|y-x|^{\alpha+2\beta}+\|Y\|_{\beta}\cdot\|X^{2}\|_{2\alpha}\cdot|y-x|^{2\alpha+\beta}.
\]
In addition, note that $R(x,y)=0$ when $x,y$ are adjacent partition
points. Hence a direct application of Lemma \ref{lem:discr-sewing}
yields 
\[
|R(x,y)|\leqslant K_{\alpha+2\beta}\|\mathcal{R}^{\mathcal{Y}_{x}}\|_{2\beta}\|X^{1}\|_{\alpha}\cdot(y-x)^{\alpha+2\beta}+K_{2\alpha+\beta}\|Y\|_{\beta}\|X^{2}\|_{2\alpha}\cdot(y-x)^{2\alpha+\beta}.
\]
Our claim (\ref{i1}) thus follows.
\end{proof}
In the random walk setting, rough integrals will be approximated by
weighted sums of trapezoidal type. Below we state a proposition bounding
this kind of sum.
\begin{prop}
\label{prop:trapezoidal} Under the same conditions as in Proposition
\ref{prop:DisRIGen}, for $(x,y)\in\Delta_{a,b}^{\mathcal{P}}$ we
set 
\begin{equation}
\mathcal{I}(x,y)\triangleq\sum_{\llbracket u,v\rrbracket\in\llbracket x,y\rrbracket}\frac{1}{2}\,\left(Y(u)+Y(v)\right)\cdot X^{1}(u,v).\label{eq:icalxy}
\end{equation}
Then $\mathcal{I}(x,y)$ enjoys the same property as $\int\mathcal{Y}(u)d{\bf X}(u)$ in Proposition \ref{prop:DisRIGen}. Namely we have
\begin{multline}
\big|\mathcal{I}(x,y)-Y(x)X^{1}(x,y)-\partial_{X}Y(x)X^{2}(x,y)\big|\\\leqslant C_{\alpha,\beta}\left(\|\partial_{X}Y\|_{\beta}\cdot\|X^{2}\|_{2\alpha}\cdot|y-x|^{2\alpha+\beta}
+\|\mathcal{R}^{\mathcal{Y}_{x}}\|_{2\beta}\cdot\|X^{1}\|_{\alpha}\cdot|y-x|^{\alpha+2\beta}\right) .\label{i2}
\end{multline}
\end{prop}

\begin{proof}
Starting from the expression (\ref{eq:icalxy}) for ${\mathcal{I}}(x,y)$,
some elementary algebraic manipulations show that 
\[
{\mathcal{I}}(x,y)=\sum_{\llbracket u,v\rrbracket\in\llbracket x,y\rrbracket}Y(u)X^{1}(u,v)+\frac{1}{2}\,\delta Y(u,v)X^{1}(u,v).
\]
Hence plugging the decomposition (\ref{eq:RYx}) in the expression
above, we get 
\[
{\mathcal{I}}(x,y)=\sum_{\llbracket u,v\rrbracket\in\llbracket x,y\rrbracket}Y(u)X^{1}(u,v)+\frac{1}{2}\,\left(\partial_{X}Y(u)X^{1}(u,v)+{\mathcal{R}}^{{\mathcal{Y}}_{x}}(u,v)\right)X^{1}(u,v).
\]
Recalling the definition (\ref{eq:disc-rough-Z}) of $\int_{x}^{y}\mathcal{Y}(u)d{\bf X}(u)$
it is thus readily seen that 
\begin{equation}
{\mathcal{I}}(x,y)=\int_{x}^{y}\mathcal{Y}(u)d{\bf X}(u)+{\mathcal{J}}(x,y),,\label{eq:caljxy}
\end{equation}
where the term ${\mathcal{J}}(x,y)$ is given by 
\[
{\mathcal{J}}(x,y):=\frac{1}{2}\,\sum_{\llbracket u,v\rrbracket\in\llbracket x,y\rrbracket}{\mathcal{R}}^{{\mathcal{Y}}_{x}}(u,v)X^{1}(u,v).
\]
Now the term $\int_{x}^{y}\mathcal{Y}(u)d{\bf X}(u)$ in (\ref{eq:caljxy})
is upper bounded thanks to (\ref{i1}). Moreover, since $\|{\mathcal{R}}^{{\mathcal{Y}}_{x}}\|_{2\beta}+\|X^{1}\|_{\alpha}<\infty$
and $2\beta+\alpha>1$, we easily get the following estimate for the
term ${\mathcal{J}}$ in (\ref{eq:caljxy}): 
\[
\big|{\mathcal{J}}(x,y)\big|\leqslant C_{\alpha,\beta}\,\|\mathcal{R}^{{\mathcal{Y}}_{x}}\|_{2\beta}\cdot\|X^{1}\|_{\alpha}\cdot|y-x|^{\alpha+2\beta}.
\]
Plugging this estimate and (\ref{i1}) into (\ref{eq:caljxy}), we
have proved our claim (\ref{i2}).
\end{proof}

\subsection{\label{sec:apriori-disc-PDE} Developing Step 1: uniform estimate
on $\tilde{v}^{\delta}$}

In this section, we will give a bound on $\tilde{v}^{\delta}$ uniformly
in $\delta$. Otherwise stated we will achieve \eqref{eq:S1Const} in Step 1 as described in Section~\ref{sec:strategy-convergence}. The key ingredient for this step is the following uniform
estimate on the discrete rough path $v^{\delta}$ with respect to
the discrete metric on the grid.
\begin{prop}
\label{prop:discrete-FPT} Consider a finite time horizon $T$ and $\te>1$. The
exponents $\alpha,\beta,\chi$ satisfy~(\ref{f1}) and we set
$\gamma=\frac{\alpha-\beta}{4}$. Recall that $\hat{U}^\delta_{1}$ is defined in Theorem \ref{thm:StrongApprox} %and $\bar{U}^\delta_{1}$ in \eqref{eq:BarUDecomp}.  
and the quantity $\kappa_{\alpha,\chi}(\hat{{\bf U}}_1^{\delta})$
is defined similarly to (\ref{f11}), albeit on a discrete grid. The
norm $\Theta^{\theta,\lambda}$ given in (\ref{eq:def_normTheta})
is also assumed to be restricted to a discrete setting. The norm $\|\cdot\|_{\mathcal{C}_L^3}$ is defined in Definition \ref{def:spaceCrL}. Let us also set 
\begin{equation}\label{eq:HatKap}
\hat{\kappa}(\omega)\triangleq\sup_{\delta\in(0,1]}\kappa_{\alpha,\chi}(\hat{\mathbf{U}}_{1}^{\delta}),
\end{equation}which is an $a.s.$ finite quenched random variable. Then there
exist  universal constants $C_1,C_2$ depending only on the exponents and $T$, such that the discrete controlled
process $v^{\delta}$ defined by (\ref{eq:disc-PDE-vdelta}) satisfies
\[
\Theta^{\theta,\lambda}(v^{\delta})\leqslant C_{1}\left(\|f_{0}\|_{\mathcal{C}_{L}^{3}}+\|g\|_{\mathcal{C}_{L}^{3}}\right),
\] provided that 
 $\lambda$ is chosen to satisfy
\begin{equation}\label{eq:LamCond}
C_{2}\lambda^{-\frac{\alpha-\beta}{4}}\hat{\kappa}(\omega)=\frac{1}{4}
\end{equation} and $\delta\in (0,\hat{\kappa}'(\omega)]$ where
\begin{equation}\label{eq:HatKap'}
\hat{\kappa}'(\omega)\triangleq\min\big\{\left(4C_{2}\right)^{-2},\left(4C_{2}\right)^{-\frac{2}{\alpha-\beta}}\cdot\hat{\kappa}(\omega)^{-\frac{1}{\alpha-\beta}}\big\}.
\end{equation}
\end{prop}

\begin{proof}
From the $\bar{U}^{\delta}$-decomposition (\ref{eq:BarUDecomp}), we can write the discrete equation
(\ref{eq:disc-PDE-vdelta}) for $v_{\delta}$ as 
\begin{equation}\label{eq:DiscPDENew}
v_{t_{k}}^{\delta}(x)=\eta_{t_{k}}^{\delta}(x)-\frac{1}{2}\left((\mathcal{M}_{1}^{\delta}+\mathcal{M}_{2}^{\delta})v^{\delta}\right)_{t_{k}}(x),
\end{equation}
where 
\begin{align}\label{eq:discrete-Mdeltavdelta}
(\mathcal{M}_{i}^{\delta}v^{\delta})_{t_{k}}(x) & \triangleq\delta^{3}\sum_{j=0}^{k-1}\sum_{y\in\delta\mathbb{Z}}\nabla_{x}^{2,\delta}\hat{p}_{t_{j}}^{\delta}(x-y)\mathcal{I}_{t_{k-1}-t_{j}}^{i,\delta}(x,y)\ \ \ (i=1,2)
\end{align}
and 
\begin{equation}\label{eq:IiDel}
\mathcal{I}_{t_{k-1}-t_{j}}^{i,\delta}(x,y)\triangleq\sum_{z=x}^{y}\frac{1}{2}\left(v_{t}^{\delta}(z)+v_{t}^{\delta}(z-\delta)\right)\cdot\bar{U}_{i}^{\delta}(z).
\end{equation}
As we explained in Lemma \ref{lem:U=U1+U2}, $\hat{U}_{1}^{\delta}$ is regarded as a
discrete rough path approximation of $W$ and $\hat{U}_{2}^{\delta}$
is a remainder. Correspondingly, the estimation of $\mathcal{M}_{1}^{\delta}v^{\delta}$
is a discrete equivalent of Proposition
\ref{prop:EstM}. Most of the computations are adaptations
of what we did in Section \ref{sec:EstMV} for the proof of Proposition \ref{prop:EstM} and we will not repeat all
of them for the sake of conciseness. We will thus focus on proper
adaptations of Lemmas \ref{eq:RIEst}, \ref{lem:KeyLem} and \ref{lem:UniformEst} (the uniform estimate). As far as the term $\mathcal{M}_{2}^{\delta}v^{\delta}$ is concerned, let us mention at this point that its estimation relies on rather trivial arguments. They hinge on the fact that 
\begin{equation}\label{eq:supUbar2}
\sup_{z\in\mathbb{Z}}|\bar{U}_{2}^{\delta}(z)|\equiv c_{\delta}
\lesssim\delta^{3/2}\,. 
\end{equation}
Therefore the quantities $\mathcal{I}_{t_{k-1}-t_{j}}^{\delta}(x,y)$ have to be considered as mere discrete Lebesgue integrals. In the sequel we will also 
consider a family of constants $c' _\delta$ given by 
\begin{equation}\label{eq:constcprimedelta}
c'_{\delta}\equiv \frac{c_{\delta}}{\delta}=O(\delta^{1/2}),
\end{equation}
where the last identity stems from \eqref{eq:supUbar2}. We now divide our proof in several steps. 

\vspace{2mm}\noindent \textit{Step 1: Extension of Lemma \ref{eq:RIEst}.}
Notice that we are considering here discrete equivalents of the quantities
$\Theta^{\theta,\lambda}$ given in Definition \ref{def:controlled-process}.
The corresponding weighted space of discrete controlled path will
be denoted by ${\mathcal{B}}^{\delta,\theta,\lambda}$. Let us now consider an element $v^{\delta}$ in ${\mathcal{B}}^{\delta,\theta,\lambda}$
and set 
\[
\kappa\triangleq\kappa_{\alpha,\chi}(\hat{{\bf U}}_1^{\delta}),\quad\Theta\triangleq\Theta^{\theta,\lambda}(v^{\delta}),\quad E\triangleq E^{\theta,\lambda}(a,t).
\]
Also recall that the quantity $D$ is defined by \eqref{eq:DFunction} and 
$\mathcal{I}_{t}^{1,\delta}(x,y)$ is given by \eqref{eq:IiDel}.
Then we claim that 
\begin{align}
\Big|\mathcal{I}_{t}^{1,\delta}(x,y)-v_{t}^{\delta}(x)X^{1}(x,y)\Big| & \leqslant C_{\beta}\,\kappa\,\Theta E\,\lambda^{\gamma}\,D(a,t,y-x)\quad\mbox{ and }\nonumber \\
\Big|\mathcal{I}_{t}^{1,\delta}(x,y)\Big| & \leqslant C_{\beta}\,\kappa\,\Theta E\,\lambda^{\gamma}\,\left(a^{\chi}|y-x|^{\alpha}+D(a,t,y-x)\right).\label{eq:discrete-RI0-RI1}
\end{align}
With Proposition \ref{prop:trapezoidal} in hand, the proof of (\ref{eq:discrete-RI0-RI1})
is identical to the proof of (\ref{eq:RI0})-(\ref{eq:RI1}) and omitted
here for the sake of conciseness.

\vspace{2mm}\noindent \textit{Step 2: Extension of Lemmas \ref{lem:KeyLem} and \ref{lem:UniformEst}}. To estimate $\mathcal{M}_{1}^{\delta}v^{\delta}$, we divide the summation (\ref{eq:discrete-Mdeltavdelta}) into the two parts: $j>0$ and $j=0$.

\vspace{2mm}\noindent (i) The $j>0$ part is defined by \[
(\mathcal{M}_{1}^{1,\delta}v^{\delta})_{t_{k}}(x)\triangleq\delta^{3}\sum_{j=1}^{k-1}\sum_{y\in\delta\mathbb{Z}}\nabla_{x}^{2,\delta}\hat{p}_{t_{j}}^{\delta}(x-y)\mathcal{I}_{t_{k-1}-t_{j}}^{1,\delta}(x,y).
\]By applying the discrete heat kernel estimate (\ref{eq:DiscUnifGau>0}) together with a change of variables $y=x+\delta\lfloor\delta^{-1}t_{j}^{1/2}w\rfloor$, we have
\begin{equation}\label{eq:j>0DiscM1}
\big|(\mathcal{M}_{1}^{1,\delta}v^{\delta})_{t_{k}}(x)\big|\leqslant C_1\delta^{2}\sum_{j=1}^{k-1}(\delta t_{j}^{-1/2})\sum_{w\in(\delta t_{j}^{-1/2})\cdot\mathbb{Z}}t_j^{-1}e^{-C_2|w|^{2}}\mathcal{A}^\delta_{x,t_{j},w,t_{k}},
\end{equation}
where 
\[
\mathcal{A}^\delta_{x,t_{j},w,t_{k}}\triangleq\big|\mathcal{I}_{t_{k-1}-t_{j}}^{1,\delta}(x,x+\delta\lfloor\delta^{-1}\sqrt{t_{j}}w\rfloor)\big|.
\]The right hand side of (\ref{eq:j>0DiscM1}) resembles its continuous counterpart (\ref{eq:Key1}). As a result, we consider
a stochastic integral term similar to (\ref{eq:Key1}) in our discrete
context. 
%Namely consider $0\leqslant\gamma_{1}\leqslant\gamma_{2}\leqslant\beta/2$,
%$\ell\geqslant1$ and $1<k_{1}\leqslant k_{2}$ such that $t_{k_{2}}\in\llbracket0,T\rrbracket$ (in the current estimate, we take $\gamma_1=\gamma_2=0$). 
For $x\in[-a,a]$ and our time horizon $T$, set $\rho\triangleq a+T^{1/2}|w|$.
Then along the same lines as for (\ref{eq:Key1sInt}) and invoking Proposition \ref{prop:trapezoidal}, we get 
\[
t_{j}^{-1}{\mathcal{A}}_{x,t_{j},w,t_{k}}^{\delta}\leqslant c_{\beta}\kappa\Theta\lambda^{\gamma}E(a,t_{k})P_{u}(w)e^{\theta(1+T)T^{1/2}|w|}e^{-(\lambda+\theta\rho)t_{j}}\phi_{1}(\rho,t_{j}),
\]
where we recall that $P_{u}(w)$ designates an arbitrary polynomial
in $w$ and where the function $\phi_{i}$, $i=1,2$, is defined by
(\ref{eq:phifunction}) (in the current estimate we simply take $\gamma_1=\gamma_2=0$ and thus $\phi_1=\phi_2$). By substituting the above inequality into (\ref{eq:j>0DiscM1}), we obtain that 
\begin{equation}\label{eq:j>0bd}
\big|(\mathcal{M}_{1}^{1,\delta}v^{\delta})_{t_{k}}(x)\big|\leqslant C_{3}\kappa\Theta\lambda^{\gamma}E(a,t_{k})\cdot\delta^{2}\sum_{j=1}^{k-1} \delta t_{j}^{-1/2}\sum_{w\in(\delta t_{j}^{-1/2})\cdot\mathbb{Z}}e^{-C_{3}|w|^{2}}e^{-(\lambda+\theta\rho)t_{j}}\phi_{1}(\rho,t_{j}).
\end{equation}Due to the explicit form of $\phi_1(\rho,t_j)$, a simple monotonicity consideration (replacing each grid point $w = (\delta t_j^{-1/2})\cdot l$ with a generic point in the interval $[(\delta t_j^{-1/2})\cdot l,(\delta t_j^{-1/2})\cdot (l+1)]$) shows that the discrete spatial $w$-integral in (\ref{eq:j>0bd}) is uniformly bounded by its continuous counterpart, namely we have
\[
\delta t_{j}^{-1/2}\sum_{w\in(\delta t_{j}^{-1/2})\cdot\mathbb{Z}}e^{-C_{3}|w|^{2}}e^{-(\lambda+\theta\rho)t_{j}}\phi_{1}(\rho,t_{j})\leqslant C_{4}\int_{\mathbb{R}}e^{-C_{3}|w|^{2}}e^{-(\lambda+\theta\rho)t_{j}}\phi_{1}(\rho,t_{j})dw.
\]It follows that \[
\big|(\mathcal{M}_{1}^{1,\delta}v^{\delta})_{t_{k}}(x)\big|\leqslant C_{5}\kappa\Theta\lambda^{\gamma}E(a,t_{k})\cdot\int_{\mathbb{R}}e^{-C_{3}|w|^{2}}dw\times\delta^{2}\sum_{j=1}^{k-1} e^{-(\lambda+\theta\rho)t_{j}}\phi_{1}(\rho,t_{j}).
\]
Note that the discrete $t_j$-integral on the right hand side is  the discrete counterpart of the quantity (\ref{eq:JancientI}) (with $\tau_1 = t_k$). Now we can perform and estimate the discrete $t_j$-integral in exactly the same way as in the proof of (\ref{eq:Key1I_12}) to conclude that \[
\delta^{2}\sum_{j=1}^{k-1}e^{-(\lambda+\theta\rho)t_{j}}\phi_{1}(\rho,t_{j})\leqslant C_{6}\lambda^{\frac{\beta-\alpha}{2}}.
\]
Therefore, we arrive at the following estimate:
\begin{equation}\label{eq:M11deltavdelta}
\big|(\mathcal{M}_{1}^{1,\delta}v^{\delta})_{t_{k}}(x)\big|\leqslant C_7\kappa\Theta E\lambda^{-\frac{\alpha-\beta}{4}}.
\end{equation}

\noindent (ii) The $j=0$ part in the summation \eqref{eq:discrete-Mdeltavdelta} is defined by \[
(\mathcal{M}_{1}^{0,\delta}v^{\delta})_{t_{k}}(x)\triangleq\delta^{3}\sum_{y\in\delta\mathbb{Z}}\nabla_{x}^{2,\delta}\hat{p}_{0}^{\delta}(x-y)\mathcal{I}_{t_{k-1}}^{1,\delta}(x,y).
\] This is essentially the quantity $\mathcal{A}_t^1(x)$ defined in (\ref{eq:A1tx}) below. By using the discrete heat kernel estimate (\ref{eq:DiscUnifGau=0}), exactly the same argument leading to the estimate (\ref{eq:A1tx-second-estimate}) gives that 
\begin{equation}\label{eq:M10deltavdelta}
\big|(\mathcal{M}_{1}^{0,\delta}v^{\delta})_{t_{k}}(x)\big|\leqslant C_8\kappa\Theta E(a,t_{k})\times\lambda^{\frac{\alpha-\beta}{4}}\delta^{\alpha-2\chi}.
\end{equation}
Gathering the estimates \eqref{eq:M11deltavdelta} and \eqref{eq:M10deltavdelta}, one can thus conclude that 
\begin{equation}\label{eq:M1deltavdelta}
\big|(\mathcal{M}_{1}v^{\delta})_{t_{k}}(x)\big|\leq C\kappa\Theta E\Big(\lambda^{-\frac{\alpha-\beta}{4}}+\lambda^{\frac{\alpha-\beta}{4}}\delta^{\alpha-2\chi}\Big).
\end{equation}

\vspace{2mm}\noindent\textit{Step 3: Estimating }$\mathcal{M}_{2}^{\delta}v^{\delta}.$
By the definition \eqref{eq:discrete-Mdeltavdelta} of $\mathcal{M}_{2}^{\delta}v^{\delta},$
we can write it as 
\begin{equation}
(\mathcal{M}_{2}^{\delta}v^{\delta})_{t_{k}}(x)=\delta^{3}\sum_{j=0}^{k-1}\sum_{y\in\delta\mathbb{Z}}\nabla_{x}^{2,\delta}\hat{p}_{t_{j}}^{\delta}(x-y)\cdot\left(\delta\sum_{z=x}^{y}\frac{v_{t_{k-1}-t_{j}}^{\delta}(z)+v_{t_{k-1}-t_{j}}^{\delta}(z-\delta)}{2}\right)\cdot c'_{\delta},\label{eq:M2Del}
\end{equation}
where we recall that  $c'_{\delta}$ is defined by \eqref{eq:constcprimedelta}, 
and that  $c'_{\delta}=O(\delta^ {1/2})$. Here we view $\mathcal{M}_{2}^{\delta}v^{\delta}$ as a discrete controlled path with respect to $\hat{{\bf U}}_1^\delta$ with zero Gubinelli derivative.
Viewed as a discrete Lebesgue integral, it is easily seen that
\begin{equation}
\left|\delta\sum_{z=x}^{y}\frac{v_{t_{k-1}-t_{j}}^{\delta}(z)+v_{t_{k-1}-t_{j}}^{\delta}(z-\delta)}{2}\right|\leqslant|y-x|\cdot\|v\|_{\infty}^{[0,t_{k}]\times[(-|x|\vee|y|,|x|\vee|y|)]}.\label{eq:DiscLebIntEst}
\end{equation}
We now estimate $(\mathcal{M}_{2}^{\delta}v^{\delta})_{t_{k}}(x)$
(for $t_{k}\in\llbracket0,T\rrbracket$ and $x\in\llbracket-a,a\rrbracket$)
by decomposing the $j$-summation (\ref{eq:M2Del}) into $j=0$ and
$1\leqslant j\leqslant k-1$ as before (since the two parts have different heat kernel
estimates). Firstly, according to the estimate (\ref{eq:DiscUnifGau=0}) for $\nabla_{x}^{2,\delta}\hat{p}_{0}^{\delta}$
as well as (\ref{eq:DiscLebIntEst}), the ``$j=0$'' term in (\ref{eq:M2Del})
is magnified by 
\[
\delta^{3}\sum_{y\in\delta\mathbb{Z}:|y-x|\leqslant\delta}\frac{C}{\delta^{3}}\times\delta\Theta E(a,t_{k})\times c'_{\delta}\leqslant C\delta^{3/2}\Theta E(a,t_{k}).
\]
Similarly, by applying (\ref{eq:DiscUnifGau>0}) together with a change of variables $y=x+\delta\lfloor\delta^{-1}t_{j}^{1/2}w\rfloor$,
the ``$j>0$'' term is estimated as 
\[
C_{1}c'_{\delta}\times\delta^{2}\sum_{j=1}^{k-1}\frac{1}{\sqrt{t_{j}}}\left(\frac{\delta}{\sqrt{t_{j}}}
\sum_{w\in( t_{j}^{-1/2}\delta)\cdot\mathbb{Z}}|w|e^{-C_{2}|w|^{2}}e^{C_{3}|w|}\right)\times\Theta E(a,t_{k}).
\]
The normalised $w$-summation inside the bracket, resembled as a discrete
approximation of the continuous integral $\int_{\mathbb{R}}|w|e^{-C_{2}|w|^{2}+C_{3}|w|}dw$,
is uniformly bounded. The normalised $j$-summation $\delta^{2}\sum_{j=1}^{k-1}1/\sqrt{t_{j}}$,
resembled as a discrete approximation of $\int_{0}^{t}\frac{ds}{\sqrt{s}}$,
is also uniformly bounded. As a consequence, we arrive at the estimate
\begin{equation}\label{eq:M2deltavdelta}
\big|(\mathcal{M}_{2}^{\delta}v^{\delta})_{t_{k}}(x)\big|\leqslant C\sqrt{\delta}\Theta E(a,t_{k}).
\end{equation}
\textit{Step 4: Conclusion.} As mentioned above,
we have focused on proper generalizations of Lemmas \ref{eq:RIEst},
\ref{lem:KeyLem} and \ref{lem:UniformEst} for the sake of conciseness.
The estimations on the time variation, space variation and remainder
terms for $\mathcal{M}_i^{\delta}v^{\delta}$ ($i=1,2$) can be established in the same lines as Lemmas \ref{lem:TimEst}, \ref{lem:SpaEst}, \ref{lem:RemEst} in the continuous case with the aid of the uniform discrete heat kernel estimates (\ref{eq:DiscUnifGau>0}). Note that the $j=0$ part is always handled separately as before by using (\ref{eq:DiscUnifGau=0}) instead. This will result in a factor of $\lambda^{\frac{\alpha-\beta}{4}}\delta^{\alpha-\beta}$ (cf. (\ref{eq:estimA1tau}), the second last term of (\ref{eq:SpVarEst}) and (\ref{eq:final-estimR0}) below, respectively). The main estimates, which are similar to \eqref{eq:M1deltavdelta}-\eqref{eq:M2deltavdelta} are summarised as follows:

\vspace{1mm}\noindent (ii) Time variation estimate:
\[
\begin{cases}
\big|(\mathcal{M}_{1}^{\delta}v^{\delta})_{t_{2}}(x)-(\mathcal{M}_{1}^{\delta}v^{\delta})_{t_{1}}(x)\big|\leqslant C\left(\lambda^{-\frac{\alpha-\beta}{4}}+\lambda^{\frac{\alpha-\beta}{4}}\delta^{\alpha-\beta}\right)\cdot\kappa\Theta E(a,t_{2})a^{\beta/2}\cdot|t_{2}-t_{1}|^{\beta/2};\\
\big|(\mathcal{M}_{2}^{\delta}v^{\delta})_{t_{2}}(x)-(\mathcal{M}_{2}^{\delta}v^{\delta})_{t_{1}}(x)\big|\leqslant C\sqrt{\delta}\cdot\Theta E(a,t_{2})\cdot|t_{2}-t_{1}|^{\beta/2},
\end{cases}
\]
(iii) Space variation estimate:
\[
\begin{cases}
\big|(\mathcal{M}_{2}^{\delta}v^{\delta})_{t}(x')-(\mathcal{M}_{2}^{\delta}v^{\delta})_{t}(x)\big|\leqslant C\left(\lambda^{-\frac{\alpha-\beta}{4}}+\lambda^{\frac{\alpha-\beta}{4}}\delta^{\alpha-\beta}\right)\cdot\kappa\Theta E(a,t)a^{\beta/2}\cdot|x'-x|^{\beta};\\
\big|(\mathcal{M}_{2}^{\delta}v^{\delta})_{t}(x')-(\mathcal{M}_{2}^{\delta}v^{\delta})_{t}(x)\big|\leqslant C\sqrt{\delta}\cdot\Theta E(a,t)\cdot|x'-x|^{\beta},
\end{cases}
\]
(iv) Remainder estimate:
\[
\begin{cases}
\big|\mathcal{R}_{t}^{\mathcal{M}_{1}^{\delta}v^{\delta}}(x,x')\big|\leqslant C\left(1+\lambda^{\frac{\alpha-\beta}{4}}\delta^{\alpha-\beta}\right)\cdot\kappa\Theta E(a,t)Q(a,t)\cdot|x'-x|^{2\beta};\\
\big|\mathcal{R}_{t}^{\mathcal{M}_{2}^{\delta}v^{\delta}}(x,x')\big|=\big|(\mathcal{M}_{2}^{\delta}v^{\delta})_{t}(x')-(\mathcal{M}_{2}^{\delta}v^{\delta})_{t}(x)-0\cdot\hat{U}^{\delta}(x,x')\big|\\
\ \ \ \ \ \ \ \ \ \ \leqslant C\sqrt{\delta}\cdot\Theta E(a,t)t^{-\beta/2}\cdot|x'-x|^{2\beta}\leqslant C\sqrt{\delta}\cdot\Theta E(a,t)Q(a,t)\cdot|x'-x|^{2\beta}.
\end{cases}
\]
As a consequence, we arrive at the desired contraction estimate: 
\[
\Theta^{\theta,\lambda}(\mathcal{M}_{1}^{\delta}v^{\delta}+{\mathcal{M}}_{2}^{\delta}v^{\delta})\leqslant C\left(\lambda^{-\frac{\alpha-\beta}{4}}\kappa+\lambda^{\frac{\alpha-\beta}{4}}\delta^{\alpha-\beta}+\sqrt{\delta}\right)\Theta^{\theta,\lambda}(v^{\delta}).
\]
Now we can choose $\lambda$ to satisfy\[
C\lambda^{-\frac{\alpha-\beta}{4}}\hat{\kappa}(\omega)=\frac{1}{4}
\]where $\hat{\kappa}(\omega)$ is defined by (\ref{eq:HatKap}) and then require that \[
0<\delta<\min\big\{\left(\frac{1}{4C}\right)^{2},\left(\frac{1}{4C}\right)^{\frac{2}{\alpha-\beta}}\hat{\kappa}(\omega)^{-\frac{1}{\alpha-\beta}}\big\}.
\]This ensures that 
\[
C\lambda^{-\frac{\alpha-\beta}{4}}\kappa+C\lambda^{\frac{\alpha-\beta}{4}}\delta^{\alpha-\beta}+C\sqrt{\delta}\leqslant\frac{3}{4},
\]which gives the desired contraction property of $\mathcal{M}^\delta$ and thus concludes our proof. \end{proof}
The result of Proposition \ref{prop:discrete-FPT} easily leads to
the following uniform estimate on the linearly interpolated path $\tilde{v}^{\delta}$ (cf. (\ref{eq:PLInterp})).
\begin{lem}
\label{lem:UniEstPLI} Let $\alpha,\beta,\chi,\theta$ be given as
in Proposition \ref{prop:discrete-FPT}, with $\te>2$. Recall that $\hat{U}_1^{\delta}$ is defined in \eqref{eq:BarUDecomp} and Theorem \ref{thm:StrongApprox}, and that  $\kappa_{\alpha,\chi}(\hat{{\bf U}}_1^{\delta})$ is the
$\alpha$-H\"older continuous rough path norm of $\hat{{\bf U}}_1^{\delta}$.
In addition, $\Theta^{\theta,\lambda}(\tilde{v}^{\delta})$ is the $\beta$-H\"older
continuous rough path norm of $\tilde{v}^{\delta}$ with respect to
$\hat{{\bf U}}_1^{\delta}$. Then there exist universal constants $C_1,C_2$ depending only on the exponents and $T$, such that
\[
\Theta^{\theta,\lambda}(\tilde{v}^{\delta})\leqslant C_1e^{\lambda\delta^{2}}\left(\|f_{0}\|_{\mathcal{C}_{L}^{3}}+\|g\|_{\mathcal{C}_{L}^{3}}\right)
\]for all $\delta\in(0,\hat{\kappa}'(\omega)]$, 
where $\lambda>0$ is chosen to satisfy (\ref{eq:LamCond}) and
$\hat{\kappa}'(\omega)$ is defined by (\ref{eq:HatKap'}).
\end{lem}

\begin{comment}
\begin{rem}
Since $\hat{U}_1^{\delta}\rightarrow W$ with respect to the rough path
metric as $\delta\rightarrow0$, it is clear that $\kappa(\hat{{\bf U}}_1^{\delta})$
is uniformly bounded. As a result, the choice of $\lambda$ can be
made independent of $\delta$.
\end{rem}
\end{comment}

\begin{proof}
According to the definition of $\Theta^{\theta,\lambda}$, we need
to estimate four terms: uniform norm, space-variation, time-variation
and the remainder. Let $c$ be a universal upper bound of the discrete
rough path norm $\Theta^{\theta/2,\lambda}(v^{\delta})$ given by
Proposition \ref{prop:discrete-FPT}.

We first consider the uniform norm estimate. 
Consider $t\in[0,T]$ and $x\in[-a,a]$. We assume that
\begin{equation}\label{i3}
t_{1}\le t< t_{2},
\quad\text{and}\quad
x_{1}\le x< x_{2},
\end{equation}
for $t_{1},t_{2}$ and $x_{1},x_{2}$ adjacent points on the grid $\delta^{2}\N\times\delta\Z$. Since we have chosen $\te>2$, one can apply Proposition~\ref{prop:discrete-FPT} with $\te:=\te/2$. This yields
\begin{equation}\label{i4}
|\tilde{v}_{t}^{\delta}(x)|\leqslant\max\big\{|v_{t_{i}}^{\delta}(x_{j})|:i,j=1,2\big\}\leqslant c\left(\|f_{0}\|_{\mathcal{C}_{L}^{3}}+\|g\|_{\mathcal{C}_{L}^{3}}\right)\cdot E^{\theta/2,\lambda}(a+\delta,t_{2}).
\end{equation}
Furthermore, since we are working with a small $\delta<1$, it is readily checked from~\eqref{eq:EQDef} that
\begin{equation}\label{eq:EAdjust}
E^{\theta/2,\lambda}(a+\delta,t_{2})=e^{\lambda t_{2}+\theta(a+\delta)/2+\theta at_{2}/2}\leqslant C_{\theta}e^{\lambda\delta^{2}}E^{\theta,\lambda}(a,t).
\end{equation}
Reporting  \eqref{eq:EAdjust} into \eqref{i4} we end up with the following uniform bound, valid for $t\in[0,T]$ and $x\in[-a,a]$:
\[
E^{\theta,\lambda}(a,t)^{-1}|\tilde{v}_{t}^{\delta}(x)|\leqslant cC_{\theta}e^{\lambda\delta^{2}}\left(\|f_{0}\|_{\mathcal{C}_{L}^{3}}+\|g\|_{\mathcal{C}_{L}^{3}}\right).
\]

Next, we consider the space variation estimate. To this aim, we pick $t\in[0,T]$ and $x,x'\in[-a,a]$. Similarly to \eqref{i3}, we assume for now that
\begin{equation*}
t_{1}\le t< t_{2},
\quad\text{and}\quad
x_{1}\le x,x'< x_{2} \, ,
\end{equation*}
where $t_{1},t_{2}$ and $x_{1},x_{2}$ are adjacent grid points. By the definition \eqref{eq:PLInterp} of
$\tilde{v}^{\delta}$ it is clear that 
\[
\tilde{v}_{t}^{\delta}(x')-\tilde{v}_{t}^{\delta}(x)=\frac{x-x'}{x_{2}-x_{1}}\left((1-\mu)v_{t_{1}}^{\delta}(x_{1},x_{2})+\mu v_{t_{2}}^{\delta}(x_{1},x_{2})\right),
\]
where $\mu\triangleq\frac{t-t_{1}}{\delta^{2}}$ and $v_{t_{i}}^{\delta}(x_{1},x_{2})\triangleq v_{t_{i}}^{\delta}(x_{2})-v_{t_{i}}^{\delta}(x_{1})$.
Along the same lines as for the uniform bound above, we invoke Proposition \ref{prop:discrete-FPT} with $\te:=\te/2$. We get 
\[
v_{t_{i}}^{\delta}(x_{1},x_{2})\leqslant c\left(\|f_{0}\|_{\mathcal{C}_{L}^{3}}+\|g\|_{\mathcal{C}_{L}^{3}}\right)\cdot E^{\theta/2,\lambda}(a+\delta,t_{i})\lambda^{\frac{\alpha-\beta}{4}}a^{\beta/2}|x_{2}-x_{1}|^{\beta}.
\]
It follows from (\ref{eq:EAdjust}) that 
\begin{align}
\big|\tilde{v}_{t}^{\delta}(x')-\tilde{v}_{t}^{\delta}(x)\big| & \lesssim \left(\|f_{0}\|_{\mathcal{C}_{L}^{3}}+\|g\|_{\mathcal{C}_{L}^{3}}\right)e^{\lambda\delta^{2}}\cdot\big|\frac{x-x'}{x_{2}-x_{1}}\big|\cdot E^{\theta,\lambda}(a,t)\cdot\lambda^{\frac{\alpha-\beta}{4}}a^{\beta/2}\cdot|x_{2}-x_{1}|^{\beta}\nonumber \\
 & \leqslant \left(\|f_{0}\|_{\mathcal{C}_{L}^{3}}+\|g\|_{\mathcal{C}_{L}^{3}}\right)e^{\lambda\delta^{2}}\cdot E^{\theta,\lambda}(a,t)\cdot\lambda^{\frac{\alpha-\beta}{4}}a^{\beta/2}\cdot|x'-x|^{\beta},\label{eq:InterSV}
\end{align}
which is the desired estimate. For general $x<x'$, let $x_{1}$ (respectively,
$x_{2}$) be the smallest (respectively, largest) grid point that
is larger than $x$ (respectively, smaller than $x'$). By considering
the decomposition 
\[
\tilde{v}_{t}^{\delta}(x')-\tilde{v}_{t}^{\delta}(x)=\tilde{v}_{t}^{\delta}(x')-\tilde{v}_{t}^{\delta}(x_{2})+\tilde{v}_{t}^{\delta}(x_{2})-\tilde{v}_{t}^{\delta}(x_{1})+\tilde{v}_{t}^{\delta}(x_{1})-\tilde{v}_{t}^{\delta}(x),
\]
we easily obtain the same type of estimate as in (\ref{eq:InterSV}).

The time variation and remainder estimates are treated in a similar
way. For the remainder, the extra point to note is that for $t\in(t_{1},t_{2}]$,
we have $t_{1}^{-\beta/2}\lesssim t^{-\beta/2}$ if $t\geqslant\delta^{2}$,
while in the case of $t<\delta^{2}$ we do not need the term $t_{1}^{-\beta/2}$
due to our definition of $\tilde{v}_{t}^{\delta}$ (cf. (\ref{eq:PLInterp})).
\end{proof}

\subsection{Developing Step 3: comparing $\cw^{\delta}$ and $\cw$}\label{subsec:S3} Recall that $\cw^{\delta}$ and $\cw$ are respectively defined by (\ref{eq:disc-wdelta})
and (\ref{eq:cont-w}). This section is devoted to establish an estimate for the distance between $\cw^{\delta}$ and $\cw$, as announced at the end of the Section \ref{sec:strategy-convergence}. First in view of the $\hat{U}^{\delta}$-decomposition introduced in \eqref{eq:BarUDecomp}-\eqref{eq:BarU1Decomp}, we write
\[
\mathcal{W}_{t}^{\delta}(x)=\mathcal{W}_{t}^{1,\delta}(x)+\mathcal{W}_{t}^{2,\delta}(x),
\]
where, for $i=1,2$ we have set  
\begin{equation}\label{eq:WiDel}
\mathcal{W}_{t}^{i,\delta}(x)\triangleq\int_{0}^{t}\int_{\mathbb{R}}\nabla_{x}^{2,\delta}\hat{p}_{\lfloor s\rfloor}^{\delta}(\lfloor x-y\rfloor)\left(\int_{x}^{y}\tilde{v}^\delta_{t-s}(z)d\hat{U}_{i}^{\delta}(z)\right)dyds.
\end{equation}
From earlier discussions, it is natural to compare $\mathcal{W}^{1,\delta}$ with $\mathcal{W}$ and view $\mathcal{W}^{2,\delta}$ as a remainder. 
Let us first handle this remainder term in the decomposition \eqref{eq:WiDel}. We label this preliminary step in a lemma.
\begin{lem}\label{lem:remainder-in-WiDel}
For $\delta>0$ let $\mathcal{W}^{2,\delta}$ be the process defined by 
\eqref{eq:WiDel}, where we recall that $\bar{U}_{2}^{\delta}$ is defined by 
\eqref{eq:BarU1Decomp} and $\hat{U}_{2}^{\delta}$ is the linear interpolation of  $\bar{U}_{2}^{\delta}$. We consider a set of parameters 
$\alpha,\beta,\chi,\theta,\lambda,\delta$ as in Proposition \ref{prop:discrete-FPT}. Then there exists a constant $C=C_{\alpha,\beta,\chi,\theta,\lambda,T}$ such 
that the norm of  $\hat{U}_{2}^{\delta}$ as a process controlled by $\hat{U}_{1}^{\delta}$ satisfies 
\begin{equation}\label{eq:controlednorm-lemWiDel}
\Theta_{\hat{U}_{1}^{\delta}}^{\theta,\lambda}(\mathcal{W}^{2,\delta})\leqslant C\left(\|f_{0}\|_{\mathcal{C}_{L}^{3}}+\|g\|_{\mathcal{C}_{L}^{3}}\right)\sqrt{\delta},
\end{equation}
where the subscript $\hat{U}_{1}^{\delta}$ above means that we consider 
$\mathcal{W}^{2,\delta}$ as a process controlled by $\hat{U}_{1}^{\delta}$.
\end{lem} 

\begin{proof}
The term $\mathcal{W}^{2,\delta}$ is very similar to $\mathcal{M}_{2}^{\delta}v^{\delta}$ in \eqref{eq:M2Del}. Hence its analysis ressembles what we did in Steps 4 and 5 for the proof of Proposition \ref{prop:discrete-FPT}. Namely we consider 
$\mathcal{W}^{2,\delta}$ as a controlled path with respect to $\hat{U}_{1}^{\delta}$, with zero Gubinelli derivative. All the integrals have to be treated as discrete Lebesgue integrals. With the aid of the uniform discrete heat kernel bounds \eqref{eq:DiscUnifGau>0}-\eqref{eq:DiscUnifGau=0}, we let the reader check the details leading 
to \eqref{eq:controlednorm-lemWiDel}.
\end{proof}

With Lemma \ref{lem:remainder-in-WiDel} in hand, in what follows we focus on developing the comparison between $\mathcal{W}^{1,\delta}$ and $\mathcal{W}$.
Similar to the strategies for Proposition \ref{prop:EstM}
and Proposition \ref{prop:discrete-FPT}, and recalling that $d_{\hat{{\bf U}}_1^{\delta},{\bf W}}(\cw^{1,\delta},\cw)$
is defined similarly to (\ref{eq:ConDistVTilV}), we will split the
study of this quantity into four parts: the uniform distance, the
time variation distance, the space variation distance, and the remainder
distance. In what follows, to simplify notation we will omit all the
super/subscripts when writing the norms and distances. For instance,
$\kappa(\hat{{\bf U}}_1^{\delta})=\kappa_{\alpha,\chi}(\hat{{\bf U}}_1^{\delta})$.
From time to time, we will use ``$\apprle$'' to denote an inequality
up to a multiplicative constant $C$ that does not depend on $\delta,a,t,x,y.$
The value of the notation $C$ (sometimes denoted as $C_{i}$) may
also differ from line to line.

\subsubsection{The uniform distance estimate}\label{subsec:S3Unif}

In this subsection we mimic Lemma \ref{lem:UniformEst} and get a
uniform estimate for $\cw^{1,\delta}-\cw$.
\begin{lem}
\label{lem:UnifWDelvsW} Recall that the exponents $\alpha,\beta,\chi,\theta,\lambda$
satisfy (\ref{f1}) and we set $\gamma=\frac{\alpha-\beta}{4}$. For
$\delta>0$, $\cw^{1,\delta}$ and $\cw$ are respectively defined by (\ref{eq:disc-wdelta})
and (\ref{eq:cont-w}). Then, for any $a\geqslant1$, $x\in[-a,a]$
and $t\in(0,T]$, we have 
\begin{multline}\label{eq:unif-estim-wdelta-w}
\big|\mathcal{W}_{t}^{1,\delta}(x)-\mathcal{W}_{t}(x)\big| \\
\leqslant C E(a,t)\lambda^{-\frac{\alpha-\beta}{4}}\left(\kappa(\hat{{\bf U}}_1^{\delta})d_{\hat{{\bf U}}_1^{\delta},{\bf W}}(\tilde{v}^{\delta},v)+\Theta(v)\rho(\hat{{\bf U}}_1^{\delta},{\bf W})
  +\kappa(\hat{{\bf U}}_1^{\delta})\Theta(\tilde{v}^{\delta})\delta^{r}\right),
\end{multline}
where $E(a,t)$ is defined by (\ref{eq:EQDef}), $\rho(\hat{{\bf U}}_1^{\delta},{\bf W})$
is defined by (\ref{eq:RhoUW}), $r$ is an arbitrary number
that is less than $\alpha-2\chi$ and $C$ is a positive constant depending only on $\alpha,\beta,T,r$.
%and $C_{2}$ depends additionally on $r$.
\end{lem}

Before we prove Lemma \ref{lem:UnifWDelvsW}, we will state some convergence estimates for the discrete and continuous heat kernels $p$. We start with a uniform bound.

\begin{lem}
\noindent \label{lem:UnifUp2ndDiff} 
Consider $\delta>0$, $s\ge \delta^{2}$ and $w\in\R$. Let 
$\hat{p}^{\delta}$ be the rescaled  kernel from \eqref{eq:hat-p-delta} and $p$ be the Gaussian kernel in \eqref{eq:GauKernelSigma}.
For notational sake we set (see Notation~\ref{h2} for our conventions on integer parts):
\begin{equation}\label{eq:k=00003D000026n}
k \triangleq \frac{\lfloor\sqrt{s}w\rfloor}{\delta}
= \left\lfloor \frac{\sqrt{s}w}{\delta}\right\rfloor_{\Z},
\quad\text{and}\quad 
n \triangleq \frac{\lfloor s\rfloor}{\delta^{2}}
= \left\lfloor \frac{s}{\delta^{2}}\right\rfloor_{\N}.
\end{equation}
Then the following
uniform upper bound holds true for
$w\in\mathbb{R}$:
\begin{equation}
\big|\nabla_{x}^{2,\delta}\hat{p}_{\delta^{2}n}^{\delta}(\delta k)
\big|\vee
\left|\frac{1}{\delta^{3}}\nabla_{k}^{2}p_{n}(k)\right|\leqslant\frac{C_{1}}{s^{3/2}}e^{-C_{2}w^{2}}.
\label{eq:estim-nabla-hat-p}
\end{equation}
In \eqref{eq:estim-nabla-hat-p}, $\nabla^{2,\delta}$ stands for the rescaled gradient given by \eqref{d10} and $\nabla_{k}^{2}p_{n}(k)$ denotes the second discrete gradient applied to the 
heat kernel $p_{t}$ in \eqref{eq:GauKernelSigma}, that is 
\begin{equation}\label{eq:discrete-second-gradient}
\nabla_{k}^{2}p_{n}(k)=p_{n}(k+1)+p_{n}(k-1)-2p_{n}(k). 
\end{equation}
\end{lem}

\begin{rem}
We are stating and proving Lemma \ref{lem:UnifUp2ndDiff} with the second order derivatives in order to handle the most challenging context of interest for
us. However, one can prove non gradient estimates in the same way. Let us label the following one, which holds true under the same conditions as for 
Lemma \ref{lem:UnifUp2ndDiff}: 
\begin{equation}\label{eq:non-gradient-estimate}
\hat{p}^{\delta}_{\delta^{2}n}(\delta k)\leqslant e^{-C_3 w^{2}}, 
\end{equation}
where $k$ and $n$ are given in \eqref{eq:k=00003D000026n}.
\end{rem}

\begin{proof}[Proof of Lemma \ref{lem:UnifUp2ndDiff}]
We start with the following elementary observation, valid for $s\ge\delta^{2}$ and $w\in\R$:
\begin{equation}
s\geqslant\delta^{2},w\in\mathbb{R}\implies\begin{cases}
C_{1}s\leqslant\lfloor s\rfloor\leqslant C_{2}s,\\
C_{3}e^{-C_{4}w^{2}}\leqslant e^{-C_{5}(\lfloor\sqrt{s}w\rfloor/\sqrt{\lfloor s\rfloor})^{2}}\leqslant C_{6}e^{-C_{7}w^{2}},
\end{cases}\label{eq:EquivIntNonInt}
\end{equation}
where the $C_{i}$'s are universal constants whose exact values are irrelevant. With~(\ref{eq:EquivIntNonInt}) in hand, plus recalling that $\lfloor\sqrt{s}w\rfloor=\delta k$ and $\lfloor s\rfloor=\delta^{2} n$, we get that the estimate on  $\nabla_{x}^{2,\delta}\hat{p}_{\delta^{2} n}^{\delta}(\delta k)$ in~\eqref{eq:estim-nabla-hat-p}
follows directly from Proposition~\ref{prop:DiscUnifGau}. Next we observe 
that according to \eqref{eq:GauKernelSigma} we have 
\begin{equation*}
p_{n}(k)=\frac{1}{(2\pi\sigma^{2}n)^{1/2}}\exp\left(-\frac{k^{2}}{2\sigma^{2}n}\right). 
\end{equation*}
Furthermore, owing to \eqref{eq:k=00003D000026n} and setting 
\begin{equation}\label{eq:u-eta}
u\triangleq\frac{\lfloor\sqrt{s}w\rfloor}{\sqrt{\lfloor s\rfloor}} , 
\quad\text{ and }\quad\eta\triangleq\frac{\delta}{\sqrt{\lfloor s\rfloor}}=\frac{1}{\sqrt{n}},
\end{equation}
we have 
\begin{equation*}
\frac{k^{2}}{n}=\frac{\lfloor\sqrt{s}w\rfloor^{2}}{\delta^{2}}
\times \frac{\delta^{2}}{\lfloor s\rfloor}
=u^{2}. 
\end{equation*}
Hence for the quantity $\nabla_{k}^{2}p_{n}(k)$ defined by \eqref{eq:discrete-second-gradient} it is readily checked that 
\begin{equation}
\frac{1}{\delta^{3}}\nabla_{k}^{2}p_{n}(k)=\frac{1}{\sqrt{2\pi\sigma^{2}}\lfloor s\rfloor^{3/2}}\cdot\frac{1}{\eta^{2}}\left(e^{-\frac{1}{2\sigma^{2}}(u+\eta)^{2}}+e^{-\frac{1}{2\sigma^{2}}(u-\eta)^{2}}-2e^{-\frac{u^{2}}{2\sigma^{2}}}\right).\label{eq:2ndDifEven}
\end{equation}
We now bound the right hand side of (\ref{eq:2ndDifEven}). To this
aim, observe that $\eta\leqslant1$ since we have assumed $s\geqslant\delta^{2}$. Hence using a second order Taylor
approximation of the function $e^{-\frac{u^{2}}{2\sigma^{2}}}$ it
is easily seen that 
\begin{equation}\label{eq:estimTaylorexp}
\left|\frac{1}{\eta^{2}}\left(e^{-\frac{1}{2\sigma^{2}}(u+\eta)^{2}}+e^{-\frac{1}{2\sigma^{2}}(u-\eta)^{2}}-2e^{-\frac{u^{2}}{2\sigma^{2}}}\right)\right|\leqslant C_{1}e^{-C_{2}u^{2}}.
\end{equation}
Plugging this estimate in (\ref{eq:2ndDifEven}) and invoking (\ref{eq:EquivIntNonInt})
our claim (\ref{eq:estim-nabla-hat-p}) is easily checked.
\end{proof}

The next lemma quantifies the convergence of derivatives of the heat kernel.

\begin{lem}
\label{lem:Unif2ndDiffvsHK}   As in Lemma \ref{lem:UnifUp2ndDiff}, we consider the Gaussian kernel given by \eqref{eq:GauKernelSigma} and the discrete derivative $\nabla_{k}^{2}$ in \eqref{eq:discrete-second-gradient}. The notation $k$ and $n$ in \eqref{eq:k=00003D000026n} prevails. Then the following bound holds true for all $s\geq\delta^{2}$ and $w\in\mathbb{R}$:  
\[
\Big|\frac{1}{\delta^{3}}\nabla_{k}^{2}p_{n}(k)-\frac{1}{s^{3/2}}\partial_{xx}^{2}p_{1}(w)\Big|\lesssim\frac{\delta}{\sqrt{s}}\cdot\frac{e^{-Cw^{2}}}{s^{3/2}}.
\]
\end{lem}

\begin{proof}
We start from the expression \eqref{eq:2ndDifEven} for 
$\nabla_{k}^{2}p_{n}(k)$, and we substract $\partial_{xx}^{2}p_{1}(u)$ 
(recall that $u$ is introduced in \eqref{eq:u-eta}). Similarly to \eqref{eq:estimTaylorexp}, we introduce a third order taylor expansion 
(as opposed to the second order expansion alluded to in Lemma \ref{lem:UnifUp2ndDiff}). We get that  
\begin{equation}\label{eq:diffTaylorkernel}
\left|\frac{1}{\delta^{3}}\nabla_{k}^{2}p_{n}(k)-\frac{1}{\lfloor s\rfloor^{3/2}}\partial_{xx}^{2}p_{1}(u)\right|\lesssim\frac{\delta}{\sqrt{s}}\frac{e^{-Cw^{2}}}{s^{3/2}} \, ,
\end{equation}
where we have used the relations in (\ref{eq:EquivIntNonInt}). Next,
we evaluate the difference 
\begin{equation}\label{eq:calQsw}
\mathcal{Q}_{s,w}=\frac{1}{s^{3/2}}\partial_{xx}^{2}p_{1}(w)
-\frac{1}{\lfloor s\rfloor^{3/2}}\partial_{xx}^{2}p_{1}(u),
\end{equation} 
where we recall again that $u=\lfloor\sqrt{s}w\rfloor/\sqrt{\lfloor s\rfloor}$. To this aim we define two interpolating paths of the form 
\[
s(r)\triangleq(1-r)s+r\lfloor s\rfloor,
\quad\text{and}\quad
x(r)\triangleq(1-r)\frac{\lfloor\sqrt{s}w\rfloor}{\sqrt{\lfloor s\rfloor}}+rw= (1-r) u+rw.
\]
Differentiating along those paths we have 
\begin{equation}\label{eq:Q=Q1+Q2}
\mathcal{Q}_{s,w}=-\mathcal{Q}^{1}_{s,w}+\mathcal{Q}^{2}_{s,w},
\end{equation}
where the two terms in \eqref{eq:Q=Q1+Q2} are respectively defined by 
\begin{eqnarray*}
\mathcal{Q}^{1}_{s,w}
&=&
\frac{3}{2}\int_{0}^{1}s(r){}^{-5/2}\partial_{xx}^{2}p_{1}(x(r))dr\cdot(s-\lfloor s\rfloor) \\
\mathcal{Q}^{2}_{s,w}&=&
\int_{0}^{1}s(r)^{-3/2}\partial_{x}^{3}p_{1}(x(r))dr\cdot\left(w-u\right).
\end{eqnarray*}
For the term $\mathcal{Q}^{1}_{s,w}$ above, we bound $s(r)^{-5/2}$ by $s^{-5/2}$, the kernel $\partial_{x}^{3}p_{1}(w)$ by $e^{-Cw^{2}}$ and $s-\lfloor s\rfloor$ by $\delta^{2}$. This yields 
\[
\mathcal{Q}^{1}_{s,w}\lesssim\frac{\delta^{2}}{s^{5/2}}e^{-Cw^{2}}.
\]
For the term $\mathcal{Q}^{2}_{s,w}$ in \eqref{eq:Q=Q1+Q2}, we first note that 
\begin{eqnarray*}
\Big|\frac{\lfloor\sqrt{s}w\rfloor}{\sqrt{\lfloor s\rfloor}}-w\Big|  
&\leqslant&
\frac{|\lfloor\sqrt{s}w\rfloor-\sqrt{s}w|}{\sqrt{\lfloor s\rfloor}}+\frac{|w|\cdot|\sqrt{s}-\sqrt{\lfloor s\rfloor}|}{\sqrt{\lfloor s\rfloor}}\\
  &\leqslant&
  \frac{\delta}{\sqrt{\lfloor s\rfloor}}+|w|\cdot\left|\left(1+\frac{s-\lfloor s\rfloor}{\lfloor s\rfloor}\right)^{1/2}-1\right|
  \leqslant
  \frac{\delta}{\sqrt{\lfloor s\rfloor}}+\frac{|w|}{2}\cdot\frac{\delta^{2}}{\lfloor s\rfloor}.
\end{eqnarray*}
As a result, bounding $s(r)^{-3/2}$ and $\partial_{x}^{3}p_{1}$ as for 
$\mathcal{Q}^{1}_{s,w}$, we get
\[
\mathcal{Q}^{2}_{s,w}
\lesssim\left(\frac{\delta}{\sqrt{\lfloor s\rfloor}}+\frac{|w|}{2}\cdot\frac{\delta^{2}}{\lfloor s\rfloor}\right)\frac{e^{-Cw^{2}}}{s^{3/2}}.
\]
Combining the estimates for $\mathcal{Q}^{1}_{s,w}$ and $\mathcal{Q}^{2}_{s,w}$ and recalling the decomposition \eqref{eq:Q=Q1+Q2}, we have thus obtained 
\begin{equation}\label{eq:estim-calQsw}
\mathcal{Q}_{s,w}
\lesssim\frac{\delta}{\sqrt{s}}\frac{e^{-Cw^{2}}}{s^{3/2}}.
\end{equation}
Eventually combining \eqref{eq:diffTaylorkernel}, \eqref{eq:calQsw} and \eqref{eq:estim-calQsw}, this finishes the proof of our lemma.
\end{proof}

We now proceed to prove Lemma \ref{lem:UnifWDelvsW}.

\begin{proof}[Proof of Lemma \ref{lem:UnifWDelvsW}]
For $t\in[0,T]$ and $x\in[-,a,a]$ decompose the difference $\mathcal{W}_{t}^{1,\delta}(x)-\mathcal{W}_{t}(x)$
as 
\begin{equation}
\mathcal{W}_{t}^{1,\delta}(x)-\mathcal{W}_{t}(x)=\mathcal{K}_{t}^{1}(x)+\mathcal{K}_{t}^{2}(x),\label{eq:wdelta-w=00003DK1+K2}
\end{equation}
where $\mathcal{K}_{t}^{1}(x)$ and $\mathcal{K}_{t}^{2}(x)$ are
given by 
\begin{equation}
\mathcal{K}_{t}^{1}(x)\triangleq\int_{0}^{t}\int_{\mathbb{R}}\left(\nabla_{x}^{2,\delta}\hat{p}_{\lfloor s\rfloor}^{\delta}(\lfloor x-y\rfloor)-\partial_{xx}^{2}p_{s}(x-y)\right)\mathcal{I}_{t-s}^{1,\delta}(x,y)dyds,\label{eq:defK1}
\end{equation}
\begin{equation}
\mathcal{K}_{t}^{2}(x)\triangleq\int_{0}^{t}\int_{\mathbb{R}}\partial_{xx}^{2}p_{s}(x-y)\left(\mathcal{I}_{t-s}^{1,\delta}(x,y)-\mathcal{I}_{t-s}(x,y)\right)dyds,\label{eq:defK2}
\end{equation}
and where we recall that 
\[
\mathcal{I}_{t-s}^{1,\delta}(x,y)\triangleq\int_{x}^{y}\tilde{v}_{t-s}^{\delta}(z)d\hat{{\bf U}}_1^{\delta}(z),
\qquad 
\mathcal{I}_{t-s}(x,y)\triangleq\int_{x}^{y}v_{t-s}(z)dW(z).
\]
With this decomposition in hand, we now divide our estimate (\ref{eq:unif-estim-wdelta-w})
in several steps.

\smallskip

\noindent\textit{Step 1: estimate for} $\mathcal{K}_{t}^{2}(x)$. For the term
$\mathcal{K}_{t}^{2}(x)$ in (\ref{eq:defK2}) we easily get that,
for $a\geqslant1$, $t\in(0,T]$ and $x\in[-,a,a]$ we have 
\begin{equation}
|\mathcal{K}_{t}^{2}(x)|\leqslant CE(a,t)\lambda^{-\frac{\alpha-\beta}{4}}\left(\kappa(\hat{{\bf U}}_1^{\delta})d_{\hat{{\bf U}}_1^{\delta},{\bf W}}(\tilde{v}^{\delta},v)+\Theta(v)\rho(\hat{{\bf U}}_1^{\delta},{\bf W})\right).\label{eq:UnifK2}
\end{equation}
Indeed, the patient reader can check that (\ref{eq:UnifK2}) is an
easy variation of Lemma \ref{lem:KeyLem} and Lemma \ref{lem:UniformEst}.

\smallskip

\noindent\textit{Step2 : decomposition for} $\mathcal{K}_{t}^{1}(x)$. In order
to estimate $\mathcal{K}_{t}^{1}(x)$ we will decompose this integral
into small and large time domain. Namely, we write 
\begin{equation}
\mathcal{K}_{t}^{1}(x)=\mathcal{K}_{t}^{11}(x)+\mathcal{K}_{t}^{12}(x),\label{eq:K1=00003DK11+K12}
\end{equation}
where 
\begin{equation}
\mathcal{K}_{t}^{11}(x)=\int_{0}^{t\wedge\delta^{2}}\int_{\mathbb{R}}\left(\nabla_{x}^{2,\delta}\hat{p}_{\lfloor s\rfloor}^{\delta}(\lfloor x-y\rfloor)-\partial_{xx}^{2}p_{s}(x-y)\right)\mathcal{I}_{t-s}^{1,\delta}(x,y)dyds,\label{eq:K11}
\end{equation}
and 
\begin{equation}
\mathcal{K}_{t}^{12}(x)=\int_{t\wedge\delta^{2}}^{t}\int_{\mathbb{R}}\left(\nabla_{x}^{2,\delta}\hat{p}_{\lfloor s\rfloor}^{\delta}(\lfloor x-y\rfloor)-\partial_{xx}^{2}p_{s}(x-y)\right)\mathcal{I}_{t-s}^{1,\delta}(x,y)dyds.\label{eq:K12}
\end{equation}
We now bound $\mathcal{K}^{11}$ and $\mathcal{K}^{12}$ in two different
ways. More specifically, for the small time regime we shall rely on
the fact that $\hat{p}^{\delta}$ and $p$ average to small quantities,
while in the large time regime we will invoke the local central limit
theorem for $\hat{p}^{\delta}$.

\smallskip

\noindent\textit{Step 3: Small time estimates}. For small times, that is when
$s\leqslant\delta^{2}$, we will not use the difference $\hat{p}^{\delta}-p$.
Therefore we simply decompose $\mathcal{K}_{t}^{11}(x)$ into 
\begin{equation}
\mathcal{K}_{t}^{11}(x)=\mathcal{A}_{t}^{1}(x)+\mathcal{A}_{t}^{2}(x),\label{eq:K11=00003DA1+A2}
\end{equation}
with 
\begin{equation}
\mathcal{A}_{t}^{1}(x)=\int_{0}^{t\wedge\delta^{2}}\int_{\mathbb{R}}\big|\nabla_{x}^{2,\delta}\hat{p}_{0}^{\delta}(\lfloor y\rfloor)\big|\cdot\big|\mathcal{I}_{t-s}^{1,\delta}(x,x+y)\big|dyds,\label{eq:A1tx}
\end{equation}
\begin{equation}
\mathcal{A}_{t}^{2}(x)=\int_{0}^{t\wedge\delta^{2}}\int_{\mathbb{R}}s^{-1}\big|\partial_{xx}^{2}p_{1}(w)\big|\cdot\big|\mathcal{I}_{t-s}^{1,\delta}(x,x+\sqrt{s}w)\big|dwds,\label{eq:A2tx}
\end{equation}
where for $\mathcal{A}_{t}^{2}(x)$ we have invoked our usual change
of variable $y=\sqrt{s}w$.

\noindent We first estimate $\mathcal{A}_{t}^{1}(x)$. This is based
on the same kind of decomposition as in the proof of Lemma \ref{lem:KeyLem}
but requires different estimates since the spatial integral is no
longer of Gaussian type. More specifically, note that $\nabla_{x}^{2,\delta}\hat{p}_{0}(\lfloor y\rfloor)$
is supported on $\{y:|y|\leqslant\delta\}$, and replace the Gaussian
bounds on the kernel $\hat{p}^{\delta}$ by the simplified version \eqref{eq:DiscUnifGau=0} of
\eqref{eq:DiscUnifGau>0}: 
\begin{equation}\label{eq:app-estim-deriv--hat-heat-kernel}
\big|\nabla_{x}^{2,\delta}\hat{p}_{0}^{\delta}(\lfloor y\rfloor)\big|\lesssim\frac{1}{\delta^{3}} \, \mathds{1}_{\{|u|\leqslant\delta\}}.
\end{equation}
We now bound the term $\mathcal{I}_{t-s}^{1,\delta}(x,x+y)$ in the
right hand side of (\ref{eq:A1tx}). Namely one gets, similarly to
(\ref{eq:discrete-RI0-RI1}), 
\begin{align}
 & \Big|\mathcal{I}_{t-s}^{1,\delta}(x,x+y)\Big|=\Big|\int_{x}^{x+y}\tilde{v}_{t-s}^{\delta}(z)d\hat{{\bf U}}_1^{\delta}(z)\Big|\nonumber \\
 & \leqslant C\kappa(\hat{{\bf U}}_1^{\delta})\Theta(\tilde{v}^{\delta})\lambda^{\frac{\alpha-\beta}{4}}E(a,t)e^{-(\lambda+\theta(a+|y|))s}\times\left((a+|y|)^{\chi}|y|^{\alpha}+(a+|y|)^{2\chi}|y|^{2\alpha}\right.\nonumber \\
 & \ \ \ \left.+(a+|y|)^{2\chi+\beta/2}|y|^{2\alpha+\beta}+(a+|y|)^{2\chi+\beta/2}|y|^{\alpha+2\beta}+(a+|y|)^{2\chi}(t-s)^{-\beta/2}|y|^{\alpha+2\beta}\right).\label{eq:UniSmallTimDisc}
\end{align}
Applying this inequality to (\ref{eq:A1tx}), we obtain a relation
of the form 
\begin{equation}
\mathcal{A}_{t}^{1}(x)\leqslant\kappa(\hat{{\bf U}}_1^{\delta})\Theta(\tilde{v}^{\delta})\lambda^{\frac{\alpha-\beta}{4}}E(a,t)\sum_{i=1}^{5}\mathcal{A}_{t}^{1i},\label{eq:A11.15}
\end{equation}
and the terms $\mathcal{A}_{t}^{1i}$ have to be estimated separately.
For the sake of conciseness we will just focus on $\mathcal{A}_{t}^{11}$
and $\mathcal{A}_{t}^{15}$ in the sequel.

With (\ref{eq:UniSmallTimDisc}) in hand, the expression we have obtained
for $\mathcal{A}_{t}^{11}$ is 
\[
\mathcal{A}_{t}^{11}=\int_{0}^{t\wedge\delta^{2}}\int_{\mathbb{R}}\big|\nabla_{x}^{2,\delta}\hat{p}_{0}^{\delta}(\lfloor y\rfloor)\big|\cdot e^{-(\lambda+\theta(a+|y|))s}(a+|y|)^{\chi}|y|^{\alpha}dyds.
\]
Hence owing to (\ref{eq:app-estim-deriv--hat-heat-kernel}), the
fact that $\hat{p}_{0}^{\delta}$ is supported in $[-\delta,\delta]$ and the assumption $t\leqslant\delta^{2}$,
we get 
\begin{equation}\label{eq:calA11estimate}
\mathcal{A}_{t}^{11}\lesssim\frac{1}{\delta^{3}}\int_{\{y:|y|\leqslant\delta\}}(a+|y|)^{\chi}|y|^{\alpha}dy\int_{0}^{\delta^{2}}e^{-(\lambda+\theta(a+|y|))s}ds.
\end{equation}
Next in order to eliminate the factor $(a+|y|)^{\chi}$ above, we
apply H\"older's inequality with $p=\chi^{-1}$ and $q=(1-\chi)^{-1}$
to the time integral. This yields 
\[
\int_{0}^{\delta^{2}}e^{-(\lambda+\theta(a+|y|))s}ds\leqslant\left(\int_{0}^{\delta^{2}}e^{-p(\lambda+\theta(a+|y|))s}ds\right)^{1/p}\cdot\left(\int_{0}^{\delta^{2}}1ds\right)^{1/q}\lesssim\frac{\delta^{2(1-\chi)}}{(\lambda+\theta(a+|y|))^{\chi}}.
\]
It follows that 
\begin{equation}
\mathcal{A}_{t}^{11}\lesssim\frac{1}{\delta^{3}}\int_{\{y:|y|\leqslant\delta\}}\frac{(a+|y|)^{\chi}|y|^{\alpha}\delta^{2(1-\chi)}}{(\lambda+\theta(a+|y|))^{\chi}}dy\lesssim\frac{1}{\delta^{3}}\cdot\delta^{2(1-\chi)}\cdot\delta^{\alpha+1}=\delta^{\alpha-2\chi}.\label{eq:A11-first-estimate}
\end{equation}
Referring to our decomposition (\ref{eq:A11.15}), let us now analyze
the term $\mathcal{A}_{t}^{15}$. According to~(\ref{eq:UniSmallTimDisc})
we have 
\[
\mathcal{A}_{t}^{15}=\int_{0}^{t\wedge\delta^{2}}\int_{\mathbb{R}}\big|\nabla_{x}^{2,\delta}\hat{p}_{0}^{\delta}(\lfloor y\rfloor)\big|e^{-(\lambda+\theta(a+|y|))s}(a+|y|)^{2\chi}(t-s)^{-\beta/2}|y|^{\alpha+2\beta}dyds.
\]
We can thus resort to inequality (\ref{eq:EleEst}) and (\ref{eq:app-estim-deriv--hat-heat-kernel})
in order to write 
\[
\mathcal{A}_{t}^{15}\lesssim\frac{1}{\delta^{3}}\int_{0}^{t\wedge\delta^{2}}\int_{\{y:|y|\leqslant\delta\}}c^{c}(\lambda+\theta(a+|y|))^{-c}s^{-c}\cdot(a+|y|)^{2\chi}(t-s)^{-\beta/2}|y|^{\alpha+2\beta}dyds.
\]
By choosing $c\triangleq2\chi$ above, we obtain 
\[
\mathcal{A}_{t}^{15}\lesssim\frac{1}{\delta^{3}}\int_{\{y:|y|\leqslant\delta\}}|y|^{\alpha+2\beta}dy\int_{0}^{t\wedge\delta^{2}}s^{-2\chi}(t-s)^{-\beta/2}ds.
\]
Therefore the elementary change of variables $s=(t\wedge\delta^{2})v$
yields 
\begin{align}
\mathcal{A}_{t}^{15} & \lesssim\frac{1}{\delta^{3}}\int_{\{y:|y|\leqslant\delta\}}|y|^{\alpha+2\beta}dy\cdot(t\wedge\delta^{2})^{1-2\chi-\beta/2}\cdot\int_{0}^{1}v^{-2\chi}(1-v)^{-\beta/2}dv\nonumber \\
 & \lesssim\frac{1}{\delta^{3}}\cdot\delta^{\alpha+2\beta+1}\cdot\delta^{2(1-2\chi-\beta/2)}=\delta^{\alpha+\beta-4\chi}.\label{eq:A15-first-estimate}
\end{align}
The estimates for the space-time integral of the other summands in
(\ref{eq:UniSmallTimDisc}) is similar to the one for $\mathcal{A}_{t}^{11}$.
We let the patient reader check that we obtain 
\begin{equation}
\mathcal{A}_{t}^{12}\lesssim\delta^{2(\alpha-2\chi)},\qquad\mathcal{A}_{t}^{13}\lesssim\delta^{2(\alpha-2\chi)},\qquad\mathcal{A}_{t}^{14}\lesssim\delta^{\alpha+\beta-4\chi},\label{eq:A12.14-first-estimate}
\end{equation}
respectively. Summarizing our considerations on the terms $\mathcal{A}_{t}^{1i}$,
we have obtained relations~(\ref{eq:A11-first-estimate}), (\ref{eq:A15-first-estimate}),
(\ref{eq:A12.14-first-estimate}). Comparing the exponents of $\delta$
in those inequalities and recalling that~(\ref{f1}) imposes $\chi<\beta/2$,
it is clear that the dominant term is given by $\delta^{\alpha-2\chi}$.
Plugging this information back into (\ref{eq:A11.15}), we end up
with 
\begin{equation}
\mathcal{A}_{t}^{1}(x)\lesssim\kappa(\hat{{\bf U}}_1^{\delta})\Theta(\tilde{v}^{\delta})\lambda^{\frac{\alpha-\beta}{4}}E(a,t)\cdot\delta^{\alpha-2\chi}.\label{eq:A1tx-second-estimate}
\end{equation}

We now turn to the term $\mathcal{A}_{t}^{2}$ defined by (\ref{eq:A2tx}).
This quantity is given by an integral involving the continuous kernel
$p_{1}$, and we will thus argue similarly to Lemma \ref{lem:KeyLem}.
Also recall that our main goal here is to extract a factor $\delta^{r}$
with a strictly positive $r$. With this objective in mind, recall
that (as in the proof of Lemma \ref{lem:KeyLem}) the quantity $P_{u}(w)$
designates an arbitrary polynomial in $w$ and that $T$ is our time
horizon. Similarly to (\ref{eq:esti}) we thus get 
\begin{multline}
\mathcal{A}_{t}^{2}(x)\lesssim\kappa(\hat{{\bf U}}_1^{\delta})\Theta(\tilde{v}^{\delta})\lambda^{\frac{\alpha-\beta}{4}}E(a,t)\\
\times\int_{\mathbb{R}}e^{\theta(1+T)\sqrt{T}|w|-Cw^{2}}P_{u}(w)dw\int_{0}^{t\wedge\delta^{2}}e^{-(\lambda+\theta\rho)s}\varphi_{\rho}(s)ds\,,\label{eq:A2tx-first-estimate}
\end{multline}
where $\rho\triangleq a+\sqrt{T}|w|$ and the function $\varphi_{\rho}$
is given by 
\begin{multline}
\varphi_{\rho}(s):=\rho^{\chi}s^{\alpha/2-1}+\rho^{2\chi}s^{\alpha-1}+\rho^{2\chi+\beta/2}s^{\alpha+\beta/2-1}\\
+\rho^{2\chi+\beta/2}s^{\alpha/2+\beta-1}+\rho^{2\chi}(t-s)^{-\beta/2}s^{\alpha/2+\beta-1}=:\sum_{i=1}^{5}\varphi_{\rho}^{i}(s).\label{eq:varphi-rho}
\end{multline}
Plugging (\ref{eq:varphi-rho}) into (\ref{eq:A2tx-first-estimate}),
let us call $\mathcal{A}_{t}^{2i}$ the term corresponding to the
time integral $\int_{0}^{t\wedge\delta^{2}}e^{-(\lambda+\theta\rho)s}\varphi_{\rho}^{i}(s)ds$.
All those quantities are treated very similarly and we only show how
to handle $\mathcal{A}_{t}^{21}$ below. Namely using the expression
for $\varphi_{\rho}^{1}$ in (\ref{eq:varphi-rho}), one can write
$\mathcal{A}_{t}^{21}$ as 
\[
\mathcal{A}_{t}^{21}\triangleq\rho^{\chi}\int_{0}^{t\wedge\delta^{2}}e^{-(\lambda+\theta\rho)s}s^{\alpha/2-1}ds.
\]
We now apply H\"older's inequality with two conjugate exponents $p,q$
in order to get 
\begin{align*}
\mathcal{A}_{t}^{21} & \leqslant\rho^{\chi}\left(\int_{0}^{t\wedge\delta^{2}}e^{-p(\lambda+\theta\rho)s}ds\right)^{1/p}\cdot\left(\int_{0}^{t\wedge\delta^{2}}s^{(\alpha/2-1)q}ds\right)^{1/q}\\
 & \leqslant\rho^{\chi}\left(\int_{0}^{\infty}e^{-p(\lambda+\theta\rho)s}ds\right)^{1/p}\cdot\left(\int_{0}^{\delta^{2}}s^{(\alpha/2-1)q}ds\right)^{1/q}.
\end{align*}
We can easily integrate the two terms in the right hand side above.
We then choose $p\triangleq\chi^{-1},q\triangleq(1-\chi)^{-1}$ in
order to get 
\begin{equation}
\mathcal{A}_{t}^{21}\lesssim\frac{\rho^{\chi}}{(\lambda+\theta\rho)^{\chi}}\,\delta^{\alpha-2+2/q}\lesssim\delta^{\alpha-2\chi}.\label{eq:A21tx-first-estimate}
\end{equation}
In a similar way, the other summands are estimated as 
\begin{equation}
\mathcal{A}_{t}^{22}\vee\mathcal{A}_{t}^{23}\lesssim\delta^{2(\alpha-2\chi)},\quad\text{ and }\quad\mathcal{A}_{t}^{24}\vee\mathcal{A}_{t}^{25}\lesssim\delta^{\alpha+\beta-4\chi}.\label{eq:A2itx-first-estimate}
\end{equation}
Therefore gathering (\ref{eq:A2itx-first-estimate}) and (\ref{eq:A21tx-first-estimate})
into (\ref{eq:A2tx-first-estimate}) we arrive at 
\begin{equation}
\mathcal{A}_{t}^{2}(x)\lesssim\kappa(\hat{{\bf U}}_1^{\delta})\Theta(\tilde{v}^{\delta})\lambda^{\frac{\alpha-\beta}{4}}E(a,t)\cdot\delta^{\alpha-2\chi}.\label{eq:A2tx-second-estimate}
\end{equation}
We can conclude our small time estimate by plugging (\ref{eq:A2tx-second-estimate})
and (\ref{eq:A1tx-second-estimate}) into (\ref{eq:K11=00003DA1+A2}).
This yields the following uper bound, valid for $t\leqslant\delta^{2}$:
\begin{equation}
\mathcal{K}_{t}^{11}(x)\lesssim\kappa(\hat{{\bf U}}_1^{\delta})\Theta(\tilde{v}^{\delta})\lambda^{\frac{\alpha-\beta}{4}}E(a,t)\cdot\delta^{\alpha-2\chi}.\label{eq:K11tx-second-estimate}
\end{equation}

\smallskip

\noindent
\textit{Step 4: Large time estimates.} This step is dedicated to handle
the integral over $s\geqslant\delta^{2}$ defining $\mathcal{K}_{t}^{12}(x)$
in (\ref{eq:K12}). More specifically, resorting to our usual change
of variable $y=x+\sqrt{s}w$, one can recast $\mathcal{K}_{t}^{12}(x)$
as 
\[
\mathcal{K}_{t}^{12}(x)=\int_{t\wedge\delta^{2}}^{t}\int_{\mathbb{R}}\left(\nabla_{x}^{2,\delta}\hat{p}_{\lfloor s\rfloor}^{\delta}(\lfloor\sqrt{s}w\rfloor)-s^{-3/2}\partial_{xx}^{2}p_{1}(w)\right)\mathcal{I}_{t-s}^{1,\delta}(x,x+\sqrt{s}w)\cdot\sqrt{s} \, dwds.
\]
We now insert the quantity $\nabla_{x}^{2}p_{n}(x)$ defined by \eqref{eq:discrete-second-gradient} in the right
hand side above (with $k,n$ defined in 
\eqref{eq:k=00003D000026n}) to get a decomposition of the form 
\begin{equation}\label{eq:upper-decom-calB1calB2}
\Big|\mathcal{K}_{t}^{12}(x)\Big|\leqslant\mathcal{B}_{t}^{1}(x)+\mathcal{B}_{t}^{2}(x),
\end{equation}
with $\mathcal{B}_{t}^{1}(x)$ and $\mathcal{B}_{t}^{2}(x)$ respectively
defined by 
\begin{align}
\mathcal{B}_{t}^{1}(x) & \triangleq\int_{t\wedge\delta^{2}}^{t}\int_{\mathbb{R}}\big|\nabla_{x}^{2,\delta}\hat{p}_{\lfloor s\rfloor}^{\delta}(\lfloor\sqrt{s}w\rfloor)-\frac{1}{\delta^{3}}\nabla_{k}^{2}p_{n}(k)\big|\cdot\big|\mathcal{I}_{t-s}^{1,\delta}(x,x+\sqrt{s}w)\big|\cdot\sqrt{s}\,dwds,\label{eq:UnifB1}\\
\mathcal{B}_{t}^{2}(x) & \triangleq\int_{t\wedge\delta^{2}}^{t}\int_{\mathbb{R}}\big|\frac{1}{\delta^{3}}\nabla_{k}^{2}p_{n}(k)-s^{-3/2}\partial_{xx}^{2}p_{1}(w)\big|\cdot\big|\mathcal{I}_{t-s}^{1,\delta}(x,x+\sqrt{s}w)\big|\cdot\sqrt{s}\,dwds.\label{eq:UnifB2}
\end{align}
Let us proceed to bound the term $\mathcal{B}_{t}^{1}(x)$ in \eqref{eq:UnifB1}. To this aim, we start by recalling the definition 
\eqref{eq:hat-p-delta} of $\hat{p}_{t}^{\delta}$ and \eqref{d10} for 
$\nabla_{x}^{2,\delta}$. This gives 
\begin{equation}\label{d10rep}
\nabla_{x}^{2,\delta} \hat{p}_{\lfloor s\rfloor}^{\delta}(x)
=\frac{1}{\delta^{3}}\nabla_{x}^{2}p_{\lfloor s\rfloor/\delta^{2}}^{d}\!\left(\frac{x}{\delta}\right),
\end{equation}
where $\nabla_{x}^{2}$ in the right hand side above stands for the discrete derivative. In addition, Theorem~\ref{thm:LCLT} states
that for $s\geqslant \delta^{2}$
\[
\Big|\nabla_{x}^{2,\delta}p_{\lfloor s\rfloor/\delta^{2}}^{d}\left(\frac{x}{\delta}\right)-\nabla_{x}^{2}p_{\lfloor s\rfloor/\delta^{2}}\left(\frac{x}{\delta}\right)\Big|\leqslant\frac{C_{1}}{(\lfloor s\rfloor/\delta^{2})^{5/2}}\lesssim\frac{\delta^{5}}{\lfloor s\rfloor^{5/2}}.
\]
Plugging this information into \eqref{d10rep} we get 
\[
\Big|\nabla_{x}^{2,\delta} \hat{p}_{\lfloor s\rfloor}^{\delta}(x)
-\frac{1}{\delta^{3}}\nabla_{x}^{2}p_{\lfloor s\rfloor/\delta^{2}}\left(\frac{x}{\delta}\right)\Big|\lesssim\frac{\delta^{2}}{\lfloor s\rfloor^{5/2}}.
\]
In particular, for $x=\lfloor\sqrt{s}w\rfloor$ and $k,n$ defined in 
\eqref{eq:k=00003D000026n} we end up with 
\begin{equation}\label{eq:appliLCLT}
\Big|\nabla_{x}^{2,\delta}\hat{p}_{\lfloor s\rfloor}^{\delta}(\lfloor\sqrt{s}w\rfloor)-\frac{1}{\delta^{3}}\nabla_{k}^{2}p_{n}(k)\Big|\lesssim\frac{\delta^{2}}{\lfloor s\rfloor^{5/2}}.
\end{equation}
Next we improve our upper bound \eqref{eq:appliLCLT} by introducing an extra interpolating parameter $b\in[0,1]$. Whenever $s\geqslant\delta^{2}$ , combining \eqref{eq:appliLCLT} and the uniform bounds in Lemma~\ref{lem:UnifUp2ndDiff}, we get the following inequality for all $w\in\mathbb{R}$:
\begin{align*}
 & \Big|\nabla_{x}^{2,\delta}\hat{p}_{\lfloor s\rfloor}^{\delta}(\lfloor\sqrt{s}w\rfloor)-\frac{1}{\delta^{3}}\nabla_{k}^{2}p_{n}(k)\Big|\\
 & \leqslant\Big|\nabla_{x}^{2,\delta}\hat{p}_{\lfloor s\rfloor}^{\delta}(\lfloor\sqrt{s}w\rfloor)-\frac{1}{\delta^{3}}\nabla_{k}^{2}p_{n}(k)\Big|^{1-b}\cdot\left(\Big|\nabla_{x}^{2,\delta}\hat{p}_{\lfloor s\rfloor}^{\delta}(\lfloor\sqrt{s}w\rfloor)\Big|^{b}+\Big|\frac{1}{\delta^{3}}\nabla_{k}^{2}p_{n}(k)\Big|^{b}\right)\\
 & \lesssim\left(\frac{\delta^{2}}{\lfloor s\rfloor^{5/2}}\right)^{1-b}\cdot\left(\frac{1}{s^{3/2}}e^{-C_{1}w^{2}}\right)^{b}
 \lesssim\frac{\delta^{2(1-b)}}{s^{5/2-b}}e^{-C_{2}w^{2}}.
\end{align*}
Reporting this inequality in \eqref{eq:UnifB1}, it follows that 
\begin{equation}\label{eq:estim-calB1}
\mathcal{B}_{t}^{1}(x)\lesssim\int_{0}^{t}\int_{\mathbb{R}}\frac{\delta^{2(1-b)}}{s^{2-b}}e^{-Cw^{2}}\cdot\big|\mathcal{I}_{t-s}^{1,\delta}(x,x+\sqrt{s}w)\big|dwds.
\end{equation}
It remains to handle the term $\mathcal{I}_{t-s}^{1,\delta}(x,x+\sqrt{s}w)$ in the right hand side of \eqref{eq:estim-calB1}. Now recall that $\mathcal{I}_{t-s}^{1,\delta}$ is the stochastic integral defined by \eqref{eq:stochintegvv}.  As in the proof of Lemma \ref{lem:KeyLem} (see also the estimates after \eqref{eq:j>0DiscM1}), we shall upper bound bound this quantity by five terms:
\begin{equation}\label{eq:fiveterms}
\big|\mathcal{I}_{t-s}^{1,\delta}(x,x+\sqrt{s}w)\big|\lesssim\kappa(\hat{{\bf U}}_1^{\delta})\Theta(\tilde{v}^{\delta})\lambda^{\frac{\alpha-\beta}{4}}E(a,t)\cdot P_{u}(w)(\mathcal{B}_{s}^{11}+\cdots+\mathcal{B}_{s}^{15}),
\end{equation}
where $P_{u}(w)$ is some polynomial in $|w|$. For simplicity, we
only discuss the first term 
\[
\mathcal{B}_{s}^{11}=e^{-(\lambda+\theta\rho)s}s^{\alpha/2}\rho^{\chi},\quad\text{where}\quad \rho\triangleq a+\sqrt{T}|w|.
\]
Reporting this definition in \eqref{eq:fiveterms} and then \eqref{eq:estim-calB1}, we get that the corresponding term in $\mathcal{B}_{t}^{1}(x)$ is 
\[
\overline{\mathcal{B}}_{t}^{11}\equiv\int_{0}^{t}\int_{\mathbb{R}}\frac{\delta^{2(1-b)}}{s^{2-b}}P_{u}(w)e^{-Cw^{2}}\mathcal{B}_{s}^{11}\,dwds.
\]
To extract a factor of $\delta^{r}$ from this term, we use the fact that 
$s\geqslant\delta^{2}$ to write 
\begin{align*}
 &\overline{\mathcal{B}}_{t}^{11} 
=\int_{\mathbb{R}}P_{u}(w)e^{-Cw^{2}}\rho^{\chi}dw\int_{0}^{t}\frac{\delta^{2(1-b)-2\gamma_{1}+2\gamma_{1}}}{s^{1-b-\gamma_{1}}}e^{-(\lambda+\theta\rho)s}\cdot s^{\alpha/2-1-\gamma_{1}}ds\\
 & \leqslant\delta^{2\gamma_{1}}\cdot\int_{\mathbb{R}}P_{u}(w)e^{-Cw^{2}}\rho^{\chi}dw\int_{0}^{t}e^{-(\lambda+\theta\rho)s}\cdot s^{\alpha/2-1-\gamma_{1}}ds,
\end{align*}
where $\gamma_{1}>0$ is some constant to be specified.
According to H\"older's inequality with $p\triangleq1/\chi$ ($1/q=1-\chi$),
we obtain that 
\begin{align*}
\int_{0}^{t}e^{-(\lambda+\theta\rho)s}\cdot s^{\alpha/2-1-\gamma_{1}}ds & \leqslant\left(\int_{0}^{t}e^{-p(\lambda+\theta\rho)s}ds\right)^{1/p}\cdot\left(\int_{0}^{t}s^{q(\alpha/2-1-\gamma_{1})}ds\right)^{1/q}\\
 & \lesssim(\lambda+\theta\rho)^{-\chi}\cdot T^{\alpha/2-\chi-\gamma_{1}}.
\end{align*}
Note that $\gamma_{1}$ needs to be less than $\alpha/2-\chi$ so
that the second time integral in the above inequality is finite. As
a consequence, the $\overline{\mathcal{B}}^{11}_{t}$-term produces an upper estimate
of the form 
\[
\overline{\mathcal{B}}^{11}_{t}\leqslant C_{\alpha,\beta,T,\gamma_{1}}\kappa(\hat{{\bf U}}_1^{\delta})\Theta(\tilde{v}^{\delta})\lambda^{\frac{\alpha-\beta}{4}}E(a,t)\cdot\delta^{2\gamma_{1}}.
\]
The other terms $\mathcal{B}^{12},\ldots,\mathcal{B}^{15}$ are discussed
in a similar way and the resulting factors of $\delta^{r'}$ are all
of higher order (i.e. $r'>2\gamma_{1}$). Therefore, we arrive at
\[
\mathcal{B}_{t}^{1}(x)\leqslant C_{\alpha,\beta,T,r}\kappa(\hat{{\bf U}}_1^{\delta})\Theta(\tilde{v}^{\delta})\lambda^{\frac{\alpha-\beta}{4}}E(a,t)\cdot\delta^{r}
\]
where $r$ is any given constant that is less than $\alpha-2\chi$.

In order to complete the large time estimates, according to our decomposition \eqref{eq:upper-decom-calB1calB2}, it remains to upper bound the term $\mathcal{B}_{t}^{2}(x)$ given by \eqref{eq:UnifB2}.
This will be an easy consequence of Lemma \ref{lem:Unif2ndDiffvsHK}. 
Namely a direct application of this lemma yields 
\[
\mathcal{B}_{t}^{2}(x)\leqslant\int_{0}^{t}\int_{\mathbb{R}}\frac{\delta}{\sqrt{s}}\cdot\frac{1}{s}e^{-Cw^{2}}\big|\mathcal{I}_{t-s}^{1,\delta}(x,x+\sqrt{s}w)\big|dwds.
\]
Hence by the same kind of analysis as in the $\mathcal{B}^{1}$-case, we
arrive at 
\[
\mathcal{B}_{t}^{2}(x)\leqslant C_{\alpha,\beta,T,r}\kappa(\hat{{\bf U}}_1^{\delta})\Theta(\tilde{v}^{\delta})\lambda^{\frac{\alpha-\beta}{4}}E(a,t)\cdot\delta^{r}
\]
where $r$ is any constant that is less than $\alpha-2\chi$. This 
is exactly the same estimate as the one for $\mathcal{B}^{1}$.

By putting (\ref{eq:UnifK2}) and the estimates of $\mathcal{A}^{1},\mathcal{A}^{2},\mathcal{B}^{1},\mathcal{B}^{2}$
together, we have thus completed the proof of Lemma \ref{lem:UnifWDelvsW}.
\end{proof}
\subsubsection{The time variation estimate.}\label{sec:TimeVarEstim}

In order to state the main result of this section, let us introduce some more notation. 
\begin{notation}\label{not:Theta-prime}
Let $\theta>2$, $a\geqslant 1$ and $\beta$ be coefficients satisfying the assumptions of Definition~\ref{def:weights} and Definition \ref{def:controlled-process}. We consider some additional parameters  $\theta'<\theta$ and $\beta'>\beta$ still satisfying the same assumptions. The norm $\Theta$ in Definition \ref{def:controlled-process} with parameters $\theta', \beta'$ will be denoted $\Theta^{\theta',\beta'}$ 
while the usual norm is written as $\Theta$. In the sequel for a function $h$ defined on $[0,T]$, the time increments of $h$ will be denoted by
\[
h_{t_{1},t_{2}}=h_{t_{2}}-h_{t_{1}},\quad\text{for}\quad 0\leqslant t_{1}<t_{2}\leqslant T.
\]
\end{notation}

As a slight elaboration of Proposition \ref{prop:discrete-FPT}, we state a lemma about the controlled norms of $\cw$ and $\cw^{1,\delta}$. Its proof is similar to that of Proposition \ref{prop:discrete-FPT} and omitted for the sake of conciseness. 

\begin{lem}\label{lem:lambda-betabetaprime}
For $\delta>0$, let $\cw, \cw^{1,\delta}$ be the processes defined by 
\eqref{eq:disc-wdelta}-\eqref{eq:cont-w}. Let $(\theta',\beta')$ be given fixed parameters satisfying the constraints specified in  Notation \ref{not:Theta-prime}. Let $\bar{\kappa}(\omega)$ be defined by (\ref{eq:BarKap}). Then there exist positive constants $C_1,C_2$ depending only on the underlying exponents, such that  
\begin{equation}\label{i5}
\Lambda^{\theta',\beta'}\equiv \sup_{\delta\geqslant 0}
\Theta^{\theta',\beta'}(\cw^{1,\delta})\vee \Theta^{\theta',\beta'}(\cw)\leqslant C_1 \bar{\kappa}(\omega) e^{C_1\lambda}\left(\|f_{0}\|_{\mathcal{C}_{L}^{3}}+\|g\|_{\mathcal{C}_{L}^{3}}\right)
\end{equation}provided that $\lambda = \lambda_\omega$ is chosen to satisfy 
\begin{equation*}
C_2\lambda^{-\frac{\alpha-\beta'}{4}}\bar{\kappa}(\omega)= \frac{1}{4}.
\end{equation*}
\end{lem}

\noindent
We now turn to the announced estimate for the time variations of $\cw^{1,\delta}-\cw$.
\begin{lem}\label{lem:TimVarWDelvsW}
Let the notation of Lemma \ref{lem:lambda-betabetaprime}, as well as Notation \ref{not:Theta-prime}, prevail. In particular the processes 
$\cw,\cw^{1,\delta}$ are introduced in \eqref{eq:disc-wdelta}-\eqref{eq:cont-w} and the exponents $\alpha,\beta$ satisfy relation 
\eqref{f1}. Then for any $a\geqslant1,$ $x\in[-a,a]$ and
$0\le t_{1}<t_{2}\le T$, we have 
\begin{align}
\big|\mathcal{W}_{t_{1},t_{2}}^{1,\delta}(x)-\mathcal{W}_{t_{1},t_{2}}(x)\big| & \leqslant CE(a,t_{2})a^{\beta/2}|t_{2}-t_{1}|^{\beta/2}\left(\lambda^{-\frac{\alpha-\beta}{4}}\left(\kappa(\hat{{\bf U}}_1^{\delta})d_{\hat{{\bf U}}_1^{\delta},{\bf W}}(\tilde{v}^{\delta},v)\right.\right.\nonumber \\
 & \ \ \ \left.\left.+\Theta(v)\rho(\hat{{\bf U}}_1^{\delta},{\bf W})\right)+\lambda^{\frac{\alpha-\beta}{4}}\kappa(\hat{{\bf U}}_1^{\delta})\Theta(\tilde{v}^{\delta})\delta^{\alpha-\beta}+\Lambda^{\theta',\beta'}\delta^{\beta'-\beta}\right),\label{eq:TimVarWDelvsW}
\end{align}
where  we recall that $\mathcal{W}_{t_{1},t_{2}}^{1,\delta}(x)\triangleq \cw_{t_{2}}^{\delta}(x)-\cw_{t_{1}}^{\delta}(x)$
and similarly for $\mathcal{W}_{t_{1},t_{2}}(x)$. We also recall that in \eqref{eq:TimVarWDelvsW} the quantity $E$ is given by \eqref{eq:EQDef},
 $\rho$ is defined by \eqref{eq:RhoUW} and $d_{\hat{{\bf U}}_1^{\delta},{\bf W}}$ is introduced in~\eqref{eq:ConDistVTilV}.
 \end{lem}

\begin{comment}
\begin{rem}
When $\beta,\alpha$ are close to $1/3$ and $\beta'$ is close to
$1/2$, the rate of convergence indicated by Lemma \ref{lem:TimVarWDelvsW}
is $\delta^{r}$ where $r$ can be made arbitrarily close to $1/6$.
\end{rem}
\end{comment}

\begin{proof}[Proof of Lemma \ref{lem:TimVarWDelvsW}]
We divide this proof again into several steps. Some technical considerations are similar to previous results and will only be sketched.

\noindent
\textit{Step 1: Decomposition of the time increments.}
We begin by decomposing the time increment into  
\begin{equation}\label{eq:incrementwwdelta}
\mathcal{W}_{t_{1},t_{2}}^{1,\delta}(x)-\mathcal{W}_{t_{1},t_{2}}(x)=\mathcal{T}^{1}_{t_{1},t_{2}}(x)+\mathcal{T}^{2}_{t_{1},t_{2}}(x),
\end{equation}
where the increments $\mathcal{T}^{1},\mathcal{T}^{2}$ are defined by
\begin{equation}\label{eq:calT12T12}
\mathcal{T}^{1}_{t_{1},t_{2}}(x)\triangleq\int_{0}^{t_{1}}\int_{\mathbb{R}}\left(Q_{s}^{\delta}(t_{2},x,y)-Q_{s}^{\delta}(t_{1},x,y)\right)dyds,\quad
\mathcal{T}^{2}_{t_{1},t_{2}}(x) \triangleq\int_{t_{1}}^{t_{2}}\int_{\mathbb{R}}Q_{s}^{\delta}(t_{2},x,y)dyds,
\end{equation}
and where we have used the notation 
\[
Q_{s}^{\delta}(t,x,y)\triangleq\nabla_{x}^{2,\delta}\hat{p}_{\lfloor t-s\rfloor}^{\delta}(\lfloor x-y\rfloor)\cdot\mathcal{I}_{s}^{\delta}(x,y)-\partial_{xx}^{2}p_{t-s}(x-y)\cdot\mathcal{I}_{s}(x,y),
\]
with $\mathcal{I}_{s}^{\delta}$, $\mathcal{I}_{s}$ given by \eqref{eq:stochintegvv}. Also recall that in the sequel we keep on using our Notation~\ref{h2} for the $\delta$-integer parts $\lfloor x\rfloor$ and $\lfloor t\rfloor$.
The remainder of the proof is dedicated to upper bound terms $\mathcal{T}^{1},\mathcal{T}^{2}$ in \eqref{eq:calT12T12}.

\noindent
\textit{Step 2: Bounding $\mathcal{T}^{2}$.}
In order to bound  $\mathcal{T}^{2}$ in  \eqref{eq:calT12T12}, we further 
decompose this term as 
\begin{equation}\label{eq:calT2112}
\mathcal{T}^{2}_{t_{1},t_{2}}(x) =\mathcal{T}^{21}_{t_{1},t_{2}}(x) +\mathcal{T}^{22}_{t_{1},t_{2}}(x) ,
\end{equation}
where $\mathcal{T}^{21}, \mathcal{T}^{22}$ are respectively defined by
\begin{equation}\label{eq:calT1t1t2}
\mathcal{T}^{21}_{t_{1},t_{2}}(x)\triangleq\int_{t_{1}}^{t_{2}}\int_{\mathbb{R}}\left(\nabla_{x}^{2,\delta}\hat{p}_{\lfloor t_{2}-s\rfloor}^{\delta}(\lfloor x-y\rfloor)-\partial_{xx}^{2}p_{t_{2}-s}(x-y)\right)\cdot\mathcal{I}_{s}^{\delta}(x,y)dyds,
\end{equation}
\begin{equation}\label{eq:T2t1t2}
\mathcal{T}^{22}_{t_{1},t_{2}}\triangleq\int_{t_{1}}^{t_{2}}\int_{\mathbb{R}}\partial_{xx}^{2}p_{t_{2}-s}(x-y)\cdot\left(\mathcal{I}_{s}^{\delta}(x,y)-\mathcal{I}_{s}(x,y)\right)dyds.
\end{equation}
The estimate of $\mathcal{T}^{22}_{t_{1},t_{2}}(x)$ above follows the same lines as in
the proof of Lemma \ref{lem:TimEst} (see \eqref{eq:uppbdT2}). For sake of conciseness, we just state the result here:
\begin{equation}\label{j1}
|\mathcal{T}^{22}_{t_{1},t_{2}}(x)|\leqslant C\lambda^{-\frac{\alpha-\beta}{4}}E(a,t_{2})a^{\beta/2}|t_{2}-t_{1}|^{\beta/2}\cdot\left(\kappa(\hat{{\bf U}}_1^{\delta})d_{\hat{{\bf U}}_1^{\delta},{\bf W}}(\tilde{v}^{\delta},v)+\Theta(v)\rho(\hat{{\bf U}}_1^{\delta},{\bf W})\right).
\end{equation} 

To estimate $\mathcal{T}^{21}_{t_{1},t_{2}}(x)$, we set $\tau\triangleq t_{2}-t_{1}$
and resort to the decomposition 
\begin{equation}\label{eq:Ttau-tau+}
\mathcal{T}^{21}_{t_{1},t_{2}}(x)=\hat{\mathcal{T}}_{\tau}^{-}(x)+
\hat{\mathcal{T}}_{\tau}^{+}(x),
\end{equation}
where 
\begin{equation}\label{eq:Ttau-}
\hat{\mathcal{T}}_{t_{1},t_{2}}^{-}(x)=\int_{0}^{\tau\wedge\delta^{2}}\int_{\mathbb{R}}\left(\nabla^{2,\delta}\hat{p}_{\lfloor r\rfloor}^{\delta}(\lfloor y\rfloor)-\partial_{xx}^{2}p_{r}(y)\right)\cdot\mathcal{I}_{t_{2}-s}^{\delta}(x,x+y)dyds
\end{equation}
\begin{equation}\label{eq:Ttau+}
\hat{\mathcal{T}}_{t_{1},t_{2}}^{+}(x)=\int_{\tau\wedge\delta^{2}}^{\tau}\int_{\mathbb{R}}\left(\nabla^{2,\delta}\hat{p}_{\lfloor r\rfloor}^{\delta}(\lfloor y\rfloor)-\partial_{xx}^{2}p_{r}(y)\right)\cdot\mathcal{I}_{t_{2}-s}^{\delta}(x,x+y)dyds.
\end{equation}
Now the arguments in order to estimate \eqref{eq:Ttau-} and \eqref{eq:Ttau+} are essentially similar to what we did in Lemma 
\ref{lem:UnifWDelvsW} for the terms $\mathcal{K}^{11}$ and $
\mathcal{K}^{12}$ (see \eqref{eq:K11}-\eqref{eq:K12}). In order to abbreviate our computations, we will only detail the bound for the small 
time integral $\hat{\mathcal{T}}^{-}$. As in Lemma \ref{lem:UnifWDelvsW}, the large time integral $\hat{\mathcal{T}}_{t_{1},t_{2}}^{+}(x)$ is treated
by using the local CLT and the Taylor approximation of the Gaussian
kernel. 

In order to handle the small time increment $\hat{\mathcal{T}}_{\tau}^{-}(x)$, we proceed exactly as in \eqref{eq:K11=00003DA1+A2}. 
That is, recalling that $\tau=t_{2}-t_{1}$, we write
\begin{equation}\label{eq:Tt1t2A1tauA2tau}
\hat{\mathcal{T}}_{\tau}^{-}(x)\leqslant
\mathcal A^{1}_{\tau}(x)+\mathcal{A}^{2}_{\tau}(x),
\end{equation}
where $\ca^{1}$ and $\ca^{2}$ are respectively given by
\begin{eqnarray}
\mathcal{A}^{1}_{\tau}(x)
&=&
\int_{0}^{\tau\wedge\delta^{2}}\int_{\mathbb{R}}
\lln\nabla^{2,\delta}\hat{p}_{\lfloor r\rfloor}^{\delta}(\lfloor y\rfloor)\rrn \cdot \lln\mathcal{I}_{t_{2}-s}^{\delta}(x,x+y)\rrn \,dyds
\label{eq:A1tau}  \\
\mathcal{A}^{2}_{\tau}(x)
&=&
\int_{0}^{\tau\wedge\delta^{2}}\int_{\mathbb{R}}
\lln\partial_{xx}^{2}p_{r}(y)\rrn \cdot \lln \mathcal{I}_{t_{2}-s}^{\delta}(x,x+y)\rrn \, dyds. 
\label{eq:A2tau}
\end{eqnarray}
Here again, for the sake of conciseness we shall only upper bound the term $\mathcal{A}^{1}_{\tau}(x)$ above. Namely recall that 
$\nabla^{2,\delta}\hat{p}_{\lfloor 0\rfloor}^{\delta}$ is supported in 
$\{y\in\mathbb{R}; |y|\leqslant\delta\}$. Moreover we have established~\eqref{eq:app-estim-deriv--hat-heat-kernel} above (which is a consequence of \eqref{eq:DiscUnifGau=0}), that is 
\begin{equation}
\big|\nabla_{x}^{2,\delta}\hat{p}^{\delta}_{0}(\lfloor y\rfloor)\big|\lesssim\frac{1}{\delta^{3}},\qquad\ \ \ \text{for }|y|\leqslant\delta.\label{eq:TimVarUni2DP0}
\end{equation}
In addition, we know from (\ref{eq:UniSmallTimDisc}) that the rough
integral $\mathcal{I}_{t_{2}-s}^{\delta}(x,x+y)$ is bounded above
by five basic terms. The first term appearing on the right hand side
of (\ref{eq:UniSmallTimDisc}) together with~(\ref{eq:TimVarUni2DP0})
yield the following integral in $\mathcal{A}^{1}_{\tau}(x)$: 
\[
\mathcal{A}^{11}_{\tau}(x)\triangleq\frac{1}{\delta^{3}}\int_{0}^{\tau\wedge\delta^{2}}\int_{\{y:|y|\leqslant\delta\}}(a+|y|)^{\chi}|y|^{\alpha}e^{-(\lambda+\theta(a+|y|))s}dyds.
\]
Owing to condition \eqref{f1}, we have  $\chi<\beta/2$. Therefore, for all $a\geqslant 1$  we get 
\begin{multline*}
|\mathcal{A}^{11}_{\tau}(x)|  \lesssim\frac{a^{\beta/2}}{\delta^{3}}\int_{0}^{\tau\wedge\delta^{2}}\int_{\{y:|y|\leqslant\delta\}}|y|^{\alpha}dyds\\
 \lesssim\frac{a^{\beta/2}}{\delta^{3}}\cdot\delta^{\alpha+1}\cdot(\tau\wedge\delta^{2})^{\beta/2+1-\beta/2}
  \leqslant\frac{a^{\beta/2}}{\delta^{3}}\cdot\delta^{\alpha+1}\cdot\tau^{\beta/2}\cdot\delta^{2-\beta}
 =a^{\beta/2}\tau^{\beta/2}\delta^{\alpha-\beta}.
\end{multline*}
In the inequality above, notice that we have extracted a power $\tau^{\beta/2}=(t_{2}-t_{1})^{\beta/2}$, which is our expected time regularity for $v$. 
This explains the appearance of the factor $\delta^{\alpha-\beta}$
in Lemma~\ref{lem:TimVarWDelvsW}. The other four terms on the right
hand side of (\ref{eq:UniSmallTimDisc}) lead to four corresponding
integrals in $\mathcal{A}^{1}_{\tau}(x)$, all of which having order $\delta^{r}$
with some $r>\alpha-\beta$. As a result, we obtain that 
\begin{equation}\label{eq:estimA1tau}
|\mathcal{A}^{1}_{\tau}(x)|\lesssim\kappa(\hat{{\bf U}}_1^{\delta})\Theta(\tilde{v}^{\delta})\lambda^{\frac{\alpha-\beta}{4}}E(a,t_{2})a^{\beta/2}|t_{2}-t_{1}|^{\beta/2}\cdot\delta^{\alpha-\beta}.
\end{equation}
A similar argument leads to exactly the same upper bound for $\mathcal{A}^{2}_{\tau}(x)$. Hence reporting~\eqref{eq:estimA1tau}
into \eqref{eq:Tt1t2A1tauA2tau}, then back into \eqref{eq:Ttau-tau+} and \eqref{eq:calT2112}, we have proved that $\mathcal{T}^{2}_{t_{1},t_{2}}(x)$ satisfies an inequality of the form \eqref{eq:TimVarWDelvsW}. This completes the estimate for $\mathcal{T}_{2}$.

\noindent
\textit{Step 3: Bounding $\mathcal{T}^{1}$.} We now consider the term $\mathcal{T}^{1}$ defined in \eqref{eq:incrementwwdelta}. The analysis is in fact similar to what we performed in \eqref{eq:decompT1T2}, and we only point out the main ingredients. In the first place we write
\begin{equation}\label{eq:decompT1=T11+T12}
\mathcal{T}^{1}_{t_{1},t_{2}}(x)=\mathcal{T}^{11}_{t_{1},t_{2}}(x)+
\mathcal{T}^{12}_{t_{1},t_{2}}(x), 
\end{equation}
where we define 
\begin{align}
\mathcal{T}^{11}_{t_{1},t_{2}}(x) & \triangleq\int_{0}^{t_{1}}\int_{\mathbb{R}}\left(J_{t_{1},t_{2}}^{\delta}(s,x-y)-J_{t_{1},t_{2}}(s,x-y)\right)\mathcal{I}_{s}^{\delta}(x,y)dyds,\label{eq:T11}\\
\mathcal{T}^{12}_{t_{1},t_{2}}(x) & \triangleq\int_{0}^{t_{1}}\int_{\mathbb{R}}J_{t_{1},t_{2}}(s,x-y)\left(\mathcal{I}_{s}^{\delta}(x,y)-\mathcal{I}_{s}(x,y)\right)dyds \, ,\label{eq:T12}
\end{align}
and where to ease notation we have set:
\begin{align}
J_{t_{1},t_{2}}^{\delta}(s,u) & \triangleq\nabla_{x}^{2,\delta}\hat{p}_{\lfloor t_{2}-s\rfloor}^{\delta}(\lfloor u\rfloor)-\nabla_{x}^{2,\delta}\hat{p}_{\lfloor t_{1}-s\rfloor}^{\delta}(\lfloor u\rfloor),\label{eq:Jt1t2delta}\\
J_{t_{1},t_{2}}(s,u) & \triangleq\partial_{xx}^{2}p_{t_{2}-s}(u)-\partial_{xx}^{2}p_{t_{1}-s}(u) .\label{eq:Jt1Jt2}
\end{align}

We now handle the terms $\mathcal{T}^{11},\mathcal{T}^{12}$ above.
The estimate of $\mathcal{T}^{12}_{t_{1},t_{2}}(x)$ follows the
same lines as in the proof of Lemma \ref{lem:TimEst}. We end up with
\begin{equation}\label{eq:estimT12}
|\mathcal{T}^{12}_{t_{1},t_{2}}(x)|\leqslant C\lambda^{-\frac{\alpha-\beta}{4}}E(a,t_{2})a^{\beta/2}|t_{2}-t_{1}|^{\beta/2}\cdot\left(\kappa(\hat{{\bf U}}_1^{\delta})d_{\hat{{\bf U}}_1^{\delta},{\bf W}}(\tilde{v}^{\delta},v)+\Theta(v)\rho(\hat{{\bf U}}_1^{\delta},{\bf W})\right).
\end{equation}
For the integral $\mathcal{T}^{11}_{t_{1},t_{2}}(x)$, we invoke the discrete heat equation \eqref{eq:discrete-heat-equation} to write, for $s_{1}<s_{2}$ in $\delta^{2}\mathbb{N}$ and $u\in\delta\mathbb{Z}$, 
\begin{equation}\label{j11}
\hat{p}_{s_{2}}^{\delta}(u)-\hat{p}_{s_{1}}^{\delta}(u)=\sum_{t_{j}=s_{1}}^{s_{2}-\delta^{2}}\left(\hat{p}_{t_{j+1}}^{\delta}(u)-\hat{p}_{t_{j}}^{\delta}(u)\right)=\frac{\sigma^{2}}{2}\delta^{2}\cdot\sum_{t_{j}=s_{1}}^{s_{2}-\delta^{2}}\nabla_{x}^{2,\delta}\hat{p}_{t_{j}}^{\delta}(u).
\end{equation}
In addition, the sum in the right hand side above can be written as a continuous time integral. We get 
\begin{equation}
\hat{p}_{s_{2}}^{\delta}(u)-\hat{p}_{s_{1}}^{\delta}(u)=\frac{\sigma^{2}}{2}\int_{s_{1}}^{s_{2}}\nabla_{x}^{2,\delta}\hat{p}_{\lfloor r\rfloor}^{\delta}(u)dr.\label{eq:diffphatHeateq}
\end{equation}
Now recall that the fourth discrete derivative $\nabla_{x}^{4,\delta}$ is defined by \eqref{d12}. Moreover, it is readily checked that $\nabla_{x}^{2,\delta}(\nabla_{x}^{2,\delta}f)=\nabla_{x}^{4,\delta}f$ for $f$ defined on the grid $\delta\mathbb{Z}$. Therefore applying $\nabla_{x}^{2,\delta}$ on both sides of \eqref{eq:diffphatHeateq} it follows that 
\[
J_{t_{1},t_{2}}^{\delta}(s,u)=\frac{\sigma^{2}}{2}\delta^{2}\cdot\sum_{t_{j}=\lfloor t_{1}-s\rfloor}^{\lfloor t_{2}-s\rfloor-\delta^{2}}\nabla^{4,\delta}\hat{p}_{t_{j}}^{\delta}(\lfloor u\rfloor) .
\]
Applying the same kind of manipulations to the continuous difference $J_{t_{1},t_{2}}(s,u)$ in \eqref{eq:Jt1Jt2} and then reporting those expressions in \eqref{eq:T11} we end up with 
\begin{align}
\mathcal{T}^{11}_{t_{1},t_{2}}(x) & =\frac{\sigma^{2}}{2}\int_{0}^{t_{1}}\int_{\mathbb{R}}\left(\int_{\lfloor t_{1}-s\rfloor}^{\lfloor t_{2}-s\rfloor}\nabla_{x}^{4,\delta}\hat{p}_{\lfloor r\rfloor}^{\delta}(\lfloor x-y\rfloor)dr-\int_{t_{1}-s}^{t_{2}-s}\partial_{x}^{4}p_{r}(x-y)dr\right)\mathcal{I}_{s}^{\delta}(x,y)dyds.\label{eq:T11intstoch}
\end{align}
Similarly to what we did for Lemma \ref{lem:UnifWDelvsW} we need to deal with the cases $|t_{2}-t_{1}|\leqslant\delta^{2}$
and $|t_{2}-t_{1}|>\delta^{2}$ separately in~\eqref{eq:T11intstoch}. Indeed, analyzing $\mathcal{T}^{11}_{t_{1},t_{2}}(x)$
in the former case will not produce a useful time-variation estimate
and we will adopt a more generic argument instead. We first discuss
the case when $|t_{2}-t_{1}|>\delta^{2}$.

\vspace{2mm}
\noindent\emph{Case I: $|t_{2}-t_{1}|>\delta^{2}$.}  
For $|t_{2}-t_{1}|>\delta^{2}$ the following holds true for all $s\in[0,t_{1}]$:
\[
\lfloor t_{1}-s\rfloor<t_{1}-s<\lfloor t_2-s\rfloor<t_{2}-s.
\]
According to this elementary fact we decompose $\mathcal{T}^{11}$ into 
\begin{equation}
\mathcal{T}^{11}_{t_{1},t_{2}}(x)
=\mathcal{T}^{111}_{t_{1},t_{2}}(x)-\mathcal{T}^{112}_{t_{1},t_{2}}(x)+\mathcal{T}^{113}_{t_{1},t_{2}}(x)\label{eq:decompT11=T111-T112+T113}
\end{equation}
where $\mathcal{T}^{111},\mathcal{T}^{112},\mathcal{T}^{113}$ are 
respectively defined by 
\begin{align}
\mathcal{T}^{111}_{t_{1},t_{2}}(x) & =\frac{\sigma^{2}}{2}\int_{0}^{t_{1}}\int_{\mathbb{R}}\int_{\lfloor t_{1}-s\rfloor}^{t_{1}-s}\nabla_{x}^{4,\delta}\hat{p}_{\lfloor r\rfloor}^{\delta}(\lfloor x-y\rfloor)
\, \mathcal{I}_{s}^{\delta}(x,y) \, drdyds,\label{eq:T111}\\
 \mathcal{T}^{112}_{t_{1},t_{2}}(x) & =\frac{\sigma^{2}}{2}\int_{0}^{t_{1}}\int_{\mathbb{R}}\int_{\lfloor t_{2}-s\rfloor}^{t_{2}-s}\nabla_{x}^{4,\delta}\hat{p}_{\lfloor r\rfloor}^{\delta}(\lfloor x-y\rfloor)
 \, \mathcal{I}_{s}^{\delta}(x,y)drdyds,\label{eq:T112}\\
 \mathcal{T}^{113}_{t_{1},t_{2}}(x) &=\frac{\sigma^{2}}{2}\int_{0}^{t_{1}}\int_{\mathbb{R}}\int_{t_{1}-s}^{t_{2}-s}\left(\nabla_{x}^{4,\delta}\hat{p}_{\lfloor r\rfloor}^{\delta}(\lfloor x-y\rfloor)-\partial_{x}^{4}p_{r}(x-y)\right)
 \, \mathcal{I}_{s}^{\delta}(x,y) \, drdyds
.\label{eq:T113}
\end{align}
The analysis of these terms is an adaptation of the
calculations developed in the proof of Lemma \ref{lem:TimEst}. For
the sake of simplicity, we only give a brief discussion on $\mathcal{T}^{112}_{t_{1},t_{2}}(x)$
and point out the main extra ingredients. The final estimate for $\mathcal{T}^{11}_{t_{1},t_{2}}(x)$
is stated in (\ref{eq:S3TimVarT11}) below.

\noindent
Let us handle the term  $\mathcal{T}^{112}_{t_{1},t_{2}}$ in \eqref{eq:T112}. We first resort to our usual change of variable 
$y=x+\sqrt{r}w$, which yields 
\[
\mathcal{T}^{112}_{t_{1},t_{2}}(x)=\frac{\sigma^{2}}{2}\int_{0}^{t_{1}}\int_{\mathbb{R}}\int_{\lfloor t_{2}-s\rfloor}^{t_{2}-s}\nabla_{x}^{4,\delta}\hat{p}_{\lfloor r\rfloor}^{\delta}(\lfloor\sqrt{r}w\rfloor)\cdot\sqrt{r}\cdot\mathcal{I}_{s}^{\delta}(x,x+\sqrt{r}w)drdwds.
\]
Moreover, since $t_{2}-t_{1}>\delta^{2}$ and $0\leqslant s\leqslant t_{1}$, we have 
\[
t_{2}-s\geqslant t_{2}-t_{1}>\delta^{2}\implies\lfloor t_{2}-s\rfloor\geqslant\delta^{2}.
\]
According to the uniform Gaussian estimate \eqref{eq:DiscUnifGau>0},
we get the following upper bound for all $r\in[\lfloor t_{2}-s\rfloor,t_{2}-s]$ and $w\in\mathbb{R}$:
\[
\big|\nabla_{x}^{4,\delta}\hat{p}_{\lfloor r\rfloor}^{\delta}(\lfloor\sqrt{r}w\rfloor)\big|\leqslant\frac{C_{1}}{r^{5/2}}e^{-C_{2}w^{2}}.
\]
Plugging this estimate into \eqref{eq:T112} and setting $u=\lfloor t_{2}-s\rfloor+r$,
it follows that 
\begin{align*}
\mathcal{T}^{112}_{t_{1},t_{2}}(x) & \lesssim\int_{\mathbb{R}}e^{-Cw^{2}}dw\int_{0}^{t_{1}}ds\int_{\lfloor t_{2}-s\rfloor}^{t_{2}-s}u^{-2}\big|\mathcal{I}_{s}^{\delta}(x,x+\sqrt{u}w)\big|du\\
 & =\int_{\mathbb{R}}e^{-Cw^{2}}dw\int_{0}^{t_{1}}ds\int_{0}^{(t_{2}-s)-\lfloor t_{2}-s\rfloor}(\lfloor t_{2}-s\rfloor+r)^{-2} \,\big|\mathcal{I}_{s}^{\delta}(x,x+\sqrt{\lfloor t_{2}-s\rfloor+r}w)\big|dr.
 \end{align*}
 We now invoke the relation $(t_{2}-s)-\lfloor t_{2}-s\rfloor\leq\delta^{2}$ again and apply Fubini' s theorem plus the change of variable $v=t_{2}-s$. This yields
 \begin{align}\label{j2}
 \mathcal{T}^{112}_{t_{1},t_{2}}(x) & \lesssim\int_{0}^{\delta^{2}}dr\int_{\mathbb{R}}e^{-Cw^{2}}dw\int_{0}^{t_{1}}\frac{ds}{(t_{2}-s+2r)^{2}}\;\big|\mathcal{I}_{s}^{\delta}(x,x+\sqrt{\lfloor t_{2}-s\rfloor+r}w)\big|\notag\\
 & \lesssim\int_{0}^{\delta^{2}}dr\int_{\mathbb{R}}e^{-Cw^{2}}dw\int_{t_{2}-t_{1}}^{t_{2}}\frac{dv}{v^{2}}\;\big|\mathcal{I}_{t_{2}-v}^{\delta}(x,x+\sqrt{\lfloor v\rfloor+r}w)\big|.
\end{align}
As we have seen several times (cf. \eqref{eq:app-estim-deriv--hat-heat-kernel}), the estimate of the rough integral
$\mathcal{I}_{t_{2}-v}^{\delta}(x,x+\sqrt{\lfloor v\rfloor+r}w)$
involves five basic terms. With the observation that 
\[
\lfloor v\rfloor+r\leqslant v+\delta^{2}\leqslant v+t_{2}-t_{1}\leqslant2v,
\]
the first term among the five leads us to the consideration of the
integral 
\[
\rho^{\chi}\int_{t_{2}-t_{1}}^{t_{2}}e^{-(\lambda+\theta\rho)v}v^{\alpha/2-2}dv\leqslant\rho^{\beta/2}\int_{t_{2}-t_{1}}^{\infty}v^{\alpha/2-2}dv\lesssim\rho^{\beta/2}|t_{2}-t_{1}|^{\alpha/2-1}
\]
where $\rho\triangleq a+\sqrt{T}|w|$. By including the outer two
integrals in \eqref{j2}, this term is further estimated as 
\begin{equation}\label{eq:exponents}
a^{\beta/2}|t_{2}-t_{1}|^{\alpha/2-1}\delta^{2}  =a^{\beta/2}\delta^{2}|t_{2}-t_{1}|^{\frac{\alpha-\beta}{2}-1}\cdot|t_{2}-t_{1}|^{\beta/2}
 \leqslant a^{\beta/2}\delta^{\alpha-\beta}\cdot|t_{2}-t_{1}|^{\beta/2},
\end{equation}
where we have used the relation $t_{2}-t_{1}\geq\delta^{2}$ for the last inequality. Notice that \eqref{eq:exponents} is compatible with our claim \eqref{eq:TimVarWDelvsW}. Next in the five terms decomposition alluded to above, one can check (as in the proof of Lemma \ref{lem:UnifWDelvsW} ) that the other four terms yield higher powers of the form $\delta^{r}$ with $r>\alpha-\beta$. Gathering those estimates we have obtained 
\begin{equation}\label{eq:estimT112}
\big|\mathcal{T}^{112}_{t_{1},t_{2}}(x)\big|\lesssim
E(a,t_{2})a^{\beta/2}|t_{2}-t_{1}|^{\beta/2}\lambda^{\frac{\alpha-\beta}{4}}\kappa(\hat{{\bf U}}_1^{\delta})\Theta(\tilde{v}^{\delta})\delta^{\alpha-\beta}.
\end{equation}
The term $\mathcal{T}^{111}_{t_{1},t_{2}}(x)$ is treated along the same lines, with a small difference. Namely one should separate the 
cases  $t_{1}\leqslant\delta^{2}$ and $t_{1}>\delta^{2}$. For the case 
$t_{1}>\delta^{2}$, one can follow exactly the same arguments as for 
\eqref{eq:estimT112}. Whenever $t_{1}\leqslant\delta^{2}$ one should rely on the fact that $\nabla_{x}^{4,\delta}\hat{p}_{0}^{\delta}$ is compactly supported, similarly to \eqref{eq:A2tx} and \eqref{eq:calA11estimate}. We leave the details to the patient reader. 

In order to handle the term $\mathcal{T}^{113}_{t_{1},t_{2}}(x)$ in 
\eqref{eq:T113}, we gather Theorem \ref{thm:LCLT} with uniform Gaussian estimates. This enables to write, for all $r\geq\delta^{2}$ and $w\in\mathbb{R}$, 
\[
\big|\nabla_{x}^{4,\delta}\hat{p}_{\lfloor r\rfloor}^{\delta}(\lfloor\sqrt{r}w\rfloor)-r^{-5/2}\partial_{x}^{4}p_{1}(w)\big|\lesssim\frac{\delta}{\sqrt{r}}\cdot\frac{e^{-Cw^{2}}}{r^{5/2}},\ \ \ \forall r\geqslant\delta^{2},w\in\mathbb{R}.
\]
Then arguing similarly to $\mathcal{T}^{111}_{t_{1},t_{2}}(x)$ and $\mathcal{T}^{112}_{t_{1},t_{2}}(x)$ we get
\begin{equation}\label{eq:estimT113}
\big|\mathcal{T}^{113}_{t_{1},t_{2}}(x)\big|\lesssim
E(a,t_{2})a^{\beta/2}|t_{2}-t_{1}|^{\beta/2}\lambda^{\frac{\alpha-\beta}{4}}\kappa(\hat{{\bf U}}_1^{\delta})\Theta(\tilde{v}^{\delta})\delta^{\alpha-\beta}.
\end{equation}
Eventually, putting together \eqref{eq:estimT112} and \eqref{eq:estimT113} and recalling the decomposition \eqref{eq:decompT11=T111-T112+T113}, we arrive at the following bound for the case $t_{2}-t_{1}\ge \delta^{2}$: 
\begin{equation}
|\mathcal{T}^{11}_{t_{1},t_{2}}(x)|\lesssim E(a,t_{2})a^{\beta/2}|t_{2}-t_{1}|^{\beta/2}\lambda^{\frac{\alpha-\beta}{4}}
\kappa(\hat{{\bf U}}_1^{\delta})\Theta(\tilde{v}^{\delta})\delta^{\alpha-\beta}.\label{eq:S3TimVarT11}
\end{equation}

\noindent \textit{Case II}: $|t_{2}-t_{1}|\leqslant\delta^{2}.$ \quad The argument for this case is generic. Recall that $\theta'<\theta$
and $\beta'>\beta$ are parameters ($\beta'$ satisfies the same constraints
as $\beta$ does) such that 
\begin{equation}
e^{(\theta'-\theta)a+(\theta'-\theta)aT}a^{\frac{\beta'-\beta}{2}}\leqslant1 \, ,
\quad \text{for all} \quad  a\geqslant1.\label{eq:ParConTV}
\end{equation}
The notation $\Theta^{\theta',\beta'}$ for the controlled path norm
indicates that the exponents $(\theta,\beta)$ are replaced by $(\theta',\beta')$
whenever applicable. We also denote $E^{\theta'}(a,t_{2})$ as the
exponential weight function $E(\cdot,\cdot)$ defined in Definition
\ref{def:weights} with $\theta$ being replaced by $\theta'$. By
treating $w^{\delta},w$ as controlled paths with respect to the $\Theta^{\theta',\beta'}$-norm,
we directly apply the triangle inequality to get 
\begin{align}
\big|\mathcal{W}_{t_{1},t_{2}}^{1,\delta}(x)-\mathcal{W}_{t_{1},t_{2}}(x)\big| & \leqslant\big|\mathcal{W}_{t_{1},t_{2}}^{1,\delta}(x)\big|+|\mathcal{W}_{t_{1},t_{2}}(x)|\nonumber \\
 & \leqslant E^{\theta'}(a,t_{2})a^{\beta'/2}\left(\Theta^{\theta',\beta'}(w^{\delta})+\Theta^{\theta',\beta'}(w)\right)|t_{2}-t_{1}|^{\beta'/2}\nonumber \\
 & \leqslant E^{\theta}(a,t_{2})a^{\beta/2}\left(\Theta^{\theta',\beta'}(w^{\delta})+\Theta^{\theta',\beta'}(w)\right)\delta^{\beta'-\beta}|t_{2}-t_{1}|^{\beta/2},\label{eq:GenericTV}
\end{align}
where the last inequality follows from the assumption that $|t_{2}-t_{1}|\leqslant\delta^{2}$.
This yields the last term appearing on the right hand side of (\ref{eq:TimVarWDelvsW}).

\noindent
\textit{Step 4: Conclusion.}
The proof of Lemma \ref{lem:TimVarWDelvsW} is now complete by reporting \eqref{eq:estimT12}, \eqref{eq:S3TimVarT11} and~\eqref{eq:GenericTV} into the decompositions \eqref{eq:decompT1=T11+T12} and \eqref{eq:incrementwwdelta}.
\end{proof}

\subsubsection{The space variation estimate}\label{sec:SpaceVarEstim}

\noindent The analysis for this part follows the same steps as in
Section \ref{sec:SpVarEstRP}, with the discussion of multiple cases
(small time versus large time, that is $t\leqslant|x'-x|^{2}$ versus $t>|x'-x|^{2}$). The separation small versus large time is also largely similar to the time variation case in Section \ref{sec:TimeVarEstim}. We only state
the final result below and leave the details to the patient reader. Notice that below we use the notation~\eqref{eq:notation-spacediff} for spatial increments.
\begin{lem}
\noindent \label{lem:SpVarWDelvsW}
As in Lemma \ref{lem:TimVarWDelvsW} we consider some parameters 
$\theta',\beta', \theta,\theta'$ according to Notation~\ref{not:Theta-prime} as well as the processes $\mathcal{W}$, $\mathcal{W}^{1,\delta}$ in \eqref{eq:disc-wdelta}-\eqref{eq:cont-w}. Following our notation~\eqref{eq:notation-spacediff}, set 
\[
\mathcal{W}_{t}(x,x')\triangleq \mathcal{W}_{t}(x')-\mathcal{W}_{t}(x),\quad\text{and}\quad
\mathcal{W}_{t}^{1,\delta}(x,x')\triangleq \mathcal{W}_{t}^{1,\delta}(x')-\mathcal{W}_{t}^{1,\delta}(x).
\]
Also recall that $\alpha,\beta,\chi$ fulfill condition \eqref{f1} and that the quantities $\Lambda$ are introduced in~\eqref{i5}. Then, for any $a\geqslant1,$ $x,x'\in[-a,a]$
and $t\in[0,T]$, we have 
\begin{align}
\big|\mathcal{W}_{t}^{1,\delta}(x,x')-\mathcal{W}_{t}(x,x')\big| 
& \leqslant 
CE(a,t)a^{\beta/2}|x'-x|^{\beta}\left(\lambda^{-\frac{\alpha-\beta}{4}}\left(\kappa(\hat{{\bf U}}_1^{\delta})d_{\hat{{\bf U}}_1^{\delta},{\bf W}}(\tilde{v}^{\delta},v)\right.\right.\nonumber\\
 & \ \ \ \left.\left.+\Theta(v)\rho(\hat{{\bf U}}_1^{\delta},{\bf W})\right)+\lambda^{\frac{\alpha-\beta}{4}}\kappa(\hat{{\bf U}}_1^{\delta})\Theta(\tilde{v}^{\delta})\delta^{\alpha-\beta}+\Lambda^{\theta',\beta'}\delta^{\beta'-\beta}\right)\label{eq:SpVarEst}.
\end{align}
\end{lem}

\subsubsection{The remainder estimate.}\label{sec:RemainEstim}

In Proposition \ref{prop:discrete-FPT} we have bounded the rough path norm of $v^{\delta}$ without  referring explicitly to its rough path decomposition. We will now give a more specific formula in this direction. 
Namely, going back to relation \eqref{eq:disc-wdelta} and since $\mathcal{W}=\mathcal{M}\mathcal{V}$, equation~\eqref{eq:deriv-Gubi-M} asserts that the derivative $\partial_{{\bf W}}\mathcal{W}_{t}(x)$ is 
\begin{equation}\label{eq:derivativeWW}
\partial_{{\bf W}}\mathcal{W}_{t}(x)=\partial_{{\bf W}}v_{t}(x)=-\frac{2}{\sigma^{2}}v_{t}(x).
\end{equation}
As in relation \eqref{eq:ReMV} we thus introduce a remainder term for 
$\mathcal{W}$, seen as a process controlled by $\bf W$:
\begin{equation}\label{eq:ReMV-again}
\mathcal{R}_{{\bf W}}^{\mathcal{W}_{t}}(x,x')=\mathcal{W}_{t}(x,x')+\frac{2}{\sigma^{2}}v_{t}(x)W^{1}(x,x').
\end{equation}
Analogously, we shall define similar quantities for $\mathcal{W}^{1,\delta}$ seen as a process controlled by $\hat{{\bf U}}_1^{\delta}$:
\begin{eqnarray}
\label{eq:derivativeUW}
\partial_{\hat{U}_1^{\delta}}\mathcal{W}_{t}^{1,\delta}(x)
&=&
-\frac{2}{\sigma^{2}}\tilde{v}_{t}^{\delta}(x),
\\
\label{eq:ReUhat}
\mathcal{R}_{\hat{{\bf U}}_1^{\delta}}^{\mathcal{W}_{t}^{1,\delta}}(x,x')
&=&
\mathcal{W}_{t}^{1,\delta}(x,x')+\frac{2}{\sigma^{2}}\tilde{v}_{t}^{\delta}(x)\hat{U}_1^{\delta}(x,x').
\end{eqnarray}
Notice that in \eqref{eq:ReMV-again} and \eqref{eq:ReUhat} we have used the notation~\eqref{eq:notation-spacediff} for spatial increments.
Notice also that \eqref{eq:ReUhat} gives a controlled path for the stochastic integral term $\mathcal{W}^{1,\delta}$ only. Taking into account the fact that one expects the deterministic terms in 
\eqref{eq:disc-PDE-vdelta}-\eqref{eq:disc-eta-delta}
to have a null Gubinelli derivative, relation \eqref{eq:ReUhat} can also 
be translated into a decomposition for the process $v^{\delta}$. We label 
this decomposition for further use:
\begin{equation}\label{eq:rewrite-derivativeUW}
\tilde{v}_{t}^{\delta}(x,x')=-\frac{2}{\sigma^{2}}\tilde{v}_{t}^{\delta}(x)\hat{U}_1^{\delta}(x,x') 
+\mathcal{R}_{\hat{{\bf U}}_1^{\delta}}^{\tilde{v}_{t}^{\delta}}(x,x').
\end{equation}
We will go back to this decomposition in Section \ref{sec:comparing-J-eta}. For now, 
the main task in this section is to compare the remainder terms $\mathcal{R}^{\mathcal{W}_{t}}$ and 
$\mathcal{R}^{\mathcal{W}_{t}^{1,\delta}}$ above. Our result is summarized in the lemma below. 

\begin{lem}\label{lem:RemEstComp} 
We use the same notation as in Lemma \ref{lem:SpVarWDelvsW}. Let $\mathcal{R}^{\mathcal{W}_{t}}$ and 
$\mathcal{R}^{\mathcal{W}_{t}^{1,\delta}}$ be the remainders respectively defined by \eqref{eq:ReMV-again} and \eqref{eq:ReUhat}. Then for every $a\geqslant1,$ $x,x'\in[-a,a]$
and $t\in[0,T]$, we have 
\begin{align}
 & \big|\mathcal{R}_{\hat{{\bf U}}_1^{\delta}}^{\mathcal{W}_{t}^{1,\delta}}(x,x')-\mathcal{R}_{{\bf W}}^{\mathcal{W}_{t}}(x,x')\big|\nonumber \\
 & \leqslant CE(a,t)a^{\chi}(a^{\beta/2}+t^{-\beta/2})|x'-x|^{2\beta}\left[\left(\kappa(\hat{{\bf U}}_1^{\delta})d_{\hat{{\bf U}}_1^{\delta},{\bf W}}(\tilde{v}^{\delta},v)+\Theta(v)\rho(\hat{{\bf U}}_1^{\delta},{\bf W})\right)\right.\nonumber \\
 &\left. \ \ \ +\lambda^{\frac{\alpha-\beta}{4}}\left(\kappa({\bf W})\Theta(v)\delta^{\alpha-\beta}+\left(\Lambda^{\theta',\beta'}+\kappa(\hat{{\bf U}}_1^{\delta})\Theta(\tilde{v}^{\delta})\right)\delta^{\frac{\alpha}{\beta'+\beta}(\beta'-\beta)}\right)\right].\label{eq:RemEstComp}
\end{align}
\end{lem}

\begin{comment}
\begin{rem}\label{rem:delta1/10}
\noindent By taking $\beta'\uparrow\alpha,$ $\alpha\uparrow\frac{1}{2}$
and $\beta\downarrow\frac{1}{3}$, the above convergence rate is arbitrarily
close to $\delta^{1/10}$. However, by letting $\beta\to 0$ and still $\beta'\to\alpha$ and $\alpha\to\frac{1}{2}$, one can achieve a convergence rate of order $\delta^{1/2}$. This would thus require to consider $\mathcal{W}_{t}$, $\mathcal{W}_{t}^{\delta}$ as controlled processes with low regularity $\beta$, leading to tedious rough paths expansions. 
\end{rem}
\end{comment}

\noindent 
Let us prepare for the proof of Lemma \ref{lem:RemEstComp}.
As in Sections \ref{sec:TimeVarEstim} and \ref{sec:SpaceVarEstim}, we only provide the details for those estimates
that do not follow from a simple adaptation of the analysis in Section
\ref{sec:EstMV}. We begin by recalling that according to \eqref{eq:W=MV} we have set $\mathcal{W}=\mathcal{M}\mathcal{V}$ and the remainder for $\mathcal{W}$ is spelled out in \eqref{eq:ReMV-again}. Otherwise stated the remainder of 
$\mathcal{W}$ as a controlled process is given by (\ref{eq:ReMVbis}), that is
\begin{equation}
\mathcal{R}_{{\bf W}}^{\mathcal{W}_{t}}(x,x')=\mathcal{R}_{t}^{0}(x,x')+\frac{2}{\sigma^{2}}v_{t}(x) P_{t}W(x,x'),
\label{eq:ConRemDecomp}
\end{equation}
where $\mathcal{R}_{t}^{0}(x,x')$ is defined by (\ref{eq:ReMV0})
and $P_{t}$ is the heat semigroup on $\mathbb{R}$ with generator
$\frac{\sigma^{2}}{2}\Delta$. To compare \eqref{eq:ConRemDecomp} with the discrete remainder,
we shall make use of a similar decomposition of $\mathcal{R}_{\hat{{\bf U}}_1^{\delta}}^{\mathcal{W}_{t}^{1,\delta}}$
given in the lemma below. 
\begin{lem}
\noindent Let $\mathcal{W}^{1,\delta}$ be the process introduced in 
\eqref{eq:WiDel}, and recall that its remainder $\mathcal{R}_{\hat{{\bf U}}_1^{\delta}}^{\mathcal{W}}$ is defined by \eqref{eq:ReUhat}. Then for all $t\in[0,T]$ and $x,x'\in\mathbb{Z}$, the following decomposition holds true: 
\begin{equation}
\mathcal{R}_{\hat{{\bf U}}_1^{\delta}}^{\mathcal{W}_{t}^{1,\delta}}(x,x')=\mathcal{R}_{t}^{0,\delta}(x,x')
+\frac{2}{\sigma^{2}}\tilde{v}_{t}^{\delta}(x) \left(\hat{P}_{t}^{\delta}\hat{U}_1^{\delta}(x,x')
+\mathcal{Q}_{t}^{\delta}(x,x')-\mathcal{K}_{t}^{\delta}(x,x')\right).\label{eq:DisRemDecomp}
\end{equation}
In \eqref{eq:DisRemDecomp} we have set 
\begin{multline}
\mathcal{R}_{t}^{0,\delta}(x,x')  \triangleq\int_{0}^{t}\int_{\mathbb{R}}\left(\nabla_{x}^{2,\delta}\hat{p}_{\lfloor s\rfloor}^{\delta}(\lfloor x'-y\rfloor-\nabla_{x}^{2,\delta}\hat{p}_{\lfloor s\rfloor}^{\delta}(\lfloor x-y\rfloor))\right) \\
\times\int_{x}^{y}\left(\tilde{v}_{t-s}^{\delta}(z)-\tilde{v}_{t}^{\delta}(x)\right)d\hat{{\bf U}}_1^{\delta}(z)dyds,\label{j3}
\end{multline}
and for $\xi\in\mathbb{R}$ we also define
\begin{align}
\hat{P}_{t}^{\delta}\hat{U}_1^{\delta}(\xi) & \triangleq\int_{\mathbb{R}}\hat{p}_{\lfloor t\rfloor}^{\delta}(\lfloor\xi-y\rfloor)\hat{U}_1^{\delta}(y)dy,\label{j4} \\
\mathcal{Q}_{t}^{\delta}(\xi) & \triangleq\frac{t-\lfloor t\rfloor}{\delta^{2}}\int_{\mathbb{R}}\left(\hat{p}_{\lfloor t\rfloor+\delta^{2}}^{\delta}(\lfloor\xi-y\rfloor)-\hat{p}_{\lfloor t\rfloor}^{\delta}(\lfloor\xi-y\rfloor)\right)\hat{U}_1^{\delta}(y)dy,\label{eq:RemEstJDel}\\
\mathcal{K}^{\delta}(\xi) & \triangleq\hat{P}_{0}^{\delta}\hat{U}_1^{\delta}(\xi)-\hat{U}_1^{\delta}(\xi).\label{eq:RemEstKDel}
\end{align}
Also recall our convention \eqref{eq:notation-spacediff} for increments of the form $ \mathcal{Q}_t^\delta(x,x')$, $ \ck_{t}^{\delta}(x,x')$.
\end{lem}

\begin{proof}
\noindent The proof is purely algebraic and reproduces some of the steps in Lemma \ref{lem:RemDec}. Starting from \eqref{eq:disc-wdelta}, some elementary manipulations yield 
\begin{equation}\label{eq:WRI}
\mathcal{W}_{t}^{1,\delta}(x,x')=\mathcal{R}_{t}^{0,\delta}(x,x')+\tilde{v}_{t}^{\delta}(x)\cdot I_{t}^{\delta}(x,x'),
\end{equation}
where $\mathcal{R}_{t}^{0,\delta}(x,x')$ is defined by (\ref{j3})
and 
\begin{equation}\label{eq:Idelta}
I_{t}^{\delta}(\xi)\triangleq\int_{0}^{t}\int_{\mathbb{R}}\nabla_{x}^{2,\delta}\hat{p}_{\lfloor s\rfloor}^{\delta}(\lfloor\xi-y\rfloor)\hat{U}_1^{\delta}(y)dyds.
\end{equation}
Similarly to what we have done in \eqref{j11}, one can switch from space to time gradients of $\hat{p}^{\delta}$ in \eqref{eq:Idelta} thanks to the discrete heat equation~\eqref{eq:discrete-heat-equation}. Namely write
\[
\frac{\hat{p}_{s+\delta^{2}}^{\delta}(u)-\hat{p}_{s}^{\delta}(u)}{\delta^{2}}=\frac{\sigma^{2}}{2}\nabla_{x}^{2,\delta}\hat{p}_{s}^{\delta}(u),\quad\text{ for all }\quad (s,u)\in\delta^{2}\mathbb{N}\times\delta\mathbb{Z}.
\]
Then we can recast \eqref{eq:Idelta} as
\[
I_{t}^{\delta}(\xi)=\frac{2}{\sigma^{2}}\int_{0}^{t}\int_{\mathbb{R}}\frac{\hat{p}_{\lfloor s\rfloor+\delta^{2}}^{\delta}(\lfloor\xi-y\rfloor)-\hat{p}_{\lfloor s\rfloor}^{\delta}(\lfloor\xi-y\rfloor)}{\delta^{2}}\hat{U}_1^{\delta}(y)dyds.
\]
By further expressing the above time integral as $\int_{0}^{t}=\int_{0}^{\lfloor t\rfloor}+\int_{\lfloor t\rfloor}^{t}$,
the first integral yields $\frac{2}{\sigma^{2}}(\hat{P}_{t}^{\delta}\hat{U}_1^{\delta}-\hat{P}_{0}^{\delta}\hat{U}_1^{\delta})(\xi)$ by expressing the time integral as a discrete sum,
while the other one yields $\frac{2}{\sigma^{2}}\mathcal{Q}_{t}^{\delta}(\xi)$. 
Summarizing, we have obtained the relation 
\begin{equation}\label{eq:I=diffPplusJ}
I_{t}^{\delta}(\xi)=\frac{2}{\sigma^{2}}
\left((\hat{P}_{t}^{\delta}\hat{U}_1^{\delta}-\hat{P}_{0}^{\delta}\hat{U}_1^{\delta})(\xi)+\mathcal{Q}_{t}^{\delta}(\xi)\right).
\end{equation}
Plugging \eqref{eq:I=diffPplusJ} into \eqref{eq:WRI} we thus obtain 
\begin{equation*}
\mathcal{W}_{t}^{1,\delta}(x,x')=\frac{2\tilde{v}_{t}^{\delta}(x)}{\sigma^{2}}\left((\hat{P}_{t}^{\delta}\hat{U}_1^{\delta}-\hat{P}_{0}^{\delta}\hat{U}_1^{\delta})(x,x')+\mathcal{Q}_{t}^{\delta}(x,x')\right)
+\mathcal{R}_{t}^{0,\delta}(x,x').
\end{equation*}
Now add and subtract the term $\frac{2\tilde{v}_{t}^{\delta}(x)}{\sigma^{2}}\hat{U}_1^{\delta}(x,x')$ to the above expression. This yields 
\begin{equation*}
\mathcal{W}_{t}^{1,\delta}(x,x')=\frac{2\tilde{v}_{t}^{\delta}(x)}{\sigma^{2}}\hat{U}_1^{\delta}(x,x')+
\frac{2\tilde{v}_{t}^{\delta}(x)}{\sigma^{2}}\left(\hat{P}_{t}^{\delta}\hat{U}_1^{\delta}(x,x')+\mathcal{Q}_{t}^{\delta}(x,x')-\mathcal{K}_{t}^{\delta}(x,x')\right)
+\mathcal{R}_{t}^{0,\delta}(x,x'),
\end{equation*}
where $\mathcal{K}_{t}^{\delta}$ is defined by \eqref{eq:RemEstKDel}. This proves our claim \eqref{eq:DisRemDecomp}.
\end{proof}
\begin{rem}
The appearance of the $\mathcal{K}^{\delta}$-term is a special feature
in the interpolated discrete equation which does not arise in the
continuous case. Indeed, it is due to the fact that $x\mapsto\hat{p}_{0}^{\delta}(\lfloor x\rfloor)$
is not the Dirac delta function on $\mathbb{R}.$ %The estimate of
%$\mathcal{K}^{\delta}(x,x')$ indeed contributes to the worst part
%of the convergence rate (it gives rise to the factor $\delta^{1/10}$
%while all other estimates produce $\delta^{1/6}$ in the worst scenario).
\end{rem}

We now turn our attention to the proof of our main estimate 
for the remainder. 

\begin{proof}[Proof of Lemma \ref{lem:RemEstComp}] We have seen that  $\mathcal{R}_{{\bf W}}^{\mathcal{W}_{t}}$ and $\mathcal{R}_{\hat{{\bf U}}_1^{\delta}}^{\mathcal{W}_{t}^{1,\delta}}$ are respectively given by \eqref{eq:ReMV-again} and \eqref{eq:ReUhat}, with decompositions \eqref{eq:ConRemDecomp} and \eqref{eq:DisRemDecomp}. Therefore our global strategy will be based on the following estimates:

\vspace{2mm}\noindent
(i) compare $\mathcal{R}_{t}^{0,\delta}(x,x')$ versus $\mathcal{R}_{t}^{0}(x,x')$;\\
(ii) compare $\tilde{v}_{t}^{\delta}(x)(\hat{P}_{t}^{\delta}\hat{U}_1^{\delta}(x')-\hat{P}_{t}^{\delta}\hat{U}_1^{\delta}(x))$ versus 
 $v_{t}(x)(P_{t}W(x')-P_{t}W(x))$;\\
(iii) show that  both $\mathcal{Q}_{t}^{\delta}(x,x'),$ $\mathcal{K}^{\delta}(x,x')$ are small terms.

\vspace{2mm}\noindent Below we divide this task into several steps.

%\smallskip

\noindent \textit{Step 1: Case $t<\delta^{2}$.} 
 Like in the proof of Lemma \ref{lem:TimVarWDelvsW}
(see \eqref{eq:GenericTV}), we will prove that all the terms in (i)-(iii) above are small individually. We summarize the basic analysis as follows.

%\smallskip

\noindent(A) {[}Estimating $\mathcal{R}^{0}${]} By adapting the proof of
Lemma \ref{lem:R0Est}, it is seen easily that 
\begin{equation}\label{eq:estimR0xxprime}
|\mathcal{R}_{t}^{0}(x,x')|\lesssim\kappa({\bf W})\Theta(v) E(a,t)a^{\chi+\beta/2}|x'-x|^{2\beta}\cdot\delta^{\alpha-\beta}.
\end{equation}

\noindent(B) {[}Estimating $P_{t}W${]} As a consequence of Lemma \ref{lem:SGEst}
(cf. (\ref{eq:estim-incr-PtW-bis}) and (\ref{eq:estim-incr-PtW-ter})), as well as the estimate on $\mathcal{V}$ in Theorem \ref{thm:FPT},
we have 
\begin{eqnarray}
\big|v_{t}(x)\cdot(P_{t}W(x')-P_{t}W(x))\big|
&\lesssim&\kappa({\bf W})\Theta(v)E(a,t)a^{\chi}t^{-\beta/2}|x'-x|^{2\beta}\cdot t^{\frac{\alpha-\beta}{2}}
\notag\\
& \leqslant&
\kappa({\bf W})\Theta(v)E(a,t)a^{\chi}t^{-\beta/2}|x'-x|^{2\beta}\cdot\delta^{\alpha-\beta}.\label{eq:estim-diff-Ptxxprime}
\end{eqnarray}

\noindent(C) {[}Estimating $\mathcal{R}^{0,\delta}${]} We start from the expression \eqref{j3} and we set $u=x-y$. This yields 
\begin{multline}
\mathcal{R}_{t}^{0,\delta}(x,x')  =\int_{0}^{t}\int_{\mathbb{R}}\left(\nabla_{x}^{2,\delta}\hat{p}_{0}^{\delta}(\lfloor x'-x+u\rfloor)-\nabla_{x}^{2,\delta}\hat{p}_{0}^{\delta}(\lfloor u\rfloor)\right) \\
 \times\int_{x}^{x-u}\left(\tilde{v}_{t-s}^{\delta}(z)-\tilde{v}_{t}^{\delta}(x)\right)d\hat{U}_1^{\delta}(z)duds.\label{eq:EstR0Del}
\end{multline}
The above expression can be estimated in the usual way as in the continuous
case (see relation~\eqref{eq:R0eL1L1pL2}) by writing out $\nabla_{x}^{2,\delta}\hat{p}_{0}^{\delta}$ explicitly
and splitting 
\[
\tilde{v}_{t-s}^{\delta}(z)-\tilde{v}_{t}^{\delta}(x)=\tilde{v}_{t-s}^{\delta}(z)-\tilde{v}_{t-s}^{\delta}(x)+\tilde{v}_{t-s}^{\delta}(x)-\tilde{v}_{t}^{\delta}(x).
\]

For simplicity, we only consider the integral corresponding to $\tilde{v}_{t-s}^{\delta}(z)-\tilde{v}_{t-s}^{\delta}(x)$ (the other one is in fact easier),
and call this integral $\mathcal{N}_{t}^{\delta}(x,x')$. Specifically, we consider 
\begin{multline}\label{j5}
\mathcal{N}_{t}^{\delta}(x,x')=\int_{0}^{t}\int_{\mathbb{R}}\left(\nabla_{x}^{2,\delta}\hat{p}_{0}^{\delta}(\lfloor x'-x+u\rfloor)-\nabla_{x}^{2,\delta}\hat{p}_{0}^{\delta}(\lfloor u\rfloor)\right) \\
\times\int_{x}^{x-u}\left(\tilde{v}_{t-s}^{\delta}(z)-\tilde{v}_{t-s}^{\delta}(x)\right)d\hat{U}_1^{\delta}(z)duds \, .
\end{multline}
With the exact expression of $\hat{p}_{0}(x)=\delta^{-1}\mathds{1}_{(x=0)}$, a simple explicit calculation shows that
\begin{equation}\label{eq:Nab2p0Del}
\nabla_{x}^{2,\delta}\hat{p}_{0}^{\delta}(u)=-\frac{1}{2\delta^{3}}{\bf 1}_{\{u=0\}}-\frac{1}{16\delta^{3}}{\bf 1}_{\{u=\pm\delta\}}+\frac{1}{8\delta^{3}}{\bf 1}_{\{u=\pm2\delta\}}+\frac{1}{16\delta^{3}}{\bf 1}_{\{u=\pm3\delta\}}.
\end{equation}
As a result, the $u$-integral in (\ref{eq:EstR0Del})-\eqref{j5}
is supported on finitely many intervals of order $\delta$. Those regions will be paired in order to take advantage of some cancellations due to our expression \eqref{eq:Nab2p0Del}.

As an example, set $\eta=x'-x$ and consider the term corresponding to the regions 
$\lfloor \eta+u\rfloor=0$ versus $\lfloor u\rfloor=0$. This yields a term
\begin{align}
\mathcal{N}^{1,\delta}(x,x')\triangleq & \int_{0}^{t}\left(\int_{\{u:\lfloor\eta+u\rfloor=0\}}\nabla_{x}^{2,\delta}\hat{p}_{0}^{\delta}(\lfloor\eta+u\rfloor)-\int_{\{u:\lfloor u\rfloor=0\}}\nabla_{x}^{2,\delta}\hat{p}_{0}^{\delta}(\lfloor\eta+u\rfloor)\right)\nonumber\\
 & \ \ \ \times\int_{x}^{x-u}\left(\tilde{v}_{t-s}^{\delta}(z)-\tilde{v}_{t-s}^{\delta}(x)\right)d\hat{U}_{1}^{\delta}duds\nonumber\\
= & -\frac{1}{2\delta^{3}}\int_{0}^{t}\left(\int_{-\eta}^{-\eta+\delta}-\int_{0}^{\delta}\right)\int_{x}^{x-u}\left(\tilde{v}_{t-s}^{\delta}(z)-\tilde{v}_{t-s}^{\delta}(x)\right)d\hat{U}_{1}^{\delta}duds.\label{eq:N1Delta}
\end{align}We first consider the case when $0<\eta<\delta$. Then (\ref{eq:N1Delta}) becomes
\begin{equation}
\mathcal{N}^{1,\delta}(x,x')=  -\frac{1}{2\delta^{3}}\int_{0}^{t}\left(\int_{-\eta}^{0}-\int_{-\eta+\delta}^{\delta}\right)\int_{x}^{x-u}\left(\tilde{v}_{t-s}^{\delta}(z)-\tilde{v}_{t-s}^{\delta}(x)\right)d\hat{U}_{1}^{\delta}duds.\label{eq:mathcalN1delta}
\end{equation}
\begin{comment}
This yields a term \hb{(details needed below)}
\begin{equation}\label{eq:mathcalN1delta}
\mathcal{N}_{t}^{1,\delta}(x,x')=\int_{0}^{t}\left(\int_{-\eta}^{0}-\int_{-\eta+\delta}^{\delta}\right)
\int_{x}^{x-u}\left(\tilde{v}_{t-s}^{\delta}(z)-\tilde{v}_{t-s}^{\delta}(x)\right)d\hat{U}_1^{\delta}(z)duds \, .
\end{equation}
\end{comment}
In order to estimate this term, we proceed as in Proposition \ref{prop:discrete-FPT} and Lemma \ref{eq:RIEst}. Namely, we bound the integral with respect to $\hat{U}_1^{\delta}$ by four terms like in \eqref{eq:estim-fourterms}. For notational sake, we will just focus on the fourth term which will be called $\mathcal{N}_{t}^{14,\delta}(x,x')$. It reads 
\begin{multline}\label{eq:mathcalN14}
\mathcal{N}_{t}^{14,\delta}(x,x')=\frac{1}{\delta^{3}}\kappa(\hat{{\bf U}}_1^{\delta})\Theta(\tilde{v}^{\delta})E(a,t)a^{2\chi}\lambda^{\frac{\alpha-\beta}{4}}\\\times\int_{0}^{t}\left(\int_{-\eta}^{0}-\int_{-\eta+\delta}^{\delta}\right)e^{-(\lambda+\theta a)s}|t-s|^{-\beta/2}|u|^{\alpha+2\beta}duds \, .
\end{multline}
The above integral is handled by elementary methods. Since $0<\eta<\delta$, the quantity $\mathcal{N}_{t}^{14,\delta}$ can be upper bounded as 
\begin{multline}\label{eq:re-mathcalN14}
\mathcal{N}_{t}^{14,\delta}(x,x')  \\
\le \frac{1}{\delta^{3}}\kappa(\hat{{\bf U}}_1^{\delta})\Theta(\tilde{v}^{\delta})E(a,t)a^{2\chi}\lambda^{\frac{\alpha-\beta}{4}}
\int_{0}^{t}\int_{[-\eta,0]\cup[-\eta+\delta,\delta]}
e^{-(\lambda+\theta a)s}|t-s|^{-\beta/2}|u|^{\alpha+2\beta}duds \, .
\end{multline}
Then an explicit computation reveals that 
\begin{align}
 & \frac{1}{\delta^{3}}\int_{0}^{t}(t-s)^{-\beta/2}ds\int_{-\eta+\delta}^{\delta}|u|^{\alpha+2\beta}du\nonumber \\
 & \lesssim\frac{1}{\delta^{3}}t^{1-\beta/2}\cdot\delta^{\alpha+2\beta}\cdot|\eta|=\frac{1}{\delta^{3}}t^{1-\beta/2}\cdot\delta^{\alpha+2\beta}\cdot|\eta|^{1-2\beta}|\eta|^{2\beta}
 \leqslant|\eta|^{2\beta}\cdot\delta^{\alpha-\beta},\label{eq:R0DelT4}
\end{align}
where we have used the assumptions $t<\delta^{2}$ and $|\eta|<\delta$
to reach the last inequality. Moreover, the same kind of inequality holds true for the integral on $[-\eta,0]$. Hence plugging relation~\eqref{eq:R0DelT4} into \eqref{eq:re-mathcalN14} we get 
\begin{eqnarray}\label{eq:last-estim-N14}
\big|\mathcal{N}_{t}^{14,\delta}(x,x')\big|
&\lesssim&
\kappa(\hat{{\bf U}}_1^{\delta})\Theta(\tilde{v}^{\delta})E(a,t)\lambda^{\frac{\alpha-\beta}{4}}
e^{-(\lambda+\theta a)s}a^{2\chi}|\eta|^{2\beta}\delta^{\alpha-\beta} \notag\\
&\lesssim& 
\kappa(\hat{{\bf U}}_1^{\delta})\Theta(\tilde{v}^{\delta})E(a,t)\lambda^{\frac{\alpha-\beta}{4}}
a^{\chi+\beta/2}|x'-x|^{2\beta}\cdot\delta^{\alpha-\beta},
\end{eqnarray}
where we have used the relation $\chi<\beta/2$, $a\geqslant 1$ and where we recall that $\eta=x'-x$ for the last inequality. We now gather 
\eqref{eq:re-mathcalN14} with similar estimates coming from the integral with respect to $\hat{U}_1^{\delta}$, and we plug those into 
\eqref{eq:mathcalN1delta}. We end up with 
\begin{equation}\label{eq:final-estim-mathcalN1}
\big|\mathcal{N}_{t}^{1,\delta}(x,x')\big|
\lesssim\kappa(\hat{{\bf U}}_1^{\delta})\Theta(\tilde{v}^{\delta})E(a,t)\lambda^{\frac{\alpha-\beta}{4}}a^{\chi+\beta/2}|x'-x|^{2\beta}\cdot\delta^{\alpha-\beta},
\end{equation}
whenever $|x'-x|<\delta$.

In order to bound the term $\mathcal{N}^{1,\delta}$ in \eqref{eq:mathcalN1delta} when $|\eta|=|x-x'|\geqslant\delta$,  we can directly apply the triangle inequality
to the $u$-integral in (\ref{eq:N1Delta}) and the resulting estimate
is the same as in~\eqref{eq:final-estim-mathcalN1}. Regions corresponding to other cases in (\ref{eq:Nab2p0Del})
lead to the same estimate as well. To conclude, we arrive at the following
estimate: 
\begin{equation}\label{eq:final-estimR0}
|\mathcal{R}_{t}^{0,\delta}(x,x')|\lesssim\kappa(\hat{{\bf U}}_1^{\delta})\Theta(\tilde{v}^{\delta})\lambda^{\frac{\alpha-\beta}{4}}E(a,t)a^{\chi+\beta/2}|x'-x|^{2\beta}\cdot\delta^{\alpha-\beta}.
\end{equation}

\noindent(D) We estimate the quantity 
\begin{equation}\label{eq:mathcalR1}
\mathcal{R}_{t}^{1,\delta}(x,x')\triangleq\hat{P}_{t}^{\delta}\hat{U}_1^{\delta}(x,x')+\mathcal{Q}_{t}^{\delta}(x,x')-\mathcal{K}^{\delta}(x,x') \, ,
\end{equation}
appearing in (\ref{eq:DisRemDecomp}), in one go. First observe that we are dealing with the case $t\leqslant\delta^{2}$, for which $\lfloor t\rfloor=0$. Hence the expressions in \eqref{j4}-\eqref{eq:RemEstJDel}-\eqref{eq:RemEstKDel} can be reduced to
\begin{align}
\mathcal{R}_{t}^{1,\delta}(x,x') & =\frac{t}{\delta^{2}}\int_{\mathbb{R}}\left(\hat{p}_{\delta^{2}}^{\delta}(\lfloor\xi-y\rfloor)-\hat{p}_{0}^{\delta}(\lfloor\xi-y\rfloor)\right)\hat{U}_1^{\delta}(y)dy\big|_{x}^{x'}+\hat{U}_1^{\delta}(x,x')\nonumber \\
 & =\frac{t}{\delta^{2}}\int_{\mathbb{R}}\left(\hat{p}_{\delta^{2}}^{\delta}(\lfloor z\rfloor)-\hat{p}_{0}^{\delta}(\lfloor z\rfloor)\right)
 \hat{U}_1^{\delta}(x-z,x'-z) \, dz
 +\hat{U}_1^{\delta}(x,x').\label{eq:R1Del}
\end{align}
To proceed further, we shall divide our discussion into two cases:
$|x'-x|^{2}\geqslant t$ or $|x'-x|^{2}<t$. The second case will be treated
by a generic argument later on (see Step 2 below). Here we only consider
the first case. Note that integrating against $\hat{p}_{0}^{\delta}$
or $\hat{p}_{\delta^{2}}^{\delta}$ is essentially averaging over
regions of order $\delta.$ By applying the triangle inequality to
(\ref{eq:R1Del}) in the obvious way and simply using the $\alpha$-H\"older
estimate for $\hat{U}_1^{\delta}$, we see that 
\begin{equation*}
|\mathcal{R}_{t}^{1,\delta}(x,x')|\lesssim\kappa(\hat{{\bf U}}_1^{\delta})a^{\chi}|x'-x|^{\alpha}=\kappa(\hat{{\bf U}}_1^{\delta})a^{\chi}t^{-\beta/2}\cdot t^{\beta/2}|x'-x|^{\alpha-2\beta}\cdot|x'-x|^{2\beta}.
\end{equation*}
In addition, recall from \eqref{f1} that $\beta>1/3$ and $\alpha<1/2$. In particular we have $\alpha<2\beta$. Hence owing to the fact that $|x-x'|^{2}\geqslant t$ we get 
\begin{align*}
|\mathcal{R}_{t}^{1,\delta}(x,x')|
 & \leqslant\kappa(\hat{{\bf U}}_1^{\delta})a^{\chi}t^{-\beta/2}\cdot t^{\beta/2}\cdot t{}^{\alpha/2-\beta}\cdot|x'-x|^{2\beta}\\
 & =\kappa(\hat{{\bf U}}_1^{\delta})a^{\chi}t^{-\beta/2}\cdot t{}^{\alpha/2-\beta/2}\cdot|x'-x|^{2\beta}\leqslant\kappa(\hat{{\bf U}}_1^{\delta})a^{\chi}t^{-\beta/2}|x'-x|^{2\beta}\cdot\delta^{\alpha-\beta},
\end{align*}
where we have invoked our standing assumption $t<\delta^{2}$ for the last inequality. With Proposition \ref{prop:discrete-FPT} in mind, we thus obtain the following estimate:
\begin{align}\label{eq:estimvtimesR1}
\big|2\tilde{v}_{t}^{\delta}(x)\cdot\mathcal{R}_{t}^{1,\delta}(x,x')\big| & \lesssim\kappa(\hat{{\bf U}}_1^{\delta})\Theta(\tilde{v}^{\delta})E(a,t)\cdot a^{\chi}t^{-\beta/2}|x'-x|^{2\beta}\cdot\delta^{\alpha-\beta}.
\end{align}

Summarizing our considerations so far, we report our bounds 
\eqref{eq:estimR0xxprime}, \eqref{eq:estim-diff-Ptxxprime}, \eqref{eq:final-estimR0} and \eqref{eq:estimvtimesR1} into \eqref{eq:ReMV-again} and~\eqref{eq:ReUhat}, we have achieved the following inequality for the case $t\leqslant\delta^{2}$:
\begin{multline}\label{j6}
\big|\mathcal{R}_{\hat{{\bf U}}_1^{\delta}}^{\mathcal{W}_{t}^{1,\delta}}(x,x')-\mathcal{R}_{{\bf W}}^{\mathcal{W}_{t}}(x,x')\big| \\
\lesssim
\left(\kappa(\hat{{\bf U}}_1^{\delta})\Theta(\tilde{v}^{\delta})+\kappa({\bf W})\Theta(v)\right)E(a,t)\lambda^{\frac{\alpha-\beta}{4}}a^{\chi}(a^{\beta}+t^{-\beta/2})|x'-x|^{2\beta}\cdot\delta^{\alpha-\beta}.
\end{multline}

\noindent
\textit{Step 2: Case $|x'-x|^{2}<t\wedge\delta^{2\eta}$.} Here we consider a parameter $\eta\in(0,1)$
 to be chosen later on, and we assume $|x'-x|^{2}<t\wedge\delta^{2\eta}$.
 In this case, we adopt a generic argument as in Case II of the time-variation
estimate. Under the same notation leading to (\ref{eq:GenericTV})
in that part, we have 
\begin{align*}
\big|\mathcal{R}_{\hat{{\bf U}}_1^{\delta}}^{\mathcal{W}_{t}^{1,\delta}}(x,x')-\mathcal{R}_{{\bf W}}^{\mathcal{W}_{t}}(x,x')\big| & \leqslant\big|\mathcal{R}_{\hat{{\bf U}}_1^{\delta}}^{\mathcal{W}_{t}^{1,\delta}}(x,x')\big|+\big|\mathcal{R}_{{\bf W}}^{\mathcal{W}_{t}}(x,x')\big|\\
 & \leqslant\Lambda^{\theta',\beta'}E^{\theta',\beta'}(a,t)\lambda^{\frac{\alpha-\beta'}{4}}\cdot a^{\chi}(a^{\beta'/2}+t^{-\beta/2})|x'-x|^{2\beta'},
\end{align*}
where $\Lambda^{\theta',\beta'}$ is the uniform upper bound on $\Theta^{\theta',\beta'}(w^{\delta})\vee\Theta^{\theta',\beta'}(w)$ introduced in~\eqref{i5}.
In view of the constraint (\ref{eq:ParConTV}) on the parameters,
we can further write 
\begin{align*}
\big|\mathcal{R}_{\hat{{\bf U}}_1^{\delta}}^{\mathcal{W}_{t}^{1,\delta}}(x,x')-\mathcal{R}_{{\bf W}}^{\mathcal{W}_{t}}(x,x')\big| & \leqslant\Lambda^{\theta',\beta'}E(a,t)\lambda^{\frac{\alpha-\beta}{4}}\cdot a^{\chi}|x'-x|^{2\beta}\\
 & \ \ \ \cdot\left(a^{\beta/2}|x'-x|^{2(\beta'-\beta)}+t^{-\beta/2}\cdot t^{-\frac{\beta'-\beta}{2}}|x'-x|^{2(\beta'-\beta)}\right).
\end{align*}
Under the current case, we assume $|x'-x|\leqslant\delta^{\eta}$. Hence 
\[
|x'-x|^{2(\beta'-\beta)}\leqslant\delta^{2\eta(\beta'-\beta)},
\]
and since we also assume $t>|x'-x|^{2}$ we get
\[
t^{-\frac{\beta'-\beta}{2}}|x'-x|^{2(\beta'-\beta)}=t^{-\frac{\beta'-\beta}{2}}|x'-x|^{\beta'-\beta}\cdot|x'-x|^{\beta'-\beta}\leqslant\delta^{\eta(\beta'-\beta)}.
\]
As a result, we arrive at the following estimate whenever $|x'-x|^{2}<t\wedge \delta^{2\eta}$: 
\begin{equation}
\big|\mathcal{R}_{\hat{{\bf U}}_1^{\delta}}^{\mathcal{W}_{t}^{1,\delta}}(x,x')-\mathcal{R}_{{\bf W}}^{\mathcal{W}_{t}}(x,x')\big|\lesssim\Lambda^{\theta',\beta'}E(a,t)\lambda^{\frac{\alpha-\beta}{4}}a^{\chi}(a^{\beta/2}+t^{-\beta/2})|x'-x|^{2\beta}\cdot\delta^{\eta(\beta'-\beta)}.\label{eq:RemEstGen}
\end{equation}

\noindent
\textit{Step 3: Case  $t\geqslant\delta^{2}$ and $|x'-x|\geqslant\delta$}.
The following estimate is obtained in the same way leading to the
time-variation estimate in Lemma \ref{lem:UnifWDelvsW} and space-variation estimates in Lemma \ref{lem:SpVarWDelvsW}, based on the local central
limit theorem. We only state the final result here. Namely we obtain that
\begin{multline}\label{eq:R0delta-R0}
\big|\mathcal{R}_{t}^{0,\delta}(x,x')-\mathcal{R}_{t}^{0}(x,x')\big|  \lesssim E(a,t)\left(\kappa(\hat{{\bf U}}_1^{\delta})d_{\hat{{\bf U}}_1^{\delta},{\bf W}}(\tilde{v}^{\delta},v)+\Theta(v)\rho(\hat{{\bf U}}_1^{\delta},{\bf W})\right)\\
 \cdot a^{\chi}(a^{\beta/2}+t^{-\beta/2})|x'-x|^{2\beta} \, .
\end{multline}
As far as the terms involving $P_{t}$ in \eqref{eq:DisRemDecomp} are concerned, we first have 
\begin{equation}\label{eq:splitA=B}
\left|\tilde{v}_{t}^{\delta}(x)\hat{P}_{t}^{\delta}\hat{U}_{1}^{\delta}(x,x')-v_{t}(x)P_{t}W(x,x')\right|\leq \mathcal{A}+\mathcal{B},
\end{equation}
where $\mathcal{A}$ and $\mathcal{B}$ are respectively defined by 
\begin{equation}\label{eq:A}
\mathcal{A}=\big|\tilde{v}_{t}^{\delta}(x)\big|\cdot\left|\hat{P}_{t}^{\delta}\hat{U}_{1}^{\delta}(x,x')-P_{t}W(x,x')\right|
\end{equation}
\begin{equation}\label{eq:B}
\mathcal{B}=\big|\tilde{v}_{t}^{\delta}(x)-v_{t}(x)\big|\cdot\big|P_{t}W(x,x')\big|.
\end{equation}
According to Lemma \ref{lem:SGEst}, the quantity $\mathcal{B}$ is bounded as 
\begin{equation}
\mathcal{B}\lesssim E(a,t)d_{\hat{\mathbf{U}}_{1}^{\delta},{\bf W}}(\tilde{v}^{\delta},v)\cdot\kappa({\bf W})a^{\chi}t^{-\beta/2}|x'-x|^{2\beta}.\label{eq:BEstPtTerm}
\end{equation}
To estimate $\mathcal{A}$, we further write 
\begin{align}
 & \left|\hat{P}_{t}^{\delta}\hat{U}_{1}^{\delta}(x,x')-P_{t}W(x,x')\right|\nonumber\\
 & \leqslant\big|\big(\hat{P}_{t}^{\delta}-P_{t}\big)\hat{U}_{1}^{\delta}(x,x')\big|+\big|P_{t}\big(\hat{U}_{1}^{\delta}-W\big)(x,x')\big|=:\mathcal{A}_{1}+\mathcal{A}_{2}.\label{eq:splitA1+A2}
\end{align}
As in Lemma \ref{lem:SGEst}, the $\mathcal{A}_{2}$-term is easily
estimated as 
\begin{equation}\label{eq:A2upper}
\big|P_{t}\big(\hat{U}_{1}^{\delta}-W\big)(x,x')\big|\lesssim\rho\big(\hat{{\bf U}}_{1}^{\delta},{\bf W}\big)a^{\chi}t^{-\beta/2}|x'-x|^{2\beta}.
\end{equation}
To estimate the $\mathcal{A}_{1}$-term, for simplicity we only consider
the case when $|x'-x|^{2}>t$ (the other case is treated in the same
way as Case (i) in the proof of Lemma \ref{lem:SGEst}). Now for $|x'-x|^{2}>t$,
we write
\begin{equation}\label{eq:A1case>t}
\big(\hat{P}_{t}^{\delta}-P_{t}\big)\hat{U}_{1}^{\delta}(x,x') =\int_{\mathbb{R}}\left(\hat{p}_{\lfloor t\rfloor}^{\delta}(\lfloor y\rfloor)-p_{t}(y)\right)\big(\hat{U}_{1}^{\delta}(x'-y)-\hat{U}_{1}^{\delta}(x-y)\big)dy.
\end{equation}
Recall from the local CLT (cf. Theorem \ref{thm:LCLT}) and the uniform
Gaussian estimates (cf. Proposition \ref{prop:DiscUnifGau}) for the
heat kernels that 
\[
\big|\hat{p}_{\lfloor t\rfloor}^{\delta}(\lfloor y\rfloor)-p_{t}(y)\big|\lesssim\frac{\delta^{2}}{t^{3/2}},
\quad\text{ and}\quad
\big|\hat{p}_{\lfloor t\rfloor}^{\delta}(\lfloor y\rfloor)\big|\vee\big|p_{t}(y)\big|\lesssim\frac{1}{\sqrt{t}}e^{-Cy^{2}/t}.
\]
Plugging those inequalities into \eqref{eq:A1case>t}, it follows that
\begin{align}
\mathcal{A}_{1} & \leqslant\int_{\mathbb{R}}\big|\hat{p}_{\lfloor t\rfloor}^{\delta}(\lfloor y\rfloor)-p_{t}(y)\big|^{1/2}\cdot\big|\hat{p}_{\lfloor t\rfloor}^{\delta}(\lfloor y\rfloor)-p_{t}(y)\big|^{1/2}
\cdot\big|\hat{U}_{1}^{\delta}(x'-y)-\hat{U}_{1}^{\delta}(x-y)\big|dy\nonumber \\
 & \lesssim\frac{\delta}{t^{3/4}}\cdot t^{1/4}\cdot\int_{\mathbb{R}}\frac{1}{\sqrt{t}}e^{-Cy^{2}/2t}\big|\hat{U}_{1}^{\delta}(x'-y)-\hat{U}_{1}^{\delta}(x-y)\big|dy\nonumber \\
 & =\frac{\delta}{\sqrt{t}}\cdot\int_{\mathbb{R}}e^{-Cw^{2}/2}\big|\hat{U}_{1}^{\delta}(x'-\sqrt{t}w)-\hat{U}_{1}^{\delta}(x-\sqrt{t}w)\big|dw\nonumber \\
 & \lesssim\frac{\delta}{\sqrt{t}}\cdot a^{\chi}\kappa(\hat{{\bf U}}_{1}^{\delta})\,|x'-x|^{\alpha}.\label{eq:HatPPDif}
\end{align}
Similar to the analysis leading to (\ref{eq:estim-incr-PtW-ter}),
we write 
\[
|x'-x|^{\alpha}=t^{-\beta/2}|x'-x|^{2\beta}\cdot t^{\beta/2}|x'-x|^{\alpha-2\beta}.
\]
Noting that $\alpha-2\beta<0$ and $|x'-x|^{2}>t$ in the current
case, simple algebra shows that 
\[
t^{\beta/2}|x'-x|^{\alpha-2\beta}<t^{\frac{\alpha-\beta}{2}}\quad\text{and}\quad t^{-\beta/2}|x-x'|^{2\beta}<t^{\beta/2}\leq T^{\beta/2}.
\]
Therefore we get
\[
|x'-x|^{\alpha}\leq C_{T}\,t^{\frac{\alpha-\beta}{2}}. 
\]
Plugging this inequality in \eqref{eq:HatPPDif} we obtain 
\begin{equation}\label{eq:A1upper}
\mathcal{A}_{1}\lesssim a^{\chi}\kappa(\hat{{\bf U}}_{1}^{\delta})\cdot\frac{\delta}{\sqrt{t}}t^{\frac{\alpha-\beta}{2}}=a^{\chi}\kappa(\hat{{\bf U}}_{1}^{\delta})\cdot\delta^{\alpha-\beta}\cdot\left(\frac{\delta}{\sqrt{t}}\right)^{1-(\alpha-\beta)}\leqslant a^{\chi}\kappa(\hat{{\bf U}}_{1}^{\delta})\cdot\delta^{\alpha-\beta},
\end{equation}
where the last inequality holds since $t\geqslant\delta^{2}$ in the
current scenario. To summarise, we report \eqref{eq:A1upper} and \eqref{eq:A2upper} into \eqref{eq:splitA1+A2}. This yields 
\[
\mathcal{A}\lesssim\Theta(v)\left(\kappa(\hat{{\bf U}}_{1}^{\delta})\delta^{\alpha-\beta}+\rho\big(\hat{{\bf U}}_{1}^{\delta},{\bf W}\big)\right)\cdot E(a,t)a^{\chi}t^{-\beta/2}|x'-x|^{2\beta}.
\]
Together with (\ref{eq:BEstPtTerm}) and recalling \eqref{eq:splitA=B}, we thus arrive at the following estimate:
\begin{align}
 & \big|\tilde{v}_{t}^{\delta}(x)\hat{P}_{t}^{\delta}\hat{U}_{1}^{\delta}(x,x')-v_{t}(x)P_{t}W(x,x')\big|\label{eq:vdeltaPhatdeltaUhatdelta-vPW}\\
 & \lesssim\left[\Theta(v)\left(\kappa(\hat{{\bf U}}_{1}^{\delta})\delta^{\alpha-\beta}+\rho(\hat{{\bf U}}_{1}^{\delta},{\bf W})\right)+\kappa({\bf W})d_{\hat{{\bf U}}_{1}^{\delta},{\bf W}}(\tilde{v}^{\delta},v)\right]\cdot E(a,t)a^{\chi}t^{-\beta/2}|x'-x|^{2\beta}.\nonumber
\end{align}

\begin{comment}
\begin{multline}\label{eq:vdeltaPhatdeltaUhatdelta-vPW}
 \big|\tilde{v}_{t}^{\delta}(x)(\hat{P}_{t}^{\delta}\hat{U}_1^{\delta})(x,x')-v_{t}(x)(P_{t}W)(x,x')\big|\\
  \lesssim E(a,t)\left(\kappa(\hat{{\bf U}}_1^{\delta})d_{\hat{{\bf U}}_1^{\delta},{\bf W}}(\tilde{v}^{\delta},v)
 +\Theta(v)(\rho(\hat{{\bf U}}_1^{\delta},{\bf W})
 +\kappa(\hat{{\bf U}}_1^{\delta})\delta^{\alpha-\beta})\right)a^{\chi}t^{-\beta/2}|x'-x|^{2\beta}.
\end{multline}
\end{comment}

Now it remains to estimate the quantities $\mathcal{Q}_{t}^{\delta}(x,x'),$
$\mathcal{K}_{t}^{\delta}(x,x')$ defined by \eqref{eq:RemEstJDel},
\eqref{eq:RemEstKDel} respectively. We will divide those estimates in two substeps.

\noindent
\textit{Step 3-1:  estimate for  $\mathcal{Q}_{t}^{\delta}$.} Our  claim here is  that whenever $t\geqslant\delta^{2}$ and $|x'-x|\geqslant\delta$ we have:
\begin{equation}\label{j7}
|\mathcal{Q}_{t}^{\delta}(x,x')|\lesssim\kappa(\hat{{\bf U}}_1^{\delta})a^{\chi}t^{-\beta/2}|x'-x|^{2\beta}\cdot\delta^{\alpha-\beta}.
\end{equation}
In order to prove \eqref{j7}, we return to the expression of $\mathcal{Q}_{t}^{\delta}(x,x')$ in
terms of $\nabla_{x}^{2,\delta}$. Namely invoking the discrete heat equation \eqref{eq:discrete-heat-equation} in \eqref{eq:RemEstJDel} we get
\begin{equation}\label{eq:calJwithNabla}
\mathcal{Q}_{t}^{\delta}(x,x')=\frac{(t-\lfloor t\rfloor)\sigma^{2}}{2}\int_{\mathbb{R}}\nabla_{x}^{2,\delta}\hat{p}_{\lfloor t\rfloor}^{\delta}(\lfloor\xi-y\rfloor)\hat{U}_1^{\delta}(y)dy\big|_{x}^{x'}.
\end{equation}
We now split the discussion in two cases.

\noindent
\textit{Case A}: $|x'-x|^{2}\geqslant t.$ In this case, set $y=\xi-\sqrt{t}w$ in \eqref{eq:calJwithNabla}. We obtain
\begin{equation*}
\mathcal{Q}_{t}^{\delta}(x,x')  =\frac{(t-\lfloor t\rfloor)\sigma^{2}}{2}\int_{\mathbb{R}}\nabla_{x}^{2,\delta}\hat{p}_{\lfloor t\rfloor}^{\delta}(\lfloor\sqrt{t}w\rfloor)
\left(\hat{U}_1^{\delta}(x'-\sqrt{t}w)-\hat{U}_1^{\delta}(x-\sqrt{t}w)\right)\sqrt{t}dw.
\end{equation*}
We resort to the fact that $t-\lfloor t\rfloor\leqslant\delta^{2}$ and to the H\"older continuity of $\hat{U}_1^{\delta}$. Taking into account the upper bound  (\ref{eq:DiscUnifGau>0}), which reads 
\[
\big|\nabla_{x}^{2,\delta}\hat{p}_{\lfloor t\rfloor}^{\delta}(\lfloor\sqrt{t}w\rfloor)\big|\lesssim\frac{1}{t^{3/2}}e^{-Cw^{2}}\quad\text{ for all }\quad t\geqslant\delta^{2},w\in\mathbb{R},
\]
we end up with
\begin{align*}
\mathcal{Q}_{t}^{\delta}(x,x') & \lesssim\delta^{2}\cdot\kappa(\hat{{\bf U}}_1^{\delta})a^{\chi}\cdot t^{-1}|x'-x|^{\alpha}\\
 & =\kappa(\hat{{\bf U}}_1^{\delta})a^{\chi}\cdot t^{-\beta/2}|x'-x|^{2\beta}\cdot\delta^{2}\cdot t^{\beta/2-1}|x'-x|^{\alpha-2\beta}.
 \end{align*}
 Moreover, since we assume $|x-x'|^2\geqslant t$ and we have 
 $\alpha<2\beta$, we discover that 
 \begin{align*}
\mathcal{Q}_{t}^{\delta}(x,x') & \lesssim\kappa(\hat{{\bf U}}_1^{\delta})a^{\chi}\cdot t^{-\beta/2}|x'-x|^{2\beta}\cdot\delta^{2}\cdot t^{\beta/2-1}t{}^{\alpha/2-\beta}\\
 & =\kappa(\hat{{\bf U}}_1^{\delta})a^{\chi}\cdot t^{-\beta/2}|x'-x|^{2\beta}\cdot\delta^{2}\cdot t{}^{\alpha/2-\beta/2-1}.
 \end{align*}
 Eventually recall that $t>\delta^{2}$ in this step, and $\alpha-\beta-2<0$. This yields 
 \begin{equation}\label{J7b}
 \mathcal{Q}_{t}^{\delta}(x,x') \lesssim
\kappa(\hat{{\bf U}}_1^{\delta})a^{\chi}t^{-\beta/2}|x'-x|^{2\beta}\cdot\delta^{\alpha-\beta},
 \end{equation}
which is compatible with our claim \eqref{eq:RemEstComp}.

\noindent
\textit{Case B:} $|x'-x|^{2}<t$. Since $\int_{\mathbb{R}}\nabla_{x}^{2,\delta}\hat{p}_{\lfloor t\rfloor}^{\delta}(y)dy=0$,
we can insert this quantity in the expression~\eqref{eq:calJwithNabla} for 
$\mathcal{Q}_{t}^{\delta}(x,x')$. This yields
\begin{align}
\mathcal{Q}_{t}^{\delta}(x,x') & =\frac{(t-\lfloor t\rfloor)\sigma^{2}}{2}\int_{\mathbb{R}}\left(\nabla_{x}^{2,\delta}\hat{p}_{\lfloor t\rfloor}^{\delta}(\lfloor x'-y\rfloor)-\nabla_{x}^{2,\delta}\hat{p}_{\lfloor t\rfloor}^{\delta}(\lfloor x-y\rfloor)\right)\left(\hat{U}_1^{\delta}(y)-\hat{U}_1^{\delta}(x)\right)dy\nonumber\\
 & =\frac{(t-\lfloor t\rfloor)\sigma^{2}}{2}\int_{\mathbb{R}}\int_{\lfloor x-y\rfloor}^{\lfloor x'-y\rfloor}\tilde{\nabla}^{3,\delta}\hat{p}_{\lfloor t\rfloor}^{\delta}(\lfloor z\rfloor)\left(\hat{U}_1^{\delta}(y)-\hat{U}_1^{\delta}(x)\right)dzdy,\label{J8}
\end{align}
where the discrete gradient $\tilde{\nabla}^{3,\delta}$ is given by \eqref{d11}:
\[
(\tilde{\nabla}^{3,\delta}f)(x)\triangleq\frac{(\nabla_{x}^{2,\delta}f)(x+\delta)-(\nabla_{x}^{2,\delta}f)(x)}{\delta}.
\]
Moreover, the  third order difference $\tilde{\nabla}^{3,\delta}\hat{p}_{\lfloor t\rfloor}^{\delta}(\lfloor z\rfloor)$
also satisfies the uniform Gaussian estimate in Proposition \ref{prop:DiscUnifGau}. Namely we have 
\begin{equation}\label{eq:estimGaussnabla3}
\big|\tilde{\nabla}^{3,\delta}\hat{p}_{\lfloor t\rfloor}^{\delta}(\lfloor z\rfloor)\big|\lesssim\frac{1}{t^{2}}e^{-Cz^{2}/t}\quad\text{ for all }\quad t\geqslant\delta^{2},\ z\in\mathbb{R}.
\end{equation}
Suppose that $\lfloor x-y\rfloor$ and $\lfloor x'-y\rfloor$ have
the same sign, say 
\[
\lfloor x-y\rfloor\geqslant\lfloor x'-y\rfloor\geqslant0.
\]
Since $|x'-x|\geqslant\delta$ under the current discussion, it is
easily seen that 
\begin{equation}\label{eq:exponential-integral-estim}
\int_{\lfloor x'-y\rfloor}^{\lfloor x-y\rfloor}e^{-Cz^{2}/t}dz\lesssim|x'-x|e^{-C(x'-y)^{2}/t}.
\end{equation}
Plugging \eqref{eq:estimGaussnabla3} and \eqref{eq:exponential-integral-estim} into \eqref{J8}, we obtain 
\begin{align}
|\mathcal{Q}_{t}^{\delta}(x,x')|& \lesssim\delta^{2}t^{-2}|x'-x|\cdot\int_{\mathbb{R}}e^{-C(x'-y)^{2}/t}\cdot\big|\hat{U}_1^{\delta}(y)-\hat{U}_1^{\delta}(x)\big|dy\nonumber\\
& =\delta^{2}t^{-3/2}|x'-x|\cdot\int_{\mathbb{R}}e^{-Cw^{2}}\big|\hat{U}_1^{\delta}(x'-\sqrt{t}w)-\hat{U}_1^{\delta}(x)\big|dw,\label{J9}
\end{align}
where we have set $y=x'-\sqrt{t}w$ for the second identity. Next, owing to the H\"older regularity of $\hat{U}_1^{\delta}$ and invoking $|x'-x|\le\sqrt{t}$ again (also noting $x,x'\in[-a,a]$ and $a\geqslant1$), we have 
\begin{align}
\big|\hat{U}_1^{\delta}(x'-\sqrt{t}w)-\hat{U}_1^{\delta}(x)\big|  & \leqslant\kappa(\hat{{\bf U}}_1^{\delta})
(|x'|+|x|+\sqrt{t}|w|)^{\chi}\cdot(|x'-x|+\sqrt{t}|w|)^{\alpha}\nonumber\\
& \leqslant \kappa(\hat{{\bf U}}_{1}^{\delta})(2a+\sqrt{T}|w|)^{\chi}\cdot(\sqrt{t}+\sqrt{t}|w|)^{\alpha}\nonumber\\
& \leqslant\kappa(\hat{{\bf U}}_1^{\delta})a^{\chi}t^{\alpha/2}
(2+2\sqrt{T}|w|)^\chi(1+|w|)^{\alpha}.\label{eq:diffUhatdelta-Uhatdelta}
\end{align}
Reporting this inequality into \eqref{J9}, we can write 
\begin{align*}
|\mathcal{Q}_{t}^{\delta}(x,x')| & \lesssim\kappa(\hat{{\bf U}}_1^{\delta})a^{\chi}\cdot\delta^{2}t^{\alpha/2-3/2}|x'-x|\\
 & =\kappa(\hat{{\bf U}}_1^{\delta})a^{\chi}t^{-\beta/2}|x'-x|^{2\beta}\cdot\delta^{2}t^{\alpha/2+\beta/2-3/2}|x'-x|^{1-2\beta}.
 \end{align*}
 Now recall once more that $|x'-x|<\sqrt{t}$ and $t\geqslant\delta^{2}$. We thus obtain 
 \begin{align}
 |\mathcal{Q}_{t}^{\delta}(x,x')|  & \lesssim \kappa(\hat{{\bf U}}_1^{\delta})a^{\chi}t^{-\beta/2}|x'-x|^{2\beta}\cdot\delta^{2}t^{\alpha/2-\beta/2-1}\nonumber\\
&\leqslant\kappa(\hat{{\bf U}}_1^{\delta})a^{\chi}t^{-\beta/2}|x'-x|^{2\beta}\cdot\delta^{\alpha-\beta}.\label{J10}
\end{align}
This fits our claim \eqref{eq:R0delta-R0} when $\lfloor x-y\rfloor$ and $\lfloor x'-y\rfloor$ have the same sign. 
If $\lfloor x-y\rfloor$ and $\lfloor x'-y\rfloor$ have different
signs, say 
\[
x-y\leqslant0\leqslant x'-y\iff x\leqslant y\leqslant x',
\]
then we simply bound (owing to \eqref{eq:estimGaussnabla3}) the quantity 
$\tilde{\nabla}^{3,\delta}\hat{p}_{\lfloor t\rfloor}^{\delta}$ by $t^{-2}$. Plugging this crude estimate in \eqref{J8} we get 
\begin{align*}
\mathcal{Q}_{t}^{\delta}(x,x') & \lesssim\delta^{2}t^{-2}|x'-x|\cdot\int_{x}^{x'}\big|\hat{U}_1^{\delta}(y)-\hat{U}_1^{\delta}(x)\big|dy\\
 & \leqslant\delta^{2}t^{-2}|x'-x|^{2}\cdot\kappa(\hat{{\bf U}}_1^{\delta})a^{\chi}|x'-x|^{\alpha}\\
 & =\kappa(\hat{{\bf U}}_1^{\delta})a^{\chi}t^{-\beta/2}|x'-x|^{2\beta}\cdot\delta^{2}t^{-2+\beta/2}|x'-x|^{2+\alpha-2\beta}.
 \end{align*}
 Thanks to the assumption $|x'-x|<\sqrt{t}$ and $t\geqslant\delta^{2}$, we get
 \begin{equation}\label{J11}
 \mathcal{Q}_{t}^{\delta}(x,x')  \lesssim \kappa(\hat{{\bf U}}_1^{\delta})a^{\chi}t^{-\beta/2}|x'-x|^{2\beta}\cdot\delta^{\alpha-\beta}.
\end{equation}
Gathering \eqref{J10}, \eqref{J11} and \eqref{J7b}, we have thus obtained an inequality for $\mathcal{Q}_{t}^{\delta}(x,x')$ which is compatible with~\eqref{eq:R0delta-R0}.

\noindent
\textit{Step 3-2: Estimate for $\mathcal{K}^{\delta}$.}
Let us focus our attention on the term $\mathcal{K}^{\delta}(x,x')$ defined by~\eqref{eq:RemEstKDel}. Here 
under the current case, if $|x'-x|\geqslant\sqrt{t},$ we assert that 
\begin{equation}
|\mathcal{K}^{\delta}(x,x')|\lesssim\kappa(\hat{{\bf U}}_1^{\delta})a^{\chi}t^{-\beta/2}|x'-x|^{2\beta}\cdot\delta^{\alpha-\beta}.\label{eq:DisHSGEstInit1}
\end{equation}
If $|x'-x|\geqslant\delta^{\tau}$ where $\tau\in(0,1)$ is a given
number to be chosen later on, we will see that 
\begin{equation}
|\mathcal{K}^{\delta}(x,x')|\lesssim\kappa(\hat{{\bf U}}_1^{\delta})a^{\chi}t^{-\beta/2}|x'-x|^{2\beta}\cdot\delta^{\alpha-2\beta\tau}.\label{eq:DisHSGEstInit2}
\end{equation}
The proof of those claims goes as follows: by the definition of $\mathcal{K}^{\delta}(x,x')$, we can write 
\begin{align*}
\mathcal{K}^{\delta}(x,x') & =\int_{\mathbb{R}}\hat{p}_{0}^{\delta}(\lfloor y\rfloor)\left(\hat{U}_1^{\delta}(x',x'-y)-\hat{U}_1^{\delta}(x,x-y)\right)dy\\
 & =\frac{1}{\delta}\left(\int_{0}^{\delta}\hat{U}_1^{\delta}(x',x'-y)dy-\int_{0}^{\delta}\hat{U}_1^{\delta}(x,x-y)dy\right).
\end{align*}
As a result, we have 
\[
|\mathcal{K}^{\delta}(x,x')|\leqslant\frac{1}{\delta}\int_{0}^{\delta}\left(|\hat{U}_1^{\delta}(x',x'-y)|+|\hat{U}_1^{\delta}(x,x-y)|\right)dy.
\]
By using the standard $\alpha$-H\"older estimate for $\hat{U}_1^{\delta}$,
we obtain that 
\[
|\mathcal{K}^{\delta}(x,x')|\lesssim\kappa(\hat{{\bf U}}_1^{\delta})a^{\chi}\delta^{\alpha}.
\]

If $|x'-x|\geqslant\sqrt{t}$, since $t\geqslant\delta^{2}$ we can
further write 
\[
\delta^{\alpha}=\delta^{\alpha}|x'-x|^{-2\beta}|x'-x|^{2\beta}\leqslant\delta^{\alpha}t^{-\beta}|x'-x|^{2\beta}\leqslant t^{-\beta/2}|x'-x|^{2\beta}\cdot\delta^{\alpha-\beta},
\]
which leads to (\ref{eq:DisHSGEstInit1}). If $|x'-x|\geqslant\delta^{\tau},$
we have 
\[
\delta^{\alpha}=\delta^{\alpha-2\beta\tau}\delta^{2\beta\tau}t^{\beta/2}t^{-\beta/2}\leqslant C_{T}|x'-x|^{2\beta}t^{-\beta/2}\cdot\delta^{\alpha-2\beta\tau},
\]
which leads to (\ref{eq:DisHSGEstInit2}). %Notice that \eqref{eq:DisHSGEstInit1} is one of the estimates which will lead to our rate of order $\delta^{1/10}$. 

\noindent
\textit{Step 4: Conclusion.}
We first remark that Steps 1, 2 and 3 above cover all possibilities. Indeed,
the complement of Steps 1, 3 is the case when $t\geqslant\delta^{2}$
and $|x'-x|<\delta$. But this situation is contained in Step 2:
\[
t\geqslant\delta^{2},\ |x'-x|<\delta\implies|x'-x|^{2}<t\wedge\delta^{2\tau}
\]
since $\tau\in(0,1)$. In view of the estimates (\ref{eq:RemEstGen})
and (\ref{eq:DisHSGEstInit2}), we now choose $\tau$ to be such that
\[
\tau(\beta'-\beta)=\alpha-2\beta\tau\iff\tau=\frac{\alpha}{\beta'+\beta}.
\]
The resulting rate is $\delta^{r}$ with $r=\frac{\alpha}{\beta'+\beta}(\beta'-\beta)$.
Combining all the ingredients obtained so far, the desired estimate
(\ref{eq:RemEstComp}) follows.
\end{proof}

\subsection{Developing Step 2: comparing $\mathcal{J}^{\delta}$ and $\eta$}\label{sec:comparing-J-eta}

In Section \ref{sec:strategy-convergence} we mentioned that the convergence of $\tilde{v}^{\delta}$ to $v$  also relied on the convergence of a series of deterministic type terms. Specifically, consider the path $\eta$ defined by \eqref{eq:EtaTerm}. This path is approximated by the discrete process $\mathcal{J}_{t}^{\delta}$ introduced in~\eqref{eq:JDel}. The second step of the strategy recalled in Section~\ref{sec:strategy-convergence} leads to inequality 
\eqref{eq:Step2Det}. We will detail this step now. 
\begin{comment}
Our informal Proposition \ref{prop:exponent-1/10} does not specify the regularity assumptions one needs on the initial data in \eqref{eq:EtaTerm}. 
We start by introducing suitable weighted norms in order to quantify this regularity. 
\begin{defn}\label{def:spaceCrL}
Let $L>0$ and $r\in\mathbb{N}$. We say that a function $f_{0}$ of the spatial variable $x$ is an element of 
$\mathcal{C}_{L}^{r}$ if 
\[
\|f_{0}\|_{\mathcal{C}_{L}^{r}}\triangleq\sup_{a\geqslant1}a^{-L}\sup_{k\leqslant r,\,x\in[-a,a]}\big|D_{x}^{k}f_{0}(x)\big|<\infty.
\]
Similarly, we say that a function $g$ of the space-time parameter $(t,x)$ is an element of $\mathcal{C}_{L}^{r}$ if, denoting by  $D_{t,x}^{ij}g$ the derivative $\frac{\partial^{i+j}g}{\partial t^{i}\partial x^{j}}$, we have
\[
\|g\|_{\mathcal{C}_{L}^{r}}\triangleq\sup_{a\geqslant1}a^{-L}\sup_{\substack{i+j\leqslant r\\
t\in[0,T],\,x\in[-a,a]
}
}\big|D_{t,x}^{ij}g(t,x)\big|<\infty.
\]
\end{defn}
\end{comment}
Our estimate on the controlled distance between $\mathcal{J}^{\delta}$
and $\eta$ is summarized in the proposition below.
\begin{prop}
\label{prop:ComDetJEta} Let $\alpha,\beta,\chi,\theta$ be exponents satisfying condition \eqref{f1}. We consider fixed coefficients $\theta'<\theta$ and $\beta'\in(\beta,\alpha)$. The upper bound $\Lambda^{\theta',\beta'}$ on the rough path norms of $\cw^{1,\delta}$ and $\cw$ is given by \eqref{i5}. Recall that $\mathcal{J}^{\delta}$ is introduced in \eqref{eq:JDel}, 
while $\eta$ is the path given by \eqref{eq:EtaTerm}. Then we have 
\begin{equation}\label{eq:S2Est}
d_{\hat{{\bf U}}_1^{\delta},{\bf W}}(\mathcal{J}^{\delta},\eta) \leqslant C\left((\|f_{0}\|_{\mathcal{C}_{L}^{3}}+\|g\|_{\mathcal{C}_{L}^{3}})\delta^{\beta} + (1+\kappa(\hat{{\bf U}}_1^{\delta}))\lambda^{\frac{\alpha-\beta}{4}}(\Lambda^{\theta',\beta'}\delta^{\frac{\beta'(\beta'-\beta)}{\beta'+\beta}}+\Lambda^{\theta',\beta}\delta^{\alpha-2\chi})\right).
\end{equation}
\end{prop}

\begin{comment}
\begin{rem}\label{rem:delta1/2}
We should be aware that $\beta'$ is taken arbitrarily close to $\alpha$,
and also $\alpha=1/2^{-},\beta=1/3^{+}.$ As a result, the above rate
is (again) arbitrarily close to $1/10$. Note that when $\beta\rightarrow0$,
the $\delta$-exponent can be made arbitrarily close to $1/2$. The
same type of estimate also appears in Lemma \ref{lem:RemEstComp}.
\end{rem}
\end{comment}

The rest of this subsection is devoted to the proof of Proposition
\ref{prop:ComDetJEta}, for which we now outline a strategy.
In view of the definition of $d_{\hat{U}^{\delta},{\bf W}}(\mathcal{J}^{\delta},\eta)$,
the comparison between $\mathcal{J}^{\delta}$ and $\eta$ boils down
to estimating four types of differences: uniform, time-variation,
space-variation and remainder. Since a substantial part of the analysis
here is a technical repetition of the previous section (indeed it
is simpler than Step 3 except for the estimation of one term which
we will point out later on), for most of the time we will only consider
uniform distance estimates. In addition, we will only consider $\mathcal{J}_{t}^{\delta}(x)-\eta_{t}(x)$
on grid points $(t,x)\in\delta^{2}\mathbb{N}\times\delta\mathbb{Z}$.
The adaptation to non-grid points is easy since $\tilde{v}_{t}^{\delta}(x)$
is defined through piecewise linear interpolation (cf. Definition \ref{def:linear-interpolation-rough-path-v}).

\subsubsection{A decomposition of $\mathcal{J}_{t}^{\delta}(x)$}\label{decompJdelta}

Recall that in \eqref{eq:JDel}, the quantity $\mathcal{J}_{t}^{\delta}$ was implicit.  As a starting point, we shall first compute $\mathcal{J}_{t}^{\delta}(x)$
explicitly. This is the content of the lemma below. 
\begin{lem}
\label{lem:JdelDecomp} Recall that $v_{t_{k}}^{\delta}(x)$
is the discrete process defined by (\ref{eq:disc-PDE-vdelta}), in which $\eta_{t_{k}}^{\delta}(x)$
and $\mathcal{I}_{t_{k}}^{\delta}(x,y)$ are given by (\ref{eq:disc-eta-delta}) and (\ref{eq:disc-Ideltaxy})
respectively. For each $(t_{k},x)\in\delta^{2}\mathbb{N}\times\delta\mathbb{Z}$,
we have 
\begin{equation}\label{eq:explicitcalJdelta}
\mathcal{J}_{t_{k}}^{\delta}(x)=\frac{1}{2}\left(\eta_{t_{k}}^{\delta}(x)+\eta_{t_{k+1}}^{\delta}(x)\right)+\frac{1}{2}\left(v_{t_{k}}^{\delta}(x)-v_{t_{k+1}}^{\delta}(x)\right)+\mathcal{E}_{t_{k}}^{\delta}(x).
\end{equation}
The function $\mathcal{E}_{t_{k}}^{\delta}(x)$ in the above equation
is an error term defined by 
\[
\mathcal{E}_{t_{k}}^{\delta}(x)=\mathcal{E}_{t_{k}}^{1,\delta}(x)+\mathcal{E}_{t_{k}}^{2,\delta}(x),
\]
where $\mathcal{E}^{1,\delta}$ and $\mathcal{E}^{2,\delta}$ are respectively given by 
\begin{eqnarray}
\mathcal{E}_{t_{k}}^{1,\delta}(x) 
&\triangleq&
-\frac{1}{4}\delta^{3}\sum_{z\in\delta\mathbb{Z}}\nabla_{x}^{2,\delta}\hat{p}_{t_{k}}^{\delta}(\lfloor x-z\rfloor)\mathcal{I}_{0}^{\delta}(x,z)\label{eq:S2ErrorE1Del}\\
\mathcal{E}_{t_{k}}^{2,\delta}(x) 
&\triangleq&
-\frac{1}{2}\int_{0}^{t_{k}}\int_{\mathbb{R}}\nabla_{x}^{2,\delta}\hat{p}_{\lfloor s\rfloor}^{\delta}(\lfloor x-y\rfloor)\left(\int_{y}^{\lfloor y\rfloor+\delta}\tilde{v}_{t_{k}-s}^{\delta}(z)d\hat{U}^{\delta}(z)\right)dyds.\label{eq:S2ErrorE2Del}
\end{eqnarray}
\end{lem}

\begin{proof}
Using the definition (\ref{eq:JDel}) of $\mathcal{J}_{t_{k}}^{\delta}(x)$,
we first write 
\begin{equation}\label{eq:first-step-explicit-Jdelta}
\mathcal{J}_{t_{k}}^{\delta}(x)=\tilde{v}_{t_{k}}^{\delta}(x)+\frac{1}{2}\int_{0}^{t_{k}}\int_{\mathbb{R}}\nabla_{x}^{2,\delta}\hat{p}_{\lfloor t_{k}-s\rfloor}^{\delta}(\lfloor x-y\rfloor)\left(\int_{x}^{\lfloor y\rfloor+\delta}\tilde{v}_{s}^{\delta}(z)d\hat{U}^{\delta}(z)\right)dyds+\mathcal{E}_{t_{k}}^{2,\delta},
\end{equation}
where $\mathcal{E}_{t_{k}}^{2,\delta}(x)$ is defined by (\ref{eq:S2ErrorE2Del}).
Next, due to the fact that both $t_{k}$ and $x$ are on the grid $\delta^{2}\mathbb{N}\times\delta\mathbb{Z}$, the integral on the right hand side of 
\eqref{eq:first-step-explicit-Jdelta} can further be expressed
as
\begin{eqnarray}\label{eq:JDelMain}
 \mathcal{G}_{t_{k}}^{\delta}(x)
 &\equiv&
 \frac{1}{2}\int_{0}^{t_{k}}\int_{\mathbb{R}}\nabla_{x}^{2,\delta}\hat{p}_{t_{k}-s}^{\delta}(\lfloor x-y\rfloor)\left(\int_{x}^{\lfloor y\rfloor+\delta}\tilde{v}_{s}^{\delta}(z)d\hat{U}^{\delta}(z)\right)dyds \notag\\
  &=&\frac{\delta}{2}
  \sum_{j=1}^{k}\sum_{z\in\delta\mathbb{Z}}
  \int_{t_{j-1}}^{t_{j}}\nabla_{x}^{2,\delta}\hat{p}_{\lfloor t_{k}-t_{j}\rfloor}^{\delta}(x-z)\left(\int_{x}^{z}\tilde{v}_{s}^{\delta}(w)d\hat{U}^{\delta}(w)\right)ds.
\end{eqnarray}
Since $s\mapsto\tilde{v}_{s}^{\delta}(z)$ is linear on $[t_{j-1},t_{j}]$
by definition, it is easily seen that 
\[
\int_{t_{j-1}}^{t_{j}}\tilde{v}_{s}^{\delta}(w)ds=\frac{\delta^{2}}{2}\left(\tilde{v}_{t_{j-1}}^{\delta}(w)+\tilde{v}_{t_{j}}^{\delta}(w)\right).
\]
As a result, the term $\mathcal{G}_{t_{k}}^{\delta}$ in  (\ref{eq:JDelMain}) satisfies
\begin{equation}
 \mathcal{G}_{t_{k}}^{\delta}(x)=\frac{\delta}{2}\sum_{j=1}^{k}\sum_{z\in\delta\mathbb{Z}}\nabla_{x}^{2,\delta}\hat{p}_{t_{k}-t_{j}}^{\delta}(x-z)\int_{x}^{z}\frac{\delta^{2}}{2}\left(\tilde{v}_{t_{j-1}}^{\delta}(w)+\tilde{v}_{t_{j}}^{\delta}(w)\right)d\hat{U}^{\delta}(w).\label{eq:JDelMain2}
\end{equation}
By using the piecewise linear construction of $\tilde{v}_{t_{j-1}}^{\delta}(w)$
and $\hat{U}^{\delta}(w)$, given respectively in \eqref{eq:vtilde-linear-interpolation} and \eqref{eq:hatUdelta-linear-interpolation}, it is easily checked that 
\begin{align*}
\int_{x}^{z}\tilde{v}_{t_{j-1}}^{\delta}(w)d\hat{U}^{\delta}(w) & =\sum_{u\in\llbracket x+\delta,z\rrbracket}\int_{u-\delta}^{u}\frac{(u-w)v_{t_{j-1}}^{\delta}(u-\delta)+(w-u+\delta)v_{t_{j-1}}^{\delta}(u)}{\delta}\times\frac{\bar{U}^{\delta}(u)}{\delta}dw\\
 & =\sum_{\substack{u\in\llbracket x+\delta,z\rrbracket}
}\frac{v_{t_{j-1}}^{\delta}(u-\delta)+v_{t_{j-1}}^{\delta}(u)}{2} \,
\bar{U}^{\delta}(u)=\mathcal{I}_{t_{j-1}}^{\delta}(x,z).
\end{align*}
where the last equality is a direct consequence of \eqref{eq:disc-Ideltaxy}.
By substituting this into the expression (\ref{eq:JDelMain2}), we
obtain that
\[
\mathcal{G}_{t_{k}}^{\delta}(x)
=
\frac{\delta^{3}}{2}\sum_{j=1}^{k}\sum_{z\in\delta\mathbb{Z}}\nabla_{x}^{2,\delta}\hat{p}_{t_{k}-t_{j}}^{\delta}(x-z)\times\frac{\mathcal{I}_{t_{j-1}}^{\delta}(x,z)
+\mathcal{I}_{t_{j}}^{\delta}(x,z)}{2}.
\]
Hence going back to \eqref{eq:JDelMain} and \eqref{eq:first-step-explicit-Jdelta}, we end up with
\begin{align}
\mathcal{J}_{t_{k}}^{\delta}(x) & =v_{t_{k}}^{\delta}(x)+\frac{\delta^{3}}{2}\sum_{j=1}^{k}\sum_{z\in\delta\mathbb{Z}}\nabla_{x}^{2,\delta}\hat{p}_{t_{k}-t_{j}}^{\delta}(x-z)\times\frac{\mathcal{I}_{t_{j-1}}^{\delta}(x,z)+\mathcal{I}_{t_{j}}^{\delta}(x,z)}{2}+\mathcal{E}_{t_{k}}^{2,\delta}(x),\label{eq:Jdelta4}
\end{align}
where we note that $\tilde{v}_{t_{k}}^{\delta}(x)=v_{t_{k}}^{\delta}(x)$ since
$(t_{k},x)$ is assumed to be a grid point.

On the other hand, by using the discrete equation (\ref{eq:disc-PDE-vdelta}) for $v_{t_{k}}^{\delta}(x)$,
it is seen that 
\begin{align}
 & \frac{\delta^{3}}{2}\sum_{j=1}^{k}\sum_{z\in\delta\mathbb{Z}}\nabla_{x}^{2,\delta}\hat{p}_{t_{k}-t_{j}}^{\delta}(x-z)\times\frac{\mathcal{I}_{t_{j-1}}^{\delta}(x,z)+\mathcal{I}_{t_{j}}^{\delta}(x,z)}{2}\nonumber \\
 & =\frac{1}{2}\left(\eta_{t_{k}}^{\delta}(x)+\eta_{t_{k+1}}^{\delta}(x)\right)-\frac{1}{2}\left(v_{t_{k}}^{\delta}(x)+v_{t_{k+1}}^{\delta}(x)\right)+\mathcal{E}_{t_{k}}^{1,\delta}(x),\label{eq:Jdelta5}
\end{align}
where $\mathcal{E}_{t_{k}}^{1,\delta}(x)$ is defined by (\ref{eq:S2ErrorE1Del}).
By substituting (\ref{eq:Jdelta5}) into (\ref{eq:Jdelta4}), the
desired decomposition \eqref{eq:explicitcalJdelta} thus follows.
\end{proof}
In view of Lemma \ref{lem:JdelDecomp}, in order to compare $\mathcal{J}^{\delta}$
with $\eta,$ it is clear that there are three main ingredients to
be developed:

\noindent (i) compare $\frac{1}{2}\left(\eta_{t}^{\delta}(x)+\eta_{t+\delta^{2}}^{\delta}(x)\right)$
with $\eta_{t}(x)$;\\
(ii) show that $v_{t+\delta^ {2}}^{\delta}(x)-v_{t}^{\delta}(x)$ is small;\\
(iii) show that the error term $\mathcal{E}_{t}^{\delta}(x)$ is small.

\noindent We now implement these three steps separately.

\subsubsection{Comparing $\frac{1}{2}\left(\eta_{t}^{\delta}(x)+\eta_{t+\delta^{2}}^{\delta}(x)\right)$
with $\eta_{t}(x)$}\label{sec:comparing-eta+etadelta2-and-2eta} 
We begin our comparison procedure by looking at the difference between $\frac{1}{2}\left(\eta_{t}^{\delta}(x)+\eta_{t+\delta^{2}}^{\delta}(x)\right)$ and  $\eta_{t}(x)$. Notice that the analysis of this term is 
easy, since it only involves the input functions $f_{0}(x),g_{t}(x)$
as well as the heat kernels. The resulting estimate is summarised
in the lemma below.
\begin{lem}
Let $f_{0},g\in\mathcal{C}_{L}^{3}$ be given functions with some
$L>0$. Recall that $\eta$ and $\eta^ {\delta}$ are respectively defined by \eqref{eq:EtaTerm} and \eqref{eq:disc-eta-delta}. The distance $d_{\hat{U}^{\delta},{\bf W}}$ is introduced in \eqref{eq:ConDistVTilV}.  Then we have 
\begin{equation}
d_{\hat{{\bf U}}_1^{\delta},{\bf W}}\left(\frac{1}{2}(\eta_{\cdot}^{\delta}(\cdot)+\eta_{\cdot+\delta^{2}}^{\delta}(\cdot)),\eta\right)\leqslant C\cdot\left(\|f_{0}\|_{\mathcal{C}_{L}^{3}}+\|g\|_{\mathcal{C}_{L}^{3}}\right)\cdot\delta^{\beta}.\label{eq:EtaComp}
\end{equation}
\end{lem}

\begin{proof}
For simplicity, we only discuss the remainder estimate for the $f_{0}$-part.
Estimates of all other parts are routine repetitions of the same kind
of analysis. More specifically, let us define 
\[
\eta_{t}^{1,\delta}(x)\triangleq\delta\sum_{y\in\delta\mathbb{Z}}\nabla_{x}^{1,\delta}\hat{p}_{t}^{\delta}(x-y)f_{0}(y),\quad\text{ and }\quad
\eta_{t}^{1}(x)\triangleq\int_{\mathbb{R}}\partial_{x}p_{t}(x-y)f_{0}(y)dy
\]
respectively. We want to estimate the difference $\eta_{t}^{1,\delta}(x,x')-\eta_{t}^{1}(x,x')$
where $t\in(0,T]\cap\delta^{2}\mathbb{N}$ and $x,x'\in[-a,a]\cap\delta\mathbb{Z}$. For this purpose, we first resort to a discrete integration by parts like in Proposition \ref{prop:discrete-mild-pde-2}. Hence $\eta^ {1,\delta}$ can be recast as 
\begin{equation}
\eta_{t}^{1,\delta}(x)=\delta\sum_{y\in\delta\mathbb{Z}}\nabla_{x}^{1,\delta}f_{0}(x-y)\hat{p}_{t}^{\delta}(y)=\int_{\mathbb{R}}\nabla^{1,\delta}f_{0}(\lfloor x-y\rfloor)\hat{p}_{t}^{\delta}(\lfloor y\rfloor)dy.\label{eq:Eta1Dis}
\end{equation}
In the same way, due to our regularity assumption on $f$, a simple integration by parts yields
\begin{equation}
\eta_{t}^{1}(x)=\int_{\mathbb{R}}\partial_{x}f_{0}(x-y)p_{t}(y)dy.\label{eq:Eta1Cts}
\end{equation}
We now express the increments $\eta^{1,\delta}_{t}(x,x')$, Thanks to 
\eqref{eq:Eta1Dis} and a change of variable $y=\sqrt{t}w$, we get
\begin{multline}\label{eq:increments-eta-delta-first}
\eta_{t}^{1,\delta}(x,x')  =\int_{\mathbb{R}}\hat{p}_{t}^{\delta}(\lfloor y\rfloor)\left(\nabla^{1,\delta}f_{0}(\lfloor x'-y\rfloor)-\nabla^{1,\delta}f_{0}(\lfloor x-y\rfloor)\right)dy\\
  =\int_{\mathbb{R}}\sqrt{t}\hat{p}_{t}^{\delta}(\lfloor\sqrt{t}w\rfloor)\left(\nabla^{1,\delta}f_{0}(\lfloor x'-\sqrt{t}w\rfloor)-\nabla^{1,\delta}f_{0}(\lfloor x-\sqrt{t}w\rfloor)\right)dy
\end{multline}
Similarly, by using (\ref{eq:Eta1Cts}) we also have 
\begin{multline}\label{eq:increments-eta-delta-second}
\eta_{t}^{1}(x,x')  =\int_{\mathbb{R}}p_{t}(y)\left(\partial_{x}f_{0}(x'-y)-\partial_{x}f_{0}(x-y)\right)dy\\
 =\int_{\mathbb{R}}p_{t}(y)\left(\int_{x}^{x'}\partial_{xx}^{2}f_{0}(z-y)dz\right)dy=\int_{\mathbb{R}}p_{1}(w)\left(\int_{x}^{x'}\partial_{xx}^{2}f_{0}(z-\sqrt{t}w)dz\right)dw.
\end{multline}
Putting together \eqref{eq:increments-eta-delta-first} and \eqref{eq:increments-eta-delta-second}, we thus obtain the identity
\[
\eta_{t}^{1,\delta}(x,x')-\eta_{t}^{1}(x,x')=A_{t}(x,x')+B_{t}(x,x'),
\]
where $A_{t}(x,x')$ and $B_{t}(x,x')$ are respectively defined by 
\begin{align}
A_{t}(x,x')&\triangleq\int_{\mathbb{R}}\left(\sqrt{t}\hat{p}_{t}^{\delta}(\lfloor\sqrt{t}w\rfloor)-p_{1}(w)\right)\left(\int_{x}^{x'}\partial_{xx}^{2}f_{0}(z-\sqrt{t}w)dz\right)dw,\label{eq:Atxxprime}
\end{align}
\begin{align}
B_{t}(x,x') &\triangleq\int_{\mathbb{R}}\sqrt{t}\hat{p}_{t}^{\delta}(\lfloor\sqrt{t}w\rfloor)\,M_{x,x'}(\sqrt{t}w)dw, \label{eq:Btxxprime}
\end{align}
where for $a\in\mathbb{R}$ we define 
\[
M_{x,x'}(a)=\nabla^{1,\delta}f_{0}(\lfloor x'-a\rfloor)-\nabla^{1,\delta}f_{0}(\lfloor x-a\rfloor)
 -\int_{x}^{x'}\partial_{xx}^{2}f_{0}(z-a)dz.
\]
We will now estimate the terms $A_{t}(x,x')$ and $B_{t}(x,x')$.
Let us start by estimating the term $A_{t}(x,x')$ introduced in \eqref{eq:Atxxprime}. In order to bound this term, we proceed as in Section \ref{subsec:S3Unif}. Namely we write $\hat{p}^{\delta}$ in terms of the discrete kernel 
$p^{d}$ and we invoke Theorem \ref{thm:LCLT}. We let the patient reader check that  
\[
\left|\sqrt{t}\hat{p}_{t}^{\delta}(\lfloor\sqrt{t}w\rfloor)-p_{1}(w)\right|\leqslant C_{1}\cdot\frac{\delta}{\sqrt{t}}e^{-C_{2}w^{2}},
\]
for all $t\geqslant\delta^{2}$ and $w\in\mathbb{R}.$ Plugging this inequality into \eqref{eq:Atxxprime},  it follows
that
\[
\big|A_{t}(x,x')\big| \leqslant C_{1}\frac{\delta}{\sqrt{t}}\int_{\mathbb{R}}e^{-C_{2}w^{2}}\left(\int_{x}^{x'}\big|\partial_{xx}^{2}f_{0}(z-\sqrt{t}w)\big|dz\right)dw.
\]
In addition, the right hand side above can be bounded invoking the fact that $f\in{\mathcal C}^{3}_{L}$ and Definition \ref{def:spaceCrL} we obtain 
\[
\big|A_{t}(x,x')\big| \leqslant C_{3}\frac{\delta}{\sqrt{t}}\cdot|x-x'|\cdot a^{L}\|f_{0}\|_{\mathcal{C}_{L}^{3}}.
\]
Whenever $t\geqslant\delta^{2}$ and $x\in[-a,a]$ we thus get 
\begin{align}
 \big|A_{t}(x,x')\big|& \leqslant C_{3}t^{-\beta/2}\delta^{\beta}|x'-x|^{2\beta}\cdot a^{L+1-2\beta}\|f_{0}\|_{\mathcal{C}_{L}^{3}}\nonumber \\
 & \leqslant C_{4}\|f_{0}\|_{\mathcal{C}_{L}^{3}}E(a,t)a^{\chi}t^{-\beta/2}|x'-x|^{2\beta}\cdot\delta^{\beta},\label{eq:EtaEst3}
\end{align}
where we have absorbed some power of $a$ into $E(a,t)$ term for the last step. 

Next, we estimate $B_{t}(x,x')$ in \eqref{eq:Btxxprime}. By the definition of $\nabla_{x}^{1,\delta}$,
it is easily seen that
\begin{align*}
 & \nabla^{1,\delta}f_{0}(\lfloor x'-\sqrt{t}w\rfloor)-\nabla^{1,\delta}f_{0}(\lfloor x-\sqrt{t}w\rfloor)\\
 & =\frac{f_{0}(x'+\delta-\lfloor\sqrt{t}w\rfloor)-f_{0}(x+\delta-\lfloor\sqrt{t}w\rfloor)}{\delta}-\frac{f_{0}(x'-\lfloor\sqrt{t}w\rfloor)-f_{0}(x-\lfloor\sqrt{t}w\rfloor)}{\delta}\\
 & =\int_{x}^{x'}\frac{\partial_{x}f_{0}(z+\delta-\lfloor\sqrt{t}w\rfloor)-\partial_{x}f_{0}(z-\lfloor\sqrt{t}w\rfloor)}{\delta}dz\\
 & =\int_{x}^{x'}\left(\int_{0}^{1}\partial_{xx}^{2}f_{0}(z+\theta\delta-\lfloor\sqrt{t}w\rfloor)d\theta\right)dz.
\end{align*}
Using again Definition \ref{def:spaceCrL} and our assumption $f_{0}\in{\mathcal C}^{3}_{L}$, it follows that
\begin{align}
\Big|M_{x,x'}(\sqrt{t}w)\Big| & =\Big|\nabla^{1,\delta}f_{0}(\lfloor x'-\sqrt{t}w\rfloor)-\nabla^{1,\delta}f_{0}(\lfloor x-\sqrt{t}w\rfloor)-\int_{x}^{x'}\partial_{xx}^{2}f_{0}(z-\sqrt{t}w)dz\Big|\nonumber \\
 & =\left|\int_{x}^{x'}\int_{0}^{1}\left(\partial_{xx}^{2}f_{0}(z+\theta\delta-\lfloor\sqrt{t}w\rfloor)-\partial_{xx}^{2}f_{0}(z-\sqrt{t}w)\right)d\theta dz\right|\nonumber \\
 & \leqslant C\|f_{0}\|_{\mathcal{C}_{L}^{3}}a^{L}\cdot|x'-x|\cdot\delta\cdot(1+|w|)^{L}.\label{eq:EtaEst1}
\end{align}
In addition, recall from (\ref{eq:non-gradient-estimate}) that 
\begin{equation}
\big|\sqrt{t}\hat{p}_{t}^{\delta}(\lfloor\sqrt{t}w\rfloor)\big|\leqslant Ce^{-C_{2}w^{2}}.\label{eq:EtaEst2}
\end{equation}
By substituting (\ref{eq:EtaEst1}) and (\ref{eq:EtaEst2}) into the
expression \eqref{eq:Btxxprime}  of $B_{t}(x,x')$, we can bound it as 
\[
\big|B_{t}(x,x') \big|
\leqslant C\|f_{0}\|_{\mathcal{C}_{L}^{3}}a^{L}|x'-x|\cdot\delta.
\]
Now one can perform the same kind of manipulation as for \eqref{eq:EtaEst3} in order to get 
\begin{align}
\big|B_{t}(x,x')\big| & \leqslant C\|f_{0}\|_{\mathcal{C}_{L}^{3}}a^{L+1-2\beta}|x'-x|^{2\beta}\cdot\delta\nonumber \\
 & \leqslant C\|f_{0}\|_{\mathcal{C}_{L}^{3}}E(a,t)a^{\chi+\beta/2}|x'-x|^{2\beta}\cdot\delta.\label{eq:EtaEst4}
\end{align}

By putting (\ref{eq:EtaEst3}) and (\ref{eq:EtaEst4}) together, we
obtain the remainder estimate for $\eta^{1,\delta}-\eta^{1}$, which
contributes to one term in the rough path distance. All other terms
(uniform distance, spatial and time variations) and also the $g(t,x)$-part
are dealt with in a similar way. The resulting estimate is given by
(\ref{eq:EtaComp}).
\end{proof}

\subsubsection{Estimating $v_{t+\delta^{2}}^{\delta}(x)-v_{t}^{\delta}(x)$}\label{sec:estim-derivative-v}

This part requires some care since we do not view $t\mapsto v_{t}^{\delta}(\cdot)$
as a H\"older-continuous path taking values in spatial rough paths.
As opposed to that, we shall consider the path $G_{t}^{\delta}=v_{t+\delta^{2}}^{\delta}(x)-v_{t}^{\delta}(x)$ as an element of the (discrete) space $\mathcal {B}^{\theta,\lambda}$ which is introduced in Definition \ref{def:controlled-process}, and we will prove that the norm of $G^{\delta}$ is small. Our estimate is summarized in the following lemma. 
\begin{lem}\label{lem:estimate-norm-Gdelta}
In the decomposition \eqref{eq:explicitcalJdelta}, let $G^{\delta}$ be the path given by 
\begin{equation}\label{eq:Gdelta-differenceV}
G^{\delta}_{t}(x)\equiv v_{t+\delta^{2}}^{\delta}(x)-v_{t}^{\delta}(x),
\end{equation}
defined for $(t,x)\in\delta^{2}\mathbb{N}\times\delta\mathbb{Z}$. We consider $G$ as a discrete path controlled by $\hat{{\bf U}}_1^{\delta}$, with a null derivative. Pick $\alpha, \beta,\chi$ satisfying \eqref{f1}, as well as $\theta' <\theta$ and $\beta'\in(\beta, \alpha)$, so that for all $a\geqslant 1$ and 
$t\in [0,T]$ we have 
\begin{equation}\label{eq:Ethetaprime-Etheta}
E^{\theta'}(a,t)a^{\beta'/2}\leqslant E^{\theta}(a,t)a^{\beta/2},
\end{equation} 
where we recall that $E^{\theta}(a,t)=e^{\lambda t+\theta a+\theta at}$. 
As usual we choose $\gamma=(\alpha-\beta)/4$ and $\lambda$ is such that 
\eqref{eq:choice-lambda-kappa} is fulfilled. The controlled norms $\Theta$ are introduced in \eqref{eq:def_normTheta}. Then we have 
\begin{equation}\label{eq:normTheta-delta-betaprime-beta}
\Theta^{\theta',\lambda}(G^{\delta})\leqslant C\Lambda^{\theta',\beta'}(1+\kappa(\hat{{\bf U}}_1^\delta))\lambda^{\frac{\alpha-\beta}{4}}\cdot\delta^{\frac{\beta'(\beta'-\beta)}{\beta'+\beta}}.
\end{equation} 
\end{lem}

\begin{proof}
As usual, based on the definition of our rough path metric, we need
to develop four types of estimates. We again use the trick of adjusting
exponents based on the uniform boundedness of $v_{t}^{\delta}(x)$
proved in Proposition \ref{prop:discrete-FPT}.

\noindent\emph{Step 1: uniform estimate.} Let $\theta'<\theta$ and $\beta'\in(\beta,\alpha)$
be fixed parameters such that \eqref{eq:Ethetaprime-Etheta} is verified. In particular for all $a\geqslant 1$ and $t\in[0,T]$ we have 
\[
E^{\theta'}(a,t)a^{\beta'/2}\leqslant E^{\theta}(a,t)\ \ \ \forall a\geqslant1,t\in[0,T].
\]
It follows from Proposition \ref{prop:discrete-FPT} that 
\begin{equation}\label{eq:first-estim-Gdelta}
\big|G_{t}^{\delta}(x)\big|\leqslant\Lambda^{\theta',\beta'}E^{\theta'}(a,t)a^{\beta'/2}\delta^{\beta'}\leqslant\Lambda^{\theta',\beta'}E^{\theta}(a,t)\delta^{\beta'},
\end{equation}
where $\Lambda^{\theta',\beta'}$ denotes a uniform upper bound of
$\Theta^{\theta',\lambda}(v^{\delta})$. This gives the uniform estimate
in the definition of the rough path norm of $(t,x)\mapsto G_{t}^{\delta}(x).$

\noindent\emph{Step 2: time variation estimate}. We simply decompose the time variation as follows for $s,t\in\delta^{2}\mathbb{N}$ and $x\in\delta\mathbb{Z}$:
\[\big|G_{t}^{\delta}(x)-G_{s}^{\delta}(x)\big|
\leqslant \big|v_{t+\delta^{2}}^{\delta}(x)-v_{t}^{\delta}(x)\big|+\big|v_{s+\delta^{2}}^{\delta}(x)-v_{s}^{\delta}(x)\big|.\]
Next we recall that  $(\theta',\beta')$ are chosen as in \eqref{eq:Ethetaprime-Etheta}. Invoking the $\delta$-uniform bound on $v^{\delta}$ obtained in Proposition~\ref{prop:discrete-FPT}, we get 
\[\big|G_{t}^{\delta}(x)-G_{s}^{\delta}(x)\big|
\leqslant 2\Lambda^{\theta',\beta'}E^{\theta}(a,t)a^{\beta/2}\delta^{\beta'}=2\Lambda^{\theta',\beta'}E^{\theta}(a,t)a^{\beta/2}\delta^{\beta}\delta^{\beta'-\beta}.
\]
Since $s,t$ on the grid $\delta^{2}\mathbb{N},$ satisfy  $t-s\geqslant\delta^{2}$, we discover that 
\begin{equation}\label{eq:second-estim-Gdelta}
\big|G_{t}^{\delta}(x)-G_{s}^{\delta}(x)\big|
\leqslant
2\Lambda^{\theta',\beta'}E^{\theta}(a,t)a^{\beta/2}|t-s|^{\beta/2}\cdot\delta^{\beta'-\beta}.
\end{equation}

\noindent\emph{Step 3: the space variation estimate}. This term is handled 
exactly like the time increment in Step 2. We let the reader check that for all 
$x,x'\in [-a,a]$ and $t\in\delta^{2}\mathbb{N}$ we have 
\begin{equation}\label{eq:third-estim-Gdelta}
\big|G_{t}^{\delta}(x,x')\big|\leqslant2\Lambda^{\theta',\beta'}E^{\theta}(a,t)a^{\beta/2}|x'-x|^{\beta}\cdot\delta^{\beta'-\beta}.
\end{equation}

\noindent\emph{Step 4: remainder term and conclusion.} Since our term ansatz stipulates that the derivative $\partial_{\hat{U}^{\delta}}$ vanishes,
we are now left with an estimate of the remainder for $G^{\delta}$. This is 
achieved by elaborating on the upper bound \eqref{eq:third-estim-Gdelta} in order to get a second order estimate in $x'-x$. This step is detailed in 
Lemma \ref{lem:second-order-xxprime} below for the sake of clarity. Gathering this lemma and \eqref{eq:first-estim-Gdelta}-\eqref{eq:second-estim-Gdelta}-\eqref{eq:third-estim-Gdelta}, this finishes the proof of our claim~\eqref{eq:normTheta-delta-betaprime-beta}.
\end{proof}

We now state and prove the announced estimate of the remainder for the controlled process $G^{\delta}$ defined above. 
\begin{lem}\label{lem:second-order-xxprime}
Under the same conditions as for Lemma \ref{lem:lambda-betabetaprime}, let 
$t\in\llparenthesis 0,T\rrbracket$ and $x,x'\in\llbracket -a,a\rrbracket$. Also recall that the path $G^{\delta}$ is defined by \eqref{eq:Gdelta-differenceV}. Then we have 
\begin{equation}
\big|G_{t}^{\delta}(x,x')\big|\leqslant C\Lambda^{\theta',\beta'}(1+\kappa(\hat{{\bf U}}_1^\delta))E^{\theta}(a,t)\lambda^{\frac{\alpha-\beta}{4}}a^{\chi}(a^{\beta/2}+t^{-\beta/2})|x'-x|^{2\beta}\cdot\delta^{\frac{\beta'(\beta'-\beta)}{\beta'+\beta}}.\label{eq:RemEstDetVPart}
\end{equation}
\end{lem}

\begin{proof}
We divide the discussion into three (not necessarily disjoint) cases.
Note that such type of consideration already appears in the proof
of Lemma \ref{lem:RemEstComp}. In the sequel $\tau$ is a parameter in $[0,1]$ which will be chosen appropriately. 

\noindent\emph{Case I: $|x'-x|\geqslant\delta^{\tau}$.} 
In this case we decompose 
$G_{t}^{\delta}(x,x')$ as 
\begin{equation}\label{eq:trianG}
\big|G_{t}^{\delta}(x,x')\big|\leqslant\big|G_{t}^{\delta}(x)\big|+
\big|G_{t}^{\delta}(x')\big|.
\end{equation}
Then for $x,x'\in\llbracket -a,a\rrbracket$ one can apply the uniform estimate 
\eqref{eq:first-estim-Gdelta} and obtain 
\begin{equation}\label{eq:forth-estim-Gdelta}
\big|G_{t}^{\delta}(x,x')\big|\leqslant C\Lambda^{\theta',\beta'}E^{\theta}(a,t)a^{\beta/2}\delta^{\beta'}.
\end{equation}
One can then conclude by invoking the fact that $\delta^{\tau}\leqslant
|x'-x|$. That is we get 
\begin{align}\label{eq:fifth-estim-Gdelta}
 \big|G_{t}^{\delta}(x,x')\big|& \leqslant 
C\Lambda^{\theta',\beta'}E^{\theta}(a,t)a^{\beta/2}\delta^{2\tau\beta}\cdot\delta^{\beta'-2\tau\beta}\nonumber \\
 & \leqslant C\Lambda^{\theta',\beta'}E^{\theta}(a,t)a^{\beta/2}|x'-x|^{2\beta}\cdot\delta^{\beta'-2\tau\beta}.
\end{align}

\noindent\emph{Case II: $|x'-x|\geqslant\sqrt{t}$.} We choose $\beta',\theta'$ 
such that \eqref{eq:Ethetaprime-Etheta} is verified. Then we start our estimation procedure like in \eqref{eq:trianG}, which yields inequality 
\eqref{eq:forth-estim-Gdelta}. Now one can decompose the right hand side of 
\eqref{eq:forth-estim-Gdelta} in order to get 
\[
\big|G_{t}^{\delta}(x,x')\big|\leqslant 
C\Lambda^{\theta',\beta'}E^{\theta}(a,t)a^{\chi}t^{-\beta/2}|x'-x|^{2\beta}\cdot\delta^{\beta'}|x'-x|^{-2\beta}t^{\beta/2}.
\]
Next we resort to the fact that $|x'-x|\geqslant\sqrt{t}$ to write 
\begin{align*}
\big|G_{t}^{\delta}(x,x')\big|\leqslant 
& C\Lambda^{\theta',\beta'}E^{\theta}(a,t)a^{\chi}t^{-\beta/2}|x'-x|^{2\beta}\cdot\delta^{\beta'}t^{-\beta/2}\\
= & C\Lambda^{\theta',\beta'}E^{\theta}(a,t)a^{\chi}t^{-\beta/2}|x'-x|^{2\beta}\cdot\delta^{\beta'-\beta}\cdot\delta^{\beta}t^{-\beta/2}.
\end{align*}
Eventually, thanks to the relation $t\geq\delta^{2}$ we obtain a relation which is similar to \eqref{eq:fifth-estim-Gdelta}:
\begin{equation}\label{eq:sixth-estim-Gdelta}
 \big|G_{t}^{\delta}(x,x')\big|\leqslant C\Lambda^{\theta',\beta'}E^{\theta}(a,t)a^{\chi}t^{-\beta/2}|x'-x|^{2\beta}\cdot\delta^{\beta'-\beta}.
\end{equation}

\noindent\emph{Case III: $|x'-x|\leqslant\sqrt{t}\wedge\delta^{\tau}$.} Recall that in \eqref{eq:rewrite-derivativeUW} we have written a decomposition for the process $\tilde{v}_{t}^{\delta}$. Specialized on the grid $\delta^{2}\mathbb{N}\times\delta\mathbb{Z}$ this can be read as 
\begin{equation}\label{eq:diff-vtdelta2-vt}
v_{t+\delta^{2}}^{\delta}(x,x')-v_{t}^{\delta}(x,x')=-2\left(v_{t+\delta^{2}}^{\delta}(x)-v_{t}^{\delta}(x)\right)\hat{U}_1^{\delta}(x,x')+\mathcal{R}^{v_{t+\delta^{2}}^{\delta}}(x,x')-\mathcal{R}^{v_{t}^{\delta}}(x,x'),
\end{equation}
where we have written $\mathcal{R}^{v_{t}^{\delta}}$ instead 
$\mathcal{R}_{\hat{{\bf U}}_1^{\delta}}^{v_{t}^{\delta}}$ in order to alleviate our notation. We now bound the two terms in the right hand side of \eqref{eq:diff-vtdelta2-vt}. 

In order to bound the derivative term in the right hand side of \eqref{eq:diff-vtdelta2-vt}, consider again $\theta'<\theta$ and $\beta'\in(\beta,\alpha)$. Then invoking Proposition \ref{prop:discrete-FPT} and the fact that $\kappa(\hat{{\bf U}}_1^{\delta})$ is 
uniformly bounded in $\delta$, we get
\begin{equation}\label{eq:derivative-term-first-step}
\big|\left(v_{t+\delta^{2}}^{\delta}(x)-v_{t}^{\delta}(x)\right)\hat{U}_1^{\delta}(x,x')\big|
\leqslant C\Lambda^{\theta',\beta'}\kappa(\hat{{\bf U}}_1^{\delta})E^{\theta'}(a,t)a^{\chi+\beta'/2}\delta^{\beta'}|x'-x|^{\alpha}.
\end{equation}
Recall that we are now dealing with the case $|x-x'|\le\delta^{\tau}$ and $x,x'\in\llbracket -a,a\rrbracket$.
However, we also have $|x'-x|\ge\delta$, due to the fact that both $x$ and $x'$ belong to the grid $\delta\Z$.
Therefore one is enabled to write
\begin{equation*}
\delta^{\beta'}|x'-x|^{\alpha}
=
\delta^{\beta'-\beta}|x'-x|^{2\beta} \lp \frac{\delta}{|x'-x|}  \rp^{\beta} |x'-x|^{\al-\beta}
\le
\delta^{\beta'-\beta}|x'-x|^{2\beta} a^{\al-\beta} \, .
\end{equation*}
Reporting this inequality in the right hand side of \eqref{eq:derivative-term-first-step} yields
\begin{equation}\label{eq:derivative-term-second-step}
  \big|\left(v_{t+\delta^{2}}^{\delta}(x)-v_{t}^{\delta}(x)\right)\hat{U}_1^{\delta}(x,x')\big|
  \leqslant C\Lambda^{\theta',\beta'}\kappa(\hat{{\bf U}}_1^{\delta})E^{\theta'}(a,t)a^{\chi+\beta'/2+\alpha-\beta}|x'-x|^{2\beta}\delta^{\beta'-\beta}.
\end{equation}
In addition, we have assumed that \eqref{eq:Ethetaprime-Etheta} is verified. Hence we obtain 
\[
E^{\theta'}(a,t)a^{\chi+\beta'/2+\alpha-\beta}\leqslant CE^{\theta}(a,t)a^{\chi+\beta/2}.
\]
Plugging this information into \eqref{eq:derivative-term-second-step}, it follows that
\begin{equation}\label{eq:derivative-term-third-step}
\big|\left(v_{t+\delta^{2}}^{\delta}(x)-v_{t}^{\delta}(x)\right)\hat{U}_1^{\delta}(x,x')\big|\leqslant C\Lambda^{\theta',\beta'}\kappa(\hat{{\bf U}}_1^{\delta})E^{\theta}(a,t)a^{\chi+\beta/2}|x'-x|^{2\beta}\cdot\delta^{\beta'-\beta}.
\end{equation}

   We now turn to an estimate of the remainder terms in \eqref{eq:diff-vtdelta2-vt}. We shall bound each term in the difference individually, 
namely write 
\[
\big|\mathcal{R}^{v_{t+\delta^{2}}^{\delta}}(x,x')-\mathcal{R}^{v_{t}^{\delta}}(x,x')\big| \leqslant\big|\mathcal{R}^{v_{t+\delta^{2}}^{\delta}}(x,x')\big|+\big|\mathcal{R}^{v_{t}^{\delta}}(x,x')\big| \, .
\]
Then we resort to Proposition \ref{prop:discrete-FPT} and the definition of $\Theta^{\theta',\lambda}$ in \eqref{eq:def_normTheta} in order to obtain 
\begin{equation}\label{eq:diff-remainder-vdelta2-v-first}
\big|\mathcal{R}^{v_{t+\delta^{2}}^{\delta}}(x,x')-\mathcal{R}^{v_{t}^{\delta}}(x,x')\big| \leqslant C\Lambda^{\theta',\beta'}E^{\theta'}(a,t)\lambda^{\gamma'}a^{\chi}(a^{\beta'/2}+t^{-\beta'/2})|x'-x|^{2\beta'}.
\end{equation}
where $\gamma'\triangleq\frac{\alpha-\beta'}{4}.$ Since $|x'-x|\leqslant\delta^{\tau},$
we have 
\[
|x'-x|^{2\beta'}=|x'-x|^{2\beta}\cdot|x'-x|^{2\beta'-2\beta}\leqslant|x'-x|^{2\beta}\cdot\delta^{2\tau(\beta'-\beta)}.
\]
In addition, since $|x'-x|\leqslant\sqrt{t},$ we also have
\begin{align*}
t^{-\beta'/2}|x'-x|^{2\beta'} & =t^{-\beta/2}|x'-x|^{2\beta}\cdot t^{-(\beta'-\beta)/2}|x'-x|^{2(\beta'-\beta)}\\
 & \leqslant t^{-\beta/2}|x'-x|^{2\beta}\cdot|x'-x|^{\beta'-\beta}\leqslant t^{-\beta/2}|x'-x|^{2\beta}\cdot\delta^{\tau(\beta'-\beta)}.
\end{align*}
Reporting those values into \eqref{eq:diff-remainder-vdelta2-v-first} and recalling that \eqref{eq:Ethetaprime-Etheta} holds true we end up with
\begin{equation}\label{eq:diff-remainder-vdelta2-v-second}
\big|\mathcal{R}^{v_{t+\delta^{2}}^{\delta}}(x,x')-\mathcal{R}^{v_{t}^{\delta}}(x,x')\big|\leqslant C\Lambda^{\theta',\beta'}\lambda^{\gamma'}E^{\theta}(a,t)a^{\chi}(a^{\beta/2}+t^{-\beta/2})|x'-x|^{2\beta}\cdot\delta^{\tau(\beta'-\beta)}.
\end{equation}
Hence gathering \eqref{eq:derivative-term-third-step} and \eqref{eq:diff-remainder-vdelta2-v-second} into \eqref{eq:diff-vtdelta2-vt}, we have obtained
\begin{equation}\label{eq:diff-vtdelta2-vt-step}
\big|v_{t+\delta^{2}}^{\delta}(x,x')-v_{t}^{\delta}(x,x')\big|
\leqslant C\Lambda^{\theta',\beta'}\kappa(\hat{{\bf U}}_1^{\delta})\lambda^{\gamma'}E^{\theta}(a,t)
a^{\chi}(a^{\beta/2}+t^{-\beta/2})|x'-x|^{2\beta}\cdot\delta^{\tau(\beta'-\beta)}.
\end{equation}
%\hre{Here maybe I missed something since $\kappa_{\alpha,\sigma}(\hat{U}^{\delta})$ does not appear in relation (5.196) in the statement of Lemma 5.33.}

\noindent
Combining cases I to III, it is not hard to see that the optimal
choice of $\tau$ is such that 
\[
\beta'-2\tau\beta=\tau(\beta'-\beta)\iff\tau=\frac{\beta'}{\beta'+\beta}.
\]
The resulting $\delta$-factors is then found to be $\delta^{\frac{\beta'(\beta'-\beta)}{\beta'+\beta}}$.
This completes the proof of the lemma.
\end{proof}

\subsubsection{Estimating $\mathcal{E}_{t}^{\delta}(x)$}\label{sec:estim-calEdelat}
The last term we have to estimate in order to quantify the convergence of 
$\mathcal{J}^{\delta}$ in \eqref{eq:S2Est} is $\mathcal{E}_{t_{k}}^{\delta}$. We label another lemma in that direction 
\begin{lem}\label{lem:estimate-norm-calEdelta}
In the decomposition \eqref{eq:explicitcalJdelta}, consider the term 
\begin{equation}\label{eq:calEdelta=calE1delta+calE2delta}
\mathcal{E}_{t}^{\delta}(x)=\mathcal{E}_{t}^{1,\delta}(x)+\mathcal{E}_{t}^{2,\delta}(x),
\end{equation}
which is defined for $(t,x)\in\delta^{2}\mathbb{N}\times\delta\mathbb{Z}$ and where $\mathcal{E}_{t}^{1,\delta},\mathcal{E}_{t}^{2,\delta}$ are introduced in  \eqref{eq:S2ErrorE1Del}- \eqref{eq:S2ErrorE2Del}. We pick $\alpha, \beta,\chi$ satisfying \eqref{f1} as well as $\theta' <\theta$ such that 
\begin{equation}\label{eq:The'Cond}
E^{\theta'}(a,t)a^{2\chi+\beta/2}\leqslant E^{\theta}(a,t)\ \ \ \forall a\geqslant1,t\in [0,T].
\end{equation} Then we have 
\begin{equation}\label{eq:norm-calEdelta}
\Theta^{\theta',\lambda}(\mathcal{E}^{\delta})
\leqslant C(1+\kappa(\hat{U}_{1}^{\delta}))\Lambda^{\theta',\beta}\lambda^{\frac{\alpha-\beta}{4}}\delta^{\alpha-2\chi}.
\end{equation}
\end{lem}
\begin{proof}
Unlike for Lemma \ref{lem:estimate-norm-Gdelta}, we will not give details for all four terms defining the norm $\Theta$. For the sake of 
conciseness, we just focus here on estimating 
\begin{equation}\label{eq:uniform-estimate-calEdelta}
\mathcal{S}(\mathcal{E}^{\delta})
\equiv\sup\big\{\mathcal{E}_{r}^{\delta};\,r\in\llbracket 0,t\rrbracket, x\in\llbracket -a,a\rrbracket\big\},
\end{equation}
which is a part of $\llbracket \mathcal{E}^{\delta}\rrbracket^{\llbracket 0,t\rrbracket\times\llbracket -a,a\rrbracket}$ in \eqref{eq:def_normTheta}. Details are left 
to the reader for the other terms. 

Let us first give an estimate for the term $\mathcal{E}_{t}^{1,\delta}$ in 
\eqref{eq:S2ErrorE1Del}. To this aim, we write the discrete integral in the right hand side of \eqref{eq:S2ErrorE1Del} as continuous one:
\[
\mathcal{E}_{t}^{1,\delta}(x)=-\frac{\delta^{2}}{4}\int_{\mathbb{R}}\nabla_{x}^{2,\delta}\hat{p}^\delta_{t}(\lfloor x-y\rfloor)\int_{x}^{\lfloor y\rfloor}\tilde{v}_{0}^{\delta}(w)d\hat{U}^{\delta}(w)dy.
\]
Thanks to this expression and the $\hat{U}^\delta$-decomposition \eqref{eq:BarU1Decomp}, one can perform the same kind of analysis as in Section \ref{subsec:S3Unif} (for the $\hat{U}_1^\delta$ part) and Proposition \ref{prop:discrete-FPT} (for the $\hat{U}_2^\delta$). This yields 
\begin{equation}\label{eq:estimate-calE1delta}
\big|\mathcal{E}_{t}^{1,\delta}(x)\big|\leqslant C(\kappa(\hat{{\bf U}}_1^{\delta})+\sqrt{\delta})E(a,t)\Theta(\tilde{v}^{\delta})\cdot\delta^{2}.
\end{equation}
Next, we handle the term $\mathcal{E}_{t}^{2,\delta}$ in \eqref{eq:S2ErrorE2Del}. Again the $\hat{U}_2^\delta$-part is (trivially) controlled by $$C\sqrt{\delta}E(a,t)\Theta(\tilde{v}^\delta).$$ We therefore focus on the $\hat{U}_1^\delta$-part and further decompose it as
%To this aim we decompose this term again into a small and large time integral. Namely we write 
$\mathcal{K}_{t}^{1,\delta}(x)+\mathcal{K}_{t}^{2,\delta}(x)$, where 
\begin{eqnarray}
\mathcal{K}_{t}^{1,\delta}(x)
&\triangleq&
-\frac{1}{2}\int_{0}^{t\wedge\delta^{2}}\int_{\mathbb{R}}\nabla_{x}^{2,\delta}\hat{p}_{0}^{\delta}(\lfloor y\rfloor)\left(\int_{x-y}^{\lfloor x-y\rfloor+\delta}\tilde{v}_{t-s}^{\delta}(z)d\hat{U}_1^{\delta}(z)\right)dyds
\label{eq:calK1delta}\\
\mathcal{K}_{t}^{2,\delta}(x)
&\triangleq&
-\frac{1}{2}\int_{t\wedge\delta^{2}}^{t}\int_{\mathbb{R}}\nabla_{x}^{2,\delta}\hat{p}_{\lfloor s\rfloor}^{\delta}(\lfloor x-y\rfloor)\left(\int_{y}^{\lfloor y\rfloor+\delta}\tilde{v}_{t-s}^{\delta}(z)d\hat{U}_1^{\delta}(z)\right)dyds.\label{eq:calK2delta}
\end{eqnarray}
We now upper bound the two terms $\mathcal{K}_{t}^{1,\delta}(x)$ and $\mathcal{K}_{t}^{2,\delta}(x)$. 
The term $\mathcal{K}_{t}^{1,\delta}(x)$ is handled as follows: according to \eqref{eq:DiscUnifGau=0} or \eqref{eq:app-estim-deriv--hat-heat-kernel} we have 
\begin{equation}
\big|\nabla_{x}^{2,\delta}\hat{p}_{0}^{\delta}(\lfloor y\rfloor)\big|\lesssim\frac{1}{\delta^{3}}\mathds{1}_{\{y:|y|\leqslant\delta\}}.\label{eq:E1Del1}
\end{equation}
In addition the rough integral estimate \eqref{eq:RI1}, properly extended to a discrete setting as in the proof of Lemma \ref{lem:UnifWDelvsW}, entails that 
\begin{multline}\label{eq:E1Del2}
\left|\int_{x-y}^{\lfloor x-y\rfloor+\delta}\tilde{v}_{t-s}^{\delta}(z)d\hat{U}_1^{\delta}(z)\right|
\leqslant 
C\kappa(\hat{{\bf U}}_1^{\delta})\Theta(\tilde{v}^{\delta})\lambda^{\frac{\alpha-\beta}{4}}E(a,t)e^{-(\lambda+\theta(a+|y|))s}
\cdot\left((a+|y|)^{\chi}\delta^{\alpha}\right.
\\
\left.+(a+|y|)^{2\chi}\delta^{2\alpha}
+(a+|y|)^{2\chi+\beta/2}(\delta^{2\alpha+\beta}+\delta^{\alpha+2\beta})
+(a+|y|)^{2\chi}(t-s)^{-\beta/2}\delta^{\alpha+2\beta}\right).
\end{multline}
By substituting (\ref{eq:E1Del1}) and (\ref{eq:E1Del2}) into the
expression \eqref{eq:calK1delta} of $\mathcal{K}_{t}^{1,\delta}(x),$ we obtain that 
\begin{equation}\label{eq:estimate-calK1delta}
\big|\mathcal{K}_{t}^{1,\delta}(x)\big|\leqslant C\kappa(\hat{{\bf U}}_1^{\delta})\Theta(\tilde{v}^{\delta})\lambda^{\frac{\alpha-\beta}{4}}E(a,t)\cdot\delta^{\alpha-2\chi}.
\end{equation}
To estimate $\mathcal{K}_{t}^{2,\delta}(x)$, we assume that $t>\delta^{2}$
(for otherwise $\mathcal{K}_{t}^{2,\delta}(x)=0$). We first apply
the change of variable $x-y=\sqrt{s}w$ to write it as 
\begin{equation}\label{eq:changevar-calK2delta}
\mathcal{K}_{t}^{2,\delta}(x)
\triangleq
-\frac{1}{2}\int_{\delta^{2}}^{t}\int_{\mathbb{R}}\sqrt{s}\nabla_{x}^{2,\delta}\hat{p}_{\lfloor s\rfloor}^{\delta}(\lfloor\sqrt{s}w\rfloor)
\left(\int_{x-\sqrt{s}w}^{\lfloor x-\sqrt{s}w\rfloor+\delta}\tilde{v}_{t-s}^{\delta}(z)d\hat{U}_1^{\delta}(z)\right) dwds.
\end{equation}
Next recall from (\ref{eq:estim-nabla-hat-p}) that 
\begin{equation}
\big|\nabla_{x}^{2,\delta}\hat{p}_{\lfloor s\rfloor}^{\delta}(\lfloor\sqrt{s}w\rfloor)\big|\leqslant\frac{C_{1}}{s^{3/2}}e^{-C_{2}w^{2}},\ \ \ \forall s\geqslant\delta^{2},w\in\mathbb{R}.\label{eq:E1Del4}
\end{equation}
In addition, invoking again a discrete version of the rough integral estimate \eqref{eq:RI1}, we have
\begin{multline}
  \Big|\int_{x-\sqrt{s}w}^{\lfloor x-\sqrt{s}w\rfloor+\delta}\tilde{v}_{t-s}^{\delta}(z)d\hat{U}_1^{\delta}(z)\Big| 
  \leqslant 
  C\kappa(\hat{{\bf U}}_1^{\delta})\Lambda^{\theta',\beta}\lambda^{\frac{\alpha-\beta}{4}}
  E^{\theta'}(a,t)
 e^{-\left(\lambda+\theta'(a+\sqrt{T}|w|)\right)s} \\ 
  \times\Big((a+\sqrt{T}|w|)^{\chi}\delta^{\alpha}
+(a+\sqrt{T}|w|)^{2\chi}\delta^{2\alpha}
+(a+\sqrt{T}|w|)^{2\chi+\beta/2}(\delta^{2\alpha+\beta}+\delta^{\alpha+2\beta}) \\
+\big(a+\sqrt{T}|w|\big)^{2\chi}(t-s)^{-\beta/2}\delta^{\alpha+2\beta}\Big),\label{eq:E1Del3}
\end{multline}
where we recall that $\theta'<\theta$ and $\Lambda^{\theta',\beta}$ denotes a uniform
upper bound of $\Theta^{\theta',\lambda}(v^{\delta})$. As before,
the role of $\theta'$ is to absorb the extra polynomial factors in
$a$, i.e. we choose $\theta'<\theta$ so that~(\ref{eq:The'Cond}) holds. By substituting (\ref{eq:E1Del4}) and (\ref{eq:E1Del3}) into the
expression \eqref{eq:calK2delta} of $\mathcal{K}_{t}^{2,\delta}(x),$ there are four individual
terms coming out (corresponding to the four terms in (\ref{eq:E1Del3})).
The first three terms lead to an $s$-integral 
\[
\int_{\delta^{2}}^{t}\frac{1}{s}ds=\log\frac{t}{\delta^{2}}=O(\log\delta).
\]
The last term leads to an $s$-integral 
\[
\int_{\delta^{2}}^{t}s^{-1}(t-s)^{-\beta/2}ds=t^{-\beta/2}\int_{\delta^{2}/t}^{1}\rho^{-1}(1-\rho)^{-\beta/2}d\rho\leqslant\delta^{-\beta}O(\log\delta).
\]
To summarize, we obtain that 
\begin{equation}\label{eq:estimate-calK2delta}
\big|\mathcal{K}_{t}^{2,\delta}(x)\big|\leqslant C\kappa(\hat{{\bf U}}_1^{\delta})\Lambda^{\theta',\beta}\lambda^{\frac{\alpha-\beta}{4}}E^\theta(a,t)\delta^{\alpha}\log\delta.
\end{equation}
Combining the above estimates (also noting that $\alpha-2\chi < \alpha < 1/2$ and $\Lambda^{\theta,\beta}\leqslant \Lambda^{\theta',\beta}$), we can now conclude our $\mathcal{E}^{2,\delta}$-estimate as 
\begin{equation}\label{eq:estimate-calE2delta}
\big|\mathcal{E}_{t}^{2,\delta}(x)\big|\leqslant C (\kappa(\hat{{\bf U}}_1^\delta)+1)\Lambda^{\theta',\beta}\lambda^{\frac{\alpha-\beta}{4}}E^\theta(a,t)\cdot \delta^{\alpha-2\chi} .
\end{equation}
Eventually putting together \eqref{eq:estimate-calE1delta} and 
\eqref{eq:estimate-calE2delta} into the decomposition \eqref{eq:calEdelta=calE1delta+calE2delta} we obtain the desired estimate \eqref{eq:norm-calEdelta}.
\end{proof}

We finish this section by proving our main result.
\begin{proof}[Proof of Proposition \ref{prop:ComDetJEta}]
Gathering our estimates in Sections \ref{sec:comparing-eta+etadelta2-and-2eta} to \ref{sec:estim-calEdelat} the proof of relation~\eqref{eq:S2Est} is achieved. 
\end{proof}

\subsection{Completing the proof of Proposition \ref{prop:MainEst}} We now put together all estimates obtained so far to establish our main estimate (\ref{eq:MainEst}). To ease notation, we will write 
\[
\kappa_{\delta}\triangleq\kappa_{\alpha,\chi}(\hat{\mathbf{U}}_{1}^{\delta}),\kappa\triangleq\kappa_{\alpha,\chi}(\mathbf{W}),\rho_{\delta}\triangleq\rho_{\alpha,\chi}(\hat{\mathbf{U}}_{1}^{\delta},\mathbf{W}),\Theta\triangleq\Theta^{\theta,\lambda}(v),\Theta_{\delta}\triangleq\Theta^{\theta,\lambda}(\tilde{v}^{\delta}),
\]
and also
\[
d(\cdot,\cdot)\triangleq d_{\hat{\mathbf{U}}_{1}^{\delta},\mathbf{W}}(\cdot,\cdot).
\]
In addition, $C_{i}$ will denote constants depending on all underlying
exponents but not on $\delta,\lambda$. With all exponents $\alpha,\beta,\chi,\theta',\beta'$
given fixed, it follows from Lemmas \ref{lem:remainder-in-WiDel},  \ref{lem:UnifWDelvsW}, \ref{lem:TimVarWDelvsW}, \ref{lem:SpVarWDelvsW}, \ref{lem:RemEstComp} and Proposition \ref{prop:ComDetJEta} that 
\begin{align}
d(\tilde{v}^{\delta},v) & \leqslant d(\mathcal{W}^{1,\delta},\mathcal{W})+\Theta_{\hat{\mathbf{U}}_{1}^{\delta}}^{\theta,\lambda}(\mathcal{W}^{2,\delta})+d(\mathcal{J}^{\delta},\eta)\nonumber \\
 & \leqslant C_{1}\left[\lambda^{-\frac{\alpha-\beta}{4}}\kappa d(\tilde{v}^{\delta},v)+(\|f_{0}\|_{\mathcal{C}_{L}^{3}}+\|g\|_{\mathcal{C}_{L}^{3}})\delta^{\beta}\right.\nonumber \\
 & \ \ \ \left.+\lambda^{\frac{\alpha-\beta}{4}}(1+\kappa_{\delta}+\kappa)(\Theta_{\delta}+\Theta+\Lambda^{\theta',\beta'})(\rho_{\delta}+\delta^{\frac{\beta'(\beta'-\beta)}{\beta'+\beta}})\right].\label{eq:dEst1}
\end{align}
We now fix a choice of $\lambda=\lambda_{\omega}$ which satisfies
\[
\lambda^{-\frac{\alpha-\beta'}{4}}\bar{\kappa}(\omega)\leqslant \frac{1}{4}
\]
so that (cf. Lemma \ref{lem:UniEstPLI} and Lemma \ref{lem:lambda-betabetaprime})
\[
\sup_{\delta\in(0,1]}\Theta_{\delta}\vee\Theta\vee\Lambda^{\theta',\beta'}\leqslant C_{2}e^{C_{3}\lambda}\left(\|f_{0}\|_{\mathcal{C}_{L}^{3}}+\|g\|_{\mathcal{C}_{L}^{3}}\right)\bar{\kappa}(\omega),
\]where we recall that $\bar{\kappa}(\omega)$ is defined by (\ref{eq:BarKap}).
It follows from (\ref{eq:dEst1}) that 
\begin{align}
d(\tilde{v}^{\delta},v) & \leqslant\lambda^{-\frac{\alpha-\beta}{4}}C_{1}\kappa d(\tilde{v}^{\delta},v)+C_{4}e^{C_{5}\lambda}\left(\|f_{0}\|_{\mathcal{C}_{L}^{3}}+\|g\|_{\mathcal{C}_{L}^{3}}\right)\left(\delta^{\beta}+(1+\bar{\kappa}(\omega))^2(\rho_{\delta}+\delta^{\frac{\beta'(\beta'-\beta)}{\beta'+\beta}})\right).\label{eq:dEst2}
\end{align}
By enlarging $\lambda$ if necessary to ensure that $\lambda^{-\frac{\alpha-\beta}{4}}C_{1}\kappa\leqslant 1/2$,
one can then move the first term on the right hand side of (\ref{eq:dEst2})
to the left and the desired estimate (\ref{eq:MainEst}) thus follows.

\section*{Acknowledgement} XG gratefully acknowledges the support by ARC grant DE210101352. MG expresses his appreciation of the support by "Défis Scientifiques" (2019, 2021, 2023) of Université de Rennes.  ST is supported for this work by NSF grant DMS-2153915.

%%%%%%%%%%%%%%%%%%%%%%%%%%%%%%%%%%%%%%

\Addresses


\begin{thebibliography}{10}
\bibitem{AKS} Anshelevich, V. V.; Khanin, K. M.; Sinai, Ya. G. Symmetric
random walks in random environments. \textit{Comm. Math. Phys.} \textbf{85}
(1982), no. 3, 449-470.

\bibitem{BT} Bally, V.; Talay, D. The law of the Euler scheme for stochastic differential equations. I: Convergence rate of the distribution function.\textit{Probab. Theory Relat. Fields} 104 (1996), 43-60.

\bibitem{BC} Bass, R.F.; Chen Z.Q. Stochastic differential equations for Dirichlet processes. \textit{Probab. Theory Relat. Fields} 121 (2001), 422-446.																																				
\bibitem{BZ} Bounebache, S. K. ; Zambotti, L. A skew stochastic
heat equation. \textit{J. Theoret. Probab.} \textbf{27} (2014), no.
1, 168-201.

\bibitem{BF} Bovier, A.; Faggionato, A. Spectral analysis
of Sinai's walk for small eigenvalues. \textit{Ann. Probab.} \textbf{36}
(2008), no. 1, 198-254.

\bibitem{Br} Brox, Th. A one-dimensional diffusion process in a Wiener
medium. \textit{Ann. Probab.} \textbf{14} (1986), no. 4, 1206-1218.

\bibitem{CG} Catellier, R.; Gubinelli, M. Averaging along irregular curves and regularisation of ODEs. \textit{Stochastic Processes Appl.} 126 (2016),  2323-2366. 

\bibitem{CLL} Cass, Th.; Litterer, C.; Lyons, T.  
Integrability and tail estimates for Gaussian rough differential equations.
\textit{Ann. Probab.} 41 (2013), 3026-3050.


\bibitem{DD} Delarue, F; Diel, R. Rough paths and 1d sde with a time
dependent distributional drift. Application to polymers. Probab. Theory
Related Fields 165 (2016), no. 1-2, 1-63.

\bibitem{EK} Ethier, S.N.; Kurtz, T.G. \textit{Markov
processes. Characterization and convergence.} John Wiley \& Sons,
1986.

\bibitem{FH14} Friz, P.; Hairer, M. \textit{A course on rough paths},
Springer, 2014.

\bibitem{FK} Fujita, T.; Kawamishi, Y. A proof of It\^o's formula using
a discrete It\^o's formula. \textit{Studia scientiarum mathematicarum
Hungarica} 45(2):125-134.

\bibitem{GH} Gradinaru, M.; Haugomat, T. Convergence and discrete
schemes applied to some L\'evy-type processes. \textit{ArXiv }1707.02889.

\bibitem{Gu} Gubinelli, M. Controlling rough paths. \textit{J. Funct.
Anal.} \textbf{216}, 86-140 (2004).

\bibitem{HS93}Hebisch, W.,  Saloff-Coste, L. Gaussian estimates for Markov chains and random walks on groups.  \textit{Ann. Probab.} \textbf{21} (1993), no. 2, 673-709.

\bibitem{HLM} Hu, Y.; L\^e, K.; Mytnik, L. Stochastic differential
equation for Brox diffusion. \textit{Stochastic Process. Appl.} 127
(2017), no. 7, 2281-2315.

\bibitem{HS} Hu, Y.; Shi, Z. The local time of simple random walk in random environment. \textit{J. Theor. Probab.} 11 (1998), 765-793. 

\bibitem{KS} Karatzas, I.; Shreve, S.E. Brownian motion and stochastic calculus. Graduate Texts in Mathematics, 113. New York etc.: Springer-Verlag, 2nd Ed. 1991.

\bibitem{Ke76} Kesten, H. The limit distribution of Sinai's random walk in random environment, \textit{Physica A} 138 (1986), 299-309

\bibitem{KLO} Komorowski, T.; Landim, C.; Olla, S. \textit{Fluctuations
in Markov processes.} Springer, 2012.

\bibitem{KMT76} Koml\'os, J; Major, P., Tusn\'ady, G. An approximation
of partial sums of independent RV's, and the sample DF. II, \textit{Z.
Wahrscheinlichkeitstheorie verw. Gebiete} 34 (1976) 33--58, 1976.

\bibitem{La} Lawler, G. \textit{Random walk and the heat equation}
AMS Student Mathematical Library, 2010.

\bibitem{LLT} Léon, J. A.; Liu, Y.; Tindel, S. Euler scheme for SDEs driven by fractional Brownian motions: integrability and convergence in law, arXiv:2307.06759 (2023)

\bibitem{Na} Nagy, K. Symmetric random walk in random environment
in one dimension. \textit{Period. Math. Hungar.} \textbf{45} (2002),
no. 1-2, 101-120.

\bibitem{LL10} Lawler, G. F.; Limic, V. Random walk: a modern introduction,
Cambridge University Press, 2010.

\bibitem{LT} Liu, Y.; Tindel, S. First-order Euler scheme for
SDEs driven by fractional Brownian motion: the rough case, \textit{Ann.
Appl. Probab.} 29 (2) (2019): 758--826.

\bibitem{Se} Seignourel, P. Discrete schemes for processes in random
media. \textit{Probab. Theory Related Fields} \textbf{118} (2000),
no. 3, 293-322.

\bibitem{Si} Sinai, Ya. G. The limiting behaviour of a one-dimensional random walk in a random medium \textit{Th. Probab. Appl.} 27 (1982), 256-268.

\bibitem{So} F. Solomon Random walks in a random environment \textit{Ann. Prob.} 3 (1975), 1-31


\bibitem{Sh} Shi, Z. Sinai's walk via stochastic calculus. Milieux
al\'eatoires, 53-74, \textit{Panor. Synth\`{e}ses} \textbf{12}, Soc.
Math. France, Paris, 2001.

\bibitem{SV} Stroock, D., Varadhan, S.R.S. Multidimensional diffusion processes. Reprint of the 2nd corrected printing. Classics in Mathematics. Berlin: Springer (2006). 
\bibitem{Ze} Zeitouni, O. Random walks in random environment,
LNM 1837, J. Picard (Ed.), pp. 189--312, 2004.
\end{thebibliography}
\end{document}